\documentclass{article}

\usepackage{amssymb,latexsym,epic,eepic,rotating}

\newcommand{\xbD}{\Delta}

\newcommand{\xbG}{\Gamma}

\newcommand{\xbO}{\Omega}
\newcommand{\xbP}{\Pi}

\newcommand{\xbS}{\Sigma}

\newcommand{\xba}{\alpha}
\newcommand{\xbb}{\beta}

\newcommand{\xbd}{\delta}
\newcommand{\xbe}{\in}
\newcommand{\xbf}{\phi}
\newcommand{\xbg}{\gamma}
\newcommand{\xbh}{\eta}

\newcommand{\xbj}{\vartheta}
\newcommand{\xbk}{\kappa}
\newcommand{\xbl}{\lambda}
\newcommand{\xbm}{\mu}

\newcommand{\xbo}{\omega}

\newcommand{\xbq}{\psi}
\newcommand{\xbr}{\rho}
\newcommand{\xbs}{\sigma}
\newcommand{\xbt}{\tau}

\newcommand{\xCK}{\times}

\newcommand{\xCN}{\neg}
\newcommand{\xCQ}{\emptyset}

\newcommand{\xCd}{\approx}

\newcommand{\xCf}{\hspace{0.1em}}

\newcommand{\xcA}{\forall}

\newcommand{\xcC}{\not\subseteq}

\newcommand{\xcE}{\exists}

\newcommand{\xcH}{\not\Rightarrow}
\newcommand{\xcI}{\not\Leftarrow}
\newcommand{\xcJ}{\not\Leftrightarrow}

\newcommand{\xcL}{\not\vdash}
\newcommand{\xcM}{\not\models}
\newcommand{\xcN}{\hspace{0.2em}\not\sim\hspace{-0.9em}\mid\hspace{0.8em}}
\newcommand{\xcO}{\bigvee}
\newcommand{\xcP}{\not\rightarrow}

\newcommand{\xcS}{\bigcap}
\newcommand{\xcT}{\bot}
\newcommand{\xcU}{\bigwedge}
\newcommand{\xcV}{\bigcup}

\newcommand{\xcX}{\Box}

\newcommand{\xca}{\infty}
\newcommand{\xcb}{\subset}
\newcommand{\xcc}{\subseteq}
\newcommand{\xcd}{\supseteq}
\newcommand{\xce}{\not\in}
\newcommand{\xcf}{\supset}
\newcommand{\xcg}{\geq}
\newcommand{\xch}{\Rightarrow}
\newcommand{\xci}{\Leftarrow}
\newcommand{\xcj}{\Leftrightarrow}
\newcommand{\xck}{\leq}
\newcommand{\xcl}{\vdash}
\newcommand{\xcm}{\models}
\newcommand{\xcn}{\hspace{0.2em}\sim\hspace{-0.9em}\mid\hspace{0.58em}}

\newcommand{\xco}{\vee}
\newcommand{\xcp}{\rightarrow}
\newcommand{\xcq}{\leftarrow}
\newcommand{\xcr}{\leftrightarrow}
\newcommand{\xcs}{\cap}
\newcommand{\xcu}{\wedge}
\newcommand{\xcv}{\cup}
\newcommand{\xcx}{\Diamond}

\newcommand{\xcz}{\Box}

\newcommand{\xDC}{\hspace{2em}}
\newcommand{\xDH}{\item }

\newcommand{\xDM}{\circ}
\newcommand{\xDN}{\ominus}

\newcommand{\xDd}{\equiv}

\newcommand{\xdA}{\mbox{\boldmath$A$}}

\newcommand{\xdC}{\mbox{\boldmath$C$}}
\newcommand{\xdD}{\mbox{\boldmath$D$}}

\newcommand{\xdL}{\mbox{\boldmath$L$}}

\newcommand{\xdO}{\mbox{\boldmath$O$}}
\newcommand{\xdP}{\mbox{\boldmath$P$}}

\newcommand{\xdR}{\Re}

\newcommand{\xda}{{\cal A}}

\newcommand{\xdc}{{\cal C}}

\newcommand{\xdf}{{\cal F}}

\newcommand{\xdi}{{\cal I}}

\newcommand{\xdl}{{\cal L}}
\newcommand{\xdm}{{\cal M}}
\newcommand{\xdn}{{\cal N}}

\newcommand{\xdp}{{\cal P}}

\newcommand{\xds}{{\cal S}}

\newcommand{\xdu}{{\cal U}}

\newcommand{\xdw}{{\cal W}}
\newcommand{\xdx}{{\cal X}}
\newcommand{\xdy}{{\cal Y}}
\newcommand{\xdz}{{\cal Z}}

\newcommand{\xEH}{ & }
\newcommand{\xEI}{\begin{itemize}}
\newcommand{\xEJ}{\end{itemize}}
\newcommand{\xEP}{ \\ }

\newcommand{\xEd}{\neq}
\newcommand{\xEh}{\begin{enumerate}}
\newcommand{\xEj}{\end{enumerate}}

\newcommand{\xeA}{\nabla}

\newcommand{\xeb}{\prec}
\newcommand{\xec}{\preceq}
\newcommand{\xed}{\succeq}
\newcommand{\xee}{\succ}

\newcommand{\xej}{\lhd}

\newcommand{\xem}{\rhd}

\newcommand{\xes}{\sqsubseteq}

\newcommand{\xex}{\lceil}

\newcommand{\xFB}{\cdots}
\newcommand{\xFO}{\parallel}

\newcommand{\xfA}{\mid}
\newcommand{\xfB}{\uparrow}

\newcommand{\xfb}{\downarrow}

\newcommand{\xfo}{\hookrightarrow}

\newcommand{\xfq}{\cdot \cdot \cdot \rightarrow}

\newcommand{\Xl}{\ldots}

\newcommand{\ol}{\overline}
\newcommand{\ul}{\underline}

\newcommand{\wt}{\overbrace}

\newcommand{\xssc}{\scriptsize}

\newcommand{\bl}{\begin{lemma} \rm}
\newcommand{\el}{\end{lemma}}
\newcommand{\br}{\begin{remark} \rm}
\newcommand{\er}{\end{remark}}
\newcommand{\be}{\begin{example} \rm}
\newcommand{\ee}{\end{example}}
\newcommand{\bco}{\begin{corollary} \rm}
\newcommand{\eco}{\end{corollary}}
\newcommand{\bc}{\begin{claim} \rm}
\newcommand{\ec}{\end{claim}}
\newcommand{\bfa}{\begin{fact} \rm}
\newcommand{\efa}{\end{fact}}
\newcommand{\bp}{\begin{proposition} \rm}
\newcommand{\ep}{\end{proposition}}
\newcommand{\bd}{\begin{definition} \rm}
\newcommand{\ed}{\end{definition}}
\newcommand{\bcs}{\begin{construction} \rm}
\newcommand{\ecs}{\end{construction}}
\newcommand{\bcd}{\begin{condition} \rm}
\newcommand{\ecd}{\end{condition}}
\newcommand{\bt}{\begin{theorem} \rm}
\newcommand{\et}{\end{theorem}}
\newcommand{\bn}{\begin{notation} \rm}
\newcommand{\en}{\end{notation}}
\newcommand{\bfi}{\begin{bild} \rm}
\newcommand{\efi}{\end{bild}}
\newcommand{\bsta}{\begin{statement} \rm}
\newcommand{\esta}{\end{statement}}
\newcommand{\bcom}{\begin{comment} \rm}
\newcommand{\ecom}{\end{comment}}
\newcommand{\bdia}{\begin{diagram} \rm}
\newcommand{\edia}{\end{diagram}}

\newcommand{\bfc}{\begin{figure}[htb] \begin{center}}
\newcommand{\efc}{\end{center} \end{figure}}

\title{
Toolbox
}

\author{
Karl Schlechta
}

\begin{document}

\newtheorem{lemma}{Lemma}[section]
\newtheorem{theorem}[lemma]{Theorem}
\newtheorem{proposition}[lemma]{Proposition}
\newtheorem{corollary}[lemma]{Corollary}
\newtheorem{claim}[lemma]{Claim}
\newtheorem{fact}[lemma]{Fact}
\newtheorem{remark}[lemma]{Remark}
\newtheorem{definition}{Definition}[section]
\newtheorem{construction}{Construction}[section]
\newtheorem{condition}{Condition}[section]
\newtheorem{example}{Example}[section]
\newtheorem{notation}{Notation}[section]
\newtheorem{bild}{Figure}[section]
\newtheorem{comment}{Comment}[section]
\newtheorem{statement}{Statement}[section]
\newtheorem{diagram}{Diagram}[section]

\renewcommand{\labelenumi}
  {(\arabic{enumi})}
\renewcommand{\labelenumii}
  {(\arabic{enumi}.\arabic{enumii})}
\renewcommand{\labelenumiii}
  {(\arabic{enumi}.\arabic{enumii}.\arabic{enumiii})}
\renewcommand{\labelenumiv}
  {(\arabic{enumi}.\arabic{enumii}.\arabic{enumiii}.\arabic{enumiv})}

\maketitle

\setcounter{secnumdepth}{4}
\setcounter{tocdepth}{4}

\tableofcontents

%
%
%

\section{
Algebraic part
}

\subsubsection{ToolBase1-Alg}

 {\LARGE karl-search= Start ToolBase1-Alg }

{\xssc LABEL: {Section Toolbase1-Alg}}
\label{Section Toolbase1-Alg}
\index{Section Toolbase1-Alg}
\subsubsection{Definition Alg-Base}

 {\LARGE karl-search= Start Definition Alg-Base }

\index{Definition Alg-Base}

\bd

$\hspace{0.01em}$

(+++ Orig. No.:  Definition Alg-Base +++)

{\xssc LABEL: {Definition Alg-Base}}
\label{Definition Alg-Base}

We use $ \xdp $ to denote the power set operator,
$ \xbP \{X_{i}:i \xbe I\}$ $:=$ $\{g:$ $g:I \xcp \xcV \{X_{i}:i \xbe I\},$
$ \xcA i \xbe I.g(i) \xbe X_{i}\}$ is the general cartesian
product, $card(X)$ shall denote the cardinality of $X,$ and $V$ the
set-theoretic
universe we work in - the class of all sets. Given a set of pairs $ \xdx
,$ and a
set $X,$ we denote by $ \xdx \xex X:=\{<x,i> \xbe \xdx:x \xbe X\}.$ When
the context is clear, we
will sometime simply write $X$ for $ \xdx \xex X.$ (The intended use is
for preferential
structures, where $x$ will be a point (intention: a classical
propositional
model), and $i$ an index, permitting copies of logically identical
points.)

$A \xcc B$ will denote that $ \xCf A$ is a subset of $B$ or equal to $B,$
and $A \xcb B$ that $ \xCf A$ is
a proper subset of $B,$ likewise for $A \xcd B$ and $A \xcf B.$

Given some fixed set $U$ we work in, and $X \xcc U,$ then $ \xdC (X):=U-X$
.

If $ \xdy \xcc \xdp (X)$ for some
$X,$ we say that $ \xdy $ satisfies

$( \xcs )$ iff it is closed under finite intersections,

$( \xcS )$ iff it is closed under arbitrary intersections,

$( \xcv )$ iff it is closed under finite unions,

$( \xcV )$ iff it is closed under arbitrary unions,

$( \xdC )$ iff it is closed under complementation,

$ \xCf (-)$ iff it is closed under set difference.

We will sometimes write $A=B \xFO C$ for: $A=B,$ or $A=C,$ or $A=B \xcv
C.$

We make ample and tacit use of the Axiom of Choice.

 karl-search= End Definition Alg-Base
\vspace{7mm}

 *************************************

\vspace{7mm}

\subsubsection{Definition Rel-Base}

 {\LARGE karl-search= Start Definition Rel-Base }

\index{Definition Rel-Base}

\ed

\bd

$\hspace{0.01em}$

(+++ Orig. No.:  Definition Rel-Base +++)

{\xssc LABEL: {Definition Rel-Base}}
\label{Definition Rel-Base}

$ \xeb^{*}$ will denote the transitive closure of the relation $ \xeb.$
If a relation $<,$
$ \xeb,$ or similar is given, $a \xcT b$ will express that a and $b$ are
$<-$ (or $ \xeb -)$
incomparable - context will tell. Given any relation $<,$ $ \xck $ will
stand for
$<$ or $=,$ conversely, given $ \xck,$ $<$ will stand for $ \xck,$ but
not $=,$ similarly
for $ \xeb $ etc.

 karl-search= End Definition Rel-Base
\vspace{7mm}

 *************************************

\vspace{7mm}

\subsubsection{Definition Tree-Base}

 {\LARGE karl-search= Start Definition Tree-Base }

\index{Definition Tree-Base}

\ed

\bd

$\hspace{0.01em}$

(+++ Orig. No.:  Definition Tree-Base +++)

{\xssc LABEL: {Definition Tree-Base}}
\label{Definition Tree-Base}

A child (or successor) of an element $x$ in a tree $t$ will be a direct
child in $t.$
A child of a child, etc. will be called an indirect child. Trees will be
supposed to grow downwards, so the root is the top element.

 karl-search= End Definition Tree-Base
\vspace{7mm}

 *************************************

\vspace{7mm}

\subsubsection{Definition Seq-Base}

 {\LARGE karl-search= Start Definition Seq-Base }

\index{Definition Seq-Base}

\ed

\bd

$\hspace{0.01em}$

(+++ Orig. No.:  Definition Seq-Base +++)

{\xssc LABEL: {Definition Seq-Base}}
\label{Definition Seq-Base}

A subsequence $ \xbs_{i}:i \xbe I \xcc \xbm $ of a sequence $ \xbs_{i}:i
\xbe \xbm $ is called cofinal, iff
for all $i \xbe \xbm $ there is $i' \xbe I$ $i \xck i'.$

Given two sequences $ \xbs_{i}$ and $ \xbt_{i}$ of the same length, then
their Hamming distance
is the quantity of $i$ where they differ.

 karl-search= End Definition Seq-Base
\vspace{7mm}

 *************************************

\vspace{7mm}

 karl-search= End ToolBase1-Alg
\vspace{7mm}

 *************************************

\vspace{7mm}

\subsubsection{Lemma Abs-Rel-Ext}

 {\LARGE karl-search= Start Lemma Abs-Rel-Ext }

\index{Lemma Abs-Rel-Ext}

\ed

We give a generalized abstract nonsense result, taken
from  \cite{LMS01}, which must be part of the folklore:

\bl

$\hspace{0.01em}$

(+++ Orig. No.:  Lemma Abs-Rel-Ext +++)

{\xssc LABEL: {Lemma Abs-Rel-Ext}}
\label{Lemma Abs-Rel-Ext}

Given a set $X$ and a binary relation $R$ on $X,$ there exists a total
preorder (i.e.
a total, reflexive, transitive relation) $S$ on $X$ that extends $R$ such
that

$ \xcA x,y \xbe X(xSy,ySx \xch xR^{*}y)$

where $R^{*}$ is the reflexive and transitive closure of $R.$

 karl-search= End Lemma Abs-Rel-Ext
\vspace{7mm}

 *************************************

\vspace{7mm}

\subsubsection{Lemma Abs-Rel-Ext Proof}

 {\LARGE karl-search= Start Lemma Abs-Rel-Ext Proof }

\index{Lemma Abs-Rel-Ext Proof}

\el

\subparagraph{
Proof
}

$\hspace{0.01em}$

(+++*** Orig.:  Proof )

Define $x \xDd y$ iff $xR^{*}y$ and $yR^{*}x.$
The relation $ \xDd $ is an equivalence relation.
Let $[x]$ be the equivalence class of $x$ under $ \xDd.$ Define $[x] \xec
[y]$ iff $xR^{*}y.$
The definition of $ \xec $ does not depend on the representatives $x$ and
$y$ chosen.
The relation $ \xec $ on equivalence classes is a partial order.
Let $ \xck $ be any total order on these equivalence classes that extends
$ \xec.$
Define xSy iff $[x] \xck [y].$
The relation $S$ is total (since $ \xck $ is total) and transitive
(since $ \xck $ is transitive) and is therefore a total preorder.
It extends $R$ by the definition of $ \xec $ and the fact that $ \xck $
extends $ \xec.$
Suppose now xSy and ySx. We have $[x] \xck [y]$ and $[y] \xck [x]$
and therefore $[x]=[y]$ by antisymmetry. Therefore $x \xDd y$ and
$xR^{*}y.$
$ \xcz $
\\[3ex]

 karl-search= End Lemma Abs-Rel-Ext Proof
\vspace{7mm}

 *************************************

\vspace{7mm}

\newpage

\section{
Logical rules
}

\subsubsection{ToolBase1-Log}

 {\LARGE karl-search= Start ToolBase1-Log }

{\xssc LABEL: {Section Toolbase1-Log}}
\label{Section Toolbase1-Log}
\index{Section Toolbase1-Log}

\subsection{
Logics: Base
}

\subsubsection{ToolBase1-Log-Base}

 {\LARGE karl-search= Start ToolBase1-Log-Base }

{\xssc LABEL: {Section Toolbase1-Log-Base}}
\label{Section Toolbase1-Log-Base}
\index{Section Toolbase1-Log-Base}
\subsubsection{Definition Log-Base}

 {\LARGE karl-search= Start Definition Log-Base }

\index{Definition Log-Base}

\bd

$\hspace{0.01em}$

(+++ Orig. No.:  Definition Log-Base +++)

{\xssc LABEL: {Definition Log-Base}}
\label{Definition Log-Base}

We work here in a classical propositional language $ \xdl,$ a theory $T$
will be an
arbitrary set of formulas. Formulas will often be named $ \xbf,$ $ \xbq
,$ etc., theories
$T,$ $S,$ etc.

$v( \xdl )$ will be the set of propositional variables of $ \xdl.$

$M_{ \xdl }$ will be the set of (classical) models for $ \xdl,$ $M(T)$ or
$M_{T}$
is the set of models of $T,$ likewise $M( \xbf )$ for a formula $ \xbf.$

$ \xdD_{ \xdl }:=\{M(T):$ $T$ a theory in $ \xdl \},$ the set of definable
model sets.

Note that, in classical propositional logic, $ \xCQ,M_{ \xdl } \xbe
\xdD_{ \xdl },$ $ \xdD_{ \xdl }$ contains
singletons, is closed under arbitrary intersections and finite unions.

An operation $f: \xdy \xcp \xdp (M_{ \xdl })$ for $ \xdy \xcc \xdp (M_{
\xdl })$ is called definability
preserving, $ \xCf (dp)$ or $( \xbm dp)$ in short, iff for all $X \xbe
\xdD_{ \xdl } \xcs \xdy $ $f(X) \xbe \xdD_{ \xdl }.$

We will also use $( \xbm dp)$ for binary functions $f: \xdy \xCK \xdy \xcp
\xdp (M_{ \xdl })$ - as needed
for theory revision - with the obvious meaning.

$ \xcl $ will be classical derivability, and

$ \ol{T}:=\{ \xbf:T \xcl \xbf \},$ the closure of $T$ under $ \xcl.$

$Con(.)$ will stand for classical consistency, so $Con( \xbf )$ will mean
that
$ \xbf $ is clasical consistent, likewise for $Con(T).$ $Con(T,T' )$ will
stand for
$Con(T \xcv T' ),$ etc.

Given a consequence relation $ \xcn,$ we define

$ \ol{ \ol{T} }:=\{ \xbf:T \xcn \xbf \}.$

(There is no fear of confusion with $ \ol{T},$ as it just is not useful to
close
twice under classical logic.)

$T \xco T':=\{ \xbf \xco \xbf ': \xbf \xbe T, \xbf ' \xbe T' \}.$

If $X \xcc M_{ \xdl },$ then $Th(X):=\{ \xbf:X \xcm \xbf \},$ likewise
for $Th(m),$ $m \xbe M_{ \xdl }.$ $( \xcm $ will
usually be classical validity.)

 karl-search= End Definition Log-Base
\vspace{7mm}

 *************************************

\vspace{7mm}

\subsubsection{Fact Log-Base}

 {\LARGE karl-search= Start Fact Log-Base }

\index{Fact Log-Base}

\ed

We recollect and note:

\bfa

$\hspace{0.01em}$

(+++ Orig. No.:  Fact Log-Base +++)

{\xssc LABEL: {Fact Log-Base}}
\label{Fact Log-Base}

Let $ \xdl $ be a fixed propositional language, $ \xdD_{ \xdl } \xcc X,$ $
\xbm:X \xcp \xdp (M_{ \xdl }),$ for a $ \xdl -$theory $T$
$ \ol{ \ol{T} }:=Th( \xbm (M_{T})),$ let $T,$ $T' $ be arbitrary theories,
then:

(1) $ \xbm (M_{T}) \xcc M_{ \ol{ \ol{T} }}$,

(2) $M_{T} \xcv M_{T' }=M_{T \xco T' }$ and $M_{T \xcv T' }=M_{T} \xcs
M_{T' }$,

(3) $ \xbm (M_{T})= \xCQ $ $ \xcr $ $ \xcT \xbe \ol{ \ol{T} }$.

If $ \xbm $ is definability preserving or $ \xbm (M_{T})$ is finite, then
the following also hold:

(4) $ \xbm (M_{T})=M_{ \ol{ \ol{T} }}$,

(5) $T' \xcl \ol{ \ol{T} }$ $ \xcr $ $M_{T' } \xcc \xbm (M_{T}),$

(6) $ \xbm (M_{T})=M_{T' }$ $ \xcr $ $ \ol{T' }= \ol{ \ol{T} }.$
$ \xcz $
\\[3ex]

 karl-search= End Fact Log-Base
\vspace{7mm}

 *************************************

\vspace{7mm}

\subsubsection{Fact Th-Union}

 {\LARGE karl-search= Start Fact Th-Union }

\index{Fact Th-Union}

\efa

\bfa

$\hspace{0.01em}$

(+++ Orig. No.:  Fact Th-Union +++)

{\xssc LABEL: {Fact Th-Union}}
\label{Fact Th-Union}

Let $A,B \xcc M_{ \xdl }.$

Then $Th(A \xcv B)$ $=$ $Th(A) \xcs Th(B).$

 karl-search= End Fact Th-Union
\vspace{7mm}

 *************************************

\vspace{7mm}

\subsubsection{Fact Th-Union Proof}

 {\LARGE karl-search= Start Fact Th-Union Proof }

\index{Fact Th-Union Proof}

\efa

\subparagraph{
Proof
}

$\hspace{0.01em}$

(+++*** Orig.:  Proof )

$ \xbf \xbe Th(A \xcv B)$ $ \xcj $ $A \xcv B \xcm \xbf $ $ \xcj $ $A \xcm
\xbf $ and $B \xcm \xbf $ $ \xcj $ $ \xbf \xbe Th(A)$ and $ \xbf \xbe
Th(B).$

$ \xcz $
\\[3ex]

 karl-search= End Fact Th-Union Proof
\vspace{7mm}

 *************************************

\vspace{7mm}

 karl-search= End Log-Base
\vspace{7mm}

 *************************************

\vspace{7mm}

\subsubsection{Fact Log-Form}

 {\LARGE karl-search= Start Fact Log-Form }

\index{Fact Log-Form}

\bfa

$\hspace{0.01em}$

(+++ Orig. No.:  Fact Log-Form +++)

{\xssc LABEL: {Fact Log-Form}}
\label{Fact Log-Form}

Let $X \xcc M_{ \xdl },$ $ \xbf, \xbq $ formulas.

(1) $X \xcs M( \xbf ) \xcm \xbq $ iff $X \xcm \xbf \xcp \xbq.$

(2) $X \xcs M( \xbf ) \xcm \xbq $ iff $M(Th(X)) \xcs M( \xbf ) \xcm \xbq
.$

(3) $Th(X \xcs M( \xbf ))= \ol{Th(X) \xcv \{ \xbf \}}$

(4) $X \xcs M( \xbf )= \xCQ $ $ \xcj $ $M(Th(X)) \xcs M( \xbf )= \xCQ $

(5) $Th(M(T) \xcs M(T' ))= \ol{T \xcv T' }.$

 karl-search= End Fact Log-Form
\vspace{7mm}

 *************************************

\vspace{7mm}

\subsubsection{Fact Log-Form Proof}

 {\LARGE karl-search= Start Fact Log-Form Proof }

\index{Fact Log-Form Proof}

\efa

\subparagraph{
Proof
}

$\hspace{0.01em}$

(+++*** Orig.:  Proof )

(1) `` $ \xch $ '': $X=(X \xcs M( \xbf )) \xcv (X \xcs M( \xCN \xbf )).$ In
both parts holds $ \xCN \xbf \xco \xbq,$ so
$X \xcm \xbf \xcp \xbq.$ `` $ \xci $ '': Trivial.

(2) $X \xcs M( \xbf ) \xcm \xbq $ (by (1)) iff $X \xcm \xbf \xcp \xbq $
iff $M(Th(X)) \xcm \xbf \xcp \xbq $ iff (again by (1))
$M(Th(X)) \xcs M( \xbf ) \xcm \xbq.$

(3) $ \xbq \xbe Th(X \xcs M( \xbf ))$ $ \xcj $ $X \xcs M( \xbf ) \xcm \xbq
$ $ \xcj_{(2)}$ $M(Th(X) \xcv \{ \xbf \})=M(Th(X)) \xcs M( \xbf ) \xcm
\xbq $ $ \xcj $
$Th(X) \xcv \{ \xbf \} \xcl \xbq.$

(4) $X \xcs M( \xbf )= \xCQ $ $ \xcj $ $X \xcm \xCN \xbf $ $ \xcj $
$M(Th(X)) \xcm \xCN \xbf $ $ \xcj $ $M(Th(X)) \xcs M( \xbf )= \xCQ.$

(5) $M(T) \xcs M(T' )=M(T \xcv T' ).$

$ \xcz $
\\[3ex]

 karl-search= End Fact Log-Form Proof
\vspace{7mm}

 *************************************

\vspace{7mm}

 karl-search= End ToolBase1-Log-Base
\vspace{7mm}

 *************************************

\vspace{7mm}

\subsection{
Logics: Definability
}

\subsubsection{ToolBase1-Log-Dp}

 {\LARGE karl-search= Start ToolBase1-Log-Dp }

{\xssc LABEL: {Section Toolbase1-Log-Dp}}
\label{Section Toolbase1-Log-Dp}
\index{Section Toolbase1-Log-Dp}
\subsubsection{Fact Dp-Base}

 {\LARGE karl-search= Start Fact Dp-Base }

\index{Fact Dp-Base}

\bfa

$\hspace{0.01em}$

(+++ Orig. No.:  Fact Dp-Base +++)

{\xssc LABEL: {Fact Dp-Base}}
\label{Fact Dp-Base}

If $X=M(T),$ then $M(Th(X))=X.$

 karl-search= End Fact Dp-Base
\vspace{7mm}

 *************************************

\vspace{7mm}

\subsubsection{Fact Dp-Base Proof}

 {\LARGE karl-search= Start Fact Dp-Base Proof }

\index{Fact Dp-Base Proof}

\efa

\subparagraph{
Proof
}

$\hspace{0.01em}$

(+++*** Orig.:  Proof )

$X \xcc M(Th(X))$ is trivial. $Th(M(T))= \ol{T}$ is trivial by classical
soundness and
completeness. So $M(Th(M(T))=M( \ol{T})=M(T)=X.$ $ \xcz $
\\[3ex]

 karl-search= End Fact Dp-Base Proof
\vspace{7mm}

 *************************************

\vspace{7mm}

\subsubsection{Example Not-Def}

 {\LARGE karl-search= Start Example Not-Def }

\index{Example Not-Def}

\be

$\hspace{0.01em}$

(+++ Orig. No.:  Example Not-Def +++)

{\xssc LABEL: {Example Not-Def}}
\label{Example Not-Def}

If $v( \xdl )$ is infinite, and $m$ any model for $ \xdl,$ then $M:=M_{
\xdl }-\{m\}$ is not definable
by any theory $T.$ (Proof: Suppose it were, and let $ \xbf $ hold in $M,$
but not in $m,$ so in $m$ $ \xCN \xbf $ holds, but as $ \xbf $ is finite,
there is a model $m' $ in
$M$ which coincides on all propositional variables of $ \xbf $ with $m,$
so in $m' $ $ \xCN \xbf $
holds, too, a contradiction.) Thus, in the infinite case, $ \xdp (M_{ \xdl
}) \xEd \xdD_{ \xdl }.$

(There is also a simple cardinality argument, which shows that almost no
model sets are definable, but it is not constructive and thus less
instructive
than above argument. We give it nonetheless: Let $ \xbk:=card(v( \xdl
)).$ Then
there are $ \xbk $ many formulas, so $2^{ \xbk }$ many theories, and thus
$2^{ \xbk }$ many
definable model sets. But there are $2^{ \xbk }$ many models, so $(2^{
\xbk })^{ \xbk }$ many model
sets.)

$ \xcz $
\\[3ex]

 karl-search= End Example Not-Def
\vspace{7mm}

 *************************************

\vspace{7mm}

\subsubsection{Definition Def-Clos}

 {\LARGE karl-search= Start Definition Def-Clos }

\index{Definition Def-Clos}

\ee

\bd

$\hspace{0.01em}$

(+++ Orig. No.:  Definition Def-Clos +++)

{\xssc LABEL: {Definition Def-Clos}}
\label{Definition Def-Clos}

Let $ \xdy \xcc \xdp (Z)$ be given and closed under arbitrary
intersections.

For $A \xcc Z,$ let $ \wt{A}$ $:=$ $ \xcS \{X \xbe \xdy:A \xcc X\}.$

Intuitively, $Z$ is the set of all models for $ \xdl,$ $ \xdy $ is $
\xdD_{ \xdl }$, and $ \wt{A}=M(Th(A)),$
this is the intended application. Note that then $ \wt{ \xCQ }= \xCQ.$

 karl-search= End Definition Def-Clos
\vspace{7mm}

 *************************************

\vspace{7mm}

\subsubsection{Fact Def-Clos}

 {\LARGE karl-search= Start Fact Def-Clos }

\index{Fact Def-Clos}

\ed

\bfa

$\hspace{0.01em}$

(+++ Orig. No.:  Fact Def-Clos +++)

{\xssc LABEL: {Fact Def-Clos}}
\label{Fact Def-Clos}

(1) If $ \xdy \xcc \xdp (Z)$ is closed under arbitrary intersections and
finite unions,
$Z \xbe \xdy,$ $X,Y \xcc Z,$ then the following hold:

$(Cl \xcv )$ $ \wt{X \xcv Y}$ $=$ $ \wt{X} \xcv \wt{Y}$

$(Cl \xcs )$ $ \wt{X \xcs Y} \xcc \wt{X} \xcs \wt{Y},$ but usually not
conversely,

$ \xCf (Cl-)$ $ \wt{A}- \wt{B} \xcc \wt{A-B},$

$(Cl=)$ $X=Y$ $ \xcp $ $ \wt{X}= \wt{Y},$ but not conversely,

$(Cl \xcc 1)$ $ \wt{X} \xcc Y$ $ \xcp $ $X \xcc Y,$ but not conversely,

$(Cl \xcc 2)$ $X \xcc \wt{Y}$ $ \xcp $ $ \wt{X} \xcc \wt{Y}.$

(2) If, in addition, $X \xbe \xdy $ and $ \xdC X:=Z-X \xbe \xdy,$ then
the following two properties
hold, too:

$(Cl \xcs +)$ $ \wt{A} \xcs X= \wt{A \xcs X},$

$(Cl-+)$ $ \wt{A}-X= \wt{A-X}.$

(3) In the intended application, i.e. $ \wt{A}=M(Th(A)),$ the following
hold:

(3.1) $Th(X)$ $=$ $Th( \wt{X}),$

(3.2) Even if $A= \wt{A},$ $B= \wt{B},$ it is not necessarily true that $
\wt{A-B} \xcc \wt{A}- \wt{B}.$

 karl-search= End Fact Def-Clos
\vspace{7mm}

 *************************************

\vspace{7mm}

\subsubsection{Fact Def-Clos Proof}

 {\LARGE karl-search= Start Fact Def-Clos Proof }

\index{Fact Def-Clos Proof}

\efa

\subparagraph{
Proof
}

$\hspace{0.01em}$

(+++*** Orig.:  Proof )

$(Cl=),$ $(Cl \xcc 1),$ $(Cl \xcc 2),$ (3.1) are trivial.

$(Cl \xcv )$ Let $ \xdy (U):=\{X \xbe \xdy:U \xcc X\}.$ If $A \xbe \xdy
(X \xcv Y),$ then $A \xbe \xdy (X)$ and $A \xbe \xdy (Y),$ so
$ \wt{X \xcv Y}$ $ \xcd $ $ \wt{X} \xcv \wt{Y}.$
If $A \xbe \xdy (X)$ and $B \xbe \xdy (Y),$ then $A \xcv B \xbe \xdy (X
\xcv Y),$ so $ \wt{X \xcv Y}$ $ \xcc $ $ \wt{X} \xcv \wt{Y}.$

$(Cl \xcs )$ Let $X',Y' \xbe \xdy,$ $X \xcc X',$ $Y \xcc Y',$ then $X
\xcs Y \xcc X' \xcs Y',$ so $ \wt{X \xcs Y} \xcc \wt{X} \xcs \wt{Y}.$
For the converse, set $X:=M_{ \xdl }-\{m\},$ $Y:=\{m\}$ in Example \ref{Example
Not-Def} (page \pageref{Example Not-Def}) .

$ \xCf (Cl-)$ Let $A-B \xcc X \xbe \xdy,$ $B \xcc Y \xbe \xdy,$ so $A
\xcc X \xcv Y \xbe \xdy.$ Let $x \xce \wt{B}$ $ \xch $ $ \xcE Y \xbe \xdy
(B \xcc Y,$ $x \xce Y),$
$x \xce \wt{A-B}$ $ \xch $ $ \xcE X \xbe \xdy (A-B \xcc X,$ $x \xce X),$
so $x \xce X \xcv Y,$ $A \xcc X \xcv Y,$ so $x \xce \wt{A}.$ Thus, $x \xce
\wt{B},$ $x \xce \wt{A-B}$ $ \xch $
$x \xce \wt{A},$ or $x \xbe \wt{A}- \wt{B}$ $ \xch $ $x \xbe \wt{A-B}.$

$(Cl \xcs +)$ $ \wt{A} \xcs X \xcd \wt{A \xcs X}$ by $(Cl \xcs ).$
For `` $ \xcc $ '': Let $A \xcs X \xcc A' \xbe \xdy,$ then by closure
under $( \xcv ),$
$A \xcc A' \xcv \xdC X \xbe \xdy,$ $(A' \xcv \xdC X) \xcs X \xcc A'.$ So
$ \wt{A} \xcs X \xcc \wt{A \xcs X}.$

$(Cl-+)$ $ \wt{A-X}= \wt{A \xcs \xdC X}= \wt{A} \xcs \xdC X= \wt{A}-X$ by
$(Cl \xcs +).$

(3.2) Set $A:=M_{ \xdl },$ $B:=\{m\}$ for $m \xbe M_{ \xdl }$ arbitrary, $
\xdl $ infinite. So $A= \wt{A},$ $B= \wt{B},$ but
$ \wt{A-B}=A \xEd A-B.$

$ \xcz $
\\[3ex]

 karl-search= End Fact Def-Clos Proof
\vspace{7mm}

 *************************************

\vspace{7mm}

\subsubsection{Fact Mod-Int}

 {\LARGE karl-search= Start Fact Mod-Int }

\index{Fact Mod-Int}

\bfa

$\hspace{0.01em}$

(+++ Orig. No.:  Fact Mod-Int +++)

{\xssc LABEL: {Fact Mod-Int}}
\label{Fact Mod-Int}

Let $X,Y,Z \xbe M_{ \xdl }.$

(1) $X \xcc Y \xcs Z$ $ \xch $ $ \ol{Th(Y) \xcv Th(Z)} \xcc Th(X)$

(2) If $X=Y \xcs Z$ and $Y=M(T),$ $Z=M(T' ),$ then $ \ol{Th(Y) \xcv
Th(Z)}=Th(X).$

 karl-search= End Fact Mod-Int
\vspace{7mm}

 *************************************

\vspace{7mm}

\subsubsection{Fact Mod-Int Proof}

 {\LARGE karl-search= Start Fact Mod-Int Proof }

\index{Fact Mod-Int Proof}

\efa

\subparagraph{
Proof
}

$\hspace{0.01em}$

(+++*** Orig.:  Proof )

Let $ \wt{X}:=M(Th(X)).$

(1) $X \xcc Y \xcs Z$ $ \xch $ $ \wt{X} \xcc \wt{Y \xcs Z} \xcc \wt{Y}
\xcs \wt{Z}$ by Fact \ref{Fact Def-Clos} (page \pageref{Fact Def-Clos})  $(Cl
\xcs ).$ So
$M(Th(X))$ $ \xcc $ $M(Th(Y)) \xcs M(Th(Z))$ $=$ $M( \ol{Th(Y) \xcv
Th(Z)}),$ so $ \ol{Th(Y) \xcv Th(Z)}$ $ \xcc $ $ \ol{Th(X)}$ $=$ $Th(X).$

(2) $ \wt{Y \xcs Z}$ $=$ $ \wt{M(T) \xcs M(T' )}$ $=$ $ \wt{M(T \xcv T'
)}$ $=$ $M(T \xcv T' )$ $=$ $M(T) \xcs M(T' )$ $=$ $ \wt{M(T)} \xcs
\wt{M(T' )}$ $=$ $ \wt{Y} \xcs \wt{Z}.$
Finish as for (1).

$ \xcz $
\\[3ex]

 karl-search= End Fact Mod-Int Proof
\vspace{7mm}

 *************************************

\vspace{7mm}

 karl-search= End ToolBase1-Log-Dp
\vspace{7mm}

 *************************************

\vspace{7mm}

\subsection{
Logics: Rules
}

\subsubsection{ToolBase1-Log-Rules}

 {\LARGE karl-search= Start ToolBase1-Log-Rules }

{\xssc LABEL: {Section Toolbase1-Log-Rules}}
\label{Section Toolbase1-Log-Rules}
\index{Section Toolbase1-Log-Rules}
\subsubsection{Definition Log-Cond}

 {\LARGE karl-search= Start Definition Log-Cond }

\index{Definition Log-Cond}

\bd

$\hspace{0.01em}$

(+++ Orig. No.:  Definition Log-Cond +++)

{\xssc LABEL: {Definition Log-Cond}}
\label{Definition Log-Cond}

We introduce here formally a list of properties of set functions on the
algebraic side, and their corresponding logical rules on the other side.

Recall that $ \ol{T}:=\{ \xbf:T \xcl \xbf \},$ $ \ol{ \ol{T} }:=\{ \xbf
:T \xcn \xbf \},$
where $ \xcl $ is classical consequence, and $ \xcn $ any other
consequence.

We show, wherever adequate, in parallel the formula version
in the left column, the theory version
in the middle column, and the semantical or algebraic
counterpart in the
right column. The algebraic counterpart gives conditions for a
function $f:\xdy\xcp\xdp (U)$, where $U$ is some set, and
$\xdy\xcc\xdp (U)$.

\ed

Precise connections between the columns are given in
Proposition \ref{Proposition Alg-Log} (page \pageref{Proposition Alg-Log})

When the formula version is not commonly used, we omit it,
as we normally work only with the theory version.

$A$ and $B$ in the right hand side column stand for
$M(\xbf)$ for some formula $\xbf$, whereas $X$, $Y$ stand for
$M(T)$ for some theory $T$.

{\footnotesize

\begin{tabular}{|c|c|c|}

\hline

\multicolumn{3}{|c|}{Basics} \xEP

\hline

$(AND)$
\xEH
$(AND)$
\xEH
Closure under
\xEP

$ \xbf \xcn \xbq,  \xbf \xcn \xbq '   \xch $
\xEH
$ T \xcn \xbq, T \xcn \xbq '   \xch $
\xEH
finite
\xEP

$ \xbf \xcn \xbq \xcu \xbq ' $
\xEH
$ T \xcn \xbq \xcu \xbq ' $
\xEH
intersection
\xEP

\hline

$(OR)$ \xEH $(OR)$ \xEH $( \xbm OR)$ \xEP

$ \xbf \xcn \xbq,  \xbf ' \xcn \xbq   \xch $ \xEH
$ \ol{\ol{T}} \xcs \ol{\ol{T'}} \xcc \ol{\ol{T \xco T'}} $ \xEH
$f(X \xcv Y) \xcc f(X) \xcv f(Y)$
\xEP

$ \xbf \xco \xbf ' \xcn \xbq $ \xEH
\xEH
\xEP

\hline

$(wOR)$
\xEH
$(wOR)$
\xEH
$( \xbm wOR)$
\xEP

$ \xbf \xcn \xbq,$ $ \xbf ' \xcl \xbq $ $ \xch $
\xEH
$ \ol{ \ol{T} } \xcs \ol{T' }$ $ \xcc $ $ \ol{ \ol{T \xco T' } }$
\xEH
$f(X \xcv Y) \xcc f(X) \xcv Y$
\xEP

$ \xbf \xco \xbf ' \xcn \xbq $
\xEH
\xEH
\xEP

\hline

$(disjOR)$
\xEH
$(disjOR)$
\xEH
$( \xbm disjOR)$
\xEP

$ \xbf \xcl \xCN \xbf ',$ $ \xbf \xcn \xbq,$
\xEH
$\xCN Con(T \xcv T') \xch$
\xEH
$X \xcs Y= \xCQ $ $ \xch $
\xEP

$ \xbf ' \xcn \xbq $ $ \xch $ $ \xbf \xco \xbf ' \xcn \xbq $
\xEH
$ \ol{ \ol{T} } \xcs \ol{ \ol{T' } } \xcc \ol{ \ol{T \xco T' } }$
\xEH
$f(X \xcv Y) \xcc f(X) \xcv f(Y)$
\xEP

\hline

$(LLE)$
\xEH
$(LLE)$
\xEH
\xEP

Left Logical Equivalence
\xEH
\xEH
\xEP

$ \xcl \xbf \xcr \xbf ',  \xbf \xcn \xbq   \xch $
\xEH
$ \ol{T}= \ol{T' }  \xch   \ol{\ol{T}} = \ol{\ol{T'}}$
\xEH
trivially true
\xEP

$ \xbf ' \xcn \xbq $ \xEH \xEH \xEP

\hline

$(RW)$ Right Weakening
\xEH
$(RW)$
\xEH
upward closure
\xEP

$ \xbf \xcn \xbq,  \xcl \xbq \xcp \xbq '   \xch $
\xEH
$ T \xcn \xbq,  \xcl \xbq \xcp \xbq '   \xch $
\xEH
\xEP

$ \xbf \xcn \xbq ' $
\xEH
$T \xcn \xbq ' $
\xEH
\xEP

\hline

$(CCL)$ Classical Closure \xEH $(CCL)$ \xEH \xEP

\xEH
$ \ol{ \ol{T} }$ is classically
\xEH
trivially true
\xEP

\xEH closed \xEH \xEP

\hline

$(SC)$ Supraclassicality \xEH $(SC)$ \xEH $( \xbm \xcc )$ \xEP

$ \xbf \xcl \xbq $ $ \xch $ $ \xbf \xcn \xbq $ \xEH $ \ol{T} \xcc \ol{
\ol{T} }$ \xEH $f(X) \xcc X$ \xEP

\cline{1-1}

$(REF)$ Reflexivity \xEH \xEH \xEP
$ T \xcv \{\xba\} \xcn \xba $ \xEH \xEH \xEP

\hline

$(CP)$ \xEH $(CP)$ \xEH $( \xbm \xCQ )$ \xEP

Consistency Preservation \xEH \xEH \xEP

$ \xbf \xcn \xcT $ $ \xch $ $ \xbf \xcl \xcT $ \xEH $T \xcn \xcT $ $ \xch
$ $T \xcl \xcT $ \xEH $f(X)= \xCQ $ $ \xch $ $X= \xCQ $ \xEP

\hline

\xEH
\xEH $( \xbm \xCQ fin)$
\xEP

\xEH
\xEH $X \xEd \xCQ $ $ \xch $ $f(X) \xEd \xCQ $
\xEP

\xEH
\xEH for finite $X$
\xEP

\hline

\xEH $(PR)$ \xEH $( \xbm PR)$ \xEP

$ \ol{ \ol{ \xbf \xcu \xbf ' } }$ $ \xcc $ $ \ol{ \ol{ \ol{ \xbf } } \xcv
\{ \xbf ' \}}$ \xEH
$ \ol{ \ol{T \xcv T' } }$ $ \xcc $ $ \ol{ \ol{ \ol{T} } \xcv T' }$ \xEH
$X \xcc Y$ $ \xch $
\xEP

\xEH \xEH $f(Y) \xcs X \xcc f(X)$
\xEP

\cline{3-3}

\xEH
\xEH
$(\xbm PR ')$
\xEP

\xEH
\xEH
$f(X) \xcs Y \xcc f(X \xcs Y)$
\xEP

\hline

$(CUT)$ \xEH $(CUT)$ \xEH $ (\xbm CUT) $ \xEP
$ T  \xcn \xba; T \xcv \{ \xba\} \xcn \xbb \xch $ \xEH
$T \xcc \ol{T' } \xcc \ol{ \ol{T} }  \xch $ \xEH
$f(X) \xcc Y \xcc X  \xch $ \xEP
$ T  \xcn \xbb $ \xEH
$ \ol{ \ol{T'} } \xcc \ol{ \ol{T} }$ \xEH
$f(X) \xcc f(Y)$
\xEP

\hline

\end{tabular}

}

{\footnotesize

\begin{tabular}{|c|c|c|}

\hline

\multicolumn{3}{|c|}{Cumulativity} \xEP

\hline

$(CM)$ Cautious Monotony \xEH $(CM)$ \xEH $ (\xbm CM) $ \xEP

$ \xbf \xcn \xbq,  \xbf \xcn \xbq '   \xch $ \xEH
$T \xcc \ol{T' } \xcc \ol{ \ol{T} }  \xch $ \xEH
$f(X) \xcc Y \xcc X  \xch $
\xEP

$ \xbf \xcu \xbq \xcn \xbq ' $ \xEH
$ \ol{ \ol{T} } \xcc \ol{ \ol{T' } }$ \xEH
$f(Y) \xcc f(X)$
\xEP

\cline{1-1}

\cline{3-3}

or $(ResM)$ Restricted Monotony \xEH \xEH $(\xbm ResM)$ \xEP
$ T  \xcn \xba, \xbb \xch T \xcv \{\xba\} \xcn \xbb $ \xEH \xEH
$ f(X) \xcc A \xcs B \xch f(X \xcs A) \xcc B $ \xEP

\hline

$(CUM)$ Cumulativity \xEH $(CUM)$ \xEH $( \xbm CUM)$ \xEP

$ \xbf \xcn \xbq   \xch $ \xEH
$T \xcc \ol{T' } \xcc \ol{ \ol{T} }  \xch $ \xEH
$f(X) \xcc Y \xcc X  \xch $
\xEP

$( \xbf \xcn \xbq '   \xcj   \xbf \xcu \xbq \xcn \xbq ' )$ \xEH
$ \ol{ \ol{T} }= \ol{ \ol{T' } }$ \xEH
$f(Y)=f(X)$ \xEP

\hline

\xEH
$ (\xcc \xcd) $
\xEH
$ (\xbm \xcc \xcd) $
\xEP
\xEH
$T \xcc \ol{\ol{T'}}, T' \xcc \ol{\ol{T}} \xch $
\xEH
$ f(X) \xcc Y, f(Y) \xcc X \xch $
\xEP
\xEH
$ \ol{\ol{T'}} = \ol{\ol{T}}$
\xEH
$ f(X)=f(Y) $
\xEP

\hline

\multicolumn{3}{|c|}{Rationality} \xEP

\hline

$(RatM)$ Rational Monotony \xEH $(RatM)$ \xEH $( \xbm RatM)$ \xEP

$ \xbf \xcn \xbq,  \xbf \xcN \xCN \xbq '   \xch $ \xEH
$Con(T \xcv \ol{\ol{T'}})$, $T \xcl T'$ $ \xch $ \xEH
$X \xcc Y, X \xcs f(Y) \xEd \xCQ   \xch $
\xEP

$ \xbf \xcu \xbq ' \xcn \xbq $ \xEH
$ \ol{\ol{T}} \xcd \ol{\ol{\ol{T'}} \xcv T} $ \xEH
$f(X) \xcc f(Y) \xcs X$ \xEP

\hline

\xEH $(RatM=)$ \xEH $( \xbm =)$ \xEP

\xEH
$Con(T \xcv \ol{\ol{T'}})$, $T \xcl T'$ $ \xch $ \xEH
$X \xcc Y, X \xcs f(Y) \xEd \xCQ   \xch $
\xEP

\xEH
$ \ol{\ol{T}} = \ol{\ol{\ol{T'}} \xcv T} $ \xEH
$f(X) = f(Y) \xcs X$ \xEP

\hline

\xEH
$(Log=' )$
\xEH $( \xbm =' )$
\xEP

\xEH
$Con( \ol{ \ol{T' } } \xcv T)$ $ \xch $
\xEH $f(Y) \xcs X \xEd \xCQ $ $ \xch $
\xEP

\xEH
$ \ol{ \ol{T \xcv T' } }= \ol{ \ol{ \ol{T' } } \xcv T}$
\xEH $f(Y \xcs X)=f(Y) \xcs X$
\xEP

\hline

\xEH
$(Log \xFO )$
\xEH $( \xbm \xFO )$
\xEP

\xEH
$ \ol{ \ol{T \xco T' } }$ is one of
\xEH $f(X \xcv Y)$ is one of
\xEP

\xEH
$\ol{\ol{T}},$ or $\ol{\ol{T'}},$ or $\ol{\ol{T}} \xcs \ol{\ol{T'}}$ (by (CCL))
\xEH $f(X),$ $f(Y)$ or $f(X) \xcv f(Y)$
\xEP

\hline

\xEH
$(Log \xcv )$
\xEH $( \xbm \xcv )$
\xEP

\xEH
$Con( \ol{ \ol{T' } } \xcv T),$ $ \xCN Con( \ol{ \ol{T' } }
\xcv \ol{ \ol{T} })$ $ \xch $
\xEH $f(Y) \xcs (X-f(X)) \xEd \xCQ $ $ \xch $
\xEP

\xEH
$ \xCN Con( \ol{ \ol{T \xco T' } } \xcv T' )$
\xEH $f(X \xcv Y) \xcs Y= \xCQ$
\xEP

\hline

\xEH
$(Log \xcv ' )$
\xEH $( \xbm \xcv ' )$
\xEP

\xEH
$Con( \ol{ \ol{T' } } \xcv T),$ $ \xCN Con( \ol{ \ol{T' }
} \xcv \ol{ \ol{T} })$ $ \xch $
\xEH $f(Y) \xcs (X-f(X)) \xEd \xCQ $ $ \xch $
\xEP

\xEH
$ \ol{ \ol{T \xco T' } }= \ol{ \ol{T} }$
\xEH $f(X \xcv Y)=f(X)$
\xEP

\hline

\xEH
\xEH $( \xbm \xbe )$
\xEP

\xEH
\xEH $a \xbe X-f(X)$ $ \xch $
\xEP

\xEH
\xEH $ \xcE b \xbe X.a \xce f(\{a,b\})$
\xEP

\hline

\end{tabular}

}

$(PR)$ is also called infinite conditionalization - we choose the name for
its central role for preferential structures $(PR)$ or $( \xbm PR).$

The system of rules $(AND)$ $(OR)$ $(LLE)$ $(RW)$ $(SC)$ $(CP)$ $(CM)$ $(CUM)$
is also called system $P$ (for preferential), adding $(RatM)$ gives the system
$R$ (for rationality or rankedness).

Roughly: Smooth preferential structures generate logics satisfying system
$P$, ranked structures logics satisfying system $R$.

A logic satisfying $(REF)$, $(ResM)$, and $(CUT)$ is called a consequence
relation.

$(LLE)$ and$(CCL)$ will hold automatically, whenever we work with model sets.

$(AND)$ is obviously closely related to filters, and corresponds to closure
under finite intersections. $(RW)$ corresponds to upward closure of filters.

More precisely, validity of both depend on the definition, and the
direction we consider.

Given $f$ and $(\xbm \xcc )$, $f(X)\xcc X$ generates a pricipal filter:
$\{X'\xcc X:f(X)\xcc X'\}$, with
the definition: If $X=M(T)$, then $T\xcn \xbf$  iff $f(X)\xcc M(\xbf )$.
Validity of $(AND)$ and
$(RW)$ are then trivial.

Conversely, we can define for $X=M(T)$

$\xdx:=\{X'\xcc X: \xcE \xbf (X'=X\xcs M(\xbf )$ and $T\xcn \xbf )\}$.

$(AND)$ then makes $\xdx$  closed under
finite intersections, $(RW)$ makes $\xdx$  upward
closed. This is in the infinite case usually not yet a filter, as not all
subsets of $X$ need to be definable this way.
In this case, we complete $\xdx$  by
adding all $X''$ such that there is $X'\xcc X''\xcc X$, $X'\xbe\xdx$.

Alternatively, we can define

$\xdx:=\{X'\xcc X: \xcS\{X \xcs M(\xbf ): T\xcn \xbf \} \xcc X' \}$.

$(SC)$ corresponds to the choice of a subset.

$(CP)$ is somewhat delicate, as it presupposes that the chosen model set is
non-empty. This might fail in the presence of ever better choices, without
ideal ones; the problem is addressed by the limit versions.

$(PR)$ is an infinitary version of one half of the deduction theorem: Let $T$
stand for $\xbf$, $T'$ for $\xbq$, and $\xbf \xcu \xbq \xcn \xbs$,
so $\xbf \xcn \xbq \xcp \xbs$, but $(\xbq \xcp \xbs )\xcu \xbq \xcl \xbs$.

$(CUM)$ (whose more interesting half in our context is $(CM)$) may best be seen
as normal use of lemmas: We have worked hard and found some lemmas. Now
we can take a rest, and come back again with our new lemmas. Adding them to the
axioms will neither add new theorems, nor prevent old ones to hold.

 karl-search= End Definition Log-Cond
\vspace{7mm}

 *************************************

\vspace{7mm}

\subsubsection{Definition Log-Cond-Ref-Size}

 {\LARGE karl-search= Start Definition Log-Cond-Ref-Size }

\index{Definition Log-Cond-Ref-Size}
{\xssc LABEL: {Definition Log-Cond-Ref-Size}}
\label{Definition Log-Cond-Ref-Size}

The numbers refer to Proposition \ref{Proposition Alg-Log} (page
\pageref{Proposition Alg-Log}) .

\begin{turn}{90}


{\xssc

\begin{tabular}{|c|c|c|c|c|c|c|}

\hline

\multicolumn{2}{|c|}{Logical rule}
\xEH
Correspondence
\xEH
Model set
\xEH
Correspondence
\xEH
Size
\xEH
Size-Rule
\xEP

\hline

\multicolumn{7}{|c|}{Basics}
\xEP

\hline

$(SC)$ Supraclassicality
\xEH
$(SC)$
\xEH
$\xch$ (4.1)
\xEH
$( \xbm \xcc )$
\xEH
\xEH
\xEH
$(Opt)$
\xEP

\cline{3-3}

$ \xbf \xcl \xbq $ $ \xch $ $ \xbf \xcn \xbq $
\xEH
$ \ol{T} \xcc \ol{ \ol{T} }$
\xEH
$\xci$ (4.2)
\xEH
$f(X) \xcc X$
\xEH
\xEH
\xEH
\xEP

\cline{1-1}

$(REF)$ Reflexivity
\xEH
\xEH
\xEH
\xEH
\xEH
\xEH
\xEP

$ T \xcv \{\xba\} \xcn \xba $
\xEH
\xEH
\xEH
\xEH
\xEH
\xEH
\xEP

\hline

$(LLE)$
\xEH
$(LLE)$
\xEH
\xEH
\xEH
\xEH
\xEH
\xEP

Left Logical Equivalence
\xEH
\xEH
\xEH
\xEH
\xEH
\xEH
\xEP

$ \xcl \xbf \xcr \xbf ',  \xbf \xcn \xbq   \xch $
\xEH
$ \ol{T}= \ol{T' }  \xch   \ol{\ol{T}} = \ol{\ol{T'}}$
\xEH
\xEH
\xEH
\xEH
\xEH
\xEP

$ \xbf ' \xcn \xbq $
\xEH
\xEH
\xEH
\xEH
\xEH
\xEH
\xEP

\hline

$(RW)$ Right Weakening
\xEH
$(RW)$
\xEH
\xEH
\xEH
\xEH
+
\xEH
$(iM)$
\xEP

$ \xbf \xcn \xbq,  \xcl \xbq \xcp \xbq '   \xch $
\xEH
$ T \xcn \xbq,  \xcl \xbq \xcp \xbq '   \xch $
\xEH
\xEH
\xEH
\xEH
\xEH
\xEP

$ \xbf \xcn \xbq ' $
\xEH
$T \xcn \xbq ' $
\xEH
\xEH
\xEH
\xEH
\xEH
\xEP

\hline

$(wOR)$
\xEH
$(wOR)$
\xEH
$\xch$ (3.1)
\xEH
$( \xbm wOR)$
\xEH
\xEH
+
\xEH
$(eM\xdi)$
\xEP

\cline{3-3}

$ \xbf \xcn \xbq,$ $ \xbf ' \xcl \xbq $ $ \xch $
\xEH
$ \ol{ \ol{T} } \xcs \ol{T' }$ $ \xcc $ $ \ol{ \ol{T \xco T' } }$
\xEH
$\xci$ (3.2)
\xEH
$f(X \xcv Y) \xcc f(X) \xcv Y$
\xEH
\xEH
\xEH
\xEP

$ \xbf \xco \xbf ' \xcn \xbq $
\xEH
\xEH
\xEH
\xEH
\xEH
\xEH
\xEP

\hline

$(disjOR)$
\xEH
$(disjOR)$
\xEH
$\xch$ (2.1)
\xEH
$( \xbm disjOR)$
\xEH
\xEH
$\xCd$
\xEH
$(I\xcv disj)$
\xEP

\cline{3-3}

$ \xbf \xcl \xCN \xbf ',$ $ \xbf \xcn \xbq,$
\xEH
$\xCN Con(T \xcv T') \xch$
\xEH
$\xci$ (2.2)
\xEH
$X \xcs Y= \xCQ $ $ \xch $
\xEH
\xEH
\xEH
\xEP

$ \xbf ' \xcn \xbq $ $ \xch $ $ \xbf \xco \xbf ' \xcn \xbq $
\xEH
$ \ol{ \ol{T} } \xcs \ol{ \ol{T' } } \xcc \ol{ \ol{T \xco T' } }$
\xEH
\xEH
$f(X \xcv Y) \xcc f(X) \xcv f(Y)$
\xEH
\xEH
\xEH
\xEP

\hline

$(CP)$
\xEH
$(CP)$
\xEH
$\xch$ (5.1)
\xEH
$( \xbm \xCQ )$
\xEH
\xEH
$1*s$
\xEH
$(I_1)$
\xEP

\cline{3-3}

Consistency Preservation
\xEH
\xEH
$\xci$ (5.2)
\xEH
\xEH
\xEH
\xEH
\xEP

$ \xbf \xcn \xcT $ $ \xch $ $ \xbf \xcl \xcT $
\xEH
$T \xcn \xcT $ $ \xch $ $T \xcl \xcT $
\xEH
\xEH
$f(X)= \xCQ $ $ \xch $ $X= \xCQ $
\xEH
\xEH
\xEH
\xEP

\hline

\xEH
\xEH
\xEH
$( \xbm \xCQ fin)$
\xEH
\xEH
$1*s$
\xEH
$(I_1)$
\xEP

\xEH
\xEH
\xEH
$X \xEd \xCQ $ $ \xch $ $f(X) \xEd \xCQ $
\xEH
\xEH
\xEH
\xEP

\xEH
\xEH
\xEH
for finite $X$
\xEH
\xEH
\xEH
\xEP

\hline

\xEH
$(AND_1)$
\xEH
\xEH
\xEH
\xEH
$2*s$
\xEH
$(I_2)$
\xEP

\xEH
$\xba\xcn\xbb \xch \xba\xcN\xCN\xbb$
\xEH
\xEH
\xEH
\xEH
\xEH
\xEP

\hline

\xEH
$(AND_n)$
\xEH
\xEH
\xEH
\xEH
$n*s$
\xEH
$(I_n)$
\xEP

\xEH
$\xba\xcn\xbb_1, \ldots, \xba\xcn\xbb_{n-1} \xch $
\xEH
\xEH
\xEH
\xEH
\xEH
\xEP

\xEH
$\xba\xcN(\xCN\xbb_1 \xco \ldots \xco \xCN\xbb_{n-1})$
\xEH
\xEH
\xEH
\xEH
\xEH
\xEP

\hline

$(AND)$
\xEH
$(AND)$
\xEH
\xEH
\xEH
\xEH
$\xbo *s$
\xEH
$(I_\xbo)$
\xEP

$ \xbf \xcn \xbq,  \xbf \xcn \xbq '   \xch $
\xEH
$ T \xcn \xbq, T \xcn \xbq '   \xch $
\xEH
\xEH
\xEH
\xEH
\xEH
\xEP

$ \xbf \xcn \xbq \xcu \xbq ' $
\xEH
$ T \xcn \xbq \xcu \xbq ' $
\xEH
\xEH
\xEH
\xEH
\xEH
\xEP

\hline

$(CCL)$ Classical Closure
\xEH
$(CCL)$
\xEH
\xEH
\xEH
\xEH
$\xbo *s$
\xEH
$(iM)+(I_\xbo)$
\xEP

\xEH
$ \ol{ \ol{T} }$ classically closed
\xEH
\xEH
\xEH
\xEH
\xEH
\xEP

\hline

$(OR)$
\xEH
$(OR)$
\xEH
$\xch$ (1.1)
\xEH
$( \xbm OR)$
\xEH
\xEH
$\xbo *s$
\xEH
$(eM\xdi)+(I_\xbo)$
\xEP

\cline{3-3}

$ \xbf \xcn \xbq,  \xbf ' \xcn \xbq   \xch $
\xEH
$ \ol{\ol{T}} \xcs \ol{\ol{T'}} \xcc \ol{\ol{T \xco T'}} $
\xEH
$\xci$ (1.2)
\xEH
$f(X \xcv Y) \xcc f(X) \xcv f(Y)$
\xEH
\xEH
\xEH
\xEP

$ \xbf \xco \xbf ' \xcn \xbq $
\xEH
\xEH
\xEH
\xEH
\xEH
\xEH
\xEP

\hline

\xEH
$(PR)$
\xEH
$\xch$ (6.1)
\xEH
$( \xbm PR)$
\xEH
\xEH
$\xbo *s$
\xEH
$(eM\xdi)+(I_\xbo)$
\xEP

\cline{3-3}

$ \ol{ \ol{ \xbf \xcu \xbf ' } }$ $ \xcc $ $ \ol{ \ol{ \ol{ \xbf } } \xcv
\{ \xbf ' \}}$
\xEH
$ \ol{ \ol{T \xcv T' } }$ $ \xcc $ $ \ol{ \ol{ \ol{T} } \xcv T' }$
\xEH
$\xci (\xbm dp)+(\xbm\xcc)$ (6.2)
\xEH
$X \xcc Y$ $ \xch $
\xEH
\xEH
\xEH
\xEP

\cline{3-3}

\xEH
\xEH
$\xcI$ without $(\xbm dp)$ (6.3)
\xEH
$f(Y) \xcs X \xcc f(X)$
\xEH
\xEH
\xEH
\xEP

\cline{3-3}

\xEH
\xEH
$\xci (\xbm\xcc)$ (6.4)
\xEH
\xEH
\xEH
\xEH
\xEP

\xEH
\xEH
$T'$ a formula
\xEH
\xEH
\xEH
\xEH
\xEP

\cline{3-4}

\xEH
\xEH
$\xci$ (6.5)
\xEH
$(\xbm PR ')$
\xEH
\xEH
\xEH
\xEP

\xEH
\xEH
$T'$ a formula
\xEH
$f(X) \xcs Y \xcc f(X \xcs Y)$
\xEH
\xEH
\xEH
\xEP

\hline

$(CUT)$
\xEH
$(CUT)$
\xEH
$\xch$ (7.1)
\xEH
$ (\xbm CUT) $
\xEH
\xEH
$\xbo *s$
\xEH
$(eM\xdi)+(I_\xbo)$
\xEP

\cline{3-3}

$ T  \xcn \xba; T \xcv \{ \xba\} \xcn \xbb \xch $
\xEH
$T \xcc \ol{T' } \xcc \ol{ \ol{T} }  \xch $
\xEH
$\xci$ (7.2)
\xEH
$f(X) \xcc Y \xcc X  \xch $
\xEH
\xEH
\xEH
\xEP

$ T  \xcn \xbb $
\xEH
$ \ol{ \ol{T'} } \xcc \ol{ \ol{T} }$
\xEH
\xEH
$f(X) \xcc f(Y)$
\xEH
\xEH
\xEH
\xEP

\hline

\end{tabular}

}

\end{turn}

\begin{turn}{90}

{\xssc

\begin{tabular}{|c|c|c|c|c|c|c|}

\hline

\multicolumn{2}{|c|}{Logical rule}
\xEH
Correspondence
\xEH
Model set
\xEH
Correspondence
\xEH
Size
\xEH
Size-Rule
\xEP

\hline

\multicolumn{7}{|c|}{Cumulativity}
\xEP

\hline

$(wCM)$
\xEH
\xEH
\xEH
\xEH
\xEH
\xEH
$(eM\xdf)$
\xEP

$\xba\xcn\xbb, \xba\xcl\xbb' \xch \xba\xcu\xbb'\xcn\xbb$
\xEH
\xEH
\xEH
\xEH
\xEH
\xEH
\xEP

\hline

$(CM_2)$
\xEH
\xEH
\xEH
\xEH
\xEH
$2*s$
\xEH
$(I_2)$
\xEP

$\xba\xcn\xbb, \xba\xcn\xbb' \xch \xba\xcu\xbb\xcL\xCN\xbb'$
\xEH
\xEH
\xEH
\xEH
\xEH
\xEH
\xEP

\hline

$(CM_n)$
\xEH
\xEH
\xEH
\xEH
\xEH
$n*s$
\xEH
$(I_n)$
\xEP

$\xba\xcn\xbb_1, \ldots, \xba\xcn\xbb_n \xch $
\xEH
\xEH
\xEH
\xEH
\xEH
\xEH
\xEP

$\xba \xcu \xbb_1 \xcu \ldots \xcu \xbb_{n-1} \xcL\xCN\xbb_n$
\xEH
\xEH
\xEH
\xEH
\xEH
\xEH
\xEP

\hline

$(CM)$ Cautious Monotony
\xEH
$(CM)$
\xEH
$\xch$ (8.1)
\xEH
$ (\xbm CM) $
\xEH
\xEH
$\xbo *s$
\xEH
$(I_\xbo)$
\xEP

\cline{3-3}

$ \xbf \xcn \xbq,  \xbf \xcn \xbq '   \xch $
\xEH
$T \xcc \ol{T' } \xcc \ol{ \ol{T} }  \xch $
\xEH
$\xci$ (8.2)
\xEH
$f(X) \xcc Y \xcc X  \xch $
\xEH
\xEH
\xEH
\xEP

$ \xbf \xcu \xbq \xcn \xbq ' $
\xEH
$ \ol{ \ol{T} } \xcc \ol{ \ol{T' } }$
\xEH
\xEH
$f(Y) \xcc f(X)$
\xEH
\xEH
\xEH
\xEP

\cline{1-1}

\cline{3-4}

or $(ResM)$ Restricted Monotony
\xEH
\xEH
$\xch$ (9.1)
\xEH
$(\xbm ResM)$
\xEH
\xEH
\xEH
\xEP

\cline{3-3}

$ T  \xcn \xba, \xbb \xch T \xcv \{\xba\} \xcn \xbb $
\xEH
\xEH
$\xci$ (9.2)
\xEH
$ f(X) \xcc A \xcs B \xch f(X \xcs A) \xcc B $
\xEH
\xEH
\xEH
\xEP

\hline

$(CUM)$ Cumulativity
\xEH
$(CUM)$
\xEH
$\xch$ (11.1)
\xEH
$( \xbm CUM)$
\xEH
\xEH
\xEH
\xEP

\cline{3-3}

$ \xbf \xcn \xbq   \xch $
\xEH
$T \xcc \ol{T' } \xcc \ol{ \ol{T} }  \xch $
\xEH
$\xci$ (11.2)
\xEH
$f(X) \xcc Y \xcc X  \xch $
\xEH
\xEH
\xEH
\xEP

$( \xbf \xcn \xbq '   \xcj   \xbf \xcu \xbq \xcn \xbq ' )$
\xEH
$ \ol{ \ol{T} }= \ol{ \ol{T' } }$
\xEH
\xEH
$f(Y)=f(X)$
\xEH
\xEH
\xEH
\xEP

\hline

\xEH
$ (\xcc \xcd) $
\xEH
$\xch$ (10.1)
\xEH
$ (\xbm \xcc \xcd) $
\xEH
\xEH
\xEH
\xEP

\cline{3-3}

\xEH
$T \xcc \ol{\ol{T'}}, T' \xcc \ol{\ol{T}} \xch $
\xEH
$\xci$ (10.2)
\xEH
$ f(X) \xcc Y, f(Y) \xcc X \xch $
\xEH
\xEH
\xEH
\xEP

\xEH
$ \ol{\ol{T'}} = \ol{\ol{T}}$
\xEH
\xEH
$ f(X)=f(Y) $
\xEH
\xEH
\xEH
\xEP

\hline

\multicolumn{7}{|c|}{Rationality}
\xEP

\hline

$(RatM)$ Rational Monotony
\xEH
$(RatM)$
\xEH
$\xch$ (12.1)
\xEH
$( \xbm RatM)$
\xEH
\xEH
\xEH
$(\xdm^{++})$
\xEP

\cline{3-3}

$ \xbf \xcn \xbq,  \xbf \xcN \xCN \xbq '   \xch $
\xEH
$Con(T \xcv \ol{\ol{T'}})$, $T \xcl T'$ $ \xch $
\xEH
$\xci$ $(\xbm dp)$ (12.2)
\xEH
$X \xcc Y, X \xcs f(Y) \xEd \xCQ   \xch $
\xEH
\xEH
\xEH
\xEP

\cline{3-3}

$ \xbf \xcu \xbq ' \xcn \xbq $
\xEH
$ \ol{\ol{T}} \xcd \ol{\ol{\ol{T'}} \xcv T} $
\xEH
$\xcI$ without $(\xbm dp)$ (12.3)
\xEH
$f(X) \xcc f(Y) \xcs X$
\xEH
\xEH
\xEH
\xEP

\cline{3-3}

\xEH
\xEH
$\xci$ $T$ a formula (12.4)
\xEH
\xEH
\xEH
\xEH
\xEP

\hline

\xEH
$(RatM=)$
\xEH
$\xch$ (13.1)
\xEH
$( \xbm =)$
\xEH
\xEH
\xEH
\xEP

\cline{3-3}

\xEH
$Con(T \xcv \ol{\ol{T'}})$, $T \xcl T'$ $ \xch $
\xEH
$\xci$ $(\xbm dp)$ (13.2)
\xEH
$X \xcc Y, X \xcs f(Y) \xEd \xCQ   \xch $
\xEH
\xEH
\xEH
\xEP

\cline{3-3}

\xEH
$ \ol{\ol{T}} = \ol{\ol{\ol{T'}} \xcv T} $
\xEH
$\xcI$ without $(\xbm dp)$ (13.3)
\xEH
$f(X) = f(Y) \xcs X$
\xEH
\xEH
\xEH
\xEP

\cline{3-3}

\xEH
\xEH
$\xci$ $T$ a formula (13.4)
\xEH
\xEH
\xEH
\xEH
\xEP

\hline

\xEH
$(Log=' )$
\xEH
$\xch$ (14.1)
\xEH
$( \xbm =' )$
\xEH
\xEH
\xEH
\xEP

\cline{3-3}

\xEH
$Con( \ol{ \ol{T' } } \xcv T)$ $ \xch $
\xEH
$\xci$ $(\xbm dp)$ (14.2)
\xEH
$f(Y) \xcs X \xEd \xCQ $ $ \xch $
\xEH
\xEH
\xEH
\xEP

\cline{3-3}

\xEH
$ \ol{ \ol{T \xcv T' } }= \ol{ \ol{ \ol{T' } } \xcv T}$
\xEH
$\xcI$ without $(\xbm dp)$ (14.3)
\xEH
$f(Y \xcs X)=f(Y) \xcs X$
\xEH
\xEH
\xEH
\xEP

\cline{3-3}

\xEH
\xEH
$\xci$ $T$ a formula (14.4)
\xEH
\xEH
\xEH
\xEH
\xEP

\hline

\xEH
$(Log \xFO )$
\xEH
$\xch$ (15.1)
\xEH
$( \xbm \xFO )$
\xEH
\xEH
\xEH
\xEP

\cline{3-3}

\xEH
$ \ol{ \ol{T \xco T' } }$ is one of
\xEH
$\xci$ (15.2)
\xEH
$f(X \xcv Y)$ is one of
\xEH
\xEH
\xEH
\xEP

\xEH
$\ol{\ol{T}},$ or $\ol{\ol{T'}},$ or $\ol{\ol{T}} \xcs \ol{\ol{T'}}$ (by (CCL))
\xEH
\xEH
$f(X),$ $f(Y)$ or $f(X) \xcv f(Y)$
\xEH
\xEH
\xEH
\xEP

\hline

\xEH
$(Log \xcv )$
\xEH
$\xch$ $(\xbm\xcc)+(\xbm=)$ (16.1)
\xEH
$( \xbm \xcv )$
\xEH
\xEH
\xEH
\xEP

\cline{3-3}

\xEH
$Con( \ol{ \ol{T' } } \xcv T),$ $ \xCN Con( \ol{ \ol{T' } }
\xcv \ol{ \ol{T} })$ $ \xch $
\xEH
$\xci$ $(\xbm dp)$ (16.2)
\xEH
$f(Y) \xcs (X-f(X)) \xEd \xCQ $ $ \xch $
\xEH
\xEH
\xEH
\xEP

\cline{3-3}

\xEH
$ \xCN Con( \ol{ \ol{T \xco T' } } \xcv T' )$
\xEH
$\xcI$ without $(\xbm dp)$ (16.3)
\xEH
$f(X \xcv Y) \xcs Y= \xCQ$
\xEH
\xEH
\xEH
\xEP

\hline

\xEH
$(Log \xcv ' )$
\xEH
$\xch$ $(\xbm\xcc)+(\xbm=)$ (17.1)
\xEH
$( \xbm \xcv ' )$
\xEH
\xEH
\xEH
\xEP

\cline{3-3}

\xEH
$Con( \ol{ \ol{T' } } \xcv T),$ $ \xCN Con( \ol{ \ol{T' }
} \xcv \ol{ \ol{T} })$ $ \xch $
\xEH
$\xci$ $(\xbm dp)$ (17.2)
\xEH
$f(Y) \xcs (X-f(X)) \xEd \xCQ $ $ \xch $
\xEH
\xEH
\xEH
\xEP

\cline{3-3}

\xEH
$ \ol{ \ol{T \xco T' } }= \ol{ \ol{T} }$
\xEH
$\xcI$ without $(\xbm dp)$ (17.3)
\xEH
$f(X \xcv Y)=f(X)$
\xEH
\xEH
\xEH
\xEP

\hline

\xEH
\xEH
\xEH
$( \xbm \xbe )$
\xEH
\xEH
\xEH
\xEP

\xEH
\xEH
\xEH
$a \xbe X-f(X)$ $ \xch $
\xEH
\xEH
\xEH
\xEP

\xEH
\xEH
\xEH
$ \xcE b \xbe X.a \xce f(\{a,b\})$
\xEH
\xEH
\xEH
\xEP

\hline

\end{tabular}

}

\end{turn}

 karl-search= End Definition Log-Cond-Ref-Size
\vspace{7mm}

 *************************************

\vspace{7mm}

\subsubsection{Definition Log-Cond-Ref}

 {\LARGE karl-search= Start Definition Log-Cond-Ref }

\index{Definition Log-Cond-Ref}
{\xssc LABEL: {Definition Log-Cond-Ref}}
\label{Definition Log-Cond-Ref}

The numbers refer to Proposition \ref{Proposition Alg-Log} (page
\pageref{Proposition Alg-Log}) .

{\footnotesize

\begin{tabular}{|c|c|c|c|}

\hline

\multicolumn{4}{|c|}{Basics}
\xEP

\hline

$(AND)$
\xEH
$(AND)$
\xEH
\xEH
Closure under
\xEP

$ \xbf \xcn \xbq,  \xbf \xcn \xbq '   \xch $
\xEH
$ T \xcn \xbq, T \xcn \xbq '   \xch $
\xEH
\xEH
finite
\xEP

$ \xbf \xcn \xbq \xcu \xbq ' $
\xEH
$ T \xcn \xbq \xcu \xbq ' $
\xEH
\xEH
intersection
\xEP

\hline

$(OR)$
\xEH
$(OR)$
\xEH
$\xch$ (1.1)
\xEH
$( \xbm OR)$
\xEP

\cline{3-3}

$ \xbf \xcn \xbq,  \xbf ' \xcn \xbq   \xch $
\xEH
$ \ol{\ol{T}} \xcs \ol{\ol{T'}} \xcc \ol{\ol{T \xco T'}} $
\xEH
$\xci$ (1.2)
\xEH
$f(X \xcv Y) \xcc f(X) \xcv f(Y)$
\xEP

$ \xbf \xco \xbf ' \xcn \xbq $
\xEH
\xEH
\xEH
\xEP

\hline

$(wOR)$
\xEH
$(wOR)$
\xEH
$\xch$ (3.1)
\xEH
$( \xbm wOR)$
\xEP

\cline{3-3}

$ \xbf \xcn \xbq,$ $ \xbf ' \xcl \xbq $ $ \xch $
\xEH
$ \ol{ \ol{T} } \xcs \ol{T' }$ $ \xcc $ $ \ol{ \ol{T \xco T' } }$
\xEH
$\xci$ (3.2)
\xEH
$f(X \xcv Y) \xcc f(X) \xcv Y$
\xEP

$ \xbf \xco \xbf ' \xcn \xbq $
\xEH
\xEH
\xEH
\xEP

\hline

$(disjOR)$
\xEH
$(disjOR)$
\xEH
$\xch$ (2.1)
\xEH
$( \xbm disjOR)$
\xEP

\cline{3-3}

$ \xbf \xcl \xCN \xbf ',$ $ \xbf \xcn \xbq,$
\xEH
$\xCN Con(T \xcv T') \xch$
\xEH
$\xci$ (2.2)
\xEH
$X \xcs Y= \xCQ $ $ \xch $
\xEP

$ \xbf ' \xcn \xbq $ $ \xch $ $ \xbf \xco \xbf ' \xcn \xbq $
\xEH
$ \ol{ \ol{T} } \xcs \ol{ \ol{T' } } \xcc \ol{ \ol{T \xco T' } }$
\xEH
\xEH
$f(X \xcv Y) \xcc f(X) \xcv f(Y)$
\xEP

\hline

$(LLE)$
\xEH
$(LLE)$
\xEH
\xEH
\xEP

Left Logical Equivalence
\xEH
\xEH
\xEH
\xEP

$ \xcl \xbf \xcr \xbf ',  \xbf \xcn \xbq   \xch $
\xEH
$ \ol{T}= \ol{T' }  \xch   \ol{\ol{T}} = \ol{\ol{T'}}$
\xEH
\xEH
trivially true
\xEP

$ \xbf ' \xcn \xbq $
\xEH
\xEH
\xEH
\xEP

\hline

$(RW)$ Right Weakening
\xEH
$(RW)$
\xEH
\xEH
upward closure
\xEP

$ \xbf \xcn \xbq,  \xcl \xbq \xcp \xbq '   \xch $
\xEH
$ T \xcn \xbq,  \xcl \xbq \xcp \xbq '   \xch $
\xEH
\xEH
\xEP

$ \xbf \xcn \xbq ' $
\xEH
$T \xcn \xbq ' $
\xEH
\xEH
\xEP

\hline

$(CCL)$ Classical Closure
\xEH
$(CCL)$
\xEH
\xEH
\xEP

\xEH
$ \ol{ \ol{T} }$ is classically
\xEH
\xEH
trivially true
\xEP

\xEH
closed
\xEH
\xEH
\xEP

\hline

$(SC)$ Supraclassicality
\xEH
$(SC)$
\xEH
$\xch$ (4.1)
\xEH
$( \xbm \xcc )$
\xEP

\cline{3-3}

$ \xbf \xcl \xbq $ $ \xch $ $ \xbf \xcn \xbq $
\xEH
$ \ol{T} \xcc \ol{ \ol{T} }$
\xEH
$\xci$ (4.2)
\xEH
$f(X) \xcc X$
\xEP

\cline{1-1}

$(REF)$ Reflexivity
\xEH
\xEH
\xEH
\xEP

$ T \xcv \{\xba\} \xcn \xba $
\xEH
\xEH
\xEH
\xEP

\hline

$(CP)$
\xEH
$(CP)$
\xEH
$\xch$ (5.1)
\xEH
$( \xbm \xCQ )$
\xEP

\cline{3-3}

Consistency Preservation
\xEH
\xEH
$\xci$ (5.2)
\xEH
\xEP

$ \xbf \xcn \xcT $ $ \xch $ $ \xbf \xcl \xcT $
\xEH
$T \xcn \xcT $ $ \xch $ $T \xcl \xcT $
\xEH
\xEH
$f(X)= \xCQ $ $ \xch $ $X= \xCQ $
\xEP

\hline

\xEH
\xEH
\xEH
$( \xbm \xCQ fin)$
\xEP

\xEH
\xEH
\xEH
$X \xEd \xCQ $ $ \xch $ $f(X) \xEd \xCQ $
\xEP

\xEH
\xEH
\xEH
for finite $X$
\xEP

\hline

\xEH
$(PR)$
\xEH
$\xch$ (6.1)
\xEH
$( \xbm PR)$
\xEP

\cline{3-3}

$ \ol{ \ol{ \xbf \xcu \xbf ' } }$ $ \xcc $ $ \ol{ \ol{ \ol{ \xbf } } \xcv
\{ \xbf ' \}}$
\xEH
$ \ol{ \ol{T \xcv T' } }$ $ \xcc $ $ \ol{ \ol{ \ol{T} } \xcv T' }$
\xEH
$\xci (\xbm dp)+(\xbm\xcc)$ (6.2)
\xEH
$X \xcc Y$ $ \xch $
\xEP

\cline{3-3}

\xEH
\xEH
$\xcI$ without $(\xbm dp)$ (6.3)
\xEH
$f(Y) \xcs X \xcc f(X)$
\xEP

\cline{3-3}

\xEH
\xEH
$\xci (\xbm\xcc)$ (6.4)
\xEH
\xEP

\xEH
\xEH
$T'$ a formula
\xEH
\xEP

\cline{3-4}

\xEH
\xEH
$\xci $ (6.5)
\xEH
$(\xbm PR ')$
\xEP

\xEH
\xEH
$T'$ a formula
\xEH
$f(X) \xcs Y \xcc f(X \xcs Y)$
\xEP

\hline

$(CUT)$
\xEH
$(CUT)$
\xEH
$\xch$ (7.1)
\xEH
$ (\xbm CUT) $
\xEP

\cline{3-3}

$ T  \xcn \xba; T \xcv \{ \xba\} \xcn \xbb \xch $
\xEH
$T \xcc \ol{T' } \xcc \ol{ \ol{T} }  \xch $
\xEH
$\xci$ (7.2)
\xEH
$f(X) \xcc Y \xcc X  \xch $
\xEP

$ T  \xcn \xbb $
\xEH
$ \ol{ \ol{T'} } \xcc \ol{ \ol{T} }$
\xEH
\xEH
$f(X) \xcc f(Y)$
\xEP

\hline

\end{tabular}

}

{\footnotesize

\begin{tabular}{|c|c|c|c|}

\hline

\multicolumn{4}{|c|}{Cumulativity} \xEP

\hline

$(CM)$ Cautious Monotony
\xEH
$(CM)$
\xEH
$\xch$ (8.1)
\xEH
$ (\xbm CM) $
\xEP

\cline{3-3}

$ \xbf \xcn \xbq,  \xbf \xcn \xbq '   \xch $
\xEH
$T \xcc \ol{T' } \xcc \ol{ \ol{T} }  \xch $
\xEH
$\xci$ (8.2)
\xEH
$f(X) \xcc Y \xcc X  \xch $
\xEP

$ \xbf \xcu \xbq \xcn \xbq ' $
\xEH
$ \ol{ \ol{T} } \xcc \ol{ \ol{T' } }$
\xEH
\xEH
$f(Y) \xcc f(X)$
\xEP

\cline{1-1}

\cline{3-4}

or $(ResM)$ Restricted Monotony
\xEH
\xEH
$\xch$ (9.1)
\xEH
$(\xbm ResM)$
\xEP

\cline{3-3}

$ T  \xcn \xba, \xbb \xch T \xcv \{\xba\} \xcn \xbb $
\xEH
\xEH
$\xci$ (9.2)
\xEH
$ f(X) \xcc A \xcs B \xch f(X \xcs A) \xcc B $
\xEP

\hline

$(CUM)$ Cumulativity
\xEH
$(CUM)$
\xEH
$\xch$ (11.1)
\xEH
$( \xbm CUM)$
\xEP

\cline{3-3}

$ \xbf \xcn \xbq   \xch $
\xEH
$T \xcc \ol{T' } \xcc \ol{ \ol{T} }  \xch $
\xEH
$\xci$ (11.2)
\xEH
$f(X) \xcc Y \xcc X  \xch $
\xEP

$( \xbf \xcn \xbq '   \xcj   \xbf \xcu \xbq \xcn \xbq ' )$
\xEH
$ \ol{ \ol{T} }= \ol{ \ol{T' } }$
\xEH
\xEH
$f(Y)=f(X)$
\xEP

\hline

\xEH
$ (\xcc \xcd) $
\xEH
$\xch$ (10.1)
\xEH
$ (\xbm \xcc \xcd) $
\xEP

\cline{3-3}

\xEH
$T \xcc \ol{\ol{T'}}, T' \xcc \ol{\ol{T}} \xch $
\xEH
$\xci$ (10.2)
\xEH
$ f(X) \xcc Y, f(Y) \xcc X \xch $
\xEP

\xEH
$ \ol{\ol{T'}} = \ol{\ol{T}}$
\xEH
\xEH
$ f(X)=f(Y) $
\xEP

\hline

\multicolumn{4}{|c|}{Rationality} \xEP

\hline

$(RatM)$ Rational Monotony
\xEH
$(RatM)$
\xEH
$\xch$ (12.1)
\xEH
$( \xbm RatM)$
\xEP

\cline{3-3}

$ \xbf \xcn \xbq,  \xbf \xcN \xCN \xbq '   \xch $
\xEH
$Con(T \xcv \ol{\ol{T'}})$, $T \xcl T'$ $ \xch $
\xEH
$\xci$ $(\xbm dp)$ (12.2)
\xEH
$X \xcc Y, X \xcs f(Y) \xEd \xCQ   \xch $
\xEP

\cline{3-3}

$ \xbf \xcu \xbq ' \xcn \xbq $
\xEH
$ \ol{\ol{T}} \xcd \ol{\ol{\ol{T'}} \xcv T} $
\xEH
$\xcI$ without $(\xbm dp)$ (12.3)
\xEH
$f(X) \xcc f(Y) \xcs X$
\xEP

\cline{3-3}

\xEH
\xEH
$\xci$ $T$ a formula (12.4)
\xEH
\xEP

\hline

\xEH
$(RatM=)$
\xEH
$\xch$ (13.1)
\xEH
$( \xbm =)$
\xEP

\cline{3-3}

\xEH
$Con(T \xcv \ol{\ol{T'}})$, $T \xcl T'$ $ \xch $
\xEH
$\xci$ $(\xbm dp)$ (13.2)
\xEH
$X \xcc Y, X \xcs f(Y) \xEd \xCQ   \xch $
\xEP

\cline{3-3}

\xEH
$ \ol{\ol{T}} = \ol{\ol{\ol{T'}} \xcv T} $
\xEH
$\xcI$ without $(\xbm dp)$ (13.3)
\xEH
$f(X) = f(Y) \xcs X$
\xEP

\cline{3-3}

\xEH
\xEH
$\xci$ $T$ a formula (13.4)
\xEH
\xEP

\hline

\xEH
$(Log=' )$
\xEH
$\xch$ (14.1)
\xEH
$( \xbm =' )$
\xEP

\cline{3-3}

\xEH
$Con( \ol{ \ol{T' } } \xcv T)$ $ \xch $
\xEH
$\xci$ $(\xbm dp)$ (14.2)
\xEH
$f(Y) \xcs X \xEd \xCQ $ $ \xch $
\xEP

\cline{3-3}

\xEH
$ \ol{ \ol{T \xcv T' } }= \ol{ \ol{ \ol{T' } } \xcv T}$
\xEH
$\xcI$ without $(\xbm dp)$ (14.3)
\xEH
$f(Y \xcs X)=f(Y) \xcs X$
\xEP

\cline{3-3}

\xEH
\xEH
$\xci$ $T$ a formula (14.4)
\xEH
\xEP

\hline

\xEH
$(Log \xFO )$
\xEH
$\xch$ (15.1)
\xEH
$( \xbm \xFO )$
\xEP

\cline{3-3}

\xEH
$ \ol{ \ol{T \xco T' } }$ is one of
\xEH
$\xci$ (15.2)
\xEH
$f(X \xcv Y)$ is one of
\xEP

\cline{3-3}

\xEH
$\ol{\ol{T}},$ or $\ol{\ol{T'}},$ or $\ol{\ol{T}} \xcs \ol{\ol{T'}}$ (by (CCL))
\xEH
\xEH
$f(X),$ $f(Y)$ or $f(X) \xcv f(Y)$
\xEP

\hline

\xEH
$(Log \xcv )$
\xEH
$\xch$ $(\xbm \xcc)+(\xbm =)$ (16.1)
\xEH
$( \xbm \xcv )$
\xEP

\cline{3-3}

\xEH
$Con( \ol{ \ol{T' } } \xcv T),$ $ \xCN Con( \ol{ \ol{T' } }
\xcv \ol{ \ol{T} })$ $ \xch $
\xEH
$\xci$ $(\xbm dp)$ (16.2)
\xEH
$f(Y) \xcs (X-f(X)) \xEd \xCQ $ $ \xch $
\xEP

\cline{3-3}

\xEH
$ \xCN Con( \ol{ \ol{T \xco T' } } \xcv T' )$
\xEH
$\xcI$ without $(\xbm dp)$ (16.3)
\xEH
$f(X \xcv Y) \xcs Y= \xCQ$
\xEP

\hline

\xEH
$(Log \xcv ' )$
\xEH
$\xch$ $(\xbm \xcc)+(\xbm =)$ (17.1)
\xEH
$( \xbm \xcv ' )$
\xEP

\cline{3-3}

\xEH
$Con( \ol{ \ol{T' } } \xcv T),$ $ \xCN Con( \ol{ \ol{T' }
} \xcv \ol{ \ol{T} })$ $ \xch $
\xEH
$\xci$ $(\xbm dp)$ (17.2)
\xEH
$f(Y) \xcs (X-f(X)) \xEd \xCQ $ $ \xch $
\xEP

\cline{3-3}

\xEH
$ \ol{ \ol{T \xco T' } }= \ol{ \ol{T} }$
\xEH
$\xcI$ without $(\xbm dp)$ (17.3)
\xEH
$f(X \xcv Y)=f(X)$
\xEP

\hline

\xEH
\xEH
\xEH
$( \xbm \xbe )$
\xEP

\xEH
\xEH
\xEH
$a \xbe X-f(X)$ $ \xch $
\xEP

\xEH
\xEH
\xEH
$ \xcE b \xbe X.a \xce f(\{a,b\})$
\xEP

\hline

\end{tabular}

}

 karl-search= End Definition Log-Cond-Ref
\vspace{7mm}

 *************************************

\vspace{7mm}

\subsubsection{Fact Mu-Base}

 {\LARGE karl-search= Start Fact Mu-Base }

\index{Fact Mu-Base}

\bfa

$\hspace{0.01em}$

(+++ Orig. No.:  Fact Mu-Base +++)

{\xssc LABEL: {Fact Mu-Base}}
\label{Fact Mu-Base}

The following table is to be read as follows: If the left hand side holds
for
some
function $f: \xdy \xcp \xdp (U),$ and the auxiliary properties noted in
the middle also
hold for $f$ or $ \xdy,$ then the right hand side will hold, too - and
conversely.

{\small

\begin{tabular}{|c|c|c|c|}

\hline

\multicolumn{4}{|c|}{Basics} \xEP

\hline

(1.1)
\xEH
$(\xbm PR)$
\xEH
$\xch$ $(\xcs)+(\xbm \xcc)$
\xEH
$(\xbm PR')$
\xEP

\cline{1-1}

\cline{3-3}

(1.2)
\xEH
\xEH
$\xci$
\xEH
\xEP

\hline

(2.1)
\xEH
$(\xbm PR)$
\xEH
$\xch$ $(\xbm \xcc)$
\xEH
$(\xbm OR)$
\xEP

\cline{1-1}

\cline{3-3}

(2.2)
\xEH
\xEH
$\xci$ $(\xbm \xcc)$ + $(-)$
\xEH
\xEP

\cline{1-1}

\cline{3-4}

(2.3)
\xEH
\xEH
$\xch$ $(\xbm \xcc)$
\xEH
$(\xbm wOR)$
\xEP

\cline{1-1}

\cline{3-3}

(2.4)
\xEH
\xEH
$\xci$ $(\xbm \xcc)$ + $(-)$
\xEH
\xEP

\hline

(3)
\xEH
$(\xbm PR)$
\xEH
$\xch$
\xEH
$( \xbm CUT)$
\xEP

\hline

(4)
\xEH
$(\xbm \xcc )+(\xbm \xcc \xcd )+(\xbm CUM)+$
\xEH
$\xcH$
\xEH
$( \xbm PR)$
\xEP

\xEH
$(\xbm RatM)+(\xcs )$
\xEH
\xEH
\xEP

\hline

\multicolumn{4}{|c|}{Cumulativity} \xEP

\hline

(5.1)
\xEH
$(\xbm CM)$
\xEH
$\xch$ $(\xcs)+(\xbm \xcc)$
\xEH
$(\xbm ResM)$
\xEP

\cline{1-1}

\cline{3-3}

(5.2)
\xEH
\xEH
$\xci$ (infin.)
\xEH
\xEP

\hline

(6)
\xEH
$(\xbm CM)+(\xbm CUT)$
\xEH
$\xcj$
\xEH
$(\xbm CUM)$
\xEP

\hline

(7)
\xEH
$( \xbm \xcc )+( \xbm \xcc \xcd )$
\xEH
$\xch$
\xEH
$( \xbm CUM)$
\xEP

\hline

(8)
\xEH
$( \xbm \xcc )+( \xbm CUM)+( \xcs )$
\xEH
$\xch$
\xEH
$( \xbm \xcc \xcd )$
\xEP

\hline

(9)
\xEH
$( \xbm \xcc )+( \xbm CUM)$
\xEH
$\xcH$
\xEH
$( \xbm \xcc \xcd )$
\xEP

\hline

\multicolumn{4}{|c|}{Rationality} \xEP

\hline

(10)
\xEH
$( \xbm RatM )+( \xbm PR )$
\xEH
$\xch$
\xEH
$( \xbm =)$
\xEP

\hline

(11)
\xEH
$( \xbm =)$
\xEH
$ \xch $
\xEH
$( \xbm PR),$
\xEP

\hline

(12.1)
\xEH
$( \xbm =)$
\xEH
$ \xch $ $(\xcs)+( \xbm \xcc )$
\xEH
$( \xbm =' ),$
\xEP
\cline{1-1}
\cline{3-3}
(12.2)
\xEH
\xEH
$ \xci $
\xEH
\xEP

\hline

(13)
\xEH
$( \xbm \xcc ),$ $( \xbm =)$
\xEH
$ \xch $ $(\xcv)$
\xEH
$( \xbm \xcv ),$
\xEP

\hline

(14)
\xEH
$( \xbm \xcc ),$ $( \xbm \xCQ ),$ $( \xbm =)$
\xEH
$ \xch $ $(\xcv)$
\xEH
$( \xbm \xFO ),$ $( \xbm \xcv ' ),$ $( \xbm CUM),$
\xEP

\hline

(15)
\xEH
$( \xbm \xcc )+( \xbm \xFO )$
\xEH
$ \xch $ $(-)$ of $\xdy$
\xEH
$( \xbm =),$
\xEP

\hline

(16)
\xEH
$( \xbm \xFO )+( \xbm \xbe )+( \xbm PR)+$
\xEH
$ \xch $ $(\xcv)$ + $\xdy$ contains singletons
\xEH
$( \xbm =),$
\xEP
\xEH
$( \xbm \xcc )$
\xEH
\xEH
\xEP

\hline

(17)
\xEH
$( \xbm CUM)+( \xbm =)$
\xEH
$ \xch $ $(\xcv)$ + $\xdy$ contains singletons
\xEH
$( \xbm \xbe ),$
\xEP

\hline

(18)
\xEH
$( \xbm CUM)+( \xbm =)+( \xbm \xcc )$
\xEH
$ \xch $ $(\xcv)$
\xEH
$( \xbm \xFO ),$
\xEP

\hline

(19)
\xEH
$( \xbm PR)+( \xbm CUM)+( \xbm \xFO )$
\xEH
$ \xch $ sufficient, e.g. true in $\xdD_{\xdl}$
\xEH
$( \xbm =)$.
\xEP

\hline

(20)
\xEH
$( \xbm \xcc )+( \xbm PR)+( \xbm =)$
\xEH
$ \xcH $
\xEH
$( \xbm \xFO ),$
\xEP

\hline

(21)
\xEH
$( \xbm \xcc )+( \xbm PR)+( \xbm \xFO )$
\xEH
$ \xcH $ (without $(-)$)
\xEH
$( \xbm =)$
\xEP

\hline

(22)
\xEH
$( \xbm \xcc )+( \xbm PR)+( \xbm \xFO )+$
\xEH
$ \xcH $
\xEH
$( \xbm \xbe )$
\xEP
\xEH
$( \xbm =)+( \xbm \xcv )$
\xEH
\xEH
(thus not representability
\xEP
\xEH
\xEH
\xEH
by ranked structures)
\xEP

\hline

\end{tabular}

}

 karl-search= End Fact Mu-Base
\vspace{7mm}

 *************************************

\vspace{7mm}

\subsubsection{Fact Mu-Base Proof}

 {\LARGE karl-search= Start Fact Mu-Base Proof }

\index{Fact Mu-Base Proof}

\efa

\subparagraph{
Proof
}

$\hspace{0.01em}$

(+++*** Orig.:  Proof )

All sets are to be in $ \xdy.$

(1.1) $( \xbm PR)+( \xcs )+( \xbm \xcc )$ $ \xch $ $( \xbm PR' ):$

By $X \xcs Y \xcc X$ and $( \xbm PR),$ $f(X) \xcs X \xcs Y \xcc f(X \xcs
Y).$ By $( \xbm \xcc )$ $f(X) \xcs Y=f(X) \xcs X \xcs Y.$

(1.2) $( \xbm PR' ) \xch ( \xbm PR):$

Let $X \xcc Y,$ so $X=X \xcs Y,$ so by $( \xbm PR' )$ $f(Y) \xcs X \xcc
f(X \xcs Y)=f(X).$

(2.1) $( \xbm PR)+( \xbm \xcc )$ $ \xch $ $( \xbm OR):$

$f(X \xcv Y) \xcc X \xcv Y$ by $( \xbm \xcc ),$ so $f(X \xcv Y)$ $=$ $(f(X
\xcv Y) \xcs X) \xcv (f(X \xcv Y) \xcs Y)$ $ \xcc $ $f(X) \xcv f(Y).$

(2.2) $( \xbm OR)$ $+$ $( \xbm \xcc )$ $+$ $ \xCf (-)$ $ \xch $ $( \xbm
PR):$

Let $X \xcc Y,$ $X':=Y-X$. $f(Y) \xcc f(X) \xcv f(X' )$ by $( \xbm OR),$
so $f(Y) \xcs X$ $ \xcc $
$(f(X) \xcs X) \xcv (f(X' ) \xcs X)$ $=_{( \xbm \xcc )}$ $f(X) \xcv \xCQ $
$=$ $f(X).$

(2.3) $( \xbm PR)+( \xbm \xcc )$ $ \xch $ $( \xbm wOR):$

Trivial by (2.1).

(2.4) $( \xbm wOR)$ $+$ $( \xbm \xcc )$ $+$ $ \xCf (-)$ $ \xch $ $( \xbm
PR):$

Let $X \xcc Y,$ $X':=Y-X$. $f(Y) \xcc f(X) \xcv X' $ by $( \xbm wOR),$
so $f(Y) \xcs X$ $ \xcc $
$(f(X) \xcs X) \xcv (X' \xcs X)$ $=_{( \xbm \xcc )}$ $f(X) \xcv \xCQ $ $=$
$f(X).$

(3) $( \xbm PR)$ $ \xch $ $( \xbm CUT):$

$f(X) \xcc Y \xcc X$ $ \xch $ $f(X) \xcc f(X) \xcs Y \xcc f(Y)$ by $( \xbm
PR).$

(4) $( \xbm \xcc )+( \xbm \xcc \xcd )+( \xbm CUM)+( \xbm RatM)+( \xcs )$ $
\xcH $ $( \xbm PR):$

This is shown in Example \ref{Example Need-Pr} (page \pageref{Example Need-Pr})
.

(5.1) $( \xbm CM)+( \xcs )+( \xbm \xcc )$ $ \xch $ $( \xbm ResM):$

Let $f(X) \xcc A \xcs B,$ so $f(X) \xcc A,$ so by $( \xbm \xcc )$ $f(X)
\xcc A \xcs X \xcc X,$
so by $( \xbm CM)$ $f(A \xcs X) \xcc f(X) \xcc B.$

(5.2) $( \xbm ResM) \xch ( \xbm CM):$

We consider here the infinitary version, where all sets can be model sets
of infinite theories.
Let $f(X) \xcc Y \xcc X,$ so $f(X) \xcc Y \xcs f(X),$ so by $( \xbm ResM)$
$f(Y)=f(X \xcs Y) \xcc f(X).$

(6) $( \xbm CM)+( \xbm CUT)$ $ \xcj $ $( \xbm CUM):$

Trivial.

(7) $( \xbm \xcc )+( \xbm \xcc \xcd )$ $ \xch $ $( \xbm CUM):$

Suppose $f(D) \xcc E \xcc D.$ So by $( \xbm \xcc )$ $f(E) \xcc E \xcc D,$
so by $( \xbm \xcc \xcd )$ $f(D)=f(E).$

(8) $( \xbm \xcc )+( \xbm CUM)+( \xcs )$ $ \xch $ $( \xbm \xcc \xcd ):$

Let $f(D) \xcc E,$ $f(E) \xcc D,$ so by $( \xbm \xcc )$ $f(D) \xcc D \xcs
E \xcc D,$ $f(E) \xcc D \xcs E \xcc E.$ As $f(D \xcs E)$
is defined, so $f(D)=f(D \xcs E)=f(E)$ by $( \xbm CUM).$

(9) $( \xbm \xcc )+( \xbm CUM)$ $ \xcH $ $( \xbm \xcc \xcd ):$

This is shown in Example \ref{Example Mu-Cum-Cd} (page \pageref{Example
Mu-Cum-Cd}) .

(10) $( \xbm RatM)+( \xbm PR)$ $ \xch $ $( \xbm =):$

Trivial.

(11) $( \xbm =)$ entails $( \xbm PR):$

Trivial.

(12.1) $( \xbm =) \xcp ( \xbm =' ):$

Let $f(Y) \xcs X \xEd \xCQ,$ we have to show $f(X \xcs Y)=f(Y) \xcs X.$
By $( \xbm \xcc )$ $f(Y) \xcc Y,$ so $f(Y) \xcs X=f(Y) \xcs (X \xcs Y),$
so by $( \xbm =)$ $f(Y) \xcs X$ $=$
$f(Y) \xcs (X \xcs Y)$ $=$ $f(X \xcs Y).$

(12.2) $( \xbm =' ) \xcp ( \xbm =):$

Let $X \xcc Y,$ $f(Y) \xcs X \xEd \xCQ,$ then $f(X)=f(Y \xcs X)=f(Y) \xcs
X.$

(13) $( \xbm \xcc ),$ $( \xbm =)$ $ \xcp $ $( \xbm \xcv ):$

If not, $f(X \xcv Y) \xcs Y \xEd \xCQ,$ but $f(Y) \xcs (X-f(X)) \xEd \xCQ
.$ By (11), $( \xbm PR)$ holds,
so $f(X \xcv Y) \xcs X \xcc f(X),$ so $ \xCQ $ $ \xEd $ $f(Y) \xcs
(X-f(X))$ $ \xcc $ $f(Y) \xcs (X-f(X \xcv Y)),$ so
$f(Y)-f(X \xcv Y) \xEd \xCQ,$ so by $( \xbm \xcc )$ $f(Y) \xcc Y$ and
$f(Y) \xEd f(X \xcv Y) \xcs Y.$
But by $( \xbm =)$ $f(Y)=f(X \xcv Y) \xcs Y,$ a contradiction.

(14)

$( \xbm \xcc ),$ $( \xbm \xCQ ),$ $( \xbm =)$ $ \xch $ $( \xbm \xFO ):$

If $X$ or $Y$ or both are empty, then this is trivial.
Assume then $X \xcv Y \xEd \xCQ,$ so by $( \xbm \xCQ )$ $f(X \xcv Y) \xEd
\xCQ.$
By $( \xbm \xcc )$ $f(X \xcv Y) \xcc X \xcv Y,$ so $f(X \xcv Y) \xcs X=
\xCQ $ and $f(X \xcv Y) \xcs Y= \xCQ $ together are
impossible.
Case 1, $f(X \xcv Y) \xcs X \xEd \xCQ $ and $f(X \xcv Y) \xcs Y \xEd \xCQ
:$ By $( \xbm =)$ $f(X \xcv Y) \xcs X=f(X)$ and
$f(X \xcv Y) \xcs Y=f(Y),$ so by $( \xbm \xcc )$ $f(X \xcv Y)=f(X) \xcv
f(Y).$
Case 2, $f(X \xcv Y) \xcs X \xEd \xCQ $ and $f(X \xcv Y) \xcs Y= \xCQ:$
So by $( \xbm =)$ $f(X \xcv Y)=f(X \xcv Y) \xcs X=f(X).$
Case 3, $f(X \xcv Y) \xcs X= \xCQ $ and $f(X \xcv Y) \xcs Y \xEd \xCQ:$
Symmetrical.

$( \xbm \xcc ),$ $( \xbm \xCQ ),$ $( \xbm =)$ $ \xch $ $( \xbm \xcv ' ):$

Let $f(Y) \xcs (X-f(X)) \xEd \xCQ.$
If $X \xcv Y= \xCQ,$ then $f(X \xcv Y)=f(X)= \xCQ $ by $( \xbm \xcc ).$
So suppose $X \xcv Y \xEd \xCQ.$ By
(13), $f(X \xcv Y) \xcs Y= \xCQ,$ so $f(X \xcv Y) \xcc X$ by $( \xbm \xcc
).$ By $( \xbm \xCQ ),$ $f(X \xcv Y) \xEd \xCQ,$ so
$f(X \xcv Y) \xcs X \xEd \xCQ,$ and $f(X \xcv Y)=f(X)$ by $( \xbm =).$

$( \xbm \xcc ),$ $( \xbm \xCQ ),$ $( \xbm =)$ $ \xch $ $( \xbm CUM):$

Let $f(Y) \xcc X \xcc Y.$
If $Y= \xCQ,$ this is trivial by $( \xbm \xcc ).$ If $Y \xEd \xCQ,$ then
by $( \xbm \xCQ )$ - which is
crucial here - $f(Y) \xEd \xCQ,$ so by $f(Y) \xcc X$ $f(Y) \xcs X \xEd
\xCQ,$ so by $( \xbm =)$
$f(Y)=f(Y) \xcs X=f(X).$

(15) $( \xbm \xcc )+( \xbm \xFO )$ $ \xcp $ $( \xbm =):$

Let $X \xcc Y,$ $X \xcs f(Y) \xEd \xCQ,$ and consider $Y=X \xcv (Y-$X).
Then $f(Y)=f(X) \xFO f(Y-$X). As
$f(Y) \xcs X \xEd \xCQ,$ $f(Y)=f(Y-$X) is impossible. Otherwise,
$f(X)=f(Y) \xcs X,$ and we are
done.

(16) $( \xbm \xFO )+( \xbm \xbe )+( \xbm PR)+( \xbm \xcc )$ $ \xcp $ $(
\xbm =):$

Suppose $X \xcc Y,$ $x \xbe f(Y) \xcs X,$ we have to show $f(Y) \xcs
X=f(X).$ `` $ \xcc $ '' is trivial
by $( \xbm PR).$ `` $ \xcd $ '': Assume $a \xce f(Y)$ (by $( \xbm \xcc )),$
but $a \xbe f(X).$ By $( \xbm \xbe )$ $ \xcE b \xbe Y.a \xce f(\{a,b\}).$
As $a \xbe f(X),$ by $( \xbm PR),$ $a \xbe f(\{a,x\}).$ By $( \xbm \xFO
),$ $f(\{a,b,x\})$ $=$ $f(\{a,x\}) \xFO f(\{b\}).$
As $a \xce f(\{a,b,x\}),$ $f(\{a,b,x\})$ $=$ $f(\{b\}),$ so $x \xce
f(\{a,b,x\}),$ contradicting $( \xbm PR),$
as $a,b,x \xbe Y.$

(17) $( \xbm CUM)+( \xbm =)$ $ \xcp $ $( \xbm \xbe ):$

Let $a \xbe X-f(X).$ If $f(X)= \xCQ,$ then $f(\{a\})= \xCQ $ by $( \xbm
CUM).$ If not: Let
$b \xbe f(X),$ then $a \xce f(\{a,b\})$ by $( \xbm =).$

(18) $( \xbm CUM)+( \xbm =)+( \xbm \xcc )$ $ \xcp $ $( \xbm \xFO ):$

By $( \xbm CUM),$ $f(X \xcv Y) \xcc X \xcc X \xcv Y$ $ \xcp $ $f(X)=f(X
\xcv Y),$ and $f(X \xcv Y) \xcc Y \xcc X \xcv Y$ $ \xcp $
$f(Y)=f(X \xcv Y).$ Thus, if $( \xbm \xFO )$ were to fail, $f(X \xcv Y)
\xcC X,$ $f(X \xcv Y) \xcC Y,$ but then
by $( \xbm \xcc )$ $f(X \xcv Y) \xcs X \xEd \xCQ,$ so $f(X)=f(X \xcv Y)
\xcs X,$ and $f(X \xcv Y) \xcs Y \xEd \xCQ,$ so
$f(Y)=f(X \xcv Y) \xcs Y$ by $( \xbm =).$ Thus, $f(X \xcv Y)$ $=$ $(f(X
\xcv Y) \xcs X) \xcv (f(X \xcv Y) \xcs Y)$ $=$
$f(X) \xcv f(Y).$

(19) $( \xbm PR)+( \xbm CUM)+( \xbm \xFO )$ $ \xcp $ $( \xbm =):$

Suppose $( \xbm =)$ does not hold. So, by $( \xbm PR),$ there are $X,Y,y$
s.t. $X \xcc Y,$ $X \xcs f(Y) \xEd \xCQ,$
$y \xbe Y-f(Y),$ $y \xbe f(X).$ Let $a \xbe X \xcs f(Y).$ If $f(Y)=\{a\},$
then by $( \xbm CUM)$ $f(Y)=f(X),$ so
there must be $b \xbe f(Y),$ $b \xEd a.$ Take now $Y',$ $Y'' $ s.t. $Y=Y'
\xcv Y'',$ $a \xbe Y',$ $a \xce Y'',$ $b \xbe Y'',$
$b \xce Y',$ $y \xbe Y' \xcs Y''.$ Assume now $( \xbm \xFO )$ to hold,
we show a contradiction.
If $y \xce f(Y'' ),$ then by $( \xbm PR)$ $y \xce f(Y'' \xcv \{a\}).$ But
$f(Y'' \xcv \{a\})$ $=$ $f(Y'' ) \xFO f(\{a,y\}),$
so $f(Y'' \xcv \{a\})=f(Y'' ),$ contradicting $a \xbe f(Y).$ If $y \xbe
f(Y'' ),$ then by $f(Y)$ $=$
$f(Y' ) \xFO f(Y'' ),$ $f(Y)=f(Y' ),$ $contradiction$ as $b \xce f(Y' ).$

(20) $( \xbm \xcc )+( \xbm PR)+( \xbm =)$ $ \xcH $ $( \xbm \xFO ):$

See Example \ref{Example Mu-Barbar} (page \pageref{Example Mu-Barbar}) .

(21) $( \xbm \xcc )+( \xbm PR)+( \xbm \xFO )$ $ \xcH $ $( \xbm =):$

See Example \ref{Example Mu-Equal} (page \pageref{Example Mu-Equal}) .

(22) $( \xbm \xcc )+( \xbm PR)+( \xbm \xFO )+( \xbm =)+( \xbm \xcv )$ $
\xcH $ $( \xbm \xbe ):$

See Example \ref{Example Mu-Epsilon} (page \pageref{Example Mu-Epsilon}) .

Thus, by Fact \ref{Fact Rank-Hold} (page \pageref{Fact Rank-Hold}) , the
conditions do not assure
representability by ranked structures.

$ \xcz $
\\[3ex]

 karl-search= End Fact Mu-Base Proof
\vspace{7mm}

 *************************************

\vspace{7mm}

\subsubsection{Example Mu-Cum-Cd}

 {\LARGE karl-search= Start Example Mu-Cum-Cd }

\index{Example Mu-Cum-Cd}

\be

$\hspace{0.01em}$

(+++ Orig. No.:  Example Mu-Cum-Cd +++)

{\xssc LABEL: {Example Mu-Cum-Cd}}
\label{Example Mu-Cum-Cd}

We show here $( \xbm \xcc )+( \xbm CUM)$ $ \xcH $ $( \xbm \xcc \xcd ).$

Consider $X:=\{a,b,c\},$ $Y:=\{a,b,d\},$ $f(X):=\{a\},$ $f(Y):=\{a,b\},$ $
\xdy:=\{X,Y\}.$
(If $f(\{a,b\})$ were defined, we would have $f(X)=f(\{a,b\})=f(Y),$
$contradiction.)$

Obviously, $( \xbm \xcc )$ and $( \xbm CUM)$ hold, but not $( \xbm \xcc
\xcd ).$

$ \xcz $
\\[3ex]

 karl-search= End Example Mu-Cum-Cd
\vspace{7mm}

 *************************************

\vspace{7mm}

\subsubsection{Example Need-Pr}

 {\LARGE karl-search= Start Example Need-Pr }

\index{Example Need-Pr}

\ee

\be

$\hspace{0.01em}$

(+++ Orig. No.:  Example Need-Pr +++)

{\xssc LABEL: {Example Need-Pr}}
\label{Example Need-Pr}

We show here $( \xbm \xcc )+( \xbm \xcc \xcd )+( \xbm CUM)+( \xbm RatM)+(
\xcs )$ $ \xcH $ $( \xbm PR).$

Let $U:=\{a,b,c\}.$ Let $ \xdy = \xdp (U).$ So $( \xcs )$ is trivially
satisfied.
Set $f(X):=X$ for all $X \xcc U$ except for $f(\{a,b\})=\{b\}.$ Obviously,
this cannot be represented by a preferential structure and $( \xbm PR)$ is
false
for $U$ and $\{a,b\}.$ But it satisfies $( \xbm \xcc ),$ $( \xbm CUM),$ $(
\xbm RatM).$ $( \xbm \xcc )$ is trivial.
$( \xbm CUM):$ Let $f(X) \xcc Y \xcc X.$ If $f(X)=X,$ we are done.
Consider $f(\{a,b\})=\{b\}.$ If
$\{b\} \xcc Y \xcc \{a,b\},$ then $f(Y)=\{b\},$ so we are done again. It
is shown in
Fact \ref{Fact Mu-Base} (page \pageref{Fact Mu-Base}) , (8) that $( \xbm \xcc
\xcd )$ follows.
$( \xbm RatM):$ Suppose $X \xcc Y,$ $X \xcs f(Y) \xEd \xCQ,$ we have to
show $f(X) \xcc f(Y) \xcs X.$ If $f(Y)=Y,$ the result holds by $X \xcc Y,$
so it does if $X=Y.$
The only remaining case is $Y=\{a,b\},$ $X=\{b\},$ and the result holds
again.

$ \xcz $
\\[3ex]

 karl-search= End Example Need-Pr
\vspace{7mm}

 *************************************

\vspace{7mm}

\subsubsection{Example Mu-Barbar}

 {\LARGE karl-search= Start Example Mu-Barbar }

\index{Example Mu-Barbar}

\ee

\be

$\hspace{0.01em}$

(+++ Orig. No.:  Example Mu-Barbar +++)

{\xssc LABEL: {Example Mu-Barbar}}
\label{Example Mu-Barbar}

The example shows that $( \xbm \xcc )+( \xbm PR)+( \xbm =)$ $ \xcH $ $(
\xbm \xFO ).$

Consider the following structure without transitivity:
$U:=\{a,b,c,d\},$ $c$ and $d$ have $ \xbo $ many copies in descending
order $c_{1} \xed c_{2}$  \Xl., etc.
$a,b$ have one single copy each. $a \xed b,$ $a \xed d_{1},$ $b \xed a,$
$b \xed c_{1}.$
$( \xbm \xFO )$ does not hold: $f(U)= \xCQ,$ but $f(\{a,c\})=\{a\},$
$f(\{b,d\})=\{b\}.$
$( \xbm PR)$ holds as in all preferential structures.
$( \xbm =)$ holds: If it were to fail, then for some $A \xcc B,$ $f(B)
\xcs A \xEd \xCQ,$ so $f(B) \xEd \xCQ.$
But the only possible cases for $B$ are now: $(a \xbe B,$ $b,d \xce B)$ or
$(b \xbe B,$ $a,c \xce B).$
Thus, $B$ can be $\{a\},$ $\{a,c\},$ $\{b\},$ $\{b,d\}$ with $f(B)=$
$\{a\},$ $\{a\},$ $\{b\},$ $\{b\}.$
If $A=B,$ then the result will hold trivially. Moreover, $ \xCf A$ has to
be $ \xEd \xCQ.$
So the remaining cases of $B$ where it might fail are $B=$ $\{a,c\}$ and
$\{b,d\},$ and
by $f(B) \xcs A \xEd \xCQ,$ the only cases of $ \xCf A$ where it might
fail, are $A=$ $\{a\}$ or $\{b\}$
respectively.
So the only cases remaining are: $B=\{a,c\},$ $A=\{a\}$ and $B=\{b,d\},$
$A=\{b\}.$
In the first case, $f(A)=f(B)=\{a\},$ in the second $f(A)=f(B)=\{b\},$ but
$( \xbm =)$
holds in both.

$ \xcz $
\\[3ex]

 karl-search= End Example Mu-Barbar
\vspace{7mm}

 *************************************

\vspace{7mm}

\subsubsection{Example Mu-Equal}

 {\LARGE karl-search= Start Example Mu-Equal }

\index{Example Mu-Equal}

\ee

\be

$\hspace{0.01em}$

(+++ Orig. No.:  Example Mu-Equal +++)

{\xssc LABEL: {Example Mu-Equal}}
\label{Example Mu-Equal}

The example shows that $( \xbm \xcc )+( \xbm PR)+( \xbm \xFO )$ $ \xcH $
$( \xbm =).$

Work in the set of theory definable model sets of an infinite
propositional
language. Note that this is not closed under set difference, and closure
properties will play a crucial role in the argumentation.
Let $U:=\{y,a,x_{i< \xbo }\},$ where $x_{i} \xcp a$ in the standard
topology. For the order,
arrange s.t. $y$ is minimized by any set iff this set contains a cofinal
subsequence of
the $x_{i},$ this can be done by the standard construction. Moreover, let
the $x_{i}$
all kill themselves, i.e. with $ \xbo $ many copies $x^{1}_{i} \xed
x^{2}_{i} \xed $  \Xl. There are no other
elements in the relation. Note that if $a \xce \xbm (X),$ then $a \xce X,$
and $X$ cannot contain
a cofinal subsequence of the $x_{i},$ as $X$ is closed in the standard
topology.
(A short argument: suppose $X$ contains such a subsequence, but $a \xce
X.$ Then the
theory of a $Th(a)$ is inconsistent with $Th(X),$ so already a finite
subset of
$Th(a)$ is inconsistent with $Th(X),$ but such a finite subset will
finally hold
in a cofinal sequence converging to a.)
Likewise, if $y \xbe \xbm (X),$ then $X$ cannot contain a cofinal
subsequence of the $x_{i}.$

Obviously, $( \xbm \xcc )$ and $( \xbm PR)$ hold, but $( \xbm =)$ does not
hold: Set $B:=U,$ $A:=\{a,y\}.$
Then $ \xbm (B)=\{a\},$ $ \xbm (A)=\{a,y\},$ contradicting $( \xbm =).$

It remains to show that $( \xbm \xFO )$ holds.

$ \xbm (X)$ can only be $ \xCQ,$ $\{a\},$ $\{y\},$ $\{a,y\}.$ As $ \xbm
(A \xcv B) \xcc \xbm (A) \xcv \xbm (B)$ by $( \xbm PR),$

Case 1, $ \xbm (A \xcv B)=\{a,y\}$ is settled.

Note that if $y \xbe X- \xbm (X),$ then $X$ will contain a cofinal
subsequence, and thus
$a \xbe \xbm (X).$

Case 2: $ \xbm (A \xcv B)=\{a\}.$

Case 2.1: $ \xbm (A)=\{a\}$ - we are done.

Case 2.2: $ \xbm (A)=\{y\}:$ $ \xCf A$ does not contain $ \xCf a,$ nor a
cofinal subsequence.
If $ \xbm (B)= \xCQ,$ then $a \xce B,$ so $a \xce A \xcv B,$ a
contradiction.
If $ \xbm (B)=\{a\},$ we are done.
If $y \xbe \xbm (B),$ then $y \xbe B,$ but $B$ does not contain a cofinal
subsequence, so
$A \xcv B$ does not either, so $y \xbe \xbm (A \xcv B),$ $contradiction.$

Case 2.3: $ \xbm (A)= \xCQ:$ $ \xCf A$ cannot contain a cofinal
subsequence.
If $ \xbm (B)=\{a\},$ we are done.
$a \xbe \xbm (B)$ does have to hold, so $ \xbm (B)=\{a,y\}$ is the only
remaining possibility.
But then $B$ does not contain a cofinal subsequence, and neither does $A
\xcv B,$ so
$y \xbe \xbm (A \xcv B),$ $contradiction.$

Case 2.4: $ \xbm (A)=\{a,y\}:$ $ \xCf A$ does not contain a cofinal
subsequence.
If $ \xbm (B)=\{a\},$ we are done.
If $ \xbm (B)= \xCQ,$ $B$ does not contain a cofinal subsequence (as $a
\xce B),$ so neither
does $A \xcv B,$ so $y \xbe \xbm (A \xcv B),$ $contradiction.$
If $y \xbe \xbm (B),$ $B$ does not contain a cofinal subsequence, and we
are done again.

Case 3: $ \xbm (A \xcv B)=\{y\}:$
To obtain a contradiction, we need $a \xbe \xbm (A)$ or $a \xbe \xbm (B).$
But in both cases
$a \xbe \xbm (A \xcv B).$

Case 4: $ \xbm (A \xcv B)= \xCQ:$
Thus, $A \xcv B$ contains no cofinal subsequence. If, e.g. $y \xbe \xbm
(A),$ then $y \xbe \xbm (A \xcv B),$
if $a \xbe \xbm (A),$ then $a \xbe \xbm (A \xcv B),$ so $ \xbm (A)= \xCQ
.$

$ \xcz $
\\[3ex]

 karl-search= End Example Mu-Equal
\vspace{7mm}

 *************************************

\vspace{7mm}

\subsubsection{Example Mu-Epsilon}

 {\LARGE karl-search= Start Example Mu-Epsilon }

\index{Example Mu-Epsilon}

\ee

\be

$\hspace{0.01em}$

(+++ Orig. No.:  Example Mu-Epsilon +++)

{\xssc LABEL: {Example Mu-Epsilon}}
\label{Example Mu-Epsilon}

The example show that $( \xbm \xcc )+( \xbm PR)+( \xbm \xFO )+( \xbm =)+(
\xbm \xcv )$ $ \xcH $ $( \xbm \xbe ).$

Let $U:=\{y,x_{i< \xbo }\},$ $x_{i}$ a sequence, each $x_{i}$ kills
itself, $x^{1}_{i} \xed x^{2}_{i} \xed  \Xl $
and $y$ is killed by all cofinal subsequences of the $x_{i}.$ Then for any
$X \xcc U$
$ \xbm (X)= \xCQ $ or $ \xbm (X)=\{y\}.$

$( \xbm \xcc )$ and $( \xbm PR)$ hold obviously.

$( \xbm \xFO ):$ Let $A \xcv B$ be given. If $y \xce X,$ then for all $Y
\xcc X$ $ \xbm (Y)= \xCQ.$
So, if $y \xce A \xcv B,$ we are done. If $y \xbe A \xcv B,$ if $ \xbm (A
\xcv B)= \xCQ,$ one of $A,B$ must contain
a cofinal sequence, it will have $ \xbm = \xCQ.$ If not, then $ \xbm (A
\xcv B)=\{y\},$ and this will
also hold for the one $y$ is in.

$( \xbm =):$ Let $A \xcc B,$ $ \xbm (B) \xcs A \xEd \xCQ,$ show $ \xbm
(A)= \xbm (B) \xcs A.$ But now $ \xbm (B)=\{y\},$ $y \xbe A,$
so $B$ does not contain a cofinal subsequence, neither does A, so $ \xbm
(A)=\{y\}.$

$( \xbm \xcv ):$ $(A- \xbm (A)) \xcs \xbm (A' ) \xEd \xCQ,$ so $ \xbm (A'
)=\{y\},$ so $ \xbm (A \xcv A' )= \xCQ,$ as $y \xbe A- \xbm (A).$

But $( \xbm \xbe )$ does not hold: $y \xbe U- \xbm (U),$ but there is no
$x$ s.t. $y \xce \xbm (\{x,y\}).$

$ \xcz $
\\[3ex]

 karl-search= End Example Mu-Epsilon
\vspace{7mm}

 *************************************

\vspace{7mm}

\subsubsection{Fact Mwor}

 {\LARGE karl-search= Start Fact Mwor }

\index{Fact Mwor}

\ee

\bfa

$\hspace{0.01em}$

(+++ Orig. No.:  Fact Mwor +++)

{\xssc LABEL: {Fact Mwor}}
\label{Fact Mwor}

$( \xbm wOR)+( \xbm \xcc )$ $ \xch $ $f(X \xcv Y) \xcc f(X) \xcv f(Y) \xcv
(X \xcs Y)$

 karl-search= End Fact Mwor
\vspace{7mm}

 *************************************

\vspace{7mm}

\subsubsection{Fact Mwor Proof}

 {\LARGE karl-search= Start Fact Mwor Proof }

\index{Fact Mwor Proof}

\efa

\subparagraph{
Proof
}

$\hspace{0.01em}$

(+++*** Orig.:  Proof )

$f(X \xcv Y) \xcc f(X) \xcv Y,$ $f(X \xcv Y) \xcc X \xcv f(Y),$ so $f(X
\xcv Y)$ $ \xcc $ $(f(X) \xcv Y) \xcs (X \xcv f(Y))$ $=$
$f(X) \xcv f(Y) \xcv (X \xcs Y)$ $ \xcz $
\\[3ex]

 karl-search= End Fact Mwor Proof
\vspace{7mm}

 *************************************

\vspace{7mm}

\subsubsection{Proposition Alg-Log}

 {\LARGE karl-search= Start Proposition Alg-Log }

\index{Proposition Alg-Log}

\bp

$\hspace{0.01em}$

(+++ Orig. No.:  Proposition Alg-Log +++)

{\xssc LABEL: {Proposition Alg-Log}}
\label{Proposition Alg-Log}

The following table is to be read as follows:

Let a logic $ \xcn $ satisfy $ \xCf (LLE)$ and $ \xCf (CCL),$ and define a
function $f: \xdD_{ \xdl } \xcp \xdD_{ \xdl }$
by $f(M(T)):=M( \ol{ \ol{T} }).$ Then $f$ is well defined, satisfies $(
\xbm dp),$ and $ \ol{ \ol{T} }=Th(f(M(T))).$

If $ \xcn $ satisfies a rule in the left hand side,
then - provided the additional properties noted in the middle for $ \xch $
hold, too -
$f$ will satisfy the property in the right hand side.

Conversely, if $f: \xdy \xcp \xdp (M_{ \xdl })$ is a function, with $
\xdD_{ \xdl } \xcc \xdy,$ and we define a logic
$ \xcn $ by $ \ol{ \ol{T} }:=Th(f(M(T))),$ then $ \xcn $ satisfies $ \xCf
(LLE)$ and $ \xCf (CCL).$
If $f$ satisfies $( \xbm dp),$ then $f(M(T))=M( \ol{ \ol{T} }).$

If $f$ satisfies a property in the right hand side,
then - provided the additional properties noted in the middle for $ \xci $
hold, too -
$ \xcn $ will satisfy the property in the left hand side.

If ``formula'' is noted in the table, this means that, if one of the
theories
(the one named the same way in Definition \ref{Definition Log-Cond} (page
\pageref{Definition Log-Cond}) )
is equivalent to a formula, we do not need $( \xbm dp).$

{\small

\begin{tabular}{|c|c|c|c|}

\hline

\multicolumn{4}{|c|}{Basics} \xEP

\hline

(1.1) \xEH $(OR)$ \xEH $\xch$ \xEH $(\xbm OR)$ \xEP

\cline{1-1}

\cline{3-3}

(1.2) \xEH \xEH $\xci$ \xEH \xEP

\hline

(2.1) \xEH $(disjOR)$ \xEH $\xch$ \xEH $(\xbm disjOR)$ \xEP

\cline{1-1}

\cline{3-3}

(2.2) \xEH \xEH $\xci$ \xEH \xEP

\hline

(3.1) \xEH $(wOR)$ \xEH $\xch$ \xEH $(\xbm wOR)$ \xEP

\cline{1-1}

\cline{3-3}

(3.2) \xEH \xEH $\xci$ \xEH \xEP

\hline

(4.1) \xEH $(SC)$ \xEH $\xch$ \xEH $(\xbm \xcc)$ \xEP

\cline{1-1}

\cline{3-3}

(4.2) \xEH \xEH $\xci$ \xEH \xEP

\hline

(5.1) \xEH $(CP)$ \xEH $\xch$ \xEH $(\xbm \xCQ)$ \xEP

\cline{1-1}

\cline{3-3}

(5.2) \xEH \xEH $\xci$ \xEH \xEP

\hline

(6.1) \xEH $(PR)$ \xEH $\xch$ \xEH $(\xbm PR)$ \xEP

\cline{1-1}

\cline{3-3}

(6.2) \xEH \xEH $\xci$ $(\xbm dp)+(\xbm \xcc)$ \xEH \xEP

\cline{1-1}

\cline{3-3}

(6.3) \xEH \xEH $\xcI$ without $(\xbm dp)$ \xEH \xEP

\cline{1-1}

\cline{3-3}

(6.4) \xEH \xEH $\xci$ $(\xbm \xcc)$ \xEH \xEP

\xEH \xEH $T'$ a formula \xEH \xEP

\hline

(6.5) \xEH $(PR)$ \xEH $\xci$ \xEH $(\xbm PR')$ \xEP

\xEH \xEH $T'$ a formula \xEH \xEP

\hline

(7.1) \xEH $(CUT)$ \xEH $\xch$ \xEH $(\xbm CUT)$ \xEP

\cline{1-1}

\cline{3-3}

(7.2) \xEH \xEH $\xci$ \xEH \xEP

\hline

\multicolumn{4}{|c|}{Cumulativity} \xEP

\hline

(8.1) \xEH $(CM)$ \xEH $\xch$ \xEH $(\xbm CM)$ \xEP

\cline{1-1}

\cline{3-3}

(8.2) \xEH \xEH $\xci$ \xEH \xEP

\hline

(9.1) \xEH $(ResM)$ \xEH $\xch$ \xEH $(\xbm ResM)$ \xEP

\cline{1-1}

\cline{3-3}

(9.2) \xEH \xEH $\xci$ \xEH \xEP

\hline

(10.1) \xEH $(\xcc \xcd)$ \xEH $\xch$ \xEH $(\xbm \xcc \xcd)$ \xEP

\cline{1-1}

\cline{3-3}

(10.2) \xEH \xEH $\xci$ \xEH \xEP

\hline

(11.1) \xEH $(CUM)$ \xEH $\xch$ \xEH $(\xbm CUM)$ \xEP

\cline{1-1}

\cline{3-3}

(11.2) \xEH \xEH $\xci$ \xEH \xEP

\hline

\multicolumn{4}{|c|}{Rationality} \xEP

\hline

(12.1) \xEH $(RatM)$ \xEH $\xch$ \xEH $(\xbm RatM)$ \xEP

\cline{1-1}

\cline{3-3}

(12.2) \xEH \xEH $\xci$ $(\xbm dp)$ \xEH \xEP

\cline{1-1}

\cline{3-3}

(12.3) \xEH \xEH $\xcI$ without $(\xbm dp)$ \xEH \xEP

\cline{1-1}

\cline{3-3}

(12.4) \xEH \xEH $\xci$ \xEH \xEP

\xEH \xEH $T$ a formula \xEH \xEP

\hline

(13.1) \xEH $(RatM=)$ \xEH $\xch$ \xEH $(\xbm =)$ \xEP

\cline{1-1}

\cline{3-3}

(13.2) \xEH \xEH $\xci$ $(\xbm dp)$ \xEH \xEP

\cline{1-1}

\cline{3-3}

(13.3) \xEH \xEH $\xcI$ without $(\xbm dp)$ \xEH \xEP

\cline{1-1}

\cline{3-3}

(13.4) \xEH \xEH $\xci$ \xEH \xEP

\xEH \xEH $T$ a formula \xEH \xEP

\hline

(14.1) \xEH $(Log = ')$ \xEH $\xch$ \xEH $(\xbm = ')$ \xEP

\cline{1-1}
\cline{3-3}

(14.2) \xEH \xEH $\xci$ $(\xbm dp)$ \xEH \xEP

\cline{1-1}
\cline{3-3}

(14.3) \xEH \xEH $\xcI$ without $(\xbm dp)$ \xEH \xEP

\cline{1-1}
\cline{3-3}

(14.4) \xEH \xEH $\xci$ $T$ a formula \xEH \xEP

\hline

(15.1) \xEH $(Log \xFO )$ \xEH $\xch$ \xEH $(\xbm \xFO )$ \xEP

\cline{1-1}
\cline{3-3}

(15.2) \xEH \xEH $\xci$ \xEH \xEP

\hline

(16.1)
\xEH
$(Log \xcv )$
\xEH
$\xch$ $(\xbm \xcc)+(\xbm =)$
\xEH
$(\xbm \xcv )$
\xEP

\cline{1-1}
\cline{3-3}

(16.2) \xEH \xEH $\xci$ $(\xbm dp)$ \xEH \xEP

\cline{1-1}
\cline{3-3}

(16.3) \xEH \xEH $\xcI$ without $(\xbm dp)$ \xEH \xEP

\hline

(17.1)
\xEH
$(Log \xcv ')$
\xEH
$\xch$ $(\xbm \xcc)+(\xbm =)$
\xEH
$(\xbm \xcv ')$
\xEP

\cline{1-1}
\cline{3-3}

(17.2) \xEH \xEH $\xci$ $(\xbm dp)$ \xEH \xEP

\cline{1-1}
\cline{3-3}

(17.3) \xEH \xEH $\xcI$ without $(\xbm dp)$ \xEH \xEP

\hline

\end{tabular}

}

 karl-search= End Proposition Alg-Log
\vspace{7mm}

 *************************************

\vspace{7mm}

\subsubsection{Proposition Alg-Log Proof}

 {\LARGE karl-search= Start Proposition Alg-Log Proof }

\index{Proposition Alg-Log Proof}

\ep

\subparagraph{
Proof
}

$\hspace{0.01em}$

(+++*** Orig.:  Proof )

Set $f(T):=f(M(T)),$ note that $f(T \xcv T' ):=f(M(T \xcv T' ))=f(M(T)
\xcs M(T' )).$

We show first the general framework.

Let $ \xcn $ satisfy $ \xCf (LLE)$ and $ \xCf (CCL).$ Let $f: \xdD_{ \xdl
} \xcp \xdD_{ \xdl }$ be defined by $f(M(T)):=M( \ol{ \ol{T} }).$
If $M(T)=M(T' ),$ then $ \ol{T}= \ol{T' },$ so by $ \xCf (LLE)$ $ \ol{
\ol{T} }= \ol{ \ol{T' } },$ so $f(M(T))=f(M(T' )),$ so $f$ is
well defined and satisfies $( \xbm dp).$ By $ \xCf (CCL)$ $Th(M( \ol{
\ol{T} }))= \ol{ \ol{T} }.$

Let $f$ be given, and $ \xcn $ be defined by $ \ol{ \ol{T}
}:=Th(f(M(T))).$ Obviously, $ \xcn $ satisfies
$ \xCf (LLE)$ and $ \xCf (CCL)$ (and thus $ \xCf (RW)).$ If $f$ satisfies
$( \xbm dp),$ then $f(M(T))=M(T' )$ for
some $T',$ and $f(M(T))=M(Th(f(M(T))))=M( \ol{ \ol{T} })$ by Fact \ref{Fact
Dp-Base} (page \pageref{Fact Dp-Base}) . (We will use
Fact \ref{Fact Dp-Base} (page \pageref{Fact Dp-Base})  now without further
mentioning.)

Next we show the following fact:

(a) If $f$ satisfies $( \xbm dp),$ or $T' $ is equivalent to a formula,
then
$Th(f(T) \xcs M(T' ))= \ol{ \ol{ \ol{T} } \xcv T' }.$

Case 1, $f$ satisfies $( \xbm dp).$ $Th(f(M(T)) \xcs M(T' ))$ $=$ $Th(M(
\ol{ \ol{T} }) \xcs M(T' )$ $=$ $ \ol{ \ol{ \ol{T} } \xcv T' }$
by Fact \ref{Fact Log-Form} (page \pageref{Fact Log-Form})  (5).

Case 2, $T' $ is equivalent to $ \xbf '.$ $Th(f(M(T)) \xcs M( \xbf ' ))$
$=$ $ \ol{Th(f(M(T))) \xcv \{ \xbf ' \}}$ $=$
$ \ol{ \ol{ \ol{T} } \xcv \{ \xbf ' \}}$ by Fact \ref{Fact Log-Form} (page
\pageref{Fact Log-Form})
(3).

We now prove the individual properties.

(1.1) $ \xCf (OR)$ $ \xch $ $( \xbm OR)$

Let $X=M(T),$ $Y=M(T' ).$ $f(X \xcv Y)$ $=$ $f(M(T) \xcv M(T' ))$ $=$
$f(M(T \xco T' ))$ $:=$ $M( \ol{ \ol{T \xco T' } })$ $ \xcc_{(OR)}$
$M( \ol{ \ol{T} } \xcs \ol{ \ol{T' } })$ $=_{(CCL)}$ $M( \ol{ \ol{T} })
\xcv M( \ol{ \ol{T' } })$ $=:$ $f(X) \xcv f(Y).$

(1.2) $( \xbm OR)$ $ \xch $ $ \xCf (OR)$

$ \ol{ \ol{T \xco T' } }$ $:=$ $Th(f(M(T \xco T' )))$ $=$ $Th(f(M(T) \xcv
M(T' )))$ $ \xcd_{( \xbm OR)}$ $Th(f(M(T)) \xcv f(M(T' )))$ $=$
(by Fact \ref{Fact Th-Union} (page \pageref{Fact Th-Union}) ) $Th(f(M(T))) \xcs
Th(f(M(T' )))$ $=:$ $
\ol{ \ol{T} } \xcs \ol{ \ol{T' } }.$

(2) By $ \xCN Con(T,T' ) \xcj M(T) \xcs M(T' )= \xCQ,$ we can use
directly the proofs for 1.

(3.1) $ \xCf (wOR)$ $ \xch $ $( \xbm wOR)$

Let $X=M(T),$ $Y=M(T' ).$ $f(X \xcv Y)$ $=$ $f(M(T) \xcv M(T' ))$ $=$
$f(M(T \xco T' ))$ $:=$ $M( \ol{ \ol{T \xco T' } })$ $ \xcc_{(wOR)}$
$M( \ol{ \ol{T} } \xcs \ol{T' })$ $=_{(CCL)}$ $M( \ol{ \ol{T} }) \xcv M(
\ol{T' })$ $=:$ $f(X) \xcv Y.$

(3.2) $( \xbm wOR)$ $ \xch $ $ \xCf (wOR)$

$ \ol{ \ol{T \xco T' } }$ $:=$ $Th(f(M(T \xco T' )))$ $=$ $Th(f(M(T) \xcv
M(T' )))$ $ \xcd_{( \xbm wOR)}$ $Th(f(M(T)) \xcv M(T' ))$ $=$
(by Fact \ref{Fact Th-Union} (page \pageref{Fact Th-Union}) ) $Th(f(M(T))) \xcs
Th(M(T' ))$ $=:$ $
\ol{ \ol{T} } \xcs \ol{T' }.$

(4.1) $ \xCf (SC)$ $ \xch $ $( \xbm \xcc )$

Trivial.

(4.2) $( \xbm \xcc )$ $ \xch $ $ \xCf (SC)$

Trivial.

(5.1) $ \xCf (CP)$ $ \xch $ $( \xbm \xCQ )$

Trivial.

(5.2) $( \xbm \xCQ )$ $ \xch $ $ \xCf (CP)$

Trivial.

(6.1) $ \xCf (PR)$ $ \xch $ $( \xbm PR):$

Suppose $X:=M(T),$ $Y:=M(T' ),$ $X \xcc Y,$ we have to show $f(Y) \xcs X
\xcc f(X).$
By prerequisite, $ \ol{T' } \xcc \ol{T},$ so $ \ol{T \xcv T' }= \ol{T},$
so $ \ol{ \ol{T \xcv T' } }= \ol{ \ol{T} }$ by $ \xCf (LLE).$ By $ \xCf
(PR)$ $ \ol{ \ol{T \xcv T' } } \xcc \ol{ \ol{ \ol{T'
} } \xcv T},$
so $f(Y) \xcs X=f(T' ) \xcs M(T)=M( \ol{ \ol{T' } } \xcv T) \xcc M( \ol{
\ol{T \xcv T' } })=M( \ol{ \ol{T} })=f(X).$

(6.2) $( \xbm PR)+( \xbm dp)+( \xbm \xcc )$ $ \xch $ $ \xCf (PR):$

$f(T) \xcs M(T' )=_{( \xbm \xcc )}f(T) \xcs M(T) \xcs M(T' )=f(T) \xcs M(T
\xcv T' ) \xcc_{( \xbm PR)}f(T \xcv T' ),$ so
$ \ol{ \ol{T \xcv T' } }=Th(f(T \xcv T' )) \xcc Th(f(T) \xcs M(T' ))= \ol{
\ol{ \ol{T} } \xcv T' }$ by (a) above and $( \xbm dp).$

(6.3) $( \xbm PR)$ $ \xcH $ $ \xCf (PR)$ without $( \xbm dp):$

$( \xbm PR)$ holds in all preferential structures
(see Definition \ref{Definition Pref-Str} (page \pageref{Definition Pref-Str}) )
by Fact \ref{Fact Pref-Sound} (page \pageref{Fact Pref-Sound}) .
Example \ref{Example Pref-Dp} (page \pageref{Example Pref-Dp})  shows that $
\xCf (DP)$ may fail in the
resulting logic.

(6.4) $( \xbm PR)+( \xbm \xcc )$ $ \xch $ $ \xCf (PR)$ if $T' $ is
classically equivalent to a formula:

It was shown in the proof of (6.2) that $f(T) \xcs M( \xbf ' ) \xcc f(T
\xcv \{ \xbf ' \}),$ so
$ \ol{ \ol{T \xcv \{ \xbf ' \}} }=Th(f(T \xcv \{ \xbf ' \})) \xcc Th(f(T)
\xcs M( \xbf ' ))= \ol{ \ol{ \ol{T} } \xcv \{ \xbf ' \}}$ by (a) above.

(6.5) $( \xbm PR' )$ $ \xch $ $ \xCf (PR),$ if $T' $ is classically
equivalent to a formula:

$f(M(T)) \xcs M( \xbf ' )$ $ \xcc_{( \xbm PR' )}$ $f(M(T) \xcs M( \xbf '
))$ $=$ $f(M(T \xcv \{ \xbf ' \})).$ So again
$ \ol{ \ol{T \xcv \{ \xbf ' \}} }=Th(f(T \xcv \{ \xbf ' \})) \xcc Th(f(T)
\xcs M( \xbf ' ))= \ol{ \ol{ \ol{T} } \xcv \{ \xbf ' \}}$ by (a) above.

(7.1) $ \xCf (CUT)$ $ \xch $ $( \xbm CUT)$

So let $X=M(T),$ $Y=M(T' ),$ and $f(T):=M( \ol{ \ol{T} }) \xcc M(T' ) \xcc
M(T)$ $ \xcp $ $ \ol{T} \xcc \ol{T' } \xcc \ol{ \ol{T} }=_{ \xCf (LLE)}
\ol{ \ol{( \ol{T})} }$ $ \xcp $
(by $ \xCf (CUT))$
$ \ol{ \ol{T} }= \ol{ \ol{( \ol{T})} } \xcd \ol{ \ol{( \ol{T' })} }= \ol{
\ol{T' } }$ $ \xcp $ $f(T)=M( \ol{ \ol{T} }) \xcc M( \ol{ \ol{T' } })=f(T'
),$ $thus$ $f(X) \xcc f(Y).$

(7.2) $( \xbm CUT)$ $ \xch $ $ \xCf (CUT)$

Let $T$ $ \xcc $ $ \ol{T' }$ $ \xcc $ $ \ol{ \ol{T} }.$ Thus $f(T) \xcc M(
\ol{ \ol{T} })$ $ \xcc $ $M(T' )$ $ \xcc $ $M(T),$ so by $( \xbm CUT)$
$f(T) \xcc f(T' ),$
so $ \ol{ \ol{T} }$ $=$ $Th(f(T))$ $ \xcd $ $Th(f(T' ))$ $=$ $ \ol{ \ol{T'
} }.$

(8.1) $ \xCf (CM)$ $ \xch $ $( \xbm CM)$

So let $X=M(T),$ $Y=M(T' ),$ and $f(T):=M( \ol{ \ol{T} }) \xcc M(T' ) \xcc
M(T)$ $ \xcp $ $ \ol{T} \xcc \ol{T' } \xcc \ol{ \ol{T} }=_{ \xCf (LLE)}
\ol{ \ol{( \ol{T})} }$ $ \xcp $
(by $ \xCf (LLE),$ $ \xCf (CM))$
$ \ol{ \ol{T} }= \ol{ \ol{( \ol{T})} } \xcc \ol{ \ol{( \ol{T' })} }= \ol{
\ol{T' } }$ $ \xcp $ $f(T)=M( \ol{ \ol{T} }) \xcd M( \ol{ \ol{T' } })=f(T'
),$ $thus$ $f(X) \xcd f(Y).$

(8.2) $( \xbm CM)$ $ \xch $ $ \xCf (CM)$

Let $T$ $ \xcc $ $ \ol{T' }$ $ \xcc $ $ \ol{ \ol{T} }.$ Thus by $( \xbm
CM)$ and $f(T) \xcc M( \ol{ \ol{T} })$ $ \xcc $ $M(T' )$ $ \xcc $ $M(T),$
so $f(T) \xcd f(T' )$
by $( \xbm CM),$ so $ \ol{ \ol{T} }$ $=$ $Th(f(T))$ $ \xcc $ $Th(f(T' ))$
$=$ $ \ol{ \ol{T' } }.$

(9.1) $ \xCf (ResM)$ $ \xch $ $( \xbm ResM)$

Let $f(X):=M( \ol{ \ol{ \xbD } }),$ $A:=M( \xba ),$ $B:=M( \xbb ).$ So
$f(X) \xcc A \xcs B$ $ \xch $ $ \xbD \xcn \xba, \xbb $ $ \xch_{ \xCf
(ResM)}$
$ \xbD, \xba \xcn \xbb $ $ \xch $ $M( \ol{ \ol{ \xbD, \xba } }) \xcc M(
\xbb )$ $ \xch $ $f(X \xcs A) \xcc B.$

(9.2) $( \xbm ResM)$ $ \xch $ $ \xCf (ResM)$

Let $f(X):=M( \ol{ \ol{ \xbD } }),$ $A:=M( \xba ),$ $B:=M( \xbb ).$ So $
\xbD \xcn \xba, \xbb $ $ \xch $ $f(X) \xcc A \xcs B$ $ \xch_{( \xbm
ResM)}$
$f(X \xcs A) \xcc B$ $ \xch $ $ \xbD, \xba \xcn \xbb.$

(10.1) $( \xcc \xcd )$ $ \xch $ $( \xbm \xcc \xcd )$

Let $f(T) \xcc M(T' ),$ $f(T' ) \xcc M(T).$
So $Th(M(T' )) \xcc Th(f(T)),$ $Th(M(T)) \xcc Th(f(T' )),$ so $T' \xcc
\ol{T' } \xcc \ol{ \ol{T} },$ $T \xcc \ol{T} \xcc \ol{ \ol{T' } },$
so by $( \xcc \xcd )$ $ \ol{ \ol{T} }= \ol{ \ol{T' } },$ so $f(T):=M( \ol{
\ol{T} })=M( \ol{ \ol{T' } })=:f(T' ).$

(10.2) $( \xbm \xcc \xcd )$ $ \xch $ $( \xcc \xcd )$

Let $T \xcc \ol{ \ol{T' } }$ and $T' \xcc \ol{ \ol{T} }.$ So by $ \xCf
(CCL)$ $Th(M(T))= \ol{T} \xcc \ol{ \ol{T' } }=Th(f(T' )).$
But $Th(M(T)) \xcc Th(X) \xch X \xcc M(T):$ $X \xcc M(Th(X)) \xcc
M(Th(M(T)))=M(T).$
So $f(T' ) \xcc M(T),$ likewise $f(T) \xcc M(T' ),$ so by $( \xbm \xcc
\xcd )$ $f(T)=f(T' ),$ so $ \ol{ \ol{T} }= \ol{ \ol{T' } }.$

(11.1) $ \xCf (CUM)$ $ \xch $ $( \xbm CUM):$

So let $X=M(T),$ $Y=M(T' ),$ and $f(T):=M( \ol{ \ol{T} }) \xcc M(T' ) \xcc
M(T)$ $ \xcp $ $ \ol{T} \xcc \ol{T' } \xcc \ol{ \ol{T} }=_{ \xCf (LLE)}
\ol{ \ol{( \ol{T})} }$ $ \xcp $
$ \ol{ \ol{T} }= \ol{ \ol{( \ol{T})} }= \ol{ \ol{( \ol{T' })} }= \ol{
\ol{T' } }$ $ \xcp $ $f(T)=M( \ol{ \ol{T} })=M( \ol{ \ol{T' } })=f(T' ),$
$thus$ $f(X)=f(Y).$

(11.2) $( \xbm CUM)$ $ \xch $ $ \xCf (CUM)$:

Let $T$ $ \xcc $ $ \ol{T' }$ $ \xcc $ $ \ol{ \ol{T} }.$ Thus by $( \xbm
CUM)$ and $f(T) \xcc M( \ol{ \ol{T} })$ $ \xcc $ $M(T' )$ $ \xcc $ $M(T),$
so $f(T)=f(T' ),$
so $ \ol{ \ol{T} }$ $=$ $Th(f(T))$ $=$ $Th(f(T' ))$ $=$ $ \ol{ \ol{T' }
}.$

(12.1) $ \xCf (RatM)$ $ \xch $ $( \xbm RatM)$

Let $X=M(T),$ $Y=M(T' ),$ and $X \xcc Y,$ $X \xcs f(Y) \xEd \xCQ,$ so $T
\xcl T' $ and $M(T) \xcs f(M(T' )) \xEd \xCQ,$
so $Con(T, \ol{ \ol{T' } }),$ so $ \ol{ \ol{ \ol{T' } } \xcv T} \xcc \ol{
\ol{T} }$ by $ \xCf (RatM),$ so $f(X)=f(M(T))=M( \ol{ \ol{T} }) \xcc M(
\ol{ \ol{T' } } \xcv T)=$
$M( \ol{ \ol{T' } }) \xcs M(T)=f(Y) \xcs X.$

(12.2) $( \xbm RatM)+( \xbm dp)$ $ \xch $ $ \xCf (RatM):$

Let $X=M(T),$ $Y=M(T' ),$ $T \xcl T',$ $Con(T, \ol{ \ol{T' } }),$ so $X
\xcc Y$ and by $( \xbm dp)$ $X \xcs f(Y) \xEd \xCQ,$
so by $( \xbm RatM)$ $f(X) \xcc f(Y) \xcs X,$ so
$ \ol{ \ol{T} }= \ol{ \ol{T \xcv T' } }=Th(f(T \xcv T' )) \xcd Th(f(T' )
\xcs M(T))= \ol{ \ol{ \ol{T' } } \xcv T}$ by (a) above and $( \xbm dp).$

(12.3) $( \xbm RatM)$ $ \xcH $ $ \xCf (RatM)$ without $( \xbm dp):$

$( \xbm RatM)$ holds in all ranked preferential structures
(see Definition \ref{Definition Rank-Rel} (page \pageref{Definition Rank-Rel}) )
by Fact \ref{Fact Rank-Hold} (page \pageref{Fact Rank-Hold}) . Example
\ref{Example Rank-Dp} (page \pageref{Example Rank-Dp})  (2)
shows
that $ \xCf (RatM)$ may fail in the resulting logic.

(12.4) $( \xbm RatM)$ $ \xch $ $ \xCf (RatM)$ if $T$ is classically
equivalent to a formula:

$ \xbf \xcl T' $ $ \xch $ $M( \xbf ) \xcc M(T' ).$ $Con( \xbf, \ol{
\ol{T' } })$ $ \xcj $ $M( \ol{ \ol{T' } }) \xcs M( \xbf ) \xEd \xCQ $ $
\xcj $ $f(T' ) \xcs M( \xbf ) \xEd \xCQ $
by Fact \ref{Fact Log-Form} (page \pageref{Fact Log-Form})  (4). Thus $f(M( \xbf
)) \xcc f(M(T' ))
\xcs M( \xbf )$ by $( \xbm RatM).$
Thus by (a) above $ \ol{ \ol{ \ol{T' } } \xcv \{ \xbf \}} \xcc \ol{ \ol{
\xbf } }.$

(13.1) $ \xCf (RatM=)$ $ \xch $ $( \xbm =)$

Let $X=M(T),$ $Y=M(T' ),$ and $X \xcc Y,$ $X \xcs f(Y) \xEd \xCQ,$ so $T
\xcl T' $ and $M(T) \xcs f(M(T' )) \xEd \xCQ,$
so $Con(T, \ol{ \ol{T' } }),$ so $ \ol{ \ol{ \ol{T' } } \xcv T}= \ol{
\ol{T} }$ by $ \xCf (RatM=),$ so $f(X)=f(M(T))=M( \ol{ \ol{T} })=M( \ol{
\ol{T' } } \xcv T)=$
$M( \ol{ \ol{T' } }) \xcs M(T)=f(Y) \xcs X.$

(13.2) $( \xbm =)+( \xbm dp)$ $ \xch $ $ \xCf (RatM=)$

Let $X=M(T),$ $Y=M(T' ),$ $T \xcl T',$ $Con(T, \ol{ \ol{T' } }),$ so $X
\xcc Y$ and by $( \xbm dp)$ $X \xcs f(Y) \xEd \xCQ,$
so by $( \xbm =)$ $f(X)=f(Y) \xcs X.$ So $ \ol{ \ol{ \ol{T' } } \xcv T}=
\ol{ \ol{T} }$ (a) above and $( \xbm dp).$

(13.3) $( \xbm =)$ $ \xcH $ $ \xCf (RatM=)$ without $( \xbm dp):$

$( \xbm =)$ holds in all ranked preferential structures
(see Definition \ref{Definition Rank-Rel} (page \pageref{Definition Rank-Rel}) )
by Fact \ref{Fact Rank-Hold} (page \pageref{Fact Rank-Hold}) . Example
\ref{Example Rank-Dp} (page \pageref{Example Rank-Dp})  (1)
shows
that $ \xCf (RatM=)$ may fail in the resulting logic.

(13.4) $( \xbm =)$ $ \xch $ $ \xCf (RatM=)$ if $T$ is classically
equivalent to a formula:

The proof is almost identical to the one for (12.4). Again, the
prerequisites of $( \xbm =)$ are satisfied, so $f(M( \xbf ))=f(M(T' ))
\xcs M( \xbf ).$
Thus, $ \ol{ \ol{ \ol{T' } } \xcv \{ \xbf \}}= \ol{ \ol{ \xbf } }$ by (a)
above.

Of the last four, we show (14), (15), (17), the proof for (16) is similar
to the
one for (17).

(14.1) $(Log=' )$ $ \xch $ $( \xbm =' ):$

$f(M(T' )) \xcs M(T) \xEd \xCQ $ $ \xch $ $Con( \ol{ \ol{T' } } \xcv T)$ $
\xch_{(Log=' )}$ $ \ol{ \ol{T \xcv T' } }= \ol{ \ol{ \ol{T' } } \xcv T}$ $
\xch $
$f(M(T \xcv T' ))=f(M(T' )) \xcs M(T).$

(14.2) $( \xbm =' )+( \xbm dp)$ $ \xch $ $(Log=' ):$

$Con( \ol{ \ol{T' } } \xcv T)$ $ \xch_{( \xbm dp)}$ $f(M(T' )) \xcs M(T)
\xEd \xCQ $ $ \xch $ $f(M(T' \xcv T))=f(M(T' ) \xcs M(T))$ $=_{( \xbm ='
)}$
$f(M(T' )) \xcs M(T),$ so $ \ol{ \ol{T' \xcv T} }$ $=$ $ \ol{ \ol{ \ol{T'
} } \xcv T}$ by (a) above and $( \xbm dp).$

(14.3) $( \xbm =' )$ $ \xcH $ $(Log=' )$ without $( \xbm dp):$

By Fact \ref{Fact Rank-Hold} (page \pageref{Fact Rank-Hold})  $( \xbm =' )$
holds in ranked
structures.
Consider Example \ref{Example Rank-Dp} (page \pageref{Example Rank-Dp})  (2).
There, $Con(T, \ol{ \ol{T' }
}),$ $T=T \xcv T',$ and
it was shown that $ \ol{ \ol{ \ol{T' } } \xcv T}$ $ \xcC $ $ \ol{ \ol{T}
}$ $=$ $ \ol{ \ol{T \xcv T' } }$

(14.4) $( \xbm =' )$ $ \xch $ $(Log=' )$ if $T$ is classically equivalent
to a formula:

$Con( \ol{ \ol{T' } } \xcv \{ \xbf \})$ $ \xch $ $ \xCQ \xEd M( \ol{
\ol{T' } }) \xcs M( \xbf )$ $ \xch $ $f(T' ) \xcs M( \xbf ) \xEd \xCQ $ by
Fact \ref{Fact Log-Form} (page \pageref{Fact Log-Form})  (4). So $f(M(T' \xcv \{
\xbf \}))$ $=$
$f(M(T' ) \xcs M( \xbf ))$ $=$
$f(M(T' )) \xcs M( \xbf )$ by $( \xbm =' ),$ so $ \ol{ \ol{T' \xcv \{ \xbf
\}} }= \ol{ \ol{ \ol{T' } } \xcv \{ \xbf \}}$ by (a) above.

(15.1) $(Log \xFO )$ $ \xch $ $( \xbm \xFO ):$

Trivial.

(15.2) $( \xbm \xFO )$ $ \xch $ $(Log \xFO ):$

Trivial.

(16) $(Log \xcv )$ $ \xcj $ $( \xbm \xcv ):$ Analogous to the proof of
(17).

(17.1) $(Log \xcv ' )+( \xbm \xcc )+( \xbm =)$ $ \xch $ $( \xbm \xcv ' ):$

$f(M(T' )) \xcs (M(T)-f(M(T))) \xEd \xCQ $ $ \xch $ (by $( \xbm \xcc ),$
$( \xbm =),$
Fact \ref{Fact Rank-Auxil} (page \pageref{Fact Rank-Auxil}) ) $f(M(T' )) \xcs
M(T) \xEd \xCQ,$ $f(M(T'
)) \xcs f(M(T))= \xCQ $ $ \xch $
$Con( \ol{ \ol{T' } },T),$ $ \xCN Con( \ol{ \ol{T' } }, \ol{ \ol{T} })$ $
\xch $ $ \ol{ \ol{T \xco T' } }= \ol{ \ol{T} }$ $ \xch $ $f(M(T))=f(M(T
\xco T' ))=f(M(T) \xcv M(T' )).$

(17.2) $( \xbm \xcv ' )+( \xbm dp)$ $ \xch $ $(Log \xcv ' ):$

$Con( \ol{ \ol{T' } } \xcv T),$ $ \xCN Con( \ol{ \ol{T' } } \xcv \ol{
\ol{T} })$ $ \xch_{( \xbm dp)}$ $f(T' ) \xcs M(T) \xEd \xCQ,$ $f(T' )
\xcs f(T)= \xCQ $ $ \xch $
$f(M(T' )) \xcs (M(T)-f(M(T))) \xEd \xCQ $ $ \xch $ $f(M(T))$ $=$ $f(M(T)
\xcv M(T' ))$ $=$ $f(M(T \xco T' )).$
So $ \ol{ \ol{T} }= \ol{ \ol{T \xco T' } }.$

(17.3) and (16.3) are solved by Example \ref{Example Rank-Dp} (page
\pageref{Example Rank-Dp})  (3).

$ \xcz $
\\[3ex]

 karl-search= End Proposition Alg-Log Proof
\vspace{7mm}

 *************************************

\vspace{7mm}

\subsubsection{Example Rank-Dp}

 {\LARGE karl-search= Start Example Rank-Dp }

\index{Example Rank-Dp}

\be

$\hspace{0.01em}$

(+++ Orig. No.:  Example Rank-Dp +++)

{\xssc LABEL: {Example Rank-Dp}}
\label{Example Rank-Dp}

(1) $( \xbm =)$ without $( \xbm dp)$ does not imply $(RatM=):$

Take $\{p_{i}:i \xbe \xbo \}$ and put $m:=m_{ \xcU p_{i}},$ the model
which makes all $p_{i}$ true, in the top
layer, all the other in the bottom layer. Let $m' \xEd m,$ $T':= \xCQ,$
$T:=Th(m,m' ).$ Then
Then $ \ol{ \ol{T' } }=T',$ so $Con( \ol{ \ol{T' } },T),$ $ \ol{ \ol{T}
}=Th(m' ),$ $ \ol{ \ol{ \ol{T' } } \xcv T}=T.$

So $(RatM=)$ fails, but $( \xbm =)$ holds in all ranked structures.

(2) $( \xbm RatM)$ without $( \xbm dp)$ does not imply (RatM):

Take $\{p_{i}:i \xbe \xbo \}$ and let $m:=m_{ \xcU p_{i}},$ the model
which makes all $p_{i}$ true.

Let $X:=M( \xCN p_{0}) \xcv \{m\}$ be the top layer, put the rest of $M_{
\xdl }$ in the bottom layer.
Let $Y:=M_{ \xdl }.$ The structure is ranked, as shown in Fact \ref{Fact
Rank-Hold} (page \pageref{Fact Rank-Hold}) ,
$( \xbm RatM)$ holds.

Let $T':= \xCQ,$ $T:=Th(X).$ We have to show that $Con(T, \ol{ \ol{T' }
}),$ $T \xcl T',$ but
$ \ol{ \ol{ \ol{T' } } \xcv T} \xcC \ol{ \ol{T} }.$ $ \ol{ \ol{T' } }$ $=$
$Th(M(p_{0})-\{m\})$ $=$ $ \ol{p_{0}}.$ $T$ $=$ $ \ol{\{ \xCN p_{0}\} \xco
Th(m)},$ $ \ol{ \ol{T} }=T.$ $So$ $Con(T, \ol{
\ol{T' } }).$
$M( \ol{ \ol{T' } })=M(p_{0}),$ $M(T)=X,$ $M( \ol{ \ol{T' } } \xcv T)=M(
\ol{ \ol{T' } }) \xcs M(T)=\{m\},$ $m \xcm p_{1},$ so $p_{1} \xbe \ol{
\ol{ \ol{T' } } \xcv T},$ but
$X \xcM p_{1}.$

(3) This example shows that we need $( \xbm dp)$ to go from $( \xbm \xcv
)$ to
$(Log \xcv )$ and from $( \xbm \xcv ' )$ to $(Log \xcv ' ).$

Let $v( \xdl ):=\{p,q\} \xcv \{p_{i}:i< \xbo \}.$ Let $m$ make all
variables true.

Put all models of $ \xCN p,$ and $m,$ in the upper layer, all other models
in the
lower layer. This is ranked, so by Fact \ref{Fact Rank-Hold} (page \pageref{Fact
Rank-Hold})  $( \xbm
\xcv )$
and $( \xbm \xcv ' )$ hold.
Set $X:=M( \xCN q) \xcv \{m\},$ $X':=M(q)-\{m\},$ $T:=Th(X)= \xCN q \xco
Th(m),$ $T':=Th(X' )= \ol{q}.$
Then $ \ol{ \ol{T} }= \ol{p \xcu \xCN q},$ $ \ol{ \ol{T' } }= \ol{p \xcu
q}.$ We have $Con( \ol{ \ol{T' } },T),$ $ \xCN Con( \ol{ \ol{T' } }, \ol{
\ol{T} }).$
But $ \ol{ \ol{T \xco T' } }= \ol{p} \xEd \ol{ \ol{T} }= \ol{p \xcu \xCN
q}$ and $Con( \ol{ \ol{T \xco T' } },T' ),$ so $(Log \xcv )$ and $(Log
\xcv ' )$ fail.

$ \xcz $
\\[3ex]

 karl-search= End Example Rank-Dp
\vspace{7mm}

 *************************************

\vspace{7mm}

\subsubsection{Fact Cut-Pr}

 {\LARGE karl-search= Start Fact Cut-Pr }

\index{Fact Cut-Pr}

\ee

\bfa

$\hspace{0.01em}$

(+++ Orig. No.:  Fact Cut-Pr +++)

{\xssc LABEL: {Fact Cut-Pr}}
\label{Fact Cut-Pr}

$ \xCf (CUT)$ $ \xcH $ $ \xCf (PR)$

 karl-search= End Fact Cut-Pr
\vspace{7mm}

 *************************************

\vspace{7mm}

\subsubsection{Fact Cut-Pr Proof}

 {\LARGE karl-search= Start Fact Cut-Pr Proof }

\index{Fact Cut-Pr Proof}

\efa

\subparagraph{
Proof
}

$\hspace{0.01em}$

(+++*** Orig.:  Proof )

We give two proofs:

(1) If $ \xCf (CUT)$ $ \xch $ $ \xCf (PR),$ then by $( \xbm PR)$ $ \xch $
(by Fact \ref{Fact Mu-Base} (page \pageref{Fact Mu-Base})  (3))
$( \xbm CUT)$ $ \xch $ (by Proposition \ref{Proposition Alg-Log} (page
\pageref{Proposition Alg-Log})  (7.2) $
\xCf (CUT)$ $ \xch $ $ \xCf (PR)$
we would have a proof of $( \xbm PR)$ $ \xch $ $ \xCf (PR)$ without $(
\xbm dp),$ which is impossible,
as shown by Example \ref{Example Pref-Dp} (page \pageref{Example Pref-Dp}) .

(2) Reconsider Example \ref{Example Need-Pr} (page \pageref{Example Need-Pr}) ,
and say $a \xcm p \xcu q,$
$b \xcm p \xcu \xCN q,$ $c \xcm \xCN p \xcu q.$
It is shown there that $( \xbm CUM)$ holds, so $( \xbm CUT)$ holds, so by
Proposition \ref{Proposition Alg-Log} (page \pageref{Proposition Alg-Log}) 
(7.2) $ \xCf (CUT)$ holds, if we
define $ \ol{ \ol{T} }:=Th(f(M(T)).$
Set $T:=\{p \xco ( \xCN p \xcu q)\},$ $T':=\{p\},$ then
$ \ol{ \ol{T \xcv T' } }= \ol{ \ol{T' } }= \ol{\{p \xcu \xCN q\}},$ $ \ol{
\ol{T} }= \ol{T},$ $ \ol{T \xcv T' }= \ol{T' }= \ol{\{p\}},$ $so$ $ \xCf
(PR)$ $fails.$

$ \xcz $
\\[3ex]

 karl-search= End Fact Cut-Pr Proof
\vspace{7mm}

 *************************************

\vspace{7mm}

 karl-search= End ToolBase1-Log-Rules
\vspace{7mm}

 *************************************

\vspace{7mm}

 karl-search= End ToolBase1-Log
\vspace{7mm}

 *************************************

\vspace{7mm}

\newpage

\section{
General and smooth preferential structures
}
{\xssc LABEL: {General and smooth preferential structures}}
\label{General and smooth preferential structures}

\subsubsection{ToolBase1-Pref}

 {\LARGE karl-search= Start ToolBase1-Pref }

{\xssc LABEL: {Section Toolbase1-Pref}}
\label{Section Toolbase1-Pref}
\index{Section Toolbase1-Pref}
\subsubsection{ToolBase1-Pref-Base}

 {\LARGE karl-search= Start ToolBase1-Pref-Base }

{\xssc LABEL: {Section Toolbase1-Pref-Base}}
\label{Section Toolbase1-Pref-Base}
\index{Section Toolbase1-Pref-Base}

\subsection{
General and smooth: Basics
}
{\xssc LABEL: {General and smooth: Basics}}
\label{General and smooth: Basics}

\subsubsection{Definition Pref-Str}

 {\LARGE karl-search= Start Definition Pref-Str }

\index{Definition Pref-Str}

\bd

$\hspace{0.01em}$

(+++ Orig. No.:  Definition Pref-Str +++)

{\xssc LABEL: {Definition Pref-Str}}
\label{Definition Pref-Str}

Fix $U \xEd \xCQ,$ and consider arbitrary $X.$
Note that this $X$ has not necessarily anything to do with $U,$ or $ \xdu
$ below.
Thus, the functions $ \xbm_{ \xdm }$ below are in principle functions from
$V$ to $V$ - where $V$
is the set theoretical universe we work in.

(A) Preferential models or structures.

(1) The version without copies:

A pair $ \xdm:=<U, \xeb >$ with $U$ an arbitrary set, and $ \xeb $ an
arbitrary binary relation
on $U$ is called a preferential model or structure.

(2) The version with copies:

A pair $ \xdm:=< \xdu, \xeb >$ with $ \xdu $ an arbitrary set of pairs,
and $ \xeb $ an arbitrary binary
relation on $ \xdu $ is called a preferential model or structure.

If $<x,i> \xbe \xdu,$ then $x$ is intended to be an element of $U,$ and
$i$ the index of the
copy.

We sometimes also need copies of the relation $ \xeb,$ we will then
replace $ \xeb $
by one or several arrows $ \xba $ attacking non-minimal elements, e.g. $x
\xeb y$ will
be written $ \xba:x \xcp y,$ $<x,i> \xeb <y,i>$ will be written $ \xba
:<x,i> \xcp <y,i>,$ and
finally we might have $< \xba,k>:x \xcp y$ and $< \xba,k>:<x,i> \xcp
<y,i>,$ etc.

(B) Minimal elements, the functions $ \xbm_{ \xdm }$

(1) The version without copies:

Let $ \xdm:=<U, \xeb >,$ and define

$ \xbm_{ \xdm }(X)$ $:=$ $\{x \xbe X:$ $x \xbe U$ $ \xcu $ $ \xCN \xcE x'
\xbe X \xcs U.x' \xeb x\}.$

$ \xbm_{ \xdm }(X)$ is called the set of minimal elements of $X$ (in $
\xdm ).$

Thus, $ \xbm_{ \xdm }(X)$ is the set of elements such that there is no
smaller one
in $X.$

(2) The version with copies:

Let $ \xdm:=< \xdu, \xeb >$ be as above. Define

$ \xbm_{ \xdm }(X)$ $:=$ $\{x \xbe X:$ $ \xcE <x,i> \xbe \xdu. \xCN \xcE
<x',i' > \xbe \xdu (x' \xbe X$ $ \xcu $ $<x',i' >' \xeb <x,i>)\}.$

Thus, $ \xbm_{ \xdm }(X)$ is the projection on the first coordinate of the
set of elements
such that there is no smaller one in $X.$

Again, by abuse of language, we say that $ \xbm_{ \xdm }(X)$ is the set of
minimal elements
of $X$ in the structure. If the context is clear, we will also write just
$ \xbm.$

We sometimes say that $<x,i>$ ``kills'' or ``minimizes'' $<y,j>$ if
$<x,i> \xeb <y,j>.$ By abuse of language we also say a set $X$ kills or
minimizes a set
$Y$ if for all $<y,j> \xbe \xdu,$ $y \xbe Y$ there is $<x,i> \xbe \xdu,$
$x \xbe X$ s.t. $<x,i> \xeb <y,j>.$

$ \xdm $ is also called injective or 1-copy, iff there is always at most
one copy
$<x,i>$ for each $x.$ Note that the existence of copies corresponds to a
non-injective labelling function - as is often used in nonclassical
logic, e.g. modal logic.

We say that $ \xdm $ is transitive, irreflexive, etc., iff $ \xeb $ is.

Note that $ \xbm (X)$ might well be empty, even if $X$ is not.

 karl-search= End Definition Pref-Str
\vspace{7mm}

 *************************************

\vspace{7mm}

\subsubsection{Definition Pref-Log}

 {\LARGE karl-search= Start Definition Pref-Log }

\index{Definition Pref-Log}

\ed

\bd

$\hspace{0.01em}$

(+++ Orig. No.:  Definition Pref-Log +++)

{\xssc LABEL: {Definition Pref-Log}}
\label{Definition Pref-Log}

We define the consequence relation of a preferential structure for a
given propositional language $ \xdl.$

(A)

(1) If $m$ is a classical model of a language $ \xdl,$ we say by abuse
of language

$<m,i> \xcm \xbf $ iff $m \xcm \xbf,$

and if $X$ is a set of such pairs, that

$X \xcm \xbf $ iff for all $<m,i> \xbe X$ $m \xcm \xbf.$

(2) If $ \xdm $ is a preferential structure, and $X$ is a set of $ \xdl
-$models for a
classical propositional language $ \xdl,$ or a set of pairs $<m,i>,$
where the $m$ are
such models, we call $ \xdm $ a classical preferential structure or model.

(B)

Validity in a preferential structure, or the semantical consequence
relation
defined by such a structure:

Let $ \xdm $ be as above.

We define:

$T \xcm_{ \xdm } \xbf $ iff $ \xbm_{ \xdm }(M(T)) \xcm \xbf,$ i.e. $
\xbm_{ \xdm }(M(T)) \xcc M( \xbf ).$

$ \xdm $ will be called definability preserving iff for all $X \xbe \xdD_{
\xdl }$ $ \xbm_{ \xdm }(X) \xbe \xdD_{ \xdl }.$

As $ \xbm_{ \xdm }$ is defined on $ \xdD_{ \xdl },$ but need by no means
always result in some new
definable set, this is (and reveals itself as a quite strong) additional
property.

 karl-search= End Definition Pref-Log
\vspace{7mm}

 *************************************

\vspace{7mm}

\subsubsection{Example NeedCopies}

 {\LARGE karl-search= Start Example NeedCopies }

\index{Example NeedCopies}

\ed

\be

$\hspace{0.01em}$

(+++ Orig. No.:  Example NeedCopies +++)

{\xssc LABEL: {Example NeedCopies}}
\label{Example NeedCopies}

This simple example illustrates the
importance of copies. Such examples seem to have appeared for the first
time
in print in  \cite{KLM90}, but can probably be
attibuted to folklore.

Consider the propositional language $ \xdl $ of two propositional
variables $p,q$, and
the classical preferential model $ \xdm $ defined by

$m \xcm p \xcu q,$ $m' \xcm p \xcu q,$ $m_{2} \xcm \xCN p \xcu q,$ $m_{3}
\xcm \xCN p \xcu \xCN q,$ with $m_{2} \xeb m$, $m_{3} \xeb m' $, and
let $ \xcm_{ \xdm }$ be its consequence relation. (m and $m' $ are
logically identical.)

Obviously, $Th(m) \xco \{ \xCN p\} \xcm_{ \xdm } \xCN p$, but there is no
complete theory $T' $ s.t.
$Th(m) \xco T' \xcm_{ \xdm } \xCN p$. (If there were one, $T' $ would
correspond to $m$, $m_{2},$ $m_{3},$
or the missing $m_{4} \xcm p \xcu \xCN q$, but we need two models to kill
all copies of $m.)$
On the other hand, if there were just one copy of $m,$ then one other
model,
i.e. a complete theory would suffice. More formally, if we admit at most
one
copy of each model in a structure $ \xdm,$ $m \xcM T,$ and $Th(m) \xco T
\xcm_{ \xdm } \xbf $ for some $ \xbf $ s.t.
$m \xcm \xCN \xbf $ - i.e. $m$ is not minimal in the models of $Th(m) \xco
T$ - then there is a
complete $T' $ with $T' \xcl T$ and $Th(m) \xco T' \xcm_{ \xdm } \xbf $,
i.e. there is $m'' $ with $m'' \xcm T' $ and
$m'' \xeb m.$ $ \xcz $
\\[3ex]

 karl-search= End Example NeedCopies
\vspace{7mm}

 *************************************

\vspace{7mm}

\subsubsection{Definition Smooth}

 {\LARGE karl-search= Start Definition Smooth }

\index{Definition Smooth}

\ee

\bd

$\hspace{0.01em}$

(+++ Orig. No.:  Definition Smooth +++)

{\xssc LABEL: {Definition Smooth}}
\label{Definition Smooth}

Let $ \xdy \xcc \xdp (U).$ (In applications to logic, $ \xdy $ will be $
\xdD_{ \xdl }.)$

A preferential structure $ \xdm $ is called $ \xdy -$smooth iff for every
$X \xbe \xdy $ every
element
$x \xbe X$ is either minimal in $X$ or above an element, which is minimal
in $X.$ More
precisely:

(1) The version without copies:

If $x \xbe X \xbe \xdy,$ then either $x \xbe \xbm (X)$ or there is $x'
\xbe \xbm (X).x' \xeb x.$

(2) The version with copies:

If $x \xbe X \xbe \xdy,$ and $<x,i> \xbe \xdu,$ then either there is no
$<x',i' > \xbe \xdu,$ $x' \xbe X,$
$<x',i' > \xeb <x,i>$ or there is $<x',i' > \xbe \xdu,$ $<x',i' > \xeb
<x,i>,$ $x' \xbe X,$ s.t. there is
no $<x'',i'' > \xbe \xdu,$ $x'' \xbe X,$ with $<x'',i'' > \xeb <x',i'
>.$

When considering the models of a language $ \xdl,$ $ \xdm $ will be
called smooth iff
it is $ \xdD_{ \xdl }-$smooth; $ \xdD_{ \xdl }$ is the default.

Obviously, the richer the set $ \xdy $ is, the stronger the condition $
\xdy -$smoothness
will be.

 karl-search= End Definition Smooth
\vspace{7mm}

 *************************************

\vspace{7mm}

\subsubsection{ToolBase1-Rank-Base}

 {\LARGE karl-search= Start ToolBase1-Rank-Base }

{\xssc LABEL: {Section Toolbase1-Rank-Base}}
\label{Section Toolbase1-Rank-Base}
\index{Section Toolbase1-Rank-Base}
\subsubsection{Fact Rank-Base}

 {\LARGE karl-search= Start Fact Rank-Base }

\index{Fact Rank-Base}

\ed

\bfa

$\hspace{0.01em}$

(+++ Orig. No.:  Fact Rank-Base +++)

{\xssc LABEL: {Fact Rank-Base}}
\label{Fact Rank-Base}

Let $ \xeb $ be an irreflexive, binary relation on $X,$ then the following
two conditions
are equivalent:

(1) There is $ \xbO $ and an irreflexive, total, binary relation $ \xeb '
$ on $ \xbO $ and a
function $f:X \xcp \xbO $ s.t. $x \xeb y$ $ \xcr $ $f(x) \xeb ' f(y)$ for
all $x,y \xbe X.$

(2) Let $x,y,z \xbe X$ and $x \xcT y$ wrt. $ \xeb $ (i.e. neither $x \xeb
y$ nor $y \xeb x),$ then $z \xeb x$ $ \xcp $ $z \xeb y$
and $x \xeb z$ $ \xcp $ $y \xeb z.$

$ \xcz $
\\[3ex]

 karl-search= End Fact Rank-Base
\vspace{7mm}

 *************************************

\vspace{7mm}

\subsubsection{Definition Rank-Rel}

 {\LARGE karl-search= Start Definition Rank-Rel }

\index{Definition Rank-Rel}

\efa

\bd

$\hspace{0.01em}$

(+++ Orig. No.:  Definition Rank-Rel +++)

{\xssc LABEL: {Definition Rank-Rel}}
\label{Definition Rank-Rel}

We call an irreflexive, binary relation $ \xeb $ on $X,$ which satisfies
(1)
(equivalently (2)) of Fact \ref{Fact Rank-Base} (page \pageref{Fact Rank-Base})
, ranked.
By abuse of language, we also call a preferential structure $<X, \xeb >$
ranked, iff
$ \xeb $ is.

 karl-search= End Definition Rank-Rel
\vspace{7mm}

 *************************************

\vspace{7mm}

 karl-search= End ToolBase1-Pref-Base
\vspace{7mm}

 *************************************

\vspace{7mm}

\newpage

\subsection{
General and smooth: Summary
}
\subsubsection{ToolBase1-Pref-ReprSumm}

 {\LARGE karl-search= Start ToolBase1-Pref-ReprSumm }

{\xssc LABEL: {Section Toolbase1-Pref-ReprSumm}}
\label{Section Toolbase1-Pref-ReprSumm}
\index{Section Toolbase1-Pref-ReprSumm}
\subsubsection{Proposition Pref-Representation-With-Ref}

 {\LARGE karl-search= Start Proposition Pref-Representation-With-Ref }

\index{Proposition Pref-Representation-With-Ref}

\ed

The following table summarizes representation by
preferential structures. The positive implications on the right are shown
in
Proposition \ref{Proposition Alg-Log} (page \pageref{Proposition Alg-Log}) 
(going via the $ \xbm -$functions),
those on the left
are shown in the respective representation theorems.

``singletons'' means that the domain must contain all singletons, ``1 copy''
or `` $ \xcg 1$ copy'' means that the structure may contain only 1 copy for
each point,
or several, `` $( \xbm \xCQ )$ '' etc. for the preferential structure mean
that the
$ \xbm -$function of the structure has to satisfy this property.
{\xssc LABEL: {Proposition Pref-Representation-With-Ref}}
\label{Proposition Pref-Representation-With-Ref}

{\scriptsize

\begin{tabular}{|c|c|c|c|c|}

\hline

$\xbm-$ function
\xEH
\xEH
Pref.Structure
\xEH
\xEH
Logic
\xEP

\hline

$(\xbm \xcc)+(\xbm PR)$
\xEH
$\xci$
\xEH
general
\xEH
$\xch$ $(\xbm dp)$
\xEH
$(LLE)+(RW)+$
\xEP

\xEH
Fact \ref{Fact Pref-Sound}
\xEH
\xEH
\xEH
$(SC)+(PR)$
\xEP

\xEH
page \pageref{Fact Pref-Sound}
\xEH
\xEH
\xEH
\xEP

\cline{2-2}
\cline{4-4}

\xEH
$\xch$
\xEH
\xEH
$\xci$
\xEH
\xEP

\xEH
Proposition \ref{Proposition Pref-Complete}
\xEH
\xEH
\xEH
\xEP

\xEH
page \pageref{Proposition Pref-Complete}
\xEH
\xEH
\xEH
\xEP

\cline{2-2}
\cline{4-4}

\xEH
\xEH
\xEH
$\xcH$ without $(\xbm dp)$
\xEH
\xEP

\xEH
\xEH
\xEH
Example \ref{Example Pref-Dp}
\xEH
\xEP

\xEH
\xEH
\xEH
page \pageref{Example Pref-Dp}
\xEH
\xEP

\cline{4-5}

\xEH
\xEH
\xEH
$\xcJ$ without $(\xbm dp)$
\xEH
any ``normal''
\xEP

\xEH
\xEH
\xEH
Proposition \ref{Proposition No-Norm} (1)
\xEH
characterization
\xEP

\xEH
\xEH
\xEH
page \pageref{Proposition No-Norm}
\xEH
of any size
\xEP

\hline

$(\xbm \xcc)+(\xbm PR)$
\xEH
$\xci$
\xEH
transitive
\xEH
$\xch$ $(\xbm dp)$
\xEH
$(LLE)+(RW)+$
\xEP

\xEH
Fact \ref{Fact Pref-Sound}
\xEH
\xEH
\xEH
$(SC)+(PR)$
\xEP

\xEH
page \pageref{Fact Pref-Sound}
\xEH
\xEH
\xEH
\xEP

\cline{2-2}
\cline{4-4}

\xEH
$\xch$
\xEH
\xEH
$\xci$
\xEH
\xEP

\xEH
Proposition \ref{Proposition Pref-Complete-Trans}
\xEH
\xEH
\xEH
\xEP

\xEH
page \pageref{Proposition Pref-Complete-Trans}
\xEH
\xEH
\xEH
\xEP

\cline{2-2}
\cline{4-4}

\xEH
\xEH
\xEH
$\xcH$ without $(\xbm dp)$
\xEH
\xEP

\xEH
\xEH
\xEH
Example \ref{Example Pref-Dp}
\xEH
\xEP

\xEH
\xEH
\xEH
page \pageref{Example Pref-Dp}
\xEH
\xEP

\cline{4-5}

\xEH
\xEH
\xEH
$\xcj$ without $(\xbm dp)$
\xEH
using ``small''
\xEP

\xEH
\xEH
\xEH
See \cite{Sch04}
\xEH
exception sets
\xEP

\hline

$(\xbm \xcc)+(\xbm PR)+(\xbm CUM)$
\xEH
$\xci$
\xEH
smooth
\xEH
$\xch$ $(\xbm dp)$
\xEH
$(LLE)+(RW)+$
\xEP

\xEH
Fact \ref{Fact Smooth-Sound}
\xEH
\xEH
\xEH
$(SC)+(PR)+$
\xEP

\xEH
page \pageref{Fact Smooth-Sound}
\xEH
\xEH
\xEH
$(CUM)$
\xEP

\cline{2-2}
\cline{4-4}

\xEH
$\xch$ $(\xcv)$
\xEH
\xEH
$\xci$ $(\xcv)$
\xEH
\xEP

\xEH
Proposition \ref{Proposition Smooth-Complete}
\xEH
\xEH
\xEH
\xEP

\xEH
page \pageref{Proposition Smooth-Complete}
\xEH
\xEH
\xEH
\xEP

\cline{2-2}
\cline{4-4}

\xEH
$\xcH$ without $(\xcv)$
\xEH
\xEH
$\xcH$ without $(\xbm dp)$
\xEH
\xEP

\xEH
See \cite{Sch04}
\xEH
\xEH
Example \ref{Example Pref-Dp}
\xEH
\xEP

\xEH
\xEH
\xEH
page \pageref{Example Pref-Dp}
\xEH
\xEP

\hline

$(\xbm \xcc)+(\xbm PR)+(\xbm CUM)$
\xEH
$\xci$
\xEH
smooth+transitive
\xEH
$\xch$ $(\xbm dp)$
\xEH
$(LLE)+(RW)+$
\xEP

\xEH
Fact \ref{Fact Smooth-Sound}
\xEH
\xEH
\xEH
$(SC)+(PR)+$
\xEP

\xEH
page \pageref{Fact Smooth-Sound}
\xEH
\xEH
\xEH
$(CUM)$
\xEP

\cline{2-2}
\cline{4-4}

\xEH
$\xch$ $(\xcv)$
\xEH
\xEH
$\xci$ $(\xcv)$
\xEH
\xEP

\xEH
Proposition \ref{Proposition Smooth-Complete-Trans}
\xEH
\xEH
\xEH
\xEP

\xEH
page \pageref{Proposition Smooth-Complete-Trans}
\xEH
\xEH
\xEH
\xEP

\cline{2-2}
\cline{4-4}

\xEH
\xEH
\xEH
$\xcH$ without $(\xbm dp)$
\xEH
\xEP

\xEH
\xEH
\xEH
Example \ref{Example Pref-Dp}
\xEH
\xEP

\xEH
\xEH
\xEH
page \pageref{Example Pref-Dp}
\xEH
\xEP

\cline{4-5}

\xEH
\xEH
\xEH
$\xcj$ without $(\xbm dp)$
\xEH
using ``small''
\xEP

\xEH
\xEH
\xEH
See \cite{Sch04}
\xEH
exception sets
\xEP

\hline

$(\xbm\xcc)+(\xbm=)+(\xbm PR)+$
\xEH
$\xci$
\xEH
ranked, $\xcg 1$ copy
\xEH
\xEH
\xEP

$(\xbm=')+(\xbm\xFO)+(\xbm\xcv)+$
\xEH
Fact \ref{Fact Rank-Hold}
\xEH
\xEH
\xEH
\xEP

$(\xbm\xcv ')+(\xbm\xbe)+(\xbm RatM)$
\xEH
page \pageref{Fact Rank-Hold}
\xEH
\xEH
\xEH
\xEP

\hline

$(\xbm\xcc)+(\xbm=)+(\xbm PR)+$
\xEH
$\xcH$
\xEH
ranked
\xEH
\xEH
\xEP

$(\xbm\xcv)+(\xbm\xbe)$
\xEH
Example \ref{Example Rank-Copies}
\xEH
\xEH
\xEH
\xEP

\xEH
page \pageref{Example Rank-Copies}
\xEH
\xEH
\xEH
\xEP

\hline

$(\xbm\xcc)+(\xbm=)+(\xbm \xCQ)$
\xEH
$\xcj$, $(\xcv)$
\xEH
ranked,
\xEH
\xEH
\xEP

\xEH
Proposition \ref{Proposition Rank-Rep1} (1)
\xEH
1 copy + $(\xbm \xCQ)$
\xEH
\xEH
\xEP

\xEH
page \pageref{Proposition Rank-Rep1}
\xEH
\xEH
\xEH
\xEP

\hline

$(\xbm\xcc)+(\xbm=)+(\xbm \xCQ)$
\xEH
$\xcj$, $(\xcv)$
\xEH
ranked, smooth,
\xEH
\xEH
\xEP

\xEH
Proposition \ref{Proposition Rank-Rep1} (1)
\xEH
1 copy + $(\xbm \xCQ)$
\xEH
\xEH
\xEP

\xEH
page \pageref{Proposition Rank-Rep1}
\xEH
\xEH
\xEH
\xEP

\hline

$(\xbm\xcc)+(\xbm=)+(\xbm \xCQ fin)+$
\xEH
$\xcj$, $(\xcv)$, singletons
\xEH
ranked, smooth,
\xEH
\xEH
\xEP

$(\xbm\xbe)$
\xEH
Proposition \ref{Proposition Rank-Rep1} (2)
\xEH
$\xcg$ 1 copy + $(\xbm \xCQ fin)$
\xEH
\xEH
\xEP

\xEH
page \pageref{Proposition Rank-Rep1}
\xEH
\xEH
\xEH
\xEP

\hline

$(\xbm\xcc)+(\xbm PR)+(\xbm \xFO)+$
\xEH
$\xcj$, $(\xcv)$, singletons
\xEH
ranked
\xEH
$\xcH$ without $(\xbm dp)$
\xEH
$(RatM), (RatM=)$,
\xEP

$(\xbm \xcv)+(\xbm\xbe)$
\xEH
Proposition \ref{Proposition Rank-Rep2}
\xEH
$\xcg$ 1 copy
\xEH
Example \ref{Example Rank-Dp}
\xEH
$(Log\xcv), (Log\xcv ')$
\xEP

\xEH
page \pageref{Proposition Rank-Rep2}
\xEH
\xEH
page \pageref{Example Rank-Dp}
\xEH
\xEP

\cline{4-5}

\xEH
\xEH
\xEH
$\xcJ$ without $(\xbm dp)$
\xEH
any ``normal''
\xEP

\xEH
\xEH
\xEH
Proposition \ref{Proposition No-Norm} (2)
\xEH
characterization
\xEP

\xEH
\xEH
\xEH
page \pageref{Proposition No-Norm}
\xEH
of any size
\xEP

\hline

\end{tabular}

}

 karl-search= End Proposition Pref-Representation-With-Ref
\vspace{7mm}

 *************************************

\vspace{7mm}

\subsubsection{Proposition Pref-Representation-Without-Ref}

 {\LARGE karl-search= Start Proposition Pref-Representation-Without-Ref }

\index{Proposition Pref-Representation-Without-Ref}

The following table summarizes representation by
preferential structures.

``singletons'' means that the domain must contain all singletons, ``1 copy''
or $'' \xcg 1$ copy" means that the structure may contain only 1 copy for
each point,
or several, $'' ( \xbm \xCQ )'' $ etc. for the preferential structure mean
that the
$ \xbm -$function of the structure has to satisfy this property.
{\xssc LABEL: {Proposition Pref-Representation-Without-Ref}}
\label{Proposition Pref-Representation-Without-Ref}

{\scriptsize

\begin{tabular}{|c|c|c|c|c|}

\hline

$\xbm-$ function
\xEH
\xEH
Pref.Structure
\xEH
\xEH
Logic
\xEP

\hline

$(\xbm \xcc)+(\xbm PR)$
\xEH
$\xci$
\xEH
general
\xEH
$\xch$ $(\xbm dp)$
\xEH
$(LLE)+(RW)+$
\xEP

\xEH
\xEH
\xEH
\xEH
$(SC)+(PR)$
\xEP

\cline{2-2}
\cline{4-4}

\xEH
$\xch$
\xEH
\xEH
$\xci$
\xEH
\xEP

\cline{2-2}
\cline{4-4}

\xEH
\xEH
\xEH
$\xcH$ without $(\xbm dp)$
\xEH
\xEP

\cline{4-5}

\xEH
\xEH
\xEH
$\xcJ$ without $(\xbm dp)$
\xEH
any ``normal''
\xEP

\xEH
\xEH
\xEH
\xEH
characterization
\xEP

\xEH
\xEH
\xEH
\xEH
of any size
\xEP

\hline

$(\xbm \xcc)+(\xbm PR)$
\xEH
$\xci$
\xEH
transitive
\xEH
$\xch$ $(\xbm dp)$
\xEH
$(LLE)+(RW)+$
\xEP

\xEH
\xEH
\xEH
\xEH
$(SC)+(PR)$
\xEP

\cline{2-2}
\cline{4-4}

\xEH
$\xch$
\xEH
\xEH
$\xci$
\xEH
\xEP

\cline{2-2}
\cline{4-4}

\xEH
\xEH
\xEH
$\xcH$ without $(\xbm dp)$
\xEH
\xEP

\cline{4-5}

\xEH
\xEH
\xEH
$\xcj$ without $(\xbm dp)$
\xEH
using ``small''
\xEP

\xEH
\xEH
\xEH
\xEH
exception sets
\xEP

\hline

$(\xbm \xcc)+(\xbm PR)+(\xbm CUM)$
\xEH
$\xci$
\xEH
smooth
\xEH
$\xch$ $(\xbm dp)$
\xEH
$(LLE)+(RW)+$
\xEP

\xEH
\xEH
\xEH
\xEH
$(SC)+(PR)+$
\xEP

\xEH
\xEH
\xEH
\xEH
$(CUM)$
\xEP

\cline{2-2}
\cline{4-4}

\xEH
$\xch$ $(\xcv)$
\xEH
\xEH
$\xci$ $(\xcv)$
\xEH
\xEP

\cline{2-2}
\cline{4-4}

\xEH
\xEH
\xEH
$\xcH$ without $(\xbm dp)$
\xEH
\xEP

\hline

$(\xbm \xcc)+(\xbm PR)+(\xbm CUM)$
\xEH
$\xci$
\xEH
smooth+transitive
\xEH
$\xch$ $(\xbm dp)$
\xEH
$(LLE)+(RW)+$
\xEP

\xEH
\xEH
\xEH
\xEH
$(SC)+(PR)+$
\xEP

\xEH
\xEH
\xEH
\xEH
$(CUM)$
\xEP

\cline{2-2}
\cline{4-4}

\xEH
$\xch$ $(\xcv)$
\xEH
\xEH
$\xci$ $(\xcv)$
\xEH
\xEP

\cline{2-2}
\cline{4-4}

\xEH
\xEH
\xEH
$\xcH$ without $(\xbm dp)$
\xEH
\xEP

\cline{4-5}

\xEH
\xEH
\xEH
$\xcj$ without $(\xbm dp)$
\xEH
using ``small''
\xEP

\xEH
\xEH
\xEH
\xEH
exception sets
\xEP

\hline

$(\xbm\xcc)+(\xbm=)+(\xbm PR)+$
\xEH
$\xci$
\xEH
ranked, $\xcg 1$ copy
\xEH
\xEH
\xEP

$(\xbm=')+(\xbm\xFO)+(\xbm\xcv)+$
\xEH
\xEH
\xEH
\xEH
\xEP

$(\xbm\xcv ')+(\xbm\xbe)+(\xbm RatM)$
\xEH
\xEH
\xEH
\xEH
\xEP

\hline

$(\xbm\xcc)+(\xbm=)+(\xbm PR)+$
\xEH
$\xcH$
\xEH
ranked
\xEH
\xEH
\xEP

$(\xbm\xcv)+(\xbm\xbe)$
\xEH
\xEH
\xEH
\xEH
\xEP

\hline

$(\xbm\xcc)+(\xbm=)+(\xbm \xCQ)$
\xEH
$\xcj$, $(\xcv)$
\xEH
ranked,
\xEH
\xEH
\xEP

\xEH
\xEH
1 copy + $(\xbm \xCQ)$
\xEH
\xEH
\xEP

\hline

$(\xbm\xcc)+(\xbm=)+(\xbm \xCQ)$
\xEH
$\xcj$, $(\xcv)$
\xEH
ranked, smooth,
\xEH
\xEH
\xEP

\xEH
\xEH
1 copy + $(\xbm \xCQ)$
\xEH
\xEH
\xEP

\hline

$(\xbm\xcc)+(\xbm=)+(\xbm \xCQ fin)+$
\xEH
$\xcj$, $(\xcv)$, singletons
\xEH
ranked, smooth,
\xEH
\xEH
\xEP

$(\xbm\xbe)$
\xEH
\xEH
$\xcg$ 1 copy + $(\xbm \xCQ fin)$
\xEH
\xEH
\xEP

\hline

$(\xbm\xcc)+(\xbm PR)+(\xbm \xFO)+$
\xEH
$\xcj$, $(\xcv)$, singletons
\xEH
ranked
\xEH
$\xcH$ without $(\xbm dp)$
\xEH
$(RatM), (RatM=)$,
\xEP

$(\xbm \xcv)+(\xbm\xbe)$
\xEH
\xEH
$\xcg$ 1 copy
\xEH
\xEH
$(Log\xcv), (Log\xcv ')$
\xEP

\cline{4-5}

\xEH
\xEH
\xEH
$\xcJ$ without $(\xbm dp)$
\xEH
any ``normal''
\xEP

\xEH
\xEH
\xEH
\xEH
characterization
\xEP

\xEH
\xEH
\xEH
\xEH
of any size
\xEP

\hline

\end{tabular}

}

 karl-search= End Proposition Pref-Representation-Without-Ref
\vspace{7mm}

 *************************************

\vspace{7mm}

\subsubsection{Fact Pref-Sound}

 {\LARGE karl-search= Start Fact Pref-Sound }

\index{Fact Pref-Sound}

\bfa

$\hspace{0.01em}$

(+++ Orig. No.:  Fact Pref-Sound +++)

{\xssc LABEL: {Fact Pref-Sound}}
\label{Fact Pref-Sound}

$( \xbm \xcc )$ and $( \xbm PR)$ hold in all preferential structures.

 karl-search= End Fact Pref-Sound
\vspace{7mm}

 *************************************

\vspace{7mm}

\subsubsection{Fact Pref-Sound Proof}

 {\LARGE karl-search= Start Fact Pref-Sound Proof }

\index{Fact Pref-Sound Proof}

\efa

\subparagraph{
Proof
}

$\hspace{0.01em}$

(+++*** Orig.:  Proof )

Trivial. The central argument is: if $x,y \xbe X \xcc Y,$ and $x \xeb y$
in $X,$ then also
$x \xeb y$ in $Y.$

$ \xcz $
\\[3ex]

 karl-search= End Fact Pref-Sound Proof
\vspace{7mm}

 *************************************

\vspace{7mm}

\subsubsection{Fact Smooth-Sound}

 {\LARGE karl-search= Start Fact Smooth-Sound }

\index{Fact Smooth-Sound}

\bfa

$\hspace{0.01em}$

(+++ Orig. No.:  Fact Smooth-Sound +++)

{\xssc LABEL: {Fact Smooth-Sound}}
\label{Fact Smooth-Sound}

$( \xbm \xcc ),$ $( \xbm PR),$ and $( \xbm CUM)$ hold in all smooth
preferential structures.

 karl-search= End Fact Smooth-Sound
\vspace{7mm}

 *************************************

\vspace{7mm}

\subsubsection{Fact Smooth-Sound Proof}

 {\LARGE karl-search= Start Fact Smooth-Sound Proof }

\index{Fact Smooth-Sound Proof}

\efa

\subparagraph{
Proof
}

$\hspace{0.01em}$

(+++*** Orig.:  Proof )

By Fact \ref{Fact Pref-Sound} (page \pageref{Fact Pref-Sound}) , we only have to
show $( \xbm CUM).$
By Fact \ref{Fact Mu-Base} (page \pageref{Fact Mu-Base}) , $( \xbm CUT)$ follows
from $( \xbm PR),$
so it remains to show
$( \xbm CM).$ So suppose $ \xbm (X) \xcc Y \xcc X,$ we have to show $ \xbm
(Y) \xcc \xbm (X).$ Let
$x \xbe X- \xbm (X),$ so there is $x' \xbe X,$ $x' \xeb x,$ by smoothness,
there must be $x'' \xbe \xbm (X),$
$x'' \xeb x,$ so $x'' \xbe Y,$ and $x \xce \xbm (Y).$ The proof for the
case with copies is
analogous.

 karl-search= End Fact Smooth-Sound Proof
\vspace{7mm}

 *************************************

\vspace{7mm}

\subsubsection{Example Pref-Dp}

 {\LARGE karl-search= Start Example Pref-Dp }

\index{Example Pref-Dp}

\be

$\hspace{0.01em}$

(+++ Orig. No.:  Example Pref-Dp +++)

{\xssc LABEL: {Example Pref-Dp}}
\label{Example Pref-Dp}

This example was first given in [Sch92]. It shows
that condition $ \xCf (PR)$ may fail in preferential structures which are
not
definability preserving.

Let $v( \xdl ):=\{p_{i}:i \xbe \xbo \},$ $n,n' \xbe M_{ \xdl }$ be defined
by $n \xcm \{p_{i}:i \xbe \xbo \},$
$n' \xcm \{ \xCN p_{0}\} \xcv \{p_{i}:0<i< \xbo \}.$

Let $ \xdm:=<M_{ \xdl }, \xeb >$ where only $n \xeb n',$ i.e. just two
models are
comparable. Note that the structure is transitive and smooth.
Thus, by Fact \ref{Fact Smooth-Sound} (page \pageref{Fact Smooth-Sound})  $(
\xbm \xcc ),$ $( \xbm PR),$
$( \xbm CUM)$ hold.

Let $ \xbm:= \xbm_{ \xdm },$ and $ \xcn $ be defined as usual by $ \xbm
.$

Set $T:= \xCQ,$ $T':=\{p_{i}:0<i< \xbo \}.$ We have $M_{T}=M_{ \xdl },$
$f(M_{T})=M_{ \xdl }-\{n' \},$ $M_{T' }=\{n,n' \},$
$f(M_{T' })=\{n\}.$ So by the result of Example \ref{Example Not-Def} (page
\pageref{Example Not-Def}) ,
$f$ is not
definability preserving, and, furthermore,
$ \ol{ \ol{T} }= \ol{T},$ $ \ol{ \ol{T' } }= \ol{\{p_{i}:i< \xbo \}},$
$so$ $p_{0} \xbe \ol{ \ol{T \xcv T' } },$ $but$ $ \ol{ \ol{ \ol{T} } \xcv
T' }= \ol{ \ol{T} \xcv T' }= \ol{T' },$ $so$ $p_{0} \xce
 \ol{ \ol{ \ol{T} } \xcv T' },$
contradicting $ \xCf (PR),$ which holds in all definability preserving
preferential structures $ \xcz $
\\[3ex]

 karl-search= End Example Pref-Dp
\vspace{7mm}

 *************************************

\vspace{7mm}

 karl-search= End ToolBase1-Pref-ReprSumm
\vspace{7mm}

 *************************************

\vspace{7mm}

\newpage
\subsection{
General: Representation
}
\subsubsection{ToolBase1-Pref-ReprGen}

 {\LARGE karl-search= Start ToolBase1-Pref-ReprGen }

{\xssc LABEL: {Section Toolbase1-Pref-ReprGen}}
\label{Section Toolbase1-Pref-ReprGen}
\index{Section Toolbase1-Pref-ReprGen}
\subsubsection{Proposition Pref-Complete}

 {\LARGE karl-search= Start Proposition Pref-Complete }

\index{Proposition Pref-Complete}

\ee

\bp

$\hspace{0.01em}$

(+++ Orig. No.:  Proposition Pref-Complete +++)

{\xssc LABEL: {Proposition Pref-Complete}}
\label{Proposition Pref-Complete}

Let $ \xbm: \xdy \xcp \xdp (U)$ satisfy $( \xbm \xcc )$ and $( \xbm PR).$
Then there is a preferential
structure $ \xdx $ s.t. $ \xbm = \xbm_{ \xdx }.$ See e.g.  \cite{Sch04}.

 karl-search= End Proposition Pref-Complete
\vspace{7mm}

 *************************************

\vspace{7mm}

\subsubsection{Proposition Pref-Complete Proof}

 {\LARGE karl-search= Start Proposition Pref-Complete Proof }

\index{Proposition Pref-Complete Proof}

\ep

\subparagraph{
Proof
}

$\hspace{0.01em}$

(+++*** Orig.:  Proof )

The preferential structure is defined in
Construction \ref{Construction Pref-Base} (page \pageref{Construction
Pref-Base}) , Claim \ref{Claim Pref-Rep-Base} (page \pageref{Claim
Pref-Rep-Base})  shows representation.
The construction is basic for much of
the rest of the material on non-ranked structures.
\subsubsection{Definition Y-Pi-x}

 {\LARGE karl-search= Start Definition Y-Pi-x }

\index{Definition Y-Pi-x}

\bd

$\hspace{0.01em}$

(+++ Orig. No.:  Definition Y-Pi-x +++)

{\xssc LABEL: {Definition Y-Pi-x}}
\label{Definition Y-Pi-x}

For $x \xbe Z,$ let $ \xdy_{x}:=\{Y \xbe \xdy $: $x \xbe Y- \xbm (Y)\},$
$ \xbP_{x}:= \xbP \xdy_{x}.$

\ed

Note that $ \xCQ \xce \xdy_{x}$, $ \xbP_{x} \xEd \xCQ,$ and that $
\xbP_{x}=\{ \xCQ \}$ iff $ \xdy_{x}= \xCQ.$

 karl-search= End Definition Y-Pi-x
\vspace{7mm}

 *************************************

\vspace{7mm}

\subsubsection{Claim Mu-f}

 {\LARGE karl-search= Start Claim Mu-f }

\index{Claim Mu-f}

\bc

$\hspace{0.01em}$

(+++ Orig. No.:  Claim Mu-f +++)

{\xssc LABEL: {Claim Mu-f}}
\label{Claim Mu-f}

Let $ \xbm: \xdy \xcp \xdp (Z)$ satisfy $( \xbm \xcc )$ and $( \xbm PR),$
and let $U \xbe \xdy.$
Then $x \xbe \xbm (U)$ $ \xcr $ $x \xbe U$ $ \xcu $ $ \xcE f \xbe
\xbP_{x}.ran(f) \xcs U= \xCQ.$

 karl-search= End Claim Mu-f
\vspace{7mm}

 *************************************

\vspace{7mm}

\subsubsection{Claim Mu-f Proof}

 {\LARGE karl-search= Start Claim Mu-f Proof }

\index{Claim Mu-f Proof}

\ec

\subparagraph{
Proof
}

$\hspace{0.01em}$

(+++*** Orig.:  Proof )

Case 1: $ \xdy_{x}= \xCQ,$ thus $ \xbP_{x}=\{ \xCQ \}.$
`` $ \xcp $ '': Take $f:= \xCQ.$
`` $ \xcq $ '': $x \xbe U \xbe \xdy,$ $ \xdy_{x}= \xCQ $ $ \xcp $ $x \xbe
\xbm (U)$ by definition of $ \xdy_{x}.$

Case 2: $ \xdy_{x} \xEd \xCQ.$
`` $ \xcp $ '': Let $x \xbe \xbm (U) \xcc U.$ It suffices to show $Y \xbe
\xdy_{x}$ $ \xcp $ $Y-U \xEd \xCQ.$ But if $Y \xcc U$ and
$Y \xbe \xdy_{x}$, then $x \xbe Y- \xbm (Y),$ contradicting $( \xbm PR).$
`` $ \xcq $ '': If $x \xbe U- \xbm (U),$ then $U \xbe \xdy_{x}$, so $
\xcA f \xbe \xbP_{x}.ran(f) \xcs U \xEd \xCQ.$
$ \xcz $
\\[3ex]

 karl-search= End Claim Mu-f Proof
\vspace{7mm}

 *************************************

\vspace{7mm}

\subsubsection{Construction Pref-Base}

 {\LARGE karl-search= Start Construction Pref-Base }

\index{Construction Pref-Base}

\bcs

$\hspace{0.01em}$

(+++ Orig. No.:  Construction Pref-Base +++)

{\xssc LABEL: {Construction Pref-Base}}
\label{Construction Pref-Base}

Let $ \xdx:=\{<x,f>:x \xbe Z$ $ \xcu $ $f \xbe \xbP_{x}\},$ and $<x',f'
> \xeb <x,f>$ $: \xcr $ $x' \xbe ran(f).$ Let $ \xdz:=< \xdx, \xeb >.$

 karl-search= End Construction Pref-Base
\vspace{7mm}

 *************************************

\vspace{7mm}

\subsubsection{Claim Pref-Rep-Base}

 {\LARGE karl-search= Start Claim Pref-Rep-Base }

\index{Claim Pref-Rep-Base}

\ecs

\bc

$\hspace{0.01em}$

(+++ Orig. No.:  Claim Pref-Rep-Base +++)

{\xssc LABEL: {Claim Pref-Rep-Base}}
\label{Claim Pref-Rep-Base}

For $U \xbe \xdy,$ $ \xbm (U)= \xbm_{ \xdz }(U).$

 karl-search= End Claim Pref-Rep-Base
\vspace{7mm}

 *************************************

\vspace{7mm}

\subsubsection{Claim Pref-Rep-Base Proof}

 {\LARGE karl-search= Start Claim Pref-Rep-Base Proof }

\index{Claim Pref-Rep-Base Proof}

\ec

\subparagraph{
Proof
}

$\hspace{0.01em}$

(+++*** Orig.:  Proof )

By Claim \ref{Claim Mu-f} (page \pageref{Claim Mu-f}) , it suffices to show that
for all $U \xbe
\xdy $
$x \xbe \xbm_{ \xdz }(U)$ $ \xcr $ $x \xbe U$ and $ \xcE f \xbe
\xbP_{x}.ran(f) \xcs U= \xCQ.$ So let $U \xbe \xdy.$
`` $ \xcp $ '': If $x \xbe \xbm_{ \xdz }(U),$ then there is $<x,f>$
minimal in $ \xdx \xex U$ (recall from
Definition \ref{Definition Alg-Base} (page \pageref{Definition Alg-Base})  that
$ \xdx \xex U:=\{<x,i> \xbe \xdx
:x \xbe U\}),$
so $x \xbe U,$ and there
is no $<x',f' > \xeb <x,f>,$ $x' \xbe U,$ so by $ \xbP_{x' } \xEd \xCQ $
there is no $x' \xbe ran(f),$ $x' \xbe U,$ but then
$ran(f) \xcs U= \xCQ.$
`` $ \xcq $ '': If $x \xbe U,$ and there is $f \xbe \xbP_{x}$, $ran(f)
\xcs U= \xCQ,$ then $<x,f>$ is minimal in
$ \xdx \xex U.$
$ \xcz $ (Claim \ref{Claim Pref-Rep-Base} (page \pageref{Claim Pref-Rep-Base}) 
and Proposition \ref{Proposition Pref-Complete} (page \pageref{Proposition
Pref-Complete}) )
\\[3ex]

 karl-search= End Claim Pref-Rep-Base Proof
\vspace{7mm}

 *************************************

\vspace{7mm}

 karl-search= End Proposition Pref-Complete Proof
\vspace{7mm}

 *************************************

\vspace{7mm}

\newpage
\subsubsection{Proposition Pref-Complete-Trans}

 {\LARGE karl-search= Start Proposition Pref-Complete-Trans }

\index{Proposition Pref-Complete-Trans}

\bp

$\hspace{0.01em}$

(+++ Orig. No.:  Proposition Pref-Complete-Trans +++)

{\xssc LABEL: {Proposition Pref-Complete-Trans}}
\label{Proposition Pref-Complete-Trans}

Let $ \xbm: \xdy \xcp \xdp (U)$ satisfy $( \xbm \xcc )$ and $( \xbm PR).$
Then there is a transitive
preferential structure $ \xdx $ s.t. $ \xbm = \xbm_{ \xdx }.$ See e.g.
 \cite{Sch04}.

 karl-search= End Proposition Pref-Complete-Trans
\vspace{7mm}

 *************************************

\vspace{7mm}

\subsubsection{Proposition Pref-Complete-Trans Proof}

 {\LARGE karl-search= Start Proposition Pref-Complete-Trans Proof }

\index{Proposition Pref-Complete-Trans Proof}

\ep

\subparagraph{
Proof
}

$\hspace{0.01em}$

(+++*** Orig.:  Proof )

\subsubsection{Discussion Pref-Trans}

 {\LARGE karl-search= Start Discussion Pref-Trans }

\index{Discussion Pref-Trans}

\paragraph{
Discussion Pref-Trans
}

$\hspace{0.01em}$

(+++*** Orig.:  Discussion Pref-Trans )

{\xssc LABEL: {Section Discussion Pref-Trans}}
\label{Section Discussion Pref-Trans}

The Construction \ref{Construction Pref-Base} (page \pageref{Construction
Pref-Base})
(also used in  \cite{Sch92}) cannot be made transitive as it is,
this
will be shown below in Example \ref{Example Trans-1} (page \pageref{Example
Trans-1}) . The second
construction
in  \cite{Sch92} is a
special one, which is transitive, but uses heavily lack of smoothness.
(For
completeness' sake, we give a similar proof
in Proposition \ref{Proposition Equiv-Trans} (page \pageref{Proposition
Equiv-Trans}) .) We present here a more
flexibel and more adequate construction, which avoids a certain excess in
the
relation $ \xeb $ of the construction in Proposition \ref{Proposition
Equiv-Trans} (page \pageref{Proposition Equiv-Trans}) :
There, too many elements $<y,g>$
are smaller than some $<x,f>,$ as the relation is independent from $g.$
This excess
prevents transitivity.

We refine now the construction of the relation, to have better control
over
successors.

Recall that a tree of height $ \xck \xbo $ seems the right way to encode
the successors of
an element, as far as transitivity is concerned (which speaks only about
finite
chains). Now, in the basic construction, different copies have different
successors, chosen by different functions (elements of the cartesian
product).
As it suffices to make one copy of the successor smaller than the element
to be
minimized, we do the following: Let $<x,g>,$ with $g \xbe \xbP \{X:x \xbe
X-f(X)\}$ be one of the
elements of the standard construction. Let $<x',g' >$ be s.t. $x' \xbe
ran(g),$ then we
make again copies $<x,g,g' >,$ etc. for each such $x' $ and $g',$ and
make only $<x',g' >,$
but not some other $<x',g'' >$ smaller than $<x,g,g' >,$ for some other
$g'' \xbe \xbP \{X':x' \xbe X' -f(X' )\}.$ Thus, we have a much more
restricted relation, and much
better control over it. More precisely, we make trees, where we mark all
direct
and indirect successors, and each time the choice is made by the
appropriate
choice functions of the cartesian product. An element with its tree is a
successor of another element with its tree, iff the former is an initial
segment of the latter - see the definition in Construction \ref{Construction
Pref-Trees} (page \pageref{Construction Pref-Trees}) .

Recall also that transitivity is for free as we can use the element itself
to
minimize it. This is made precise by the use of the trees $tf_{x}$ for a
given element $x$ and choice function $f_{x}.$ But they also serve another
purpose.
The trees $tf_{x}$ are constructed as follows: The root is $x,$ the first
branching is
done according to $f_{x},$ and then we continue with constant choice. Let,
e.g.
$x' \xbe ran(f_{x}),$ we can now always choose $x',$ as it will be a
legal successor of $x' $
itself, being present in all $X' $ s.t. $x' \xbe X' -f(X' ).$ So we have a
tree which
branches once, directly above the root, and is then constant without
branching.
Obviously, this is essentially equivalent to the old construction in the
not
necessarily transitive case. This shows two things: first, the
construction
with trees gives the same $ \xbm $ as the old construction with simple
choice
functions. Second, even if we consider successors of successors, nothing
changes: we are still with the old $x'.$ Consequently, considering the
transitive
closure will not change matters, an element $<x,tf_{x}>$ will be minimized
by its
direct successors iff it will be minimized by direct and indirect
successors.
If you like, the trees $tf_{x}$ are the mathematical construction
expressing the
intuition that we know so little about minimization that we have to
consider
suicide a serious possibility - the intuitive reason why transitivity
imposes
no new conditions.

To summarize: Trees seem the right way to encode all the information
needed for
full control over successors for the transitive case. The special trees
$tf_{x}$
show that we have not changed things substantially, i.e. the new $ \xbm
-$functions in
the simple case and for the transitive closure stay the same. We hope that
this
construction will show its usefulness in other contexts, its naturalness
and
generality seem to be a good promise.

We give below the example which shows that the old construction is too
brutal
for transitivity to hold.

Recall that transitivity permits substitution in the following sense:
If (the two copies of) $x$ is killed by $y_{1}$ and $y_{2}$ together, and
$y_{1}$ is killed by
$z_{1}$ and $z_{2}$ together, then $x$ should be killed by $z_{1},$
$z_{2},$ and $y_{2}$ together.

But the old construction substitutes too much: In the old construction,
we considered elements $<x,f>,$ where $f \xbe \xbP_{x}$, with $<y,g> \xeb
<x,f>$ iff $y \xbe ran(f),$
independent of $g.$ This construction can, in general, not be made
transitive,
as Example \ref{Example Trans-1} (page \pageref{Example Trans-1})  below shows.

The new construction avoids this, as it ``looks ahead'', and not all
elements
$<y_{1},t_{y_{1}}>$ are smaller than $<x,t_{x}>,$ where $y_{1}$ is a child
of $x$ in $t_{x}$ (or $y_{1} \xbe ran(f)).$
The new construction is basically the same as Construction \ref{Construction
Pref-Base} (page \pageref{Construction Pref-Base}) ,
but avoids to make too many copies smaller than the copy to be killed.

Recall that
we need no new properties of $ \xbm $ to achieve transitivity here, as a
killed
element $x$ might (partially) ``commit suicide'', i.e. for some $i,$ $i' $
$<x,i> \xeb <x,i' >,$
so we cannot substitute $x$ by any set which does not contain $x:$ In this
simple
situation, if $x \xbe X- \xbm (X),$ we cannot find out whether all copies
of $x$ are killed
by some $y \xEd x,$ $y \xbe X.$ We can assume without loss of generality
that there is an
infinite descending chain of $x-$copies, which are not
killed by other elements. Thus, we cannot replace any $y_{i}$ as above by
any set
which does not contain $y_{i}$, but then substitution becomes trivial, as
any set
substituting $y_{i}$ has to contain $y_{i}$. Thus, we need no new
properties to
achieve transitivity.

 karl-search= End Discussion Pref-Trans
\vspace{7mm}

 *************************************

\vspace{7mm}

\subsubsection{Example Trans-1}

 {\LARGE karl-search= Start Example Trans-1 }

\index{Example Trans-1}

\be

$\hspace{0.01em}$

(+++ Orig. No.:  Example Trans-1 +++)

{\xssc LABEL: {Example Trans-1}}
\label{Example Trans-1}

As we consider only one set in each case, we can index with elements,
instead
of with functions.
So suppose $x,y_{1},y_{2} \xbe X,$ $y_{1},z_{1},z_{2} \xbe Y,$ $x \xce
\xbm (X),$ $y_{1} \xce \xbm (Y),$ and that we
need $y_{1}$ and $y_{2}$ to minimize $x,$ so there are two copies
$<x,y_{1}>,$ $<x,y_{2}>,$ likewise
we need $z_{1}$ and $z_{2}$ to minimize $y_{1},$ thus we have $<x,y_{1}>
\xee <y_{1},z_{1}>,$ $<x,y_{1}> \xee <y_{1},z_{2}>,$
$<x,y_{2}> \xee y_{2},$ $<y_{1},z_{1}> \xee z_{1},$ $<y_{1},z_{2}> \xee
z_{2}$ (the $z_{i}$ and $y_{2}$ are not killed).
If we take the transitive closure, we have $<x,y_{1}> \xee z_{k}$ for any
$i,k,$ so for any $z_{k}$
$\{z_{k},y_{2}\}$ will minimize all of $x,$ which is not intended. $ \xcz
$
\\[3ex]

 karl-search= End Example Trans-1
\vspace{7mm}

 *************************************

\vspace{7mm}

\ee

The preferential structure is defined in
Construction \ref{Construction Pref-Trees} (page \pageref{Construction
Pref-Trees}) , Claim \ref{Claim Tree-Repres-1} (page \pageref{Claim
Tree-Repres-1})  shows representation
for the simple structure, Claim \ref{Claim Tree-Repres-2} (page \pageref{Claim
Tree-Repres-2})
representation
for the transitive closure of the structure.

The main idea is to use the trees $tf_{x}$, whose elements are exactly
the elements
of the range of the choice function $f.$ This makes
Construction \ref{Construction Pref-Base} (page \pageref{Construction
Pref-Base})  and
Construction \ref{Construction Pref-Trees} (page \pageref{Construction
Pref-Trees})
basically equivalent,
and shows that the transitive case is
characterized by the same conditions as the general case. These trees are
defined below in Fact \ref{Fact Pref-Trees} (page \pageref{Fact Pref-Trees}) ,
(3),
and used in the proofs of Claim \ref{Claim Tree-Repres-1} (page \pageref{Claim
Tree-Repres-1})  and
Claim \ref{Claim Tree-Repres-2} (page \pageref{Claim Tree-Repres-2}) .

Again, Construction \ref{Construction Pref-Trees} (page \pageref{Construction
Pref-Trees})  contains the basic idea for
the
treatment of the transitive case. It can certainly be re-used in other
contexts.
\subsubsection{Construction Pref-Trees}

 {\LARGE karl-search= Start Construction Pref-Trees }

\index{Construction Pref-Trees}

\bcs

$\hspace{0.01em}$

(+++ Orig. No.:  Construction Pref-Trees +++)

{\xssc LABEL: {Construction Pref-Trees}}
\label{Construction Pref-Trees}

(1) For $x \xbe Z,$ let $T_{x}$ be the set of trees $t_{x}$ s.t.

(a) all nodes are elements of $Z,$

(b) the root of $t_{x}$ is $x,$

(c) $height(t_{x}) \xck \xbo,$

(d) if $y$ is an element in $t_{x}$, then there is $f \xbe \xbP_{y}:=
\xbP \{Y \xbe \xdy $: $y \xbe Y- \xbm (Y)\}$
s.t. the set of children of $y$ is $ran(f).$

(2) For $x,y \xbe Z,$ $t_{x} \xbe T_{x}$, $t_{y} \xbe T_{y}$, set $t_{x}
\xem t_{y}$ iff $y$ is a (direct) child
of the root $x$ in $t_{x}$, and $t_{y}$ is the subtree of $t_{x}$
beginning at $y.$

(3) Let $ \xdz $ $:=$ $<$ $\{<x,t_{x}>:$ $x \xbe Z,$ $t_{x} \xbe T_{x}\}$
, $<x,t_{x}> \xee <y,t_{y}>$ iff $t_{x} \xem t_{y}$ $>.$

 karl-search= End Construction Pref-Trees
\vspace{7mm}

 *************************************

\vspace{7mm}

\subsubsection{Fact Pref-Trees}

 {\LARGE karl-search= Start Fact Pref-Trees }

\index{Fact Pref-Trees}

\ecs

\bfa

$\hspace{0.01em}$

(+++ Orig. No.:  Fact Pref-Trees +++)

{\xssc LABEL: {Fact Pref-Trees}}
\label{Fact Pref-Trees}

(1) The construction ends at some $y$ iff $ \xdy_{y}= \xCQ,$ consequently
$T_{x}=\{x\}$ iff $ \xdy_{x}= \xCQ.$ (We identify the tree of height 1
with its root.)

(2) If $ \xdy_{x} \xEd \xCQ,$ $tc_{x}$, the totally ordered tree of
height $ \xbo,$ branching with $card=1,$
and with all elements equal to $x$ is an element of $T_{x}.$ Thus, with
(1), $T_{x} \xEd \xCQ $
for any $x.$

(3) If $f \xbe \xbP_{x}$, $f \xEd \xCQ,$ then the tree $tf_{x}$ with
root $x$ and otherwise
composed of the subtrees $t_{y}$ for $y \xbe ran(f),$ where $t_{y}:=y$ iff
$ \xdy_{y}= \xCQ,$
and $t_{y}:=tc_{y}$ iff $ \xdy_{y} \xEd \xCQ,$ is an element of $T_{x}$.
(Level 0 of $tf_{x}$ has $x$ as
element, the $t_{y}' s$ begin at level 1.)

(4) If $y$ is an element in $t_{x}$ and $t_{y}$ the subtree of $t_{x}$
starting at
$y,$ then $t_{y} \xbe T_{y}$.

(5) $<x,t_{x}> \xee <y,t_{y}>$ implies $y \xbe ran(f)$ for some $f \xbe
\xbP_{x}.$
$ \xcz $
\\[3ex]

 karl-search= End Fact Pref-Trees
\vspace{7mm}

 *************************************

\vspace{7mm}

\subsubsection{Claim Tree-Repres-1}

 {\LARGE karl-search= Start Claim Tree-Repres-1 }

\index{Claim Tree-Repres-1}

\efa

Claim \ref{Claim Tree-Repres-1} (page \pageref{Claim Tree-Repres-1})  shows
basic representation.

\bc

$\hspace{0.01em}$

(+++ Orig. No.:  Claim Tree-Repres-1 +++)

{\xssc LABEL: {Claim Tree-Repres-1}}
\label{Claim Tree-Repres-1}

$ \xcA U \xbe \xdy. \xbm (U)= \xbm_{ \xdz }(U)$

 karl-search= End Claim Tree-Repres-1
\vspace{7mm}

 *************************************

\vspace{7mm}

\subsubsection{Claim Tree-Repres-1 Proof}

 {\LARGE karl-search= Start Claim Tree-Repres-1 Proof }

\index{Claim Tree-Repres-1 Proof}

\ec

\subparagraph{
Proof
}

$\hspace{0.01em}$

(+++*** Orig.:  Proof )

By Claim \ref{Claim Mu-f} (page \pageref{Claim Mu-f}) , it suffices to show that
for all $U \xbe
\xdy $
$x \xbe \xbm_{ \xdz }(U)$ $ \xcr $ $x \xbe U$ $ \xcu $ $ \xcE f \xbe
\xbP_{x}.ran(f) \xcs U= \xCQ.$
Fix $U \xbe \xdy.$
`` $ \xcp $ '': $x \xbe \xbm_{ \xdz }(U)$ $ \xcp $ ex. $<x,t_{x}>$ minimal
in $ \xdz \xex U,$ thus $x \xbe U$ and there is no $<y,t_{y}> \xbe \xdz,$
$<y,t_{y}> \xeb <x,t_{x}>,$ $y \xbe U.$ Let $f$ define the set of children
of the root
$x$ in $t_{x}$. If $ran(f) \xcs U \xEd \xCQ,$ if $y \xbe U$ is a child
of $x$ in $t_{x}$, and if $t_{y}$ is the subtree
of $t_{x}$ starting at $y,$ then $t_{y} \xbe T_{y}$ and $<y,t_{y}> \xeb
<x,t_{x}>,$ contradicting minimality of
$<x,t_{x}>$ in $ \xdz \xex U.$ So $ran(f) \xcs U= \xCQ.$
`` $ \xcq $ '': Let $x \xbe U.$ If $ \xdy_{x}= \xCQ,$ then the tree $x$
has no $ \xem -$successors, and $<x,x>$ is
$ \xee -$minimal in $ \xdz.$ If $ \xdy_{x} \xEd \xCQ $ and $f \xbe
\xbP_{x}$ s.t. $ran(f) \xcs U= \xCQ,$ then $<x,tf_{x}>$ is $ \xee
-$minimal
in $ \xdz \xex U.$
$ \xcz $
\\[3ex]

 karl-search= End Claim Tree-Repres-1 Proof
\vspace{7mm}

 *************************************

\vspace{7mm}

\subsubsection{Claim Tree-Repres-2}

 {\LARGE karl-search= Start Claim Tree-Repres-2 }

\index{Claim Tree-Repres-2}

We consider now the transitive closure of $ \xdz.$ (Recall that $
\xeb^{*}$ denotes the
transitive closure of $ \xeb.)$ Claim \ref{Claim Tree-Repres-2} (page
\pageref{Claim Tree-Repres-2})
shows that transitivity does not
destroy what we have achieved. The trees $tf_{x}$ will play a crucial role
in the
demonstration.

\bc

$\hspace{0.01em}$

(+++ Orig. No.:  Claim Tree-Repres-2 +++)

{\xssc LABEL: {Claim Tree-Repres-2}}
\label{Claim Tree-Repres-2}

Let $ \xdz ' $ $:=$ $<$ $\{<x,t_{x}>:$ $x \xbe Z,$ $t_{x} \xbe T_{x}\}$,
$<x,t_{x}> \xee <y,t_{y}>$ iff $t_{x} \xem^{*}t_{y}$ $>.$ Then $ \xbm_{
\xdz }= \xbm_{ \xdz ' }.$

 karl-search= End Claim Tree-Repres-2
\vspace{7mm}

 *************************************

\vspace{7mm}

\subsubsection{Claim Tree-Repres-2 Proof}

 {\LARGE karl-search= Start Claim Tree-Repres-2 Proof }

\index{Claim Tree-Repres-2 Proof}

\ec

\subparagraph{
Proof
}

$\hspace{0.01em}$

(+++*** Orig.:  Proof )

Suppose there is $U \xbe \xdy,$ $x \xbe U,$ $x \xbe \xbm_{ \xdz }(U),$ $x
\xce \xbm_{ \xdz ' }(U).$
Then there must be an element $<x,t_{x}> \xbe \xdz $ with no $<x,t_{x}>
\xee <y,t_{y}>$ for any $y \xbe U.$
Let $f \xbe \xbP_{x}$ determine the set of children of $x$ in $t_{x}$,
then $ran(f) \xcs U= \xCQ,$
consider $tf_{x}.$ As all elements $ \xEd x$ of $tf_{x}$ are already in
$ran(f),$ no element of $tf_{x}$
is in $U.$ Thus there is no $<z,t_{z}> \xeb^{*}<x,tf_{x}>$ in $ \xdz $
with $z \xbe U,$ so $<x,tf_{x}>$ is minimal
in $ \xdz ' \xex U,$ contradiction.
$ \xcz $ (Claim \ref{Claim Tree-Repres-2} (page \pageref{Claim Tree-Repres-2}) 
and Proposition \ref{Proposition Pref-Complete-Trans} (page \pageref{Proposition
Pref-Complete-Trans}) )
\\[3ex]

 karl-search= End Claim Tree-Repres-2 Proof
\vspace{7mm}

 *************************************

\vspace{7mm}

 karl-search= End Proposition Pref-Complete-Trans Proof
\vspace{7mm}

 *************************************

\vspace{7mm}

\subsubsection{Proposition Equiv-Trans}

 {\LARGE karl-search= Start Proposition Equiv-Trans }

\index{Proposition Equiv-Trans}

We give now the direct proof, which we cannot adapt to the smooth case.
Such
easy results must be part of the folklore, but we give them for
completeness'
sake.

\bp

$\hspace{0.01em}$

(+++ Orig. No.:  Proposition Equiv-Trans +++)

{\xssc LABEL: {Proposition Equiv-Trans}}
\label{Proposition Equiv-Trans}

In the general case, every preferential structure is equivalent to a
transitive one - i.e. they have the same $ \xbm -$functions.

 karl-search= End Proposition Equiv-Trans
\vspace{7mm}

 *************************************

\vspace{7mm}

\subsubsection{Proposition Equiv-Trans Proof}

 {\LARGE karl-search= Start Proposition Equiv-Trans Proof }

\index{Proposition Equiv-Trans Proof}

\ep

\subparagraph{
Proof
}

$\hspace{0.01em}$

(+++*** Orig.:  Proof )

If $<a,i> \xee <b,j>,$ we create
an infinite descending chain of new copies $<b,<j,a,i,n>>,$ $n \xbe \xbo
,$ where
$<b,<j,a,i,n>> \xee <b,<j,a,i,n' >>$ if $n' >n,$ and make $<a,i> \xee
<b,<j,a,i,n>>$ for all
$n \xbe \xbo,$ but cancel the pair $<a,i> \xee <b,j>$ from the relation
(otherwise, we would
not have achieved anything), but $<b,j>$ stays as element in the set.
Now, the relation is trivially transitive, and all these $<b,<j,a,i,n>>$
just kill themselves, there is no need to minimize them by anything else.
We just continued $<a,i> \xee <b,j>$ in a way it cannot bother us. For the
$<b,j>,$ we
do of course the same thing again. So, we have full equivalence, i.e. the
$ \xbm -$functions of both structures are identical (this is trivial to
see). $ \xcz $
\\[3ex]

 karl-search= End Proposition Equiv-Trans Proof
\vspace{7mm}

 *************************************

\vspace{7mm}

 karl-search= End ToolBase1-Pref-ReprGen
\vspace{7mm}

 *************************************

\vspace{7mm}

\newpage

\subsection{
Smooth: Representation
}
\subsubsection{ToolBase1-Pref-ReprSmooth}

 {\LARGE karl-search= Start ToolBase1-Pref-ReprSmooth }

{\xssc LABEL: {Section Toolbase1-Pref-ReprSmooth}}
\label{Section Toolbase1-Pref-ReprSmooth}
\index{Section Toolbase1-Pref-ReprSmooth}
\subsubsection{Proposition Smooth-Complete}

 {\LARGE karl-search= Start Proposition Smooth-Complete }

\index{Proposition Smooth-Complete}

\bp

$\hspace{0.01em}$

(+++ Orig. No.:  Proposition Smooth-Complete +++)

{\xssc LABEL: {Proposition Smooth-Complete}}
\label{Proposition Smooth-Complete}

Let $ \xbm: \xdy \xcp \xdp (U)$ satisfy $( \xbm \xcc ),$ $( \xbm PR),$
and $( \xbm CUM),$ and the domain $ \xdy $ $( \xcv ).$

Then there is a $ \xdy -$smooth preferential structure $ \xdx $ s.t. $
\xbm = \xbm_{ \xdx }.$
See e.g.  \cite{Sch04}.

 karl-search= End Proposition Smooth-Complete
\vspace{7mm}

 *************************************

\vspace{7mm}

\subsubsection{Proposition Smooth-Complete Proof}

 {\LARGE karl-search= Start Proposition Smooth-Complete Proof }

\index{Proposition Smooth-Complete Proof}

\ep

\subparagraph{
Proof
}

$\hspace{0.01em}$

(+++*** Orig.:  Proof )

\subsubsection{Comment Smooth-Complete Proof}

 {\LARGE karl-search= Start Comment Smooth-Complete Proof }

\index{Comment Smooth-Complete Proof}

Outline: We first define a structure $ \xdz $ (in a way very similar to
Construction \ref{Construction Pref-Base} (page \pageref{Construction
Pref-Base}) )
which represents $ \xbm,$ but is not
necessarily $ \xdy -$smooth, refine
it to $ \xdz ' $ and show that $ \xdz ' $ represents $ \xbm $ too, and
that $ \xdz ' $ is $ \xdy -$smooth.

In the structure $ \xdz ',$ all pairs destroying smoothness in $ \xdz $
are successively
repaired, by adding minimal elements: If $<y,j>$ is not minimal, and has
no minimal
$<x,i>$ below it, we just add one such $<x,i>.$ As the repair process
might itself
generate such ``bad'' pairs, the process may have to be repeated infinitely
often.
Of course, one has to take care that the representation property is
preserved.

The proof given is close to the minimum one has to show (except that we
avoid
$H(U),$ instead of $U$ - as was done in the old proof of [Sch96-1]). We
could simplify
further, we do not, in order to stay closer to the construction that is
really
needed. The reader will find the simplification as building block of the
proof
of the transitive case. (In the simplified proof, we would consider for
$x,U$ s.t.
$x \xbe \xbm (U)$
the pairs $<x,g_{U}>$ with $g_{U} \xbe \xbP \{ \xbm (U \xcv Y):x \xbe Y
\xcC H(U)\},$ giving minimal elements.
For the $U$ s.t. $x \xbe U- \xbm (U),$ we would choose $<x,g>$ s.t. $g
\xbe \xbP \{ \xbm (Y):x \xbe Y \xbe \xdy \}$ with
$<x',g'_{U}> \xeb <x,g>$ for $<x',g'_{U}>$ as above.)

Construction \ref{Construction Smooth-Base} (page \pageref{Construction
Smooth-Base})  represents $ \xbm.$ The
structure will not yet
be smooth, we will mend it afterwards in Construction \ref{Construction
Smooth-Admiss} (page \pageref{Construction Smooth-Admiss}) .

 karl-search= End Comment Smooth-Complete Proof
\vspace{7mm}

 *************************************

\vspace{7mm}

\subsubsection{Comment HU}

 {\LARGE karl-search= Start Comment HU }

\index{Comment HU}

\paragraph{
The constructions
}

$\hspace{0.01em}$

(+++*** Orig.:  The constructions )

{\xssc LABEL: {Section The constructions}}
\label{Section The constructions}

$ \xdy $ will be closed under finite unions
throughout this Section. We first define $H(U),$ and show some facts
about it. $H(U)$ has an important role, for the
following
reason: If $u \xbe \xbm (U),$ but $u \xbe X- \xbm (X),$ then there is $x
\xbe \xbm (X)-H(U).$ Consequently,
to kill minimality of $u$ in $X,$ we can choose $x \xbe \xbm (X)-H(U),$ $x
\xeb u,$ without
interfering with u's minimality in $U.$ Moreover, if $x \xbe Y- \xbm (Y),$
then, by $x \xce H(U),$
$ \xbm (Y) \xcC H(U),$ so we can kill minimality of $x$ in $Y$ by choosing
some $y \xce H(U).$
Thus, even in the transitive case, we can leave $U$ to destroy minimality
of $u$ in
some $X,$ without ever having to come back into $U,$ it suffices to choose
sufficiently far from $U,$ i.e. outside $H(U).$ $H(U)$ is the right notion
of
``neighborhood''.

Note: Not all $z \xbe Z$ have to occur in our structure, therefore it is
quite
possible that $X \xbe \xdy,$ $X \xEd \xCQ,$ but $ \xbm_{ \xdz }(X)= \xCQ
.$ This is why we have introduced
the set $K$ in Definition \ref{Definition K} (page \pageref{Definition K})  and
such $X$ will be subsets
of $Z-$K.

Let now $ \xbm: \xdy \xcp \xdp (Z).$

 karl-search= End Comment HU
\vspace{7mm}

 *************************************

\vspace{7mm}

\subsubsection{Definition HU}

 {\LARGE karl-search= Start Definition HU }

\index{Definition HU}

\bd

$\hspace{0.01em}$

(+++ Orig. No.:  Definition HU +++)

{\xssc LABEL: {Definition HU}}
\label{Definition HU}

Define $H(U)$ $:=$ $ \xcV \{X: \xbm (X) \xcc U\}.$

 karl-search= End Definition HU
\vspace{7mm}

 *************************************

\vspace{7mm}

\subsubsection{Definition K}

 {\LARGE karl-search= Start Definition K }

\index{Definition K}

\ed

\bd

$\hspace{0.01em}$

(+++ Orig. No.:  Definition K +++)

{\xssc LABEL: {Definition K}}
\label{Definition K}

Let $K:=\{x \xbe Z:$ $ \xcE X \xbe \xdy.x \xbe \xbm (X)\}$

 karl-search= End Definition K
\vspace{7mm}

 *************************************

\vspace{7mm}

\subsubsection{Fact HU-1}

 {\LARGE karl-search= Start Fact HU-1 }

\index{Fact HU-1}

\ed

\bfa

$\hspace{0.01em}$

(+++ Orig. No.:  Fact HU-1 +++)

{\xssc LABEL: {Fact HU-1}}
\label{Fact HU-1}

$( \xbm \xcc )$ $+$ $( \xbm PR)$ $+$ $( \xbm CUM)$ $+$ $( \xcv )$ entail:

(1) $ \xbm (A) \xcc B$ $ \xcp $ $ \xbm (A \xcv B)= \xbm (B)$

(2) $ \xbm (X) \xcc U,$ $U \xcc Y$ $ \xcp $ $ \xbm (Y \xcv X)= \xbm (Y)$

(3) $ \xbm (X) \xcc U,$ $U \xcc Y$ $ \xcp $ $ \xbm (Y) \xcs X \xcc \xbm
(U)$

(4) $ \xbm (X) \xcc U$ $ \xcp $ $ \xbm (U) \xcs X \xcc \xbm (X)$

(5) $U \xcc A,$ $ \xbm (A) \xcc H(U)$ $ \xcp $ $ \xbm (A) \xcc U$

(6) Let $x \xbe K,$ $Y \xbe \xdy,$ $x \xbe Y- \xbm (Y),$ then $ \xbm (Y)
\xEd \xCQ.$

 karl-search= End Fact HU-1
\vspace{7mm}

 *************************************

\vspace{7mm}

\subsubsection{Fact HU-1 Proof}

 {\LARGE karl-search= Start Fact HU-1 Proof }

\index{Fact HU-1 Proof}

\efa

\subparagraph{
Proof
}

$\hspace{0.01em}$

(+++*** Orig.:  Proof )

(1) $ \xbm (A) \xcc B$ $ \xcp $ $ \xbm (A \xcv B) \xcc \xbm (A) \xcv \xbm
(B) \xcc B$ $ \xcp_{( \xbm CUM)}$ $ \xbm (B)= \xbm (A \xcv B).$

(2) trivial by (1).

(3) $ \xbm (Y) \xcs X$ $=$ (by (2)) $ \xbm (Y \xcv X) \xcs X$ $ \xcc $ $
\xbm (Y \xcv X) \xcs (X \xcv U)$ $ \xcc $ (by $( \xbm PR))$
$ \xbm (X \xcv U)$ $=$ (by (1)) $ \xbm (U).$

(4) $ \xbm (U) \xcs X$ $=$ $ \xbm (X \xcv U) \xcs X$ by (1) $ \xcc $ $
\xbm (X)$ by $( \xbm PR)$

(5) Let $U \xcc A,$ $ \xbm (A) \xcc H(U).$ So $ \xbm (A)$ $=$ $ \xcV \{
\xbm (A) \xcs Y: \xbm (Y) \xcc U\}$ $ \xcc $
$ \xbm (U)$ $ \xcc $ $U$ by (3).

(6) Suppose $x \xbe \xbm (X),$ $ \xbm (Y)= \xCQ $ $ \xcp $ $ \xbm (Y) \xcc
X,$ so by (4) $Y \xcs \xbm (X) \xcc \xbm (Y),$ so
$x \xbe \xbm (Y).$

$ \xcz $
\\[3ex]

 karl-search= End Fact HU-1 Proof
\vspace{7mm}

 *************************************

\vspace{7mm}

\subsubsection{Fact HU-2}

 {\LARGE karl-search= Start Fact HU-2 }

\index{Fact HU-2}

The following Fact \ref{Fact HU-2} (page \pageref{Fact HU-2})  contains the
basic properties
of $ \xbm $ and $H(U)$ which we will need for the representation
construction.

\bfa

$\hspace{0.01em}$

(+++ Orig. No.:  Fact HU-2 +++)

{\xssc LABEL: {Fact HU-2}}
\label{Fact HU-2}

Let A, $U,$ $U',$ $Y$ and all $A_{i}$ be in $ \xdy.$ Let $( \xbm \xcc )$
$+$ $( \xbm PR)$ $+$ $( \xcv )$ hold.

(1) $A= \xcV \{A_{i}:i \xbe I\}$ $ \xcp $ $ \xbm (A) \xcc \xcV \{ \xbm
(A_{i}):i \xbe I\},$

(2) $U \xcc H(U),$ and $U \xcc U' \xcp H(U) \xcc H(U' ),$

(3) $ \xbm (U \xcv Y)-H(U) \xcc \xbm (Y).$

If, in addition, $( \xbm CUM)$ holds, then we also have:

(4) $U \xcc A,$ $ \xbm (A) \xcc H(U)$ $ \xcp $ $ \xbm (A) \xcc U,$

(5) $ \xbm (Y) \xcc H(U)$ $ \xcp $ $Y \xcc H(U)$ and $ \xbm (U \xcv Y)=
\xbm (U),$

(6) $x \xbe \xbm (U),$ $x \xbe Y- \xbm (Y)$ $ \xcp $ $Y \xcC H(U),$

(7) $Y \xcC H(U)$ $ \xcp $ $ \xbm (U \xcv Y) \xcC H(U).$

 karl-search= End Fact HU-2
\vspace{7mm}

 *************************************

\vspace{7mm}

\subsubsection{Fact HU-2 Proof}

 {\LARGE karl-search= Start Fact HU-2 Proof }

\index{Fact HU-2 Proof}

\efa

\subparagraph{
Proof
}

$\hspace{0.01em}$

(+++*** Orig.:  Proof )

(1) $ \xbm (A) \xcs A_{j} \xcc \xbm (A_{j}) \xcc \xcV \xbm (A_{i}),$ so by
$ \xbm (A) \xcc A= \xcV A_{i}$ $ \xbm (A) \xcc \xcV \xbm (A_{i}).$

(2) trivial.

(3) $ \xbm (U \xcv Y)-H(U)$ $ \xcc_{(2)}$ $ \xbm (U \xcv Y)-U$ $ \xcc_{(
\xbm \xcc )}$ $ \xbm (U \xcv Y) \xcs Y$ $ \xcc_{( \xbm PR)}$ $ \xbm (Y).$

(4) This is Fact \ref{Fact HU-1} (page \pageref{Fact HU-1})  (5).

(5) Let $ \xbm (Y) \xcc H(U),$ then by $ \xbm (U) \xcc H(U)$ and (1) $
\xbm (U \xcv Y) \xcc \xbm (U) \xcv \xbm (Y) \xcc H(U),$
so by (4) $ \xbm (U \xcv Y) \xcc U$ and $U \xcv Y \xcc H(U).$ Moreover, $
\xbm (U \xcv Y) \xcc U \xcc U \xcv Y$ $ \xcp_{( \xbm CUM)}$ $ \xbm (U \xcv
Y)= \xbm (U).$

(6) If not, $Y \xcc H(U),$ so $ \xbm (Y) \xcc H(U),$ so $ \xbm (U \xcv Y)=
\xbm (U)$ by (5), but $x \xbe Y- \xbm (Y)$ $ \xcp_{( \xbm PR)}$
$x \xce \xbm (U \xcv Y)= \xbm (U),$ $contradiction.$

(7) $ \xbm (U \xcv Y) \xcc H(U)$ $ \xcp_{(5)}$ $U \xcv Y \xcc H(U).$
$ \xcz $
\\[3ex]

 karl-search= End Fact HU-2 Proof
\vspace{7mm}

 *************************************

\vspace{7mm}

\subsubsection{Definition Gamma-x}

 {\LARGE karl-search= Start Definition Gamma-x }

\index{Definition Gamma-x}

\bd

$\hspace{0.01em}$

(+++ Orig. No.:  Definition Gamma-x +++)

{\xssc LABEL: {Definition Gamma-x}}
\label{Definition Gamma-x}

For $x \xbe Z,$ let $ \xdw_{x}:=\{ \xbm (Y)$: $Y \xbe \xdy $ $ \xcu $ $x
\xbe Y- \xbm (Y)\},$ $ \xbG_{x}:= \xbP \xdw_{x}$, and $K:=\{x \xbe Z$: $
\xcE X \xbe \xdy.x \xbe \xbm (X)\}.$

\ed

Note that we consider here now $ \xbm (Y)$ in $ \xdw_{x}$, and not $Y$ as
in $ \xdy_{x}$
in Definition \ref{Definition Y-Pi-x} (page \pageref{Definition Y-Pi-x}) .

 karl-search= End Definition Gamma-x
\vspace{7mm}

 *************************************

\vspace{7mm}

\subsubsection{Remark Gamma-x}

 {\LARGE karl-search= Start Remark Gamma-x }

\index{Remark Gamma-x}

\br

$\hspace{0.01em}$

(+++ Orig. No.:  Remark Gamma-x +++)

{\xssc LABEL: {Remark Gamma-x}}
\label{Remark Gamma-x}

Assume now $( \xbm \xcc ),$ $( \xbm PR),$ $( \xbm CUM),$ $( \xcv )$ to
hold.

(1) $x \xbe K$ $ \xcp $ $ \xbG_{x} \xEd \xCQ,$

(2) $g \xbe \xbG_{x}$ $ \xcp $ $ran(g) \xcc K.$

 karl-search= End Remark Gamma-x
\vspace{7mm}

 *************************************

\vspace{7mm}

\subsubsection{Remark Gamma-x Proof}

 {\LARGE karl-search= Start Remark Gamma-x Proof }

\index{Remark Gamma-x Proof}

\er

\subparagraph{
Proof
}

$\hspace{0.01em}$

(+++*** Orig.:  Proof )

(1) We have to show that $Y \xbe \xdy,$ $x \xbe Y- \xbm (Y)$ $ \xcp $ $
\xbm (Y) \xEd \xCQ.$ This was shown
in Fact \ref{Fact HU-1} (page \pageref{Fact HU-1})  (6).

(2) By definition, $ \xbm (Y) \xcc K$ for all $Y \xbe \xdy.$
$ \xcz $
\\[3ex]

 karl-search= End Remark Gamma-x Proof
\vspace{7mm}

 *************************************

\vspace{7mm}

\subsubsection{Claim Cum-Mu-f}

 {\LARGE karl-search= Start Claim Cum-Mu-f }

\index{Claim Cum-Mu-f}

The following claim is the analogue of Claim \ref{Claim Mu-f} (page
\pageref{Claim Mu-f})  above.

\bc

$\hspace{0.01em}$

(+++ Orig. No.:  Claim Cum-Mu-f +++)

{\xssc LABEL: {Claim Cum-Mu-f}}
\label{Claim Cum-Mu-f}

Assume now $( \xbm \xcc ),$ $( \xbm PR),$ $( \xbm CUM),$ $( \xcv )$ to
hold.

Let $U \xbe \xdy,$ $x \xbe K.$ Then

(1) $x \xbe \xbm (U)$ $ \xcr $ $x \xbe U$ $ \xcu $ $ \xcE f \xbe
\xbG_{x}.ran(f) \xcs U= \xCQ,$

(2) $x \xbe \xbm (U)$ $ \xcr $ $x \xbe U$ $ \xcu $ $ \xcE f \xbe
\xbG_{x}.ran(f) \xcs H(U)= \xCQ.$

 karl-search= End Claim Cum-Mu-f
\vspace{7mm}

 *************************************

\vspace{7mm}

\subsubsection{Claim Cum-Mu-f Proof}

 {\LARGE karl-search= Start Claim Cum-Mu-f Proof }

\index{Claim Cum-Mu-f Proof}

\ec

\subparagraph{
Proof
}

$\hspace{0.01em}$

(+++*** Orig.:  Proof )

(1)
Case 1: $ \xdw_{x}= \xCQ,$ thus $ \xbG_{x}=\{ \xCQ \}.$
`` $ \xcp $ '': Take $f:= \xCQ.$
`` $ \xcq $ '': $x \xbe U \xbe \xdy,$ $ \xdw_{x}= \xCQ $ $ \xcp $ $x \xbe
\xbm (U)$ by definition of $ \xdw_{x}.$

Case 2: $ \xdw_{x} \xEd \xCQ.$
`` $ \xcp $ '': Let $x \xbe \xbm (U) \xcc U.$ It suffices to show $Y \xbe
\xdw_{x}$ $ \xcp $ $ \xbm (Y)-H(U) \xEd \xCQ.$
But $Y \xbe \xdw_{x}$ $ \xcp $ $x \xbe Y- \xbm (Y)$ $ \xcp $ (by Fact \ref{Fact
HU-2} (page \pageref{Fact HU-2}) , (6))
$Y \xcC H(U)$ $ \xcp $ (by Fact \ref{Fact HU-2} (page \pageref{Fact HU-2}) ,
(5)) $ \xbm (Y) \xcC
H(U).$
`` $ \xcq $ '': If $x \xbe U- \xbm (U),$ $U \xbe \xdw_{x}$, moreover $
\xbG_{x} \xEd \xCQ $ by
Remark \ref{Remark Gamma-x} (page \pageref{Remark Gamma-x}) , (1) and thus (or
by the same argument) $ \xbm (U) \xEd \xCQ,$ so $ \xcA f \xbe
\xbG_{x}.ran(f) \xcs U \xEd \xCQ.$

(2): The proof is verbatim the same as for (1).
$ \xcz $
\\[3ex]

 karl-search= End Claim Cum-Mu-f Proof
\vspace{7mm}

 *************************************

\vspace{7mm}

\subparagraph{
Proof: (Prop. 6.14)
}

$\hspace{0.01em}$

(+++*** Orig.:  Proof: (Prop. 6.14) )

{\xssc LABEL: {Section Proof: (Prop. 6.14)}}
\label{Section Proof: (Prop. 6.14)}

\subsubsection{Construction Smooth-Base}

 {\LARGE karl-search= Start Construction Smooth-Base }

\index{Construction Smooth-Base}

\bcs

$\hspace{0.01em}$

(+++ Orig. No.:  Construction Smooth-Base +++)

{\xssc LABEL: {Construction Smooth-Base}}
\label{Construction Smooth-Base}

(Construction of $ \xdz )$
Let $ \xdx $ $:=$ $\{<x,g>$: $x \xbe K,$ $g \xbe \xbG_{x}\},$ $<x',g' >
\xeb <x,g>$ $: \xcr $ $x' \xbe ran(g),$

$ \xdz:=< \xdx, \xeb >.$

 karl-search= End Construction Smooth-Base
\vspace{7mm}

 *************************************

\vspace{7mm}

\subsubsection{Claim Smooth-Base}

 {\LARGE karl-search= Start Claim Smooth-Base }

\index{Claim Smooth-Base}

\ecs

\bc

$\hspace{0.01em}$

(+++ Orig. No.:  Claim Smooth-Base +++)

{\xssc LABEL: {Claim Smooth-Base}}
\label{Claim Smooth-Base}

$ \xcA U \xbe \xdy. \xbm (U)= \xbm_{ \xdz }(U)$

 karl-search= End Claim Smooth-Base
\vspace{7mm}

 *************************************

\vspace{7mm}

\subsubsection{Claim Smooth-Base Proof}

 {\LARGE karl-search= Start Claim Smooth-Base Proof }

\index{Claim Smooth-Base Proof}

\ec

\subparagraph{
Proof
}

$\hspace{0.01em}$

(+++*** Orig.:  Proof )

Case 1: $x \xce K.$ Then $x \xce \xbm (U)$ and $x \xce \xbm_{ \xdz }(U).$

Case 2: $x \xbe K.$
By Claim \ref{Claim Cum-Mu-f} (page \pageref{Claim Cum-Mu-f}) , (1) it suffices
to show that for all $U
\xbe \xdy $
$x \xbe \xbm_{ \xdz }(U)$ $ \xcr $ $x \xbe U$ $ \xcu $ $ \xcE f \xbe
\xbG_{x}.ran(f) \xcs U= \xCQ.$
Fix $U \xbe \xdy.$
`` $ \xcp $ '': $x \xbe \xbm_{ \xdz }(U)$ $ \xcp $ ex. $<x,f>$ minimal in
$ \xdx \xex U,$ thus $x \xbe U$ and there is no
$<x',f' > \xeb <x,f>,$ $x' \xbe U,$ $x' \xbe K.$ But if $x' \xbe K,$ then
by
Remark \ref{Remark Gamma-x} (page \pageref{Remark Gamma-x}) , (1), $ \xbG_{x' }
\xEd \xCQ,$
so we find suitable $f'.$ Thus, $ \xcA x' \xbe ran(f).x' \xce U$ or $x'
\xce K.$ But $ran(f) \xcc K,$ so
$ran(f) \xcs U= \xCQ.$
`` $ \xcq $ '': If $x \xbe U,$ $f \xbe \xbG_{x}$ s.t. $ran(f) \xcs U= \xCQ
,$ then $<x,f>$ is minimal in $ \xdx \xex U.$

$ \xcz $
\\[3ex]

 karl-search= End Claim Smooth-Base Proof
\vspace{7mm}

 *************************************

\vspace{7mm}

\subsubsection{Construction Smooth-Admiss}

 {\LARGE karl-search= Start Construction Smooth-Admiss }

\index{Construction Smooth-Admiss}

We now construct the refined structure $ \xdz '.$

\bcs

$\hspace{0.01em}$

(+++ Orig. No.:  Construction Smooth-Admiss +++)

{\xssc LABEL: {Construction Smooth-Admiss}}
\label{Construction Smooth-Admiss}

(Construction of $ \xdz ' )$

$ \xbs $ is called $x-$admissible sequence iff

1. $ \xbs $ is a sequence of length $ \xck \xbo,$ $ \xbs =\{ \xbs_{i}:i
\xbe \xbo \},$

2. $ \xbs_{o} \xbe \xbP \{ \xbm (Y)$: $Y \xbe \xdy $ $ \xcu $ $x \xbe Y-
\xbm (Y)\},$

3. $ \xbs_{i+1} \xbe \xbP \{ \xbm (X)$: $X \xbe \xdy $ $ \xcu $ $x \xbe
\xbm (X)$ $ \xcu $ $ran( \xbs_{i}) \xcs X \xEd \xCQ \}.$

By 2., $ \xbs_{0}$ minimizes $x,$ and by 3., if $x \xbe \xbm (X),$ and
$ran( \xbs_{i}) \xcs X \xEd \xCQ,$ i.e. we
have destroyed minimality of $x$ in $X,$ $x$ will be above some $y$
minimal in $X$ to
preserve smoothness.

Let $ \xbS_{x}$ be the set of $x-$admissible sequences, for $ \xbs \xbe
\xbS_{x}$ let $ \wt{ \xbs }:= \xcV \{ran( \xbs_{i}):i \xbe \xbo \}.$
Note that by the argument in the proof of
Remark \ref{Remark Gamma-x} (page \pageref{Remark Gamma-x}) , (1),
$ \xbS_{x} \xEd \xCQ,$ if $x \xbe K.$

Let $ \xdx ' $ $:=$ $\{<x, \xbs >$: $x \xbe K$ $ \xcu $ $ \xbs \xbe
\xbS_{x}\}$ and $<x', \xbs ' > \xeb ' <x, \xbs >$ $: \xcr $ $x' \xbe \wt{
\xbs }$.
Finally, let $ \xdz ':=< \xdx ', \xeb ' >,$ and $ \xbm ':= \xbm_{ \xdz
' }.$

\ecs

It is now easy to show that $ \xdz ' $ represents $ \xbm,$ and that $
\xdz ' $ is smooth.
For $x \xbe \xbm (U),$ we construct a special $x-$admissible sequence $
\xbs^{x,U}$ using the
properties of $H(U).$

 karl-search= End Construction Smooth-Admiss
\vspace{7mm}

 *************************************

\vspace{7mm}

\subsubsection{Claim Smooth-Admiss-1}

 {\LARGE karl-search= Start Claim Smooth-Admiss-1 }

\index{Claim Smooth-Admiss-1}

\bc

$\hspace{0.01em}$

(+++ Orig. No.:  Claim Smooth-Admiss-1 +++)

{\xssc LABEL: {Claim Smooth-Admiss-1}}
\label{Claim Smooth-Admiss-1}

For all $U \xbe \xdy $ $ \xbm (U)= \xbm_{ \xdz }(U)= \xbm ' (U).$

 karl-search= End Claim Smooth-Admiss-1
\vspace{7mm}

 *************************************

\vspace{7mm}

\subsubsection{Claim Smooth-Admiss-1 Proof}

 {\LARGE karl-search= Start Claim Smooth-Admiss-1 Proof }

\index{Claim Smooth-Admiss-1 Proof}

\ec

\subparagraph{
Proof
}

$\hspace{0.01em}$

(+++*** Orig.:  Proof )

If $x \xce K,$ then $x \xce \xbm_{ \xdz }(U),$ and $x \xce \xbm ' (U)$ for
any $U.$ So assume $x \xbe K.$ If $x \xbe U$ and
$x \xce \xbm_{ \xdz }(U),$ then for all $<x,f> \xbe \xdx,$ there is $<x'
,f' > \xbe \xdx $ with $<x',f' > \xeb <x,f>$ and
$x' \xbe U.$ Let now $<x, \xbs > \xbe \xdx ',$ then $<x, \xbs_{0}> \xbe
\xdx,$ and let $<x',f' > \xeb <x, \xbs_{0}>$ in $ \xdz $ with
$x' \xbe U.$ As $x' \xbe K,$ $ \xbS_{x' } \xEd \xCQ,$ let $ \xbs ' \xbe
\xbS_{x' }$. Then $<x', \xbs ' > \xeb ' <x, \xbs >$ in $ \xdz '.$ Thus
$x \xce \xbm ' (U).$
Thus, for all $U \xbe \xdy,$ $ \xbm ' (U) \xcc \xbm_{ \xdz }(U)= \xbm
(U).$

It remains to show $x \xbe \xbm (U) \xcp x \xbe \xbm ' (U).$

Assume $x \xbe \xbm (U)$ (so $x \xbe K),$ $U \xbe \xdy,$ we will
construct minimal $ \xbs,$ i.e. show that
there is $ \xbs^{x,U} \xbe \xbS_{x}$ s.t. $ \wt{ \xbs^{x,U}} \xcs U= \xCQ
.$ We construct this $ \xbs^{x,U}$ inductively, with the
stronger property that $ran( \xbs^{x,U}_{i}) \xcs H(U)= \xCQ $ for all $i
\xbe \xbo.$

$ \ul{ \xbs^{x,U}_{0}:}$
$x \xbe \xbm (U),$ $x \xbe Y- \xbm (Y)$ $ \xcp $ $ \xbm (Y)-H(U) \xEd \xCQ
$ by Fact \ref{Fact HU-2} (page \pageref{Fact HU-2}) , $(6)+(5).$
Let $ \xbs^{x,U}_{0}$ $ \xbe $ $ \xbP \{ \xbm (Y)-H(U):$ $Y \xbe \xdy,$
$x \xbe Y- \xbm (Y)\},$ so $ran( \xbs^{x,U}_{0}) \xcs H(U)= \xCQ.$

$ \ul{ \xbs^{x,U}_{i} \xcp \xbs^{x,U}_{i+1}:}$
By induction hypothesis, $ran( \xbs^{x,U}_{i}) \xcs H(U)= \xCQ.$ Let $X
\xbe \xdy $ be s.t. $x \xbe \xbm (X),$
$ran( \xbs^{x,U}_{i}) \xcs X \xEd \xCQ.$ Thus $X \xcC H(U),$ so $ \xbm (U
\xcv X)-H(U) \xEd \xCQ $ by Fact \ref{Fact HU-2} (page \pageref{Fact HU-2}) ,
(7).
Let $ \xbs^{x,U}_{i+1}$ $ \xbe $ $ \xbP \{ \xbm (U \xcv X)-H(U):$ $X \xbe
\xdy,$ $x \xbe \xbm (X),$ $ran( \xbs^{x,U}_{i}) \xcs X \xEd \xCQ \},$ so
$ran( \xbs^{x,U}_{i+1}) \xcs H(U)= \xCQ.$
As $ \xbm (U \xcv X)-H(U) \xcc \xbm (X)$ by Fact \ref{Fact HU-2} (page
\pageref{Fact HU-2}) ,
(3), the construction satisfies the
$x-$admissibility condition.
$ \xcz $
\\[3ex]

 karl-search= End Claim Smooth-Admiss-1 Proof
\vspace{7mm}

 *************************************

\vspace{7mm}

\subsubsection{Claim Smooth-Admiss-2}

 {\LARGE karl-search= Start Claim Smooth-Admiss-2 }

\index{Claim Smooth-Admiss-2}

We now show:

\bc

$\hspace{0.01em}$

(+++ Orig. No.:  Claim Smooth-Admiss-2 +++)

{\xssc LABEL: {Claim Smooth-Admiss-2}}
\label{Claim Smooth-Admiss-2}

$ \xdz ' $ is $ \xdy -$smooth.

 karl-search= End Claim Smooth-Admiss-2
\vspace{7mm}

 *************************************

\vspace{7mm}

\subsubsection{Claim Smooth-Admiss-2 Proof}

 {\LARGE karl-search= Start Claim Smooth-Admiss-2 Proof }

\index{Claim Smooth-Admiss-2 Proof}

\ec

\subparagraph{
Proof
}

$\hspace{0.01em}$

(+++*** Orig.:  Proof )

Let $X \xbe \xdy,$ $<x, \xbs > \xbe \xdx ' \xex X.$

Case 1, $x \xbe X- \xbm (X):$ Then $ran( \xbs_{0}) \xcs \xbm (X) \xEd \xCQ
,$ let $x' \xbe ran( \xbs_{0}) \xcs \xbm (X).$ Moreover,
$ \xbm (X) \xcc K.$ Then for all $<x', \xbs ' > \xbe \xdx ' $ $<x', \xbs
' > \xeb <x, \xbs >.$ But $<x', \xbs^{x',X}>$ as
constructed in the proof of Claim \ref{Claim Smooth-Admiss-1} (page
\pageref{Claim Smooth-Admiss-1})
is minimal in $ \xdx ' \xex X.$

Case 2, $x \xbe \xbm (X)= \xbm_{ \xdz }(X)= \xbm ' (X):$ If $<x, \xbs >$
is minimal in $ \xdx ' \xex X,$ we are done.
So suppose there is $<x', \xbs ' > \xeb <x, \xbs >,$ $x' \xbe X.$ Thus
$x' \xbe \wt{ \xbs }.$ Let
$x' \xbe ran( \xbs_{i}).$ So $x \xbe \xbm (X)$ and $ran( \xbs_{i}) \xcs X
\xEd \xCQ.$ But
$ \xbs_{i+1} \xbe \xbP \{ \xbm (X' )$: $X' \xbe \xdy $ $ \xcu $ $x \xbe
\xbm (X' )$ $ \xcu $ $ran( \xbs_{i}) \xcs X' \xEd \xCQ \},$ so $X$ is one
of the $X',$
moreover $ \xbm (X) \xcc K,$ so there is $x'' \xbe \xbm (X) \xcs ran(
\xbs_{i+1}) \xcs K,$ so for all $<x'', \xbs '' > \xbe \xdx ' $
$<x'', \xbs '' > \xeb <x, \xbs >.$ But again $<x'', \xbs^{x'',X}>$ as
constructed in the proof of
Claim \ref{Claim Smooth-Admiss-1} (page \pageref{Claim Smooth-Admiss-1})  is
minimal in $ \xdx ' \xex X.$

$ \xcz $
\\[3ex]

 karl-search= End Claim Smooth-Admiss-2 Proof
\vspace{7mm}

 *************************************

\vspace{7mm}

 karl-search= End Proposition Smooth-Complete Proof
\vspace{7mm}

 *************************************

\vspace{7mm}

\newpage
\subsubsection{Proposition Smooth-Complete-Trans}

 {\LARGE karl-search= Start Proposition Smooth-Complete-Trans }

\index{Proposition Smooth-Complete-Trans}

\bp

$\hspace{0.01em}$

(+++ Orig. No.:  Proposition Smooth-Complete-Trans +++)

{\xssc LABEL: {Proposition Smooth-Complete-Trans}}
\label{Proposition Smooth-Complete-Trans}

Let $ \xbm: \xdy \xcp \xdp (U)$ satisfy $( \xbm \xcc ),$ $( \xbm PR),$
and $( \xbm CUM),$ and the domain $ \xdy $ $( \xcv ).$

Then there is a transitive $ \xdy -$smooth preferential structure $ \xdx $
s.t. $ \xbm = \xbm_{ \xdx }.$
See e.g.  \cite{Sch04}.

 karl-search= End Proposition Smooth-Complete-Trans
\vspace{7mm}

 *************************************

\vspace{7mm}

\subsubsection{Proposition Smooth-Complete-Trans Proof}

 {\LARGE karl-search= Start Proposition Smooth-Complete-Trans Proof }

\index{Proposition Smooth-Complete-Trans Proof}

\ep

\subparagraph{
Proof
}

$\hspace{0.01em}$

(+++*** Orig.:  Proof )

\subsubsection{Discussion Smooth-Trans}

 {\LARGE karl-search= Start Discussion Smooth-Trans }

\index{Discussion Smooth-Trans}

\paragraph{
Discussion Smooth-Trans
}

$\hspace{0.01em}$

(+++*** Orig.:  Discussion Smooth-Trans )

{\xssc LABEL: {Section Discussion Smooth-Trans}}
\label{Section Discussion Smooth-Trans}

In a certain way, it is not surprising that transitivity does not impose
stronger conditions in the smooth case either. Smoothness is itself a weak
kind of
transitivity: If an element is not minimal, then there is a minimal
element
below it, i.e., $x \xee y$ with $y$ not minimal is possible, there is $z'
\xeb y,$ but
then there is $z$ minimal with $x \xee z.$ This is ``almost'' $x \xee z',$
transitivity.

To obtain representation,
we will combine here the ideas of the smooth, but not necessarily
transitive
case with those of the general transitive case - as the reader will have
suspected. Thus, we will index again with trees, and work with (suitably
adapted) admissible sequences for the construction of the trees. In the
construction of the admissible sequences, we were careful to repair all
damage
done in previous steps. We have to add now reparation of all damage done
by
using transitivity, i.e., the transitivity of the relation might destroy
minimality, and we have to construct minimal elements below all elements
for
which we thus destroyed minimality. Both cases are combined by considering
immediately all $Y$ s.t. $x \xbe Y-H(U).$ Of course, the properties
described in
Fact \ref{Fact HU-2} (page \pageref{Fact HU-2})  play again a central role.

The (somewhat complicated) construction will be commented on in more
detail
below.

Note that even beyond Fact \ref{Fact HU-2} (page \pageref{Fact HU-2}) , closure
of the domain
under finite
unions is used in the construction of the trees. This - or something like
it -
is necessary, as we have to respect the hulls of all elements treated so
far
(the predecessors), and not only of the first element, because of
transitivity.
For the same reason, we need more bookkeeping, to annotate all the hulls
(or
the union of the respective $U' $s) of all predecessors to be respected.

To summarize: we combine the ideas from the transitive general case and
the
simple smooth case, using the crucial Fact \ref{Fact HU-2} (page \pageref{Fact
HU-2})
to show that the construction
goes through. The construction leaves still some freedom, and
modifications
are possible as indicated below in the course of the proof.

Recall that $ \xdy $ will be closed under finite unions
in this section, and let again $ \xbm: \xdy \xcp \xdp (Z).$

We have to adapt Construction \ref{Construction Smooth-Admiss} (page
\pageref{Construction Smooth-Admiss})
(x-admissible sequences) to the
transitive situation, and to our construction with trees. If $< \xCQ,x>$
is the root,
$ \xbs_{0} \xbe \xbP \{ \xbm (Y):x \xbe Y- \xbm (Y)\}$ determines some
children of the root.
To preserve smoothness, we have to compensate and add
other children by the $ \xbs_{i+1}:$ $ \xbs_{i+1} \xbe \xbP \{ \xbm (X):x
\xbe \xbm (X),$ $ran( \xbs_{i}) \xcs X \xEd \xCQ \}.$
On the other hand, we have to pursue the same construction for the
children so
constructed. Moreover, these indirect children have to be added to those
children
of the root, which have to be compensated (as the first children are
compensated
by $ \xbs_{1})$ to preserve smoothness. Thus, we build the tree in a
simultaneous vertical
and horizontal induction.

This construction can be simplified, by considering immediately all $Y
\xbe \xdy $ s.t.
$x \xbe Y \xcC H(U)$ - independent of whether $x \xce \xbm (Y)$ (as done
in $ \xbs_{0}),$ or whether
$x \xbe \xbm (Y),$ and some child $y$ constructed before is in $Y$ (as
done in the $ \xbs_{i+1}),$ or
whether $x \xbe \xbm (Y),$ and some indirect child $y$ of $x$ is in $Y$
(to take care of
transitivity, as indicated above). We make this simplified construction.

There are two ways to proceed. First, we can take as $ \xej^{*}$ in the
trees
the transitive closure of $ \xej.$ Second, we can deviate from the idea
that
children are chosen by selection functions $f,$ and take nonempty subsets
of
elements instead, making more elements children than in the first case. We
take
the first alternative, as it is more in the spirit of the construction.

We will suppose for simplicity that $Z=K$ - the general case in easy to
obtain,
but complicates the picture.

For each $x \xbe Z,$ we construct trees $t_{x}$, which will be used to
index
different copies of $x,$ and control the relation $ \xeb.$

These trees $t_{x}$ will have the following form:

(a) the root of $t$ is $< \xCQ,x>$ or $<U,x>$ with $U \xbe \xdy $ and $x
\xbe \xbm (U),$

(b) all other nodes are pairs $<Y,y>,$ $Y \xbe \xdy,$ $y \xbe \xbm (Y),$

(c) $ht(t) \xck \xbo,$

(d) if $<Y,y>$ is an element in $t_{x},$ then there is some $ \xdy (y)
\xcc \{W \xbe \xdy:y \xbe W\},$ and
$f \xbe \xbP \{ \xbm (W):W \xbe \xdy (y)\}$ s.t. the set of children of
$<Y,y>$ is $\{<Y \xcv W,f(W)>:$ $W \xbe \xdy (y)\}.$

The first coordinate is used for bookkeeping when constructing children,
in
particular for condition (d).

The relation $ \xeb $ will essentially be determined by the subtree
relation.

We first construct the trees $t_{x}$ for those sets $U$ where $x \xbe \xbm
(U),$ and then take
care of the others. In the construction for the minimal elements,
at each level $n>0,$ we may have several ways to choose a selection
function $f_{n}$,
and each such choice leads to the construction of a different tree - we
construct all these trees. (We could also construct only one tree, but
then
the choice would have to be made coherently for different $x,U.$ It is
simpler to
construct more trees than necessary.)

We control the relation by indexing with trees, just as it was done in the
not
necessarily smooth case before.

 karl-search= End Discussion Smooth-Trans
\vspace{7mm}

 *************************************

\vspace{7mm}

\subsubsection{Definition Tree-TC}

 {\LARGE karl-search= Start Definition Tree-TC }

\index{Definition Tree-TC}

\bd

$\hspace{0.01em}$

(+++ Orig. No.:  Definition Tree-TC +++)

{\xssc LABEL: {Definition Tree-TC}}
\label{Definition Tree-TC}

If $t$ is a tree with root $<a,b>,$ then t/c will be the same tree,
only with the root $<c,b>.$

 karl-search= End Definition Tree-TC
\vspace{7mm}

 *************************************

\vspace{7mm}

\subsubsection{Construction Smooth-Tree}

 {\LARGE karl-search= Start Construction Smooth-Tree }

\index{Construction Smooth-Tree}

\ed

\bcs

$\hspace{0.01em}$

(+++ Orig. No.:  Construction Smooth-Tree +++)

{\xssc LABEL: {Construction Smooth-Tree}}
\label{Construction Smooth-Tree}

(A) The set $T_{x}$ of trees $t$ for fixed $x$:

(1) Construction of the set $T \xbm_{x}$ of trees for those sets $U \xbe
\xdy,$ where $x \xbe \xbm (U):$

Let $U \xbe \xdy,$ $x \xbe \xbm (U).$ The trees $t_{U,x} \xbe T \xbm_{x}$
are constructed inductively,
observing simultaneously:

If $<U_{n+1},x_{n+1}>$ is a child of $<U_{n},x_{n}>,$ then
(a) $x_{n+1} \xbe \xbm (U_{n+1})-H(U_{n}),$
and
(b) $U_{n} \xcc U_{n+1}$.

Set $U_{0}:=U,$ $x_{0}:=x.$

Level 0: $<U_{0},x_{0}>.$

Level $n \xcp n+1$:
Let $<U_{n},x_{n}>$ be in level $n.$
Suppose $Y_{n+1} \xbe \xdy,$ $x_{n} \xbe Y_{n+1},$ and $Y_{n+1} \xcC
H(U_{n}).$ Note that $ \xbm (U_{n} \xcv Y_{n+1})-H(U_{n}) \xEd \xCQ $ by
Fact \ref{Fact HU-2} (page \pageref{Fact HU-2}) , (7), and $ \xbm (U_{n} \xcv
Y_{n+1})-H(U_{n})
\xcc \xbm (Y_{n+1})$
by Fact \ref{Fact HU-2} (page \pageref{Fact HU-2}) , (3).
Choose $f_{n+1} \xbe \xbP \{ \xbm (U_{n} \xcv Y_{n+1})-H(U_{n}):$ $Y_{n+1}
\xbe \xdy,$ $x_{n} \xbe Y_{n+1} \xcC H(U_{n})\}$ (for the construction
of this tree, at this element), and let the set of
children of $<U_{n},x_{n}>$ be $\{<U_{n} \xcv Y_{n+1},f_{n+1}(Y_{n+1})>:$
$Y_{n+1} \xbe \xdy,$ $x_{n} \xbe Y_{n+1} \xcC H(U_{n})\}.$
(If there is no such $Y_{n+1}$, $<U_{n},x_{n}>$ has no children.)
Obviously, (a) and (b) hold.

We call such trees $U,x-$trees.

(2) Construction of the set $T'_{x}$ of trees for the nonminimal elements.
Let $x \xbe Z.$ Construct the tree $t_{x}$ as follows (here, one tree per
$x$ suffices for
all $U)$:

Level 0: $< \xCQ,x>$

Level 1:
Choose arbitrary $f \xbe \xbP \{ \xbm (U):x \xbe U \xbe \xdy \}.$ Note
that $U \xEd \xCQ \xcp \xbm (U) \xEd \xCQ $ by $Z=K$ (by
Remark \ref{Remark Gamma-x} (page \pageref{Remark Gamma-x}) , (1)). Let
$\{<U,f(U)>:x \xbe U \xbe \xdy
\}$ be
the set of children of $< \xCQ,x>.$
This assures that the element will be nonminimal.

Level $>1$:
Let $<U,f(U)>$ be an element of level 1, as $f(U) \xbe \xbm (U),$ there is
a $t_{U,f(U)} \xbe T \xbm_{f(U)}.$
Graft one of these trees $t_{U,f(U)} \xbe T \xbm_{f(U)}$ at $<U,f(U)>$ on
the level 1.
This assures that a minimal element will be below it to guarantee
smoothness.

Finally, let $T_{x}:=T \xbm_{x} \xcv T'_{x}.$

(B) The relation $ \xej $ between trees:
For $x,y \xbe Z,$ $t \xbe T_{x}$, $t' \xbe T_{y}$, set $t \xem t' $ iff
for some $Y$ $<Y,y>$ is a child of
the root $<X,x>$ in $t,$ and $t' $ is the subtree of $t$ beginning at this
$<Y,y>.$

(C) The structure $ \xdz $:
Let $ \xdz $ $:=$ $<$ $\{<x,t_{x}>:$ $x \xbe Z,$ $t_{x} \xbe T_{x}\}$,
$<x,t_{x}> \xee <y,t_{y}>$ iff $t_{x} \xem^{*}t_{y}$ $>.$

 karl-search= End Construction Smooth-Tree
\vspace{7mm}

 *************************************

\vspace{7mm}

\subsubsection{Fact Smooth-Tree}

 {\LARGE karl-search= Start Fact Smooth-Tree }

\index{Fact Smooth-Tree}

\ecs

The rest of the proof are simple observations.

\bfa

$\hspace{0.01em}$

(+++ Orig. No.:  Fact Smooth-Tree +++)

{\xssc LABEL: {Fact Smooth-Tree}}
\label{Fact Smooth-Tree}

(1) If $t_{U,x}$ is an $U,x-$tree, $<U_{n},x_{n}>$ an element of $t_{U,x}$
, $<U_{m},x_{m}>$ a direct or
indirect child of $<U_{n},x_{n}>,$ then $x_{m} \xce H(U_{n}).$

(2) Let $<Y_{n},y_{n}>$ be an element in $t_{U,x} \xbe T \xbm_{x}$, $t' $
the subtree starting at $<Y_{n},y_{n}>,$
then $t' $ is a $Y_{n},y_{n}-tree.$

(3) $ \xeb $ is free from cycles.

(4) If $t_{U,x}$ is an $U,x-$tree, then $<x,t_{U,x}>$ is $ \xeb -$minimal
in $ \xdz \xex U.$

(5) No $<x,t_{x}>,$ $t_{x} \xbe T'_{x}$ is minimal in any $ \xdz \xex U,$
$U \xbe \xdy.$

(6) Smoothness is respected for the elements of the form $<x,t_{U,x}>.$

(7) Smoothness is respected for the elements of the form $<x,t_{x}>$ with
$t_{x} \xbe T'_{x}.$

(8) $ \xbm = \xbm_{ \xdz }.$

 karl-search= End Fact Smooth-Tree
\vspace{7mm}

 *************************************

\vspace{7mm}

\subsubsection{Fact Smooth-Tree Proof}

 {\LARGE karl-search= Start Fact Smooth-Tree Proof }

\index{Fact Smooth-Tree Proof}

\efa

\subparagraph{
Proof
}

$\hspace{0.01em}$

(+++*** Orig.:  Proof )

(1) trivial by (a) and (b).

(2) trivial by (a).

(3) Note that no $<x,t_{x}>$ $t_{x} \xbe T'_{x}$ can be smaller than any
other element (smaller
elements require $U \xEd \xCQ $ at the root). So no cycle involves any
such $<x,t_{x}>.$
Consider now $<x,t_{U,x}>,$ $t_{U,x} \xbe T \xbm_{x}$. For any
$<y,t_{V,y}> \xeb <x,t_{U,x}>,$ $y \xce H(U)$ by (1),
but $x \xbe \xbm (U) \xcc H(U),$ so $x \xEd y.$

(4) This is trivial by (1).

(5) Let $x \xbe U \xbe \xdy,$ then $f$ as used in the construction of
level 1 of $t_{x}$ chooses
$y \xbe \xbm (U) \xEd \xCQ,$ and some $<y,t_{U,y}>$ is in $ \xdz \xex U$
and below $<x,t_{x}>.$

(6) Let $x \xbe A \xbe \xdy,$ we have to show that either $<x,t_{U,x}>$
is minimal in $ \xdz \xex A,$ or that
there is $<y,t_{y}> \xeb <x,t_{U,x}>$ minimal in $ \xdz \xex A.$
Case 1, $A \xcc H(U)$: Then $<x,t_{U,x}>$ is minimal in $ \xdz \xex A,$
again by (1).
Case 2, $A \xcC H(U)$: Then A is one of the $Y_{1}$ considered for level
1. So there is
$<U \xcv A,f_{1}(A)>$ in level 1 with $f_{1}(A) \xbe \xbm (A) \xcc A$ by
Fact \ref{Fact HU-2} (page \pageref{Fact HU-2}) , (3).
But note that by (1)
all elements below $<U \xcv A,f_{1}(A)>$ avoid $H(U \xcv A).$ Let $t$ be
the subtree of $t_{U,x}$
beginning at $<U \xcv A,f_{1}(A)>,$ then by (2) $t$ is one of the $U \xcv
A,f_{1}(A)-trees,$ and
$<f_{1}(A),t>$ is minimal in $ \xdz \xex U \xcv A$ by (4), so in $ \xdz
\xex A,$ and $<f_{1}(A),t> \xeb <x,t_{U,x}>.$

(7) Let $x \xbe A \xbe \xdy,$ $<x,t_{x}>,$ $t_{x} \xbe T_{x}',$ and
consider the subtree $t$ beginning at $<A,f(A)>,$
then $t$ is one of the $A,f(A)-$trees, and $<f(A),t>$ is minimal in $ \xdz
\xex A$ by (4).

(8) Let $x \xbe \xbm (U).$ Then any $<x,t_{U,x}>$ is $ \xeb -$minimal in $
\xdz \xex U$ by (4), so $x \xbe \xbm_{ \xdz }(U).$
Conversely, let $x \xbe U- \xbm (U).$ By (5), no $<x,t_{x}>$ is minimal in
$U.$ Consider now some
$<x,t_{V,x}> \xbe \xdz,$ so $x \xbe \xbm (V).$ As $x \xbe U- \xbm (U),$
$U \xcC H(V)$ by Fact \ref{Fact HU-2} (page \pageref{Fact HU-2}) , (6).
Thus $U$ was
considered in the construction of level 1 of $t_{V,x}.$ Let $t$ be the
subtree of $t_{V,x}$
beginning at $<V \xcv U,f_{1}(U)>,$ by $ \xbm (V \xcv U)-H(V) \xcc \xbm
(U)$ (Fact \ref{Fact HU-2} (page \pageref{Fact HU-2}) , (3)),
$f_{1}(U) \xbe \xbm (U) \xcc U,$ and $<f_{1}(U),t> \xeb <x,t_{V,x}>.$

$ \xcz $
\\[3ex]

 karl-search= End Fact Smooth-Tree Proof
\vspace{7mm}

 *************************************

\vspace{7mm}

 karl-search= End Proposition Smooth-Complete-Trans Proof
\vspace{7mm}

 *************************************

\vspace{7mm}

\subsubsection{Proposition No-Norm}

 {\LARGE karl-search= Start Proposition No-Norm }

\index{Proposition No-Norm}

\bp

$\hspace{0.01em}$

(+++ Orig. No.:  Proposition No-Norm +++)

{\xssc LABEL: {Proposition No-Norm}}
\label{Proposition No-Norm}

We call a characterization ``normal'' iff it is a universally quantified
boolean combination (of any fixed, but perhaps infinite, length) of rules
of the usual form. (We do not go into details here.)

(1) There is no ``normal'' characterization of any fixed size of not
necessarily
definability preserving preferential structures.

(2) There is no ``normal'' characterization of any fixed size of not
necessarily
definability preserving ranked preferential structures.

See  \cite{Sch04}.

 karl-search= End Proposition No-Norm
\vspace{7mm}

 *************************************

\vspace{7mm}

 karl-search= End ToolBase1-Pref-ReprSmooth
\vspace{7mm}

 *************************************

\vspace{7mm}

 karl-search= End ToolBase1-Pref
\vspace{7mm}

 *************************************

\vspace{7mm}

\newpage

\section{
Ranked structures
}

\subsubsection{ToolBase1-Rank}

 {\LARGE karl-search= Start ToolBase1-Rank }

{\xssc LABEL: {Section Toolbase1-Rank}}
\label{Section Toolbase1-Rank}
\index{Section Toolbase1-Rank}

\subsection{
Ranked: Basics
}
\subsubsection{Fact Rank-Trans}

 {\LARGE karl-search= Start Fact Rank-Trans }

\index{Fact Rank-Trans}

\ep

\bfa

$\hspace{0.01em}$

(+++ Orig. No.:  Fact Rank-Trans +++)

{\xssc LABEL: {Fact Rank-Trans}}
\label{Fact Rank-Trans}

If $ \xeb $ on $X$ is ranked, and free of cycles, then $ \xeb $ is
transitive.

 karl-search= End Fact Rank-Trans
\vspace{7mm}

 *************************************

\vspace{7mm}

\subsubsection{Fact Rank-Trans Proof}

 {\LARGE karl-search= Start Fact Rank-Trans Proof }

\index{Fact Rank-Trans Proof}

\efa

\subparagraph{
Proof
}

$\hspace{0.01em}$

(+++*** Orig.:  Proof )

Let $x \xeb y \xeb z.$ If $x \xcT z,$ then $y \xee z,$ resulting in a
cycle of length 2. If $z \xeb x,$ then
we have a cycle of length 3. So $x \xeb z.$ $ \xcz $
\\[3ex]

 karl-search= End Fact Rank-Trans Proof
\vspace{7mm}

 *************************************

\vspace{7mm}

\subsubsection{Fact Rank-Auxil}

 {\LARGE karl-search= Start Fact Rank-Auxil }

\index{Fact Rank-Auxil}

\bfa

$\hspace{0.01em}$

(+++ Orig. No.:  Fact Rank-Auxil +++)

{\xssc LABEL: {Fact Rank-Auxil}}
\label{Fact Rank-Auxil}

$M(T)-M(T' )$ is normally not definable.

In the presence of $( \xbm =)$ and $( \xbm \xcc ),$ $f(Y) \xcs (X-f(X))
\xEd \xCQ $ is equivalent to
$f(Y) \xcs X \xEd \xCQ $ and $f(Y) \xcs f(X)= \xCQ.$

 karl-search= End Fact Rank-Auxil
\vspace{7mm}

 *************************************

\vspace{7mm}

\subsubsection{Fact Rank-Auxil Proof}

 {\LARGE karl-search= Start Fact Rank-Auxil Proof }

\index{Fact Rank-Auxil Proof}

\efa

\subparagraph{
Proof
}

$\hspace{0.01em}$

(+++*** Orig.:  Proof )

$f(Y) \xcs (X-f(X))$ $=$ $(f(Y) \xcs X)-(f(Y) \xcs f(X)).$

`` $ \xci $ '': Let $f(Y) \xcs X \xEd \xCQ,$ $f(Y) \xcs f(X)= \xCQ,$ so
$f(Y) \xcs (X-f(X)) \xEd \xCQ.$

`` $ \xch $ '': Suppose $f(Y) \xcs (X-f(X)) \xEd \xCQ,$ so $f(Y) \xcs X
\xEd \xCQ.$ Suppose $f(Y) \xcs f(X) \xEd \xCQ,$ so
by $( \xbm \xcc )$ $f(Y) \xcs X \xcs Y \xEd \xCQ,$ so
by $( \xbm =)$ $f(Y) \xcs X \xcs Y=f(X \xcs Y),$ and $f(X) \xcs X \xcs Y
\xEd \xCQ,$ so by $( \xbm =)$
$f(X) \xcs X \xcs Y=f(X \xcs Y),$ so $f(X) \xcs Y=f(Y) \xcs X$ and $f(Y)
\xcs (X-f(X))= \xCQ.$

$ \xcz $
\\[3ex]

 karl-search= End Fact Rank-Auxil Proof
\vspace{7mm}

 *************************************

\vspace{7mm}

\subsubsection{Remark RatM=}

 {\LARGE karl-search= Start Remark RatM= }

\index{Remark RatM=}

\br

$\hspace{0.01em}$

(+++ Orig. No.:  Remark RatM= +++)

{\xssc LABEL: {Remark RatM=}}
\label{Remark RatM=}

Note that $( \xbm =' )$ is very close to $ \xCf (RatM):$ $ \xCf (RatM)$
says:
$ \xba \xcn \xbb,$ $ \xba \xcN \xCN \xbg $ $ \xch $ $ \xba \xcu \xbg \xcn
\xbb.$ Or, $f(A) \xcc B,$ $f(A) \xcs C \xEd \xCQ $ $ \xcp $
$f(A \xcs C) \xcc B$ for all $A,B,C.$ This is not quite, but almost: $f(A
\xcs C) \xcc f(A) \xcs C$
(it depends how many $B$ there are, if $f(A)$ is some such $B,$ the fit is
perfect).

 karl-search= End Remark RatM=
\vspace{7mm}

 *************************************

\vspace{7mm}

\subsubsection{Fact Rank-Hold}

 {\LARGE karl-search= Start Fact Rank-Hold }

\index{Fact Rank-Hold}

\er

\bfa

$\hspace{0.01em}$

(+++ Orig. No.:  Fact Rank-Hold +++)

{\xssc LABEL: {Fact Rank-Hold}}
\label{Fact Rank-Hold}

In all ranked structures, $( \xbm \xcc ),$ $( \xbm =),$ $( \xbm PR),$ $(
\xbm =' ),$ $( \xbm \xFO ),$ $( \xbm \xcv ),$ $( \xbm \xcv ' ),$
$( \xbm \xbe ),$ $( \xbm RatM)$ will hold, if the corresponding closure
conditions are
satisfied.

 karl-search= End Fact Rank-Hold
\vspace{7mm}

 *************************************

\vspace{7mm}

\subsubsection{Fact Rank-Hold Proof}

 {\LARGE karl-search= Start Fact Rank-Hold Proof }

\index{Fact Rank-Hold Proof}

\efa

\subparagraph{
Proof
}

$\hspace{0.01em}$

(+++*** Orig.:  Proof )

$( \xbm \xcc )$ and $( \xbm PR)$ hold in all preferential structures.

$( \xbm =)$ and $( \xbm =' )$ are trivial.

$( \xbm \xcv )$ and $( \xbm \xcv ' ):$ All minimal copies of elements in
$f(Y)$ have the same rank.
If some $y \xbe f(Y)$ has all its minimal copies killed by an element $x
\xbe X,$ by
rankedness, $x$ kills the rest, too.

$( \xbm \xbe ):$ If $f(\{a\})= \xCQ,$ we are done. Take the minimal
copies of a in $\{a\},$ they are
all killed by one element in $X.$

$( \xbm \xFO ):$ Case $f(X)= \xCQ:$ If below every copy of $y \xbe Y$
there is a copy of some $x \xbe X,$
then $f(X \xcv Y)= \xCQ.$ Otherwise $f(X \xcv Y)=f(Y).$ Suppose now $f(X)
\xEd \xCQ,$ $f(Y) \xEd \xCQ,$ then
the minimal ranks decide: if they are equal, $f(X \xcv Y)=f(X) \xcv f(Y),$
etc.

$( \xbm RatM):$ Let $X \xcc Y,$ $y \xbe X \xcs f(Y) \xEd \xCQ,$ $x \xbe
f(X).$ By rankedness, $y \xeb x,$ or
$y \xcT x,$ $y \xeb x$ is impossible, as $y \xbe X,$ so $y \xcT x,$ and $x
\xbe f(Y).$

$ \xcz $
\\[3ex]

 karl-search= End Fact Rank-Hold Proof
\vspace{7mm}

 *************************************

\vspace{7mm}

 karl-search= End ToolBase1-Rank-Base
\vspace{7mm}

 *************************************

\vspace{7mm}

\newpage

\subsection{
Ranked: Representation
}
\subsubsection{ToolBase1-Rank-Repr}

 {\LARGE karl-search= Start ToolBase1-Rank-Repr }

{\xssc LABEL: {Section Toolbase1-Rank-Repr}}
\label{Section Toolbase1-Rank-Repr}
\index{Section Toolbase1-Rank-Repr}

\paragraph{
(1) Results for structures without copies
}

$\hspace{0.01em}$

(+++*** Orig.:  (1) Results for structures without copies )

{\xssc LABEL: {Section (1) Results for structures without copies}}
\label{Section (1) Results for structures without copies}

\subsubsection{Proposition Rank-Rep1}

 {\LARGE karl-search= Start Proposition Rank-Rep1 }

\index{Proposition Rank-Rep1}

\bp

$\hspace{0.01em}$

(+++ Orig. No.:  Proposition Rank-Rep1 +++)

{\xssc LABEL: {Proposition Rank-Rep1}}
\label{Proposition Rank-Rep1}

The first result applies for structures without copies of elements.

(1) Let $ \xdy \xcc \xdp (U)$ be closed under finite unions.
Then $( \xbm \xcc ),$ $( \xbm \xCQ ),$ $( \xbm =)$ characterize ranked
structures for which for all
$X \xbe \xdy $ $X \xEd \xCQ $ $ \xcp $ $ \xbm_{<}(X) \xEd \xCQ $ hold,
i.e. $( \xbm \xcc ),$ $( \xbm \xCQ ),$ $( \xbm =)$ hold in such
structures for $ \xbm_{<},$ and if they hold for some $ \xbm,$ we can
find a ranked relation
$<$ on $U$ s.t. $ \xbm = \xbm_{<}.$ Moreover, the structure can be choosen
$ \xdy -$smooth.

(2) Let $ \xdy \xcc \xdp (U)$ be closed under finite unions, and contain
singletons.
Then $( \xbm \xcc ),$ $( \xbm \xCQ fin),$ $( \xbm =),$ $( \xbm \xbe )$
characterize ranked structures for which
for all finite $X \xbe \xdy $ $X \xEd \xCQ $ $ \xcp $ $ \xbm_{<}(X) \xEd
\xCQ $ hold, i.e. $( \xbm \xcc ),$ $( \xbm \xCQ fin),$ $( \xbm =),$ $(
\xbm \xbe )$
hold in such structures for $ \xbm_{<},$ and if they hold for some $ \xbm
,$ we can find
a ranked relation $<$ on $U$ s.t. $ \xbm = \xbm_{<}.$

Note that the prerequisites of (2) hold in particular in the case
of ranked structures without copies, where all elements of $U$ are present
in the
structure - we need infinite descending chains to have $ \xbm (X)= \xCQ $
for $X \xEd \xCQ.$

See  \cite{Sch04}.

 karl-search= End Proposition Rank-Rep1
\vspace{7mm}

 *************************************

\vspace{7mm}

\ep

\paragraph{
(2) Results for structures possibly with copies
}

$\hspace{0.01em}$

(+++*** Orig.:  (2) Results for structures possibly with copies )

{\xssc LABEL: {Section (2) Results for structures possibly with copies}}
\label{Section (2) Results for structures possibly with copies}

\subsubsection{Definition 1-infin}

 {\LARGE karl-search= Start Definition 1-infin }

\index{Definition 1-infin}

\bd

$\hspace{0.01em}$

(+++ Orig. No.:  Definition 1-infin +++)

{\xssc LABEL: {Definition 1-infin}}
\label{Definition 1-infin}

Let $ \xdz =< \xdx, \xeb >$ be a preferential structure. Call $ \xdz $
$1- \xca $ over $Z,$
iff for all $x \xbe Z$ there are exactly one or infinitely many copies of
$x,$ i.e.
for all $x \xbe Z$ $\{u \xbe \xdx:$ $u=<x,i>$ for some $i\}$ has
cardinality 1 or $ \xcg \xbo.$

 karl-search= End Definition 1-infin
\vspace{7mm}

 *************************************

\vspace{7mm}

\subsubsection{Lemma 1-infin}

 {\LARGE karl-search= Start Lemma 1-infin }

\index{Lemma 1-infin}

\ed

\bl

$\hspace{0.01em}$

(+++ Orig. No.:  Lemma 1-infin +++)

{\xssc LABEL: {Lemma 1-infin}}
\label{Lemma 1-infin}

Let $ \xdz =< \xdx, \xeb >$ be a preferential structure and
$f: \xdy \xcp \xdp (Z)$ with $ \xdy \xcc \xdp (Z)$ be represented by $
\xdz,$ i.e. for $X \xbe \xdy $ $f(X)= \xbm_{ \xdz }(X),$
and $ \xdz $ be ranked and free of cycles. Then there is a structure $
\xdz ' $, $1- \xca $ over
$Z,$ ranked and free of cycles, which also represents $f.$

 karl-search= End Lemma 1-infin
\vspace{7mm}

 *************************************

\vspace{7mm}

\subsubsection{Lemma 1-infin Proof}

 {\LARGE karl-search= Start Lemma 1-infin Proof }

\index{Lemma 1-infin Proof}

\el

\subparagraph{
Proof
}

$\hspace{0.01em}$

(+++*** Orig.:  Proof )

We construct $ \xdz ' =< \xdx ', \xeb ' >.$

Let $A:=\{x \xbe Z$: there is some $<x,i> \xbe \xdx,$ but for all $<x,i>
\xbe \xdx $ there is
$<x,j> \xbe \xdx $ with $<x,j> \xeb <x,i>\},$

let $B:=\{x \xbe Z$: there is some $<x,i> \xbe \xdx,$ s.t. for no $<x,j>
\xbe \xdx $ $<x,j> \xeb <x,i>\},$

let $C:=\{x \xbe Z$: there is no $<x,i> \xbe \xdx \}.$

Let $c_{i}:i< \xbk $ be an enumeration of $C.$ We introduce for each such
$c_{i}$ $ \xbo $ many
copies $<c_{i},n>:n< \xbo $ into $ \xdx ',$ put all $<c_{i},n>$ above all
elements in $ \xdx,$ and order
the $<c_{i},n>$ by $<c_{i},n> \xeb ' <c_{i' },n' >$ $: \xcr $ $(i=i' $ and
$n>n' )$ or $i>i'.$ Thus, all $<c_{i},n>$ are
comparable.

If $a \xbe A,$ then there are infinitely many copies of a in $ \xdx,$ as
$ \xdx $ was
cycle-free, we put them all into $ \xdx '.$
If $b \xbe B,$ we choose exactly one such minimal element $<b,m>$ (i.e.
there
is no $<b,n> \xeb <b,m>)$ into $ \xdx ',$ and omit all other
elements. (For definiteness, assume in all applications $m=0.)$
For all elements from A and $B,$ we take the restriction of the order $
\xeb $ of $ \xdx.$
This is the new structure $ \xdz '.$

Obviously, adding the $<c_{i},n>$ does not introduce cycles, irreflexivity
and
rankedness are preserved. Moreover, any substructure of a cycle-free,
irreflexive,
ranked structure also has these properties, so $ \xdz ' $ is $1- \xca $
over $Z,$ ranked and
free of cycles.

We show that $ \xdz $ and $ \xdz ' $ are equivalent. Let then $X \xcc Z,$
we have to prove
$ \xbm (X)= \xbm ' (X)$ $( \xbm:= \xbm_{ \xdz }$, $ \xbm ':= \xbm_{
\xdz ' }).$

Let $z \xbe X- \xbm (X).$ If $z \xbe C$ or $z \xbe A,$ then $z \xce \xbm '
(X).$ If $z \xbe B,$
let $<z,m>$ be the chosen element. As $z \xce \xbm (X),$ there is $x \xbe
X$ s.t. some $<x,j> \xeb <z,m>.$
$x$ cannot be in $C.$ If $x \xbe A,$ then also $<x,j> \xeb ' <z,m>$. If
$x \xbe B,$ then there is some
$<x,k>$ also in $ \xdx '.$ $<x,j> \xeb <x,k>$ is impossible. If $<x,k>
\xeb <x,j>,$ then $<z,m> \xee <x,k>$
by transitivity. If $<x,k> \xcT <x,j>$, then also $<z,m> \xee <x,k>$ by
rankedness. In any
case, $<z,m> \xee ' <x,k>,$ and thus $z \xce \xbm ' (X).$

Let $z \xbe X- \xbm ' (X).$ If $z \xbe C$ or $z \xbe A,$ then $z \xce \xbm
(X).$ Let $z \xbe B,$ and some $<x,j> \xeb ' <z,m>.$
$x$ cannot be in $C,$ as they were sorted on top, so $<x,j>$ exists in $
\xdx $ too and
$<x,j> \xeb <z,m>.$ But if any other $<z,i>$ is also minimal in $ \xdz $
among the $<z,k>,$
then by rankedness also $<x,j> \xeb <z,i>,$ as $<z,i> \xcT <z,m>,$ so $z
\xce \xbm (X).$ $ \xcz $
\\[3ex]

 karl-search= End Lemma 1-infin Proof
\vspace{7mm}

 *************************************

\vspace{7mm}

\subsubsection{Fact Rank-No-Rep}

 {\LARGE karl-search= Start Fact Rank-No-Rep }

\index{Fact Rank-No-Rep}

\bfa

$\hspace{0.01em}$

(+++ Orig. No.:  Fact Rank-No-Rep +++)

{\xssc LABEL: {Fact Rank-No-Rep}}
\label{Fact Rank-No-Rep}

$( \xbm \xcc )+( \xbm PR)+( \xbm =)+( \xbm \xcv )+( \xbm \xbe )$ do not
imply representation by a ranked
structure.

 karl-search= End Fact Rank-No-Rep
\vspace{7mm}

 *************************************

\vspace{7mm}

\subsubsection{Fact Rank-No-Rep Proof}

 {\LARGE karl-search= Start Fact Rank-No-Rep Proof }

\index{Fact Rank-No-Rep Proof}

\efa

\subparagraph{
Proof
}

$\hspace{0.01em}$

(+++*** Orig.:  Proof )

See Example \ref{Example Rank-Copies} (page \pageref{Example Rank-Copies}) . $
\xcz $
\\[3ex]

 karl-search= End Fact Rank-No-Rep Proof
\vspace{7mm}

 *************************************

\vspace{7mm}

\subsubsection{Example Rank-Copies}

 {\LARGE karl-search= Start Example Rank-Copies }

\index{Example Rank-Copies}

\be

$\hspace{0.01em}$

(+++ Orig. No.:  Example Rank-Copies +++)

{\xssc LABEL: {Example Rank-Copies}}
\label{Example Rank-Copies}

This example shows that the conditions $( \xbm \xcc )+( \xbm PR)+( \xbm
=)+( \xbm \xcv )+( \xbm \xbe )$
can be satisfied, and still representation by a ranked structure
is impossible.

Consider $ \xbm (\{a,b\})= \xCQ,$ $ \xbm (\{a\})=\{a\},$ $ \xbm
(\{b\})=\{b\}.$ The conditions
$( \xbm \xcc )+( \xbm PR)+( \xbm =)+( \xbm \xcv )+( \xbm \xbe )$
hold trivially. This is representable, e.g. by $a_{1} \xed b_{1} \xed
a_{2} \xed b_{2} \Xl $ without
transitivity. (Note that rankedness implies transitivity,
$a \xec b \xec c,$ but not for $a=c.)$ But this cannot be represented by a
ranked
structure: As $ \xbm (\{a\}) \xEd \xCQ,$ there must be a copy $a_{i}$ of
minimal rank, likewise for
$b$ and some $b_{i}.$ If they have the same rank, $ \xbm
(\{a,b\})=\{a,b\},$ otherwise it will be
$\{a\}$ or $\{b\}.$

$ \xcz $
\\[3ex]

 karl-search= End Example Rank-Copies
\vspace{7mm}

 *************************************

\vspace{7mm}

\subsubsection{Proposition Rank-Rep2}

 {\LARGE karl-search= Start Proposition Rank-Rep2 }

\index{Proposition Rank-Rep2}

\ee

\bp

$\hspace{0.01em}$

(+++ Orig. No.:  Proposition Rank-Rep2 +++)

{\xssc LABEL: {Proposition Rank-Rep2}}
\label{Proposition Rank-Rep2}

Let $ \xdy $ be closed under finite unions and contain singletons. Then
$( \xbm \xcc )+( \xbm PR)+( \xbm \xFO )+( \xbm \xcv )+( \xbm \xbe )$
characterize ranked structures,
where elements may appear in several copies.

See  \cite{Sch04}.

 karl-search= End Proposition Rank-Rep2
\vspace{7mm}

 *************************************

\vspace{7mm}

 karl-search= End ToolBase1-Rank-Repr
\vspace{7mm}

 *************************************

\vspace{7mm}

 karl-search= End ToolBase1-Rank
\vspace{7mm}

 *************************************

\vspace{7mm}

\subsubsection{Proposition Rank-Rep3}

 {\LARGE karl-search= Start Proposition Rank-Rep3 }

\index{Proposition Rank-Rep3}

\ep

\bp

$\hspace{0.01em}$

(+++ Orig. No.:  Proposition Rank-Rep3 +++)

{\xssc LABEL: {Proposition Rank-Rep3}}
\label{Proposition Rank-Rep3}

Let $ \xdy \xcc \xdp (U)$ be closed under finite unions.
Then $( \xbm \xcc ),$ $( \xbm \xCQ ),$ $( \xbm =)$ characterize ranked
structures for which for all
$X \xbe \xdy $ $X \xEd \xCQ $ $ \xcp $ $ \xbm_{<}(X) \xEd \xCQ $ hold,
i.e. $( \xbm \xcc ),$ $( \xbm \xCQ ),$ $( \xbm =)$ hold in such
structures for $ \xbm_{<},$ and if they hold for some $ \xbm,$ we can
find a ranked relation
$<$ on $U$ s.t. $ \xbm = \xbm_{<}.$ Moreover, the structure can be choosen
$ \xdy -$smooth.

 karl-search= End Proposition Rank-Rep3
\vspace{7mm}

 *************************************

\vspace{7mm}

\subsubsection{Proposition Rank-Rep3 Proof}

 {\LARGE karl-search= Start Proposition Rank-Rep3 Proof }

\index{Proposition Rank-Rep3 Proof}

\ep

\subparagraph{
Proof
}

$\hspace{0.01em}$

(+++*** Orig.:  Proof )

Completeness:

Note that by Fact \ref{Fact Mu-Base} (page \pageref{Fact Mu-Base})  $(3)+(4)$
$( \xbm \xFO ),$ $( \xbm \xcv ),$ $( \xbm \xcv ' )$ hold.

Define aRb iff $ \xcE A \xbe \xdy (a \xbe \xbm (A),b \xbe A)$ or $a=b.$
$R$ is reflexive and transitive:
Suppose aRb, bRc, let $a \xbe \xbm (A),$ $b \xbe A,$ $b \xbe \xbm (B),$ $c
\xbe B.$ We show $a \xbe \xbm (A \xcv B).$ By
$( \xbm \xFO )$ $a \xbe \xbm (A \xcv B)$ or $b \xbe \xbm (A \xcv B).$
Suppose $b \xbe \xbm (A \xcv B),$ then $ \xbm (A \xcv B) \xcs A \xEd \xCQ
,$
so by $( \xbm =)$ $ \xbm (A \xcv B) \xcs A= \xbm (A),$ so $a \xbe \xbm (A
\xcv B).$

Moreover, $a \xbe \xbm (A),$ $b \xbe A- \xbm (A)$ $ \xcp $ $ \xCN (bRa):$
Suppose there is $B$ s.t. $b \xbe \xbm (B),$
$a \xbe B.$ Then by $( \xbm \xcv )$ $ \xbm (A \xcv B) \xcs B= \xCQ,$ and
by $( \xbm \xcv ' )$ $ \xbm (A \xcv B)= \xbm (A),$ but
$a \xbe \xbm (A) \xcs B,$ $contradiction.$

Let by Lemma \ref{Lemma Abs-Rel-Ext} (page \pageref{Lemma Abs-Rel-Ext})  $S$ be
a total, transitive, reflexive relation on $U$ which extends
$R$ s.t. $xSy,ySx$ $ \xcp $ xRy (recall that $R$ is transitive and
reflexive). Define $a<b$
iff aSb, but not bSa. If $a \xcT b$ (i.e. neither $a<b$ nor $b<a),$ then,
by totality of
$S,$ aSb and bSa. $<$ is ranked: If $c<a \xcT b,$ then by transitivity of
$S$ cSb, but
if bSc, then again by transitivity of $S$ aSc. Similarly for $c>a \xcT b.$

$<$ represents $ \xbm $ and is $ \xdy -$smooth: Let $a \xbe A- \xbm (A).$
By $( \xbm \xCQ ),$ $ \xcE b \xbe \xbm (A),$ so bRa,
but (by above argument) not aRb, so bSa, but not aSb, so $b<a,$ so $a \xbe
A- \xbm_{<}(A),$
and, as $b$ will then be $<-$minimal (see the next sentence), $<$ is $
\xdy -$smooth. Let
$a \xbe \xbm (A),$ then for all $a' \xbe A$ aRa', so aSa', so there is
no $a' \xbe A$ $a' <a,$ so $a \xbe \xbm_{<}(A).$

$ \xcz $
\\[3ex]

 karl-search= End Proposition Rank-Rep3 Proof
\vspace{7mm}

 *************************************

\vspace{7mm}

\newpage

\section{
Theory revision
}

\subsubsection{ToolBase1-TR}

 {\LARGE karl-search= Start ToolBase1-TR }

{\xssc LABEL: {Section Toolbase1-TR}}
\label{Section Toolbase1-TR}
\index{Section Toolbase1-TR}

\subsection{
AGM revision
}
\subsubsection{ToolBase1-TR-AGM}

 {\LARGE karl-search= Start ToolBase1-TR-AGM }

{\xssc LABEL: {Section Toolbase1-TR-AGM}}
\label{Section Toolbase1-TR-AGM}
\index{Section Toolbase1-TR-AGM}
{\xssc LABEL: {Section AGM-revision}}
\label{Section AGM-revision}

All material in this Section \ref{Section AGM-revision} (page \pageref{Section
AGM-revision})  is due verbatim
or in essence
to AGM - AGM for Alchourron, Gardenfors, Makinson, see e.g.  \cite{AGM85}.
\subsubsection{Definition AGM}

 {\LARGE karl-search= Start Definition AGM }

\index{Definition AGM}

\bd

$\hspace{0.01em}$

(+++ Orig. No.:  Definition AGM +++)

{\xssc LABEL: {Definition AGM}}
\label{Definition AGM}

We present in parallel the logical
and the semantic (or purely algebraic) side. For the latter, we work in
some
fixed universe $U,$ and the intuition is $U=M_{ \xdl },$ $X=M(K),$ etc.,
so, e.g. $A \xbe K$
becomes $X \xcc B,$ etc.

(For reasons of readability, we omit most caveats about definability.)

$K_{ \xcT }$ will denote the inconsistent theory.

We consider two functions, - and $*,$ taking a deductively closed theory
and a
formula as arguments, and returning a (deductively closed) theory on the
logics
side. The algebraic counterparts work on definable model sets. It is
obvious
that $ \xCf (K-1),$ $(K*1),$ $ \xCf (K-6),$ $(K*6)$ have vacuously true
counterparts on the
semantical side. Note that $K$ $ \xCf (X)$ will never change, everything
is relative
to fixed $K$ $ \xCf (X).$ $K* \xbf $ is the result of revising $K$ with $
\xbf.$ $K- \xbf $ is the result of
subtracting enough from $K$ to be able to add $ \xCN \xbf $ in a
reasonable way, called
contraction.

Moreover,
let $ \xck_{K}$ be a relation on the formulas relative to a deductively
closed theory $K$
on the formulas of $ \xdl,$ and $ \xck_{X}$ a relation on $ \xdp (U)$ or
a suitable subset of $ \xdp (U)$
relative to fixed $X.$ When the context is clear, we simply write $ \xck
.$
$ \xck_{K}$ $( \xck_{X})$ is called a relation of epistemic entrenchment
for $K$ $ \xCf (X).$

The following table presents the ``rationality postulates'' for contraction
(-),
revision $(*)$ and epistemic entrenchment. In AGM tradition, $K$ will be a
deductively closed theory, $ \xbf, \xbq $ formulas. Accordingly, $X$ will
be the set of
models of a theory, $A,B$ the model sets of formulas.

In the further development, formulas $ \xbf $ etc. may sometimes also be
full
theories. As the transcription to this case is evident, we will not go
into
details.

\renewcommand{\arraystretch}{1.2}

{\small
\begin{tabular}{|c|c|c|c|}

\hline

\multicolumn{4}{|c|} {Contraction, $K-\xbf $} \xEP

\hline

$(K-1)$ \xEH $K-\xbf $ is deductively closed \xEH \xEH \xEP

\hline

$(K-2)$ \xEH $K-\xbf $ $ \xcc $ $K$ \xEH $(X \xDN 2)$ \xEH $X \xcc X \xDN A$
\xEP

\hline

$(K-3)$ \xEH $\xbf  \xce K$ $ \xch $ $K-\xbf =K$ \xEH $(X \xDN 3)$ \xEH $X \xcC
A$ $
\xch $ $X \xDN A=X$ \xEP

\hline

$(K-4)$ \xEH $ \xcL \xbf $ $ \xch $ $\xbf  \xce K-\xbf $ \xEH $(X \xDN 4)$ \xEH
$A \xEd
U$ $ \xch $ $X \xDN A \xcC A$ \xEP

\hline

$(K-5)$ \xEH $K \xcc \ol{(K-\xbf ) \xcv \{\xbf \}}$ \xEH $(X \xDN 5)$ \xEH $(X
\xDN
A) \xcs A$ $ \xcc $ $X$ \xEP

\hline

$(K-6)$ \xEH $ \xcl \xbf  \xcr \xbq $ $ \xch $ $K-\xbf =K-\xbq $ \xEH \xEH \xEP

\hline

$(K-7)$ \xEH $(K-\xbf ) \xcs (K-\xbq )  \xcc  $ \xEH
$(X \xDN 7)$ \xEH $X \xDN (A \xcs B)  \xcc  $ \xEP

\xEH $K-(\xbf  \xcu \xbq ) $ \xEH
\xEH $(X \xDN A) \xcv (X \xDN B)$ \xEP

\hline

$(K-8)$ \xEH $\xbf  \xce K-(\xbf  \xcu \xbq )  \xch  $ \xEH
$(X \xDN 8)$ \xEH $X \xDN (A \xcs B) \xcC A  \xch  $ \xEP

\xEH $K-(\xbf  \xcu \xbq ) \xcc K-\xbf $ \xEH
\xEH $X \xDN A \xcc X \xDN (A \xcs B)$ \xEP

\hline
\hline

\multicolumn{4}{|c|} {Revision, $K*\xbf $} \xEP

\hline

$(K*1)$ \xEH $K*\xbf $ is deductively closed \xEH - \xEH \xEP

\hline

$(K*2)$ \xEH $\xbf  \xbe K*\xbf $ \xEH $(X \xfA 2)$ \xEH $X \xfA A \xcc A$
\xEP

\hline

$(K*3)$ \xEH $K*\xbf $ $ \xcc $ $ \ol{K \xcv \{\xbf \}}$ \xEH $(X \xfA 3)$ \xEH
$X
\xcs A \xcc X \xfA A$ \xEP

\hline

$(K*4)$ \xEH $ \xCN \xbf  \xce K  \xch $ \xEH
$(X \xfA 4)$ \xEH $X \xcs A \xEd \xCQ   \xch $ \xEP

\xEH $\ol{K \xcv \{\xbf \}}  \xcc  K*\xbf $ \xEH
\xEH $X \xfA A \xcc X \xcs A$ \xEP

\hline

$(K*5)$ \xEH $K*\xbf =K_{ \xcT }$ $ \xch $ $ \xcl \xCN \xbf $ \xEH $(X \xfA 5)$
\xEH $X \xfA A= \xCQ $ $ \xch $ $A= \xCQ $ \xEP

\hline

$(K*6)$ \xEH $ \xcl \xbf  \xcr \xbq $ $ \xch $ $K*\xbf =K*\xbq $ \xEH - \xEH
\xEP

\hline

$(K*7)$ \xEH $K*(\xbf  \xcu \xbq )  \xcc $ \xEH
$(X \xfA 7)$ \xEH $(X \xfA A) \xcs B  \xcc  $ \xEP

\xEH $\ol{(K*\xbf ) \xcv \{\xbq \}}$ \xEH
\xEH $X \xfA (A \xcs B)$ \xEP

\hline

$(K*8)$ \xEH $ \xCN \xbq  \xce K*\xbf  \xch $ \xEH
$(X \xfA 8)$ \xEH $(X \xfA A) \xcs B \xEd \xCQ \xch $ \xEP

\xEH $\ol{(K*\xbf ) \xcv \{\xbq \}} \xcc K*(\xbf  \xcu \xbq )$ \xEH
\xEH $ X \xfA (A \xcs B) \xcc (X \xfA A) \xcs B$ \xEP

\hline
\hline

\multicolumn{4}{|c|} {Epistemic entrenchment} \xEP

\hline

$(EE1)$ \xEH $ \xck_{K}$ is transitive \xEH
$(EE1)$ \xEH $ \xck_{X}$ is transitive \xEP

\hline

$(EE2)$ \xEH $\xbf  \xcl \xbq   \xch  \xbf  \xck_{K}\xbq $ \xEH
$(EE2)$ \xEH $A \xcc B  \xch  A \xck_{X}B$ \xEP

\hline

$(EE3)$ \xEH $ \xcA  \xbf,\xbq  $ \xEH
$(EE3)$ \xEH $ \xcA A,B $ \xEP
\xEH $ (\xbf  \xck_{K}\xbf  \xcu \xbq $ or $\xbq  \xck_{K}\xbf  \xcu \xbq )$
\xEH
\xEH $ (A \xck_{X}A \xcs B$ or $B \xck_{X}A \xcs B)$ \xEP

\hline

$(EE4)$ \xEH $K \xEd K_{ \xcT }  \xch  $ \xEH
$(EE4)$ \xEH $X \xEd \xCQ   \xch  $ \xEP
\xEH $(\xbf  \xce K$ iff $ \xcA  \xbq.\xbf  \xck_{K}\xbq )$ \xEH
\xEH $(X \xcC A$ iff $ \xcA B.A \xck_{X}B)$ \xEP

\hline

$(EE5)$ \xEH $ \xcA \xbq.\xbq  \xck_{K}\xbf   \xch   \xcl \xbf $ \xEH
$(EE5)$ \xEH $ \xcA B.B \xck_{X}A  \xch  A=U$ \xEP

\hline

\end{tabular}
}
\\

 karl-search= End Definition AGM
\vspace{7mm}

 *************************************

\vspace{7mm}

\subsubsection{Remark TR-Rank}

 {\LARGE karl-search= Start Remark TR-Rank }

\index{Remark TR-Rank}

\ed

\br

$\hspace{0.01em}$

(+++ Orig. No.:  Remark TR-Rank +++)

{\xssc LABEL: {Remark TR-Rank}}
\label{Remark TR-Rank}

(1) Note that $(X \xfA 7)$ and $(X \xfA 8)$ express a central condition
for ranked
structures, see Section 3.10: If we note $X \xfA.$ by $f_{X}(.),$ we then
have:
$f_{X}(A) \xcs B \xEd \xCQ $ $ \xch $ $f_{X}(A \xcs B)=f_{X}(A) \xcs B.$

(2) It is trivial to see that AGM revision cannot be defined by an
individual
distance (see Definition 2.3.5 below):
Suppose $X \xfA Y$ $:=$ $\{y \xbe Y:$ $ \xcE x_{y} \xbe X( \xcA y' \xbe
Y.d(x_{y},y) \xck d(x_{y},y' ))\}.$
Consider $a,b,c.$ $\{a,b\} \xfA \{b,c\}=\{b\}$ by $(X \xfA 3)$ and $(X
\xfA 4),$ so $d(a,b)<d(a,c).$
But on the other hand $\{a,c\} \xfA \{b,c\}=\{c\},$ so $d(a,b)>d(a,c),$
$contradiction.$

 karl-search= End Remark TR-Rank
\vspace{7mm}

 *************************************

\vspace{7mm}

\subsubsection{Proposition AGM-Equiv}

 {\LARGE karl-search= Start Proposition AGM-Equiv }

\index{Proposition AGM-Equiv}

\er

\bp

$\hspace{0.01em}$

(+++ Orig. No.:  Proposition AGM-Equiv +++)

{\xssc LABEL: {Proposition AGM-Equiv}}
\label{Proposition AGM-Equiv}

Contraction, revision, and epistemic entrenchment are interdefinable by
the
following equations, i.e., if the defining side has the respective
properties,
so will the defined side.

\renewcommand{\arraystretch}{1.5}

{\scriptsize
\begin{tabular}{|c|c|}

\hline

$K*\xbf:= \ol{(K- \xCN \xbf )} \xcv {\xbf }$  \xEH  $X \xfA A:= (X \xDN  \xdC A)
\xcs A$ \xEP

\hline

$K-\xbf:= K \xcs (K* \xCN \xbf )$  \xEH  $X \xDN A:= X \xcv (X \xfA  \xdC A)$
\xEP

\hline

$K-\xbf:=\{\xbq  \xbe K:$ $(\xbf <_{K}\xbf  \xco \xbq $ or $ \xcl \xbf )\}$ \xEH

$
X \xDN A:=
\left\{
\begin{array}{rcl}
X & iff & A=U, \\
 \xcS \{B: X \xcc B \xcc U, A<_{X}A \xcv B\} & & otherwise \\
\end{array}
\right. $
\xEP

\hline

$
\xbf  \xck_{K}\xbq: \xcr
\left\{
\begin{array}{l}
\xcl \xbf  \xcu \xbq  \\
or \\
\xbf  \xce K-(\xbf  \xcu \xbq ) \\
\end{array}
\right. $
\xEH

$
A \xck_{X}B: \xcr
\left\{
\begin{array}{l}
A,B=U  \\
or \\
X \xDN (A \xcs B) \xcC A \\
\end{array}
\right. $
\xEP

\hline

\end{tabular}
}

 karl-search= End Proposition AGM-Equiv
\vspace{7mm}

 *************************************

\vspace{7mm}

\subsubsection{Intuit-Entrench}

 {\LARGE karl-search= Start Intuit-Entrench }

\index{Intuit-Entrench}
\paragraph{A remark on intuition}

\ep

The idea of epistemic entrenchment is that $ \xbf $ is more entrenched
than $ \xbq $
(relative to $K)$ iff $M( \xCN \xbq )$ is closer to $M(K)$ than $M( \xCN
\xbf )$ is to $M(K).$ In
shorthand, the more we can twiggle $K$ without reaching $ \xCN \xbf,$ the
more $ \xbf $ is
entrenched. Truth is maximally entrenched - no twiggling whatever will
reach
falsity. The more $ \xbf $ is entrenched,
the more we are certain about it. Seen this way, the properties of
epistemic
entrenchment relations are very natural (and trivial): As only the closest
points of $M( \xCN \xbf )$ count (seen from $M(K)),$ $ \xbf $ or $ \xbq $
will be as entrenched as
$ \xbf \xcu \xbq,$ and there is a logically strongest $ \xbf ' $ which is
as entrenched as $ \xbf $ -
this is just the sphere around $M(K)$ with radius $d(M(K),M( \xCN \xbf
)).$

 karl-search= End Intuit-Entrench
\vspace{7mm}

 *************************************

\vspace{7mm}

 karl-search= End ToolBase1-TR-AGM
\vspace{7mm}

 *************************************

\vspace{7mm}

\newpage

\subsection{
Distance based revision: Basics
}

\subsubsection{ToolBase1-TR-DistBase}

 {\LARGE karl-search= Start ToolBase1-TR-DistBase }

{\xssc LABEL: {Section Toolbase1-TR-DistBase}}
\label{Section Toolbase1-TR-DistBase}
\index{Section Toolbase1-TR-DistBase}
\subsubsection{Definition Distance}

 {\LARGE karl-search= Start Definition Distance }

\index{Definition Distance}

\bd

$\hspace{0.01em}$

(+++ Orig. No.:  Definition Distance +++)

{\xssc LABEL: {Definition Distance}}
\label{Definition Distance}

$d:U \xCK U \xcp Z$ is called a pseudo-distance on $U$ iff (d1) holds:

(d1) $Z$ is totally ordered by a relation $<.$

If, in addition, $Z$ has a $<-$smallest element 0, and (d2) holds, we say
that $d$
respects identity:

(d2) $d(a,b)=0$ iff $a=b.$

If, in addition, (d3) holds, then $d$ is called symmetric:

(d3) $d(a,b)=d(b,a).$

(For any $a,b \xbe U.)$

Note that we can force the triangle inequality to hold trivially (if we
can
choose the values in the real numbers): It suffices to choose the values
in
the set $\{0\} \xcv [0.5,1],$ i.e. in the interval from 0.5 to 1, or as 0.

 karl-search= End Definition Distance
\vspace{7mm}

 *************************************

\vspace{7mm}

\subsubsection{Definition Dist-Indiv-Coll}

 {\LARGE karl-search= Start Definition Dist-Indiv-Coll }

\index{Definition Dist-Indiv-Coll}

\ed

\bd

$\hspace{0.01em}$

(+++ Orig. No.:  Definition Dist-Indiv-Coll +++)

{\xssc LABEL: {Definition Dist-Indiv-Coll}}
\label{Definition Dist-Indiv-Coll}

We define the collective and the individual variant of choosing the
closest
elements in the second operand by two operators,
$ \xfA, \xfB: \xdp (U) \xCK \xdp (U) \xcp \xdp (U):$

Let $d$ be a distance or pseudo-distance.

$X \xfA Y$ $:=$ $\{y \xbe Y:$ $ \xcE x_{y} \xbe X. \xcA x' \xbe X, \xcA y'
\xbe Y(d(x_{y},y) \xck d(x',y' )\}$

(the collective variant, used in theory revision)

and

$X \xfB Y$ $:=$ $\{y \xbe Y:$ $ \xcE x_{y} \xbe X. \xcA y' \xbe
Y(d(x_{y},y) \xck d(x_{y},y' )\}$

(the individual variant, used for counterfactual conditionals and theory
update).

Thus, $A \xfA_{d}B$ is the subset of $B$ consisting of all $b \xbe B$ that
are closest to A.
Note that, if $ \xCf A$ or $B$ is infinite, $A \xfA_{d}B$ may be empty,
even if $ \xCf A$ and $B$ are not
empty. A condition assuring nonemptiness will be imposed when necessary.

 karl-search= End Definition Dist-Indiv-Coll
\vspace{7mm}

 *************************************

\vspace{7mm}

\subsubsection{Definition Dist-Repr}

 {\LARGE karl-search= Start Definition Dist-Repr }

\index{Definition Dist-Repr}

\ed

\bd

$\hspace{0.01em}$

(+++ Orig. No.:  Definition Dist-Repr +++)

{\xssc LABEL: {Definition Dist-Repr}}
\label{Definition Dist-Repr}

An operation $ \xfA: \xdp (U) \xCK \xdp (U) \xcp \xdp (U)$ is
representable iff there is a
pseudo-distance $d:U \xCK U \xcp Z$ such that

$A \xfA B$ $=$ $A \xfA_{d}B$ $:=$ $\{b \xbe B:$ $ \xcE a_{b} \xbe A \xcA
a' \xbe A \xcA b' \xbe B(d(a_{b},b) \xck d(a',b' ))\}.$

 karl-search= End Definition Dist-Repr
\vspace{7mm}

 *************************************

\vspace{7mm}

\subsubsection{Definition TR*d}

 {\LARGE karl-search= Start Definition TR*d }

\index{Definition TR*d}

\ed

The following is the central definition, it describes the way a revision
$*_{d}$ is
attached to a pseudo-distance $d$ on the set of models.

\bd

$\hspace{0.01em}$

(+++ Orig. No.:  Definition TR*d +++)

{\xssc LABEL: {Definition TR*d}}
\label{Definition TR*d}

$T*_{d}T' $ $:=$ $Th(M(T) \xfA_{d}M(T' )).$

$*$ is called representable iff there is a pseudo-distance $d$ on the set
of models
s.t. $T*T' =Th(M(T) \xfA_{d}M(T' )).$

 karl-search= End Definition TR*d
\vspace{7mm}

 *************************************

\vspace{7mm}

 karl-search= End ToolBase1-TR-DistBase
\vspace{7mm}

 *************************************

\vspace{7mm}


\subsection{
Distance based revision: Representation
}

\subsubsection{ToolBase1-TR-DistRepr}

 {\LARGE karl-search= Start ToolBase1-TR-DistRepr }

{\xssc LABEL: {Section Toolbase1-TR-DistRepr}}
\label{Section Toolbase1-TR-DistRepr}
\index{Section Toolbase1-TR-DistRepr}
\subsubsection{Fact AGM-In-Dist}

 {\LARGE karl-search= Start Fact AGM-In-Dist }

\index{Fact AGM-In-Dist}

\ed

\bfa

$\hspace{0.01em}$

(+++ Orig. No.:  Fact AGM-In-Dist +++)

{\xssc LABEL: {Fact AGM-In-Dist}}
\label{Fact AGM-In-Dist}

A distance based revision satisfies the AGM postulates provided:

(1) it respects identity, i.e. $d(a,a)<d(a,b)$ for all $a \xEd b,$

(2) it satisfies a limit condition: minima exist,

(3) it is definability preserving.

(It is trivial to see that the first two are necessary,
and Example \ref{Example TR-Dp} (page \pageref{Example TR-Dp})  (2)
below shows the necessity of (3). In particular, (2) and (3) will hold for
finite languages.)

 karl-search= End Fact AGM-In-Dist
\vspace{7mm}

 *************************************

\vspace{7mm}

\subsubsection{Fact AGM-In-Dist Proof}

 {\LARGE karl-search= Start Fact AGM-In-Dist Proof }

\index{Fact AGM-In-Dist Proof}

\efa

\subparagraph{
Proof
}

$\hspace{0.01em}$

(+++*** Orig.:  Proof )

We use $ \xfA $ to abbreviate $ \xfA_{d}.$ As a matter of fact, we show
slightly more, as
we admit also full theories on the right of $*.$

$(K*1),$ $(K*2),$ $(K*6)$ hold by definition, $(K*3)$ and $(K*4)$ as $d$
respects
identity, $(K*5)$ by existence of minima.

It remains to show $(K*7)$ and $(K*8),$ we do them together, and show:
If $T*T' $ is consistent with $T'',$ then $T*(T' \xcv T'' )$ $=$ $
\ol{(T*T' ) \xcv T'' }.$

Note that $M(S \xcv S' )=M(S) \xcs M(S' ),$ and that $M(S*S' )=M(S) \xfA
M(S' ).$ (The latter
is only true if $ \xfA $ is definability preserving.)
By prerequisite, $M(T*T' ) \xcs M(T'' ) \xEd \xCQ,$ so $(M(T) \xfA M(T'
)) \xcs M(T'' ) \xEd \xCQ.$
Let $A:=M(T),$ $B:=M(T' ),$ $C:=M(T'' ).$ `` $ \xcc $ '': Let $b \xbe A
\xfA (B \xcs C).$
By prerequisite, there is $b' \xbe (A \xfA B) \xcs C.$ Thus $d(A,b' ) \xcg
d(A,B \xcs C)=d(A,b).$
As $b \xbe B,$ $b \xbe A \xfA B,$ but $b \xbe C,$ too. `` $ \xcd $ '': Let
$b' \xbe (A \xfA B) \xcs C.$ Thus $d(A,b' )=$
$d(A,B) \xck d(A,B \xcs C),$ so by $b' \xbe B \xcs C$ $b' \xbe A \xfA (B
\xcs C).$
We conclude $M(T) \xfA (M(T' ) \xcs M(T'' ))$ $=$ $(M(T) \xfA M(T' )) \xcs
M(T'' ),$ thus that
$T*(T' \xcv T'' )= \ol{(T*T' ) \xcv T'' }.$

$ \xcz $
\\[3ex]

 karl-search= End Fact AGM-In-Dist Proof
\vspace{7mm}

 *************************************

\vspace{7mm}

\subsubsection{Definition TR-Umgeb}

 {\LARGE karl-search= Start Definition TR-Umgeb }

\index{Definition TR-Umgeb}

\bd

$\hspace{0.01em}$

(+++ Orig. No.:  Definition TR-Umgeb +++)

{\xssc LABEL: {Definition TR-Umgeb}}
\label{Definition TR-Umgeb}

For $X,Y \xEd \xCQ,$ set $U_{Y}(X):=\{z:d(X,z) \xck d(X,Y)\}.$

 karl-search= End Definition TR-Umgeb
\vspace{7mm}

 *************************************

\vspace{7mm}

\subsubsection{Fact TR-Umgeb}

 {\LARGE karl-search= Start Fact TR-Umgeb }

\index{Fact TR-Umgeb}

\ed

\bfa

$\hspace{0.01em}$

(+++ Orig. No.:  Fact Tr-Umgeb +++)

{\xssc LABEL: {Fact Tr-Umgeb}}
\label{Fact Tr-Umgeb}

Let $X,Y,Z \xEd \xCQ.$ Then

(1) $U_{Y}(X) \xcs Z \xEd \xCQ $ iff $(X \xfA (Y \xcv Z)) \xcs Z \xEd \xCQ
,$

(2) $U_{Y}(X) \xcs Z \xEd \xCQ $ iff $ \xdC Z \xck_{X} \xdC Y$ - where $
\xck_{X}$ is epistemic entrenchement relative
to $X.$

 karl-search= End Fact TR-Umgeb
\vspace{7mm}

 *************************************

\vspace{7mm}

\subsubsection{Fact TR-Umgeb Proof}

 {\LARGE karl-search= Start Fact TR-Umgeb Proof }

\index{Fact TR-Umgeb Proof}

\efa

\subparagraph{
Proof
}

$\hspace{0.01em}$

(+++*** Orig.:  Proof )

(1) Trivial.

(2) $ \xdC Z \xck_{X} \xdC Y$ iff $X \xDN ( \xdC Z \xcs \xdC Y) \xcC \xdC
Z.$ $X \xDN ( \xdC Z \xcs \xdC Y)$ $=$ $X \xcv (X \xfA \xdC ( \xdC Z \xcs
\xdC Y))$ $=$
$X \xcv (X \xfA (Z \xcv Y)).$ So $X \xDN ( \xdC Z \xcs \xdC Y) \xcC \xdC
Z$ $ \xcj $ $(X \xcv (X \xfA (Z \xcv Y))) \xcs Z \xEd \xCQ $ $ \xcj $
$X \xcs Z \xEd \xCQ $ or $(X \xfA (Z \xcv Y)) \xcs Z \xEd \xCQ $ $ \xcj $
$d(X,Z) \xck d(X,Y).$

$ \xcz $
\\[3ex]

 karl-search= End Fact TR-Umgeb Proof
\vspace{7mm}

 *************************************

\vspace{7mm}

\subsubsection{Definition TR-Dist}

 {\LARGE karl-search= Start Definition TR-Dist }

\index{Definition TR-Dist}

\bd

$\hspace{0.01em}$

(+++ Orig. No.:  Definition TR-Dist +++)

{\xssc LABEL: {Definition TR-Dist}}
\label{Definition TR-Dist}

Let $U \xEd \xCQ,$ $ \xdy \xcc \xdp (U)$ satisfy $( \xcs ),$ $( \xcv ),$
$ \xCQ \xce \xdy.$

Let $A,B,X_{i} \xbe \xdy,$ $ \xfA: \xdy \xCK \xdy \xcp \xdp (U).$

Let $*$ be a revision function defined for
arbitrary consistent theories on both sides. (This is thus a slight
extension of
the AGM framework, as AGM work with formulas only on the right of $*.)$

{\scriptsize

\begin{tabular}{|c|c|c|}

\hline

\xEH
\xEH
$(*Equiv)$
\xEP
\xEH
\xEH
$ \xcm T \xcr S,$ $ \xcm T' \xcr S',$ $\xch$ $T*T' =S*S',$
\xEP

\hline

\xEH
\xEH
$(*CCL)$
\xEP
\xEH
\xEH
$T*T' $ is a consistent, deductively closed theory,
\xEP

\hline

\xEH
$( \xfA Succ)$
\xEH
$(*Succ)$
\xEP
\xEH
$A \xfA B \xcc B$
\xEH
$T' \xcc T*T',$
\xEP

\hline

\xEH
$( \xfA Con)$
\xEH
$(*Con)$
\xEP
\xEH
$A \xcs B \xEd \xCQ $ $ \xch $ $A \xfA B=A \xcs B$
\xEH
$Con(T \xcv T') $ $\xch$ $T*T' = \ol{T \xcv T' },$
\xEP

\hline

Intuitively,
\xEH
$( \xfA Loop)$
\xEH
$(*Loop)$
\xEP
Using symmetry
\xEH
\xEH
\xEP
$d(X_{0},X_{1}) \xck d(X_{1},X_{2}),$
\xEH
$(X_{1} \xfA (X_{0} \xcv X_{2})) \xcs X_{0} \xEd \xCQ,$
\xEH
$Con(T_{0},T_{1}*(T_{0} \xco T_{2})),$
\xEP
$d(X_{1},X_{2}) \xck d(X_{2},X_{3}),$
\xEH
$(X_{2} \xfA (X_{1} \xcv X_{3})) \xcs X_{1} \xEd \xCQ,$
\xEH
$Con(T_{1},T_{2}*(T_{1} \xco T_{3})),$
\xEP
$d(X_{2},X_{3}) \xck d(X_{3},X_{4})$
\xEH
$(X_{3} \xfA (X_{2} \xcv X_{4})) \xcs X_{2} \xEd \xCQ,$
\xEH
$Con(T_{2},T_{3}*(T_{2} \xco T_{4}))$
\xEP
\Xl
\xEH
\Xl
\xEH
\Xl
\xEP
$d(X_{k-1},X_{k}) \xck d(X_{0},X_{k})$
\xEH
$(X_{k} \xfA (X_{k-1} \xcv X_{0})) \xcs X_{k-1} \xEd \xCQ $
\xEH
$Con(T_{k-1},T_{k}*(T_{k-1} \xco T_{0}))$
\xEP
$\xch$
\xEH
$\xch$
\xEH
$\xch$
\xEP
$d(X_{0},X_{1}) \xck d(X_{0},X_{k}),$
\xEH
$(X_{0} \xfA (X_{k} \xcv X_{1})) \xcs X_{1} \xEd \xCQ$
\xEH
$Con(T_{1},T_{0}*(T_{k} \xco T_{1}))$
\xEP

i.e. transitivity, or absence of
\xEH
\xEH
\xEP

loops involving $<$
\xEH
\xEH
\xEP

\hline

\end{tabular}

}

 karl-search= End Definition TR-Dist
\vspace{7mm}

 *************************************

\vspace{7mm}

\subsubsection{Definition TR-Dist-Rotate}

 {\LARGE karl-search= Start Definition TR-Dist-Rotate }

\index{Definition TR-Dist-Rotate}

\ed

\bd

$\hspace{0.01em}$

(+++ Orig. No.:  Definition TR-Dist-Rotate +++)

{\xssc LABEL: {Definition TR-Dist-Rotate}}
\label{Definition TR-Dist-Rotate}

Let $U \xEd \xCQ,$ $ \xdy \xcc \xdp (U)$ satisfy $( \xcs ),$ $( \xcv ),$
$ \xCQ \xce \xdy.$

Let $A,B,X_{i} \xbe \xdy,$ $ \xfA: \xdy \xCK \xdy \xcp \xdp (U).$

Let $*$ be a revision function defined for
arbitrary consistent theories on both sides. (This is thus a slight
extension of
the AGM framework, as AGM work with formulas only on the right of $*.)$

\begin{turn}{90}

{\scriptsize

\begin{tabular}{|c|c|c|}

\hline

\xEH
\xEH
$(*Equiv)$
\xEP
\xEH
\xEH
$ \xcm T \xcr S,$ $ \xcm T' \xcr S',$ $\xch$ $T*T' =S*S',$
\xEP

\hline

\xEH
\xEH
$(*CCL)$
\xEP
\xEH
\xEH
$T*T' $ is a consistent, deductively closed theory,
\xEP

\hline

\xEH
$( \xfA Succ)$
\xEH
$(*Succ)$
\xEP
\xEH
$A \xfA B \xcc B$
\xEH
$T' \xcc T*T',$
\xEP

\hline

\xEH
$( \xfA Con)$
\xEH
$(*Con)$
\xEP
\xEH
$A \xcs B \xEd \xCQ $ $ \xch $ $A \xfA B=A \xcs B$
\xEH
$Con(T \xcv T') $ $\xch$ $T*T' = \ol{T \xcv T' },$
\xEP

\hline

Intuitively,
\xEH
$( \xfA Loop)$
\xEH
$(*Loop)$
\xEP
Using symmetry
\xEH
\xEH
\xEP
$d(X_{0},X_{1}) \xck d(X_{1},X_{2}),$
\xEH
$(X_{1} \xfA (X_{0} \xcv X_{2})) \xcs X_{0} \xEd \xCQ,$
\xEH
$Con(T_{0},T_{1}*(T_{0} \xco T_{2})),$
\xEP
$d(X_{1},X_{2}) \xck d(X_{2},X_{3}),$
\xEH
$(X_{2} \xfA (X_{1} \xcv X_{3})) \xcs X_{1} \xEd \xCQ,$
\xEH
$Con(T_{1},T_{2}*(T_{1} \xco T_{3})),$
\xEP
$d(X_{2},X_{3}) \xck d(X_{3},X_{4})$
\xEH
$(X_{3} \xfA (X_{2} \xcv X_{4})) \xcs X_{2} \xEd \xCQ,$
\xEH
$Con(T_{2},T_{3}*(T_{2} \xco T_{4}))$
\xEP
\Xl
\xEH
\Xl
\xEH
\Xl
\xEP
$d(X_{k-1},X_{k}) \xck d(X_{0},X_{k})$
\xEH
$(X_{k} \xfA (X_{k-1} \xcv X_{0})) \xcs X_{k-1} \xEd \xCQ $
\xEH
$Con(T_{k-1},T_{k}*(T_{k-1} \xco T_{0}))$
\xEP
$\xch$
\xEH
$\xch$
\xEH
$\xch$
\xEP
$d(X_{0},X_{1}) \xck d(X_{0},X_{k}),$
\xEH
$(X_{0} \xfA (X_{k} \xcv X_{1})) \xcs X_{1} \xEd \xCQ$
\xEH
$Con(T_{1},T_{0}*(T_{k} \xco T_{1}))$
\xEP

i.e. transitivity, or absence of
\xEH
\xEH
\xEP

loops involving $<$
\xEH
\xEH
\xEP

\hline

\end{tabular}

}

\end{turn}

 karl-search= End Definition TR-Dist-Rotate
\vspace{7mm}

 *************************************

\vspace{7mm}

\subsubsection{Proposition TR-Alg-Log}

 {\LARGE karl-search= Start Proposition TR-Alg-Log }

\index{Proposition TR-Alg-Log}

\ed

\bp

$\hspace{0.01em}$

(+++ Orig. No.:  Proposition TR-Alg-Log +++)

{\xssc LABEL: {Proposition TR-Alg-Log}}
\label{Proposition TR-Alg-Log}

The following connections between the logical and the algebraic side might
be the most interesting ones. We will consider in all cases also the
variant
with full theories.

Given $*$ which respects logical equivalence, let $M(T) \xfA M(T'
):=M(T*T' ),$
conversely, given $ \xfA,$ let $T*T':=Th(M(T) \xfA M(T' )).$ We then
have:

{\small

\begin{tabular}{|c|c|c|c|}

\hline

(1.1)
\xEH
$(K*7)$
\xEH
$\xch$
\xEH
$(X \xfA 7)$
\xEP

\cline{1-1}
\cline{3-3}

(1.2)
\xEH
\xEH
$\xci$ $(\xbm dp)$
\xEH
\xEP

\cline{1-1}
\cline{3-3}

(1.3)
\xEH
\xEH
$\xci$ B is the model set for some $\xbf$
\xEH
\xEP

\cline{1-1}
\cline{3-3}

(1.4)
\xEH
\xEH
$\xcI$ in general
\xEH
\xEP

\hline

(2.1)
\xEH
$(*Loop)$
\xEH
$\xch$
\xEH
$(\xfA Loop)$
\xEP

\cline{1-1}
\cline{3-3}

(2.2)
\xEH
\xEH
$\xci$ $(\xbm dp)$
\xEH
\xEP

\cline{1-1}
\cline{3-3}

(2.3)
\xEH
\xEH
$\xci$ all $X_i$ are the model sets for some $\xbf_i$
\xEH
\xEP

\cline{1-1}
\cline{3-3}

(2.4)
\xEH
\xEH
$\xcI$ in general
\xEH
\xEP

\hline

\end{tabular}

}

 karl-search= End Proposition TR-Alg-Log
\vspace{7mm}

 *************************************

\vspace{7mm}

\subsubsection{Proposition TR-Alg-Log Proof}

 {\LARGE karl-search= Start Proposition TR-Alg-Log Proof }

\index{Proposition TR-Alg-Log Proof}

\ep

\subparagraph{
Proof
}

$\hspace{0.01em}$

(+++*** Orig.:  Proof )

(1)

We consider the equivalence of $T*(T' \xcv T'' ) \xcc \ol{(T*T' ) \xcv T''
}$ and
$(M(T) \xfA M(T' )) \xcs M(T'' ) \xcc M(T) \xfA (M(T' ) \xcs M(T'' )).$

(1.1)

$(M(T) \xfA M(T' )) \xcs M(T'' )$ $=$ $M(T*T' ) \xcs M(T'' )$ $=$ $M((T*T'
) \xcv T'' )$ $ \xcc_{(K*7)}$
$M(T*(T' \xcv T'' ))$ $=$ $M(T) \xfA M(T' \xcv T'' )$ $=$ $M(T) \xfA (M(T'
) \xcs M(T'' )).$

(1.2)

$T*(T' \xcv T'' )$ $=$ $Th(M(T) \xfA M(T' \xcv T'' ))$ $=$ $Th(M(T) \xfA
(M(T' ) \xcs M(T'' )))$ $ \xcc_{(X \xfA 7)}$
$Th((M(T) \xfA M(T' )) \xcs M(T'' ))$ $=_{( \xbm dp)}$
$ \ol{Th(M(T) \xfA M(T' )) \xcv T'' }$ $=$ $ \ol{Th(M(T*T' ) \xcv T'' }$
$=$ $ \ol{(T*T' ) \xcv T'' }.$

(1.3)

Let $T'' $ be equivalent to $ \xbf ''.$ We can then replace the use of $(
\xbm dp)$
in the proof of (1.2) by Fact \ref{Fact Log-Form} (page \pageref{Fact Log-Form})
 (3).

(1.4)

By Example \ref{Example TR-Dp} (page \pageref{Example TR-Dp})  (2), $(K*7)$ may
fail, though $(X \xfA
7)$ holds.

(2.1) and (2.2):

$Con(T_{0},T_{1}*(T_{0} \xco T_{2}))$ $ \xcj $ $M(T_{0}) \xcs
M(T_{1}*(T_{0} \xco T_{2})) \xEd \xCQ.$

$M(T_{1}*(T_{0} \xco T_{2}))$ $=$ $M(Th(M(T_{1}) \xfA M(T_{0} \xco
T_{2})))$ $=$ $M(Th(M(T_{1}) \xfA (M(T_{0}) \xcv M(T_{2}))))$ $=_{( \xbm
dp)}$
$M(T_{1}) \xfA (M(T_{0}) \xcv (T_{2})),$ so
$Con(T_{0},T_{1}*(T_{0} \xco T_{2}))$ $ \xcj $ $M(T_{0}) \xcs (M(T_{1})
\xfA (M(T_{0}) \xcv (T_{2}))) \xEd \xCQ.$

Thus, all conditions translate one-to-one, and we use $( \xfA Loop)$ and
$(*Loop)$
to go back and forth.

(2.3):

Let $A:=M(Th(M(T_{1}) \xfA (M(T_{0}) \xcv M(T_{2})))),$ $A':=M(T_{1})
\xfA (M(T_{0}) \xcv (T_{2})),$ then we do
not need $A=A',$ it suffices to have $M(T_{0}) \xcs A \xEd \xCQ $ $ \xcj
$ $M(T_{0}) \xcs A' \xEd \xCQ.$ $A= \wt{A' },$ so
we can use Fact \ref{Fact Log-Form} (page \pageref{Fact Log-Form})  (4), if
$T_{0}$ is equivalent to
some $ \xbf_{0}.$

This has to hold for all $T_{i},$ so all $T_{i}$ have to be equivalent to
some $ \xbf_{i}.$

(2.4):

By Proposition \ref{Proposition TR-Alg-Repr} (page \pageref{Proposition
TR-Alg-Repr}) , all distance defined $ \xfA $
satisfy
$( \xfA Loop).$ By Example \ref{Example TR-Dp} (page \pageref{Example TR-Dp}) 
(1), $(*Loop)$ may fail.

$ \xcz $
\\[3ex]

 karl-search= End Proposition TR-Alg-Log Proof
\vspace{7mm}

 *************************************

\vspace{7mm}

\subsubsection{Proposition TR-Representation-With-Ref}

 {\LARGE karl-search= Start Proposition TR-Representation-With-Ref }

\index{Proposition TR-Representation-With-Ref}

The following table summarizes representation of theory revision
functions by structures with a distance.

By ``pseudo-distance'' we mean here a pseudo-distance which respects
identity, and is symmetrical.

$( \xfA \xCQ )$ means that if $X,Y \xEd \xCQ,$ then $X \xfA_{d}Y \xEd
\xCQ.$
{\xssc LABEL: {Proposition TR-Representation-With-Ref}}
\label{Proposition TR-Representation-With-Ref}

{\scriptsize

\begin{tabular}{|c|c|c|c|c|}

\hline

$\xfA-$ function
\xEH
\xEH
Distance Structure
\xEH
\xEH
$*-$ function
\xEP

\hline

$(\xfA Succ)+(\xfA Con)+$
\xEH
$\xcj$, $(\xcv)+(\xcs)$
\xEH
pseudo-distance
\xEH
$\xcj$ $(\xbm dp)+(\xfA\xCQ)$
\xEH
$(*Equiv)+(*CCL)+(*Succ)+$
\xEP

$(\xfA Loop)$
\xEH
Proposition \ref{Proposition TR-Alg-Repr}
\xEH
\xEH
Proposition \ref{Proposition TR-Log-Repr}
\xEH
$(*Con)+(*Loop)$
\xEP

\xEH
page \pageref{Proposition TR-Alg-Repr}
\xEH
\xEH
page \pageref{Proposition TR-Log-Repr}
\xEH
\xEP

\cline{1-2}
\cline{4-4}

any finite
\xEH
$\xcJ$
\xEH
\xEH
$\xcH$ without $(\xbm dp)$
\xEH
\xEP

characterization
\xEH
Proposition \ref{Proposition Hamster}
\xEH
\xEH
Example \ref{Example TR-Dp}
\xEH
\xEP

\xEH
page \pageref{Proposition Hamster}
\xEH
\xEH
page \pageref{Example TR-Dp}
\xEH
\xEP

\hline

\end{tabular}

}

 karl-search= End Proposition TR-Representation-With-Ref
\vspace{7mm}

 *************************************

\vspace{7mm}

\subsubsection{Example TR-Dp}

 {\LARGE karl-search= Start Example TR-Dp }

\index{Example TR-Dp}

The following Example \ref{Example TR-Dp} (page \pageref{Example TR-Dp})  shows
that, in general, a
revision operation
defined on models via a pseudo-distance by $T*T':=Th(M(T) \xfA_{d}M(T'
))$ might not
satisfy $(*Loop)$ or $(K*7),$ unless we require $ \xfA_{d}$ to preserve
definability.

\be

$\hspace{0.01em}$

(+++ Orig. No.:  Example TR-Dp +++)

{\xssc LABEL: {Example TR-Dp}}
\label{Example TR-Dp}

Consider an infinite propositional language $ \xdl.$

Let $X$ be an infinite set of models, $m,$ $m_{1},$ $m_{2}$ be models for
$ \xdl.$
Arrange the models of $ \xdl $ in the real plane s.t. all $x \xbe X$ have
the same
distance $<2$ (in the real plane) from $m,$ $m_{2}$ has distance 2 from
$m,$ and $m_{1}$ has
distance 3 from $m.$

Let $T,$ $T_{1},$ $T_{2}$ be complete (consistent) theories, $T' $ a
theory with infinitely
many models, $M(T)=\{m\},$ $M(T_{1})=\{m_{1}\},$ $M(T_{2})=\{m_{2}\}.$ The
two variants diverge now
slightly:

(1) $M(T' )=X \xcv \{m_{1}\}.$ $T,T',T_{2}$ will be pairwise
inconsistent.

(2) $M(T' )=X \xcv \{m_{1},m_{2}\},$ $M(T'' )=\{m_{1},m_{2}\}.$

Assume in both cases $Th(X)=T',$ so $X$ will not be definable by a
theory.

Now for the results:

Then $M(T) \xfA M(T' )=X,$ but $T*T' =Th(X)=T'.$

(1) We easily verify
$Con(T,T_{2}*(T \xco T)),$ $Con(T_{2},T*(T_{2} \xco T_{1})),$
$Con(T,T_{1}*(T \xco T)),$ $Con(T_{1},T*(T_{1} \xco T' )),$
$Con(T,T' *(T \xco T)),$ and conclude by Loop (i.e. $(*Loop))$
$Con(T_{2},T*(T' \xco T_{2})),$ which
is wrong.

(2) So $T*T' $ is consistent with $T'',$ and
$ \ol{(T*T' ) \xcv T'' }=T''.$ But $T' \xcv T'' =T'',$ and $T*(T' \xcv
T'' )=T_{2} \xEd T'',$ contradicting $(K*7).$

$ \xcz $
\\[3ex]

 karl-search= End Example TR-Dp
\vspace{7mm}

 *************************************

\vspace{7mm}

\subsubsection{Proposition TR-Alg-Repr}

 {\LARGE karl-search= Start Proposition TR-Alg-Repr }

\index{Proposition TR-Alg-Repr}

\ee

\bp

$\hspace{0.01em}$

(+++ Orig. No.:  Proposition TR-Alg-Repr +++)

{\xssc LABEL: {Proposition TR-Alg-Repr}}
\label{Proposition TR-Alg-Repr}

Let $U \xEd \xCQ,$ $ \xdy \xcc \xdp (U)$ be closed under finite $ \xcs $
and finite $ \xcv,$ $ \xCQ \xce \xdy.$

(a) $ \xfA $ is representable by a symmetric pseudo-distance $d:U \xCK U
\xcp Z$ iff $ \xfA $
satisfies $( \xfA Succ)$ and $( \xfA Loop)$ in Definition \ref{Definition
TR-Dist} (page \pageref{Definition TR-Dist}) .

(b) $ \xfA $ is representable by an identity respecting symmetric
pseudo-distance
$d:U \xCK U \xcp Z$ iff $ \xfA $ satisfies $( \xfA Succ),$ $( \xfA Con),$
and $( \xfA Loop)$
in Definition \ref{Definition TR-Dist} (page \pageref{Definition TR-Dist}) .

See  \cite{LMS01} or  \cite{Sch04}.

 karl-search= End Proposition TR-Alg-Repr
\vspace{7mm}

 *************************************

\vspace{7mm}

\subsubsection{Proposition TR-Log-Repr}

 {\LARGE karl-search= Start Proposition TR-Log-Repr }

\index{Proposition TR-Log-Repr}

\ep

\bp

$\hspace{0.01em}$

(+++ Orig. No.:  Proposition TR-Log-Repr +++)

{\xssc LABEL: {Proposition TR-Log-Repr}}
\label{Proposition TR-Log-Repr}

Let $ \xdl $ be a propositional language.

(a) A revision operation $*$ is representable by a symmetric consistency
and
definability preserving pseudo-distance iff $*$ satisfies $(*Equiv),$
$(*CCL),$ $(*Succ),$ $(*Loop).$

(b) A revision operation $*$ is representable by a symmetric consistency
and
definability preserving, identity respecting pseudo-distance iff $*$
satisfies $(*Equiv),$ $(*CCL),$ $(*Succ),$ $(*Con),$ $(*Loop).$

See  \cite{LMS01} or  \cite{Sch04}.

 karl-search= End Proposition TR-Log-Repr
\vspace{7mm}

 *************************************

\vspace{7mm}

\subsubsection{Example WeakTR}

 {\LARGE karl-search= Start Example WeakTR }

\index{Example WeakTR}

\ep

\be

$\hspace{0.01em}$

(+++ Orig. No.:  Example WeakTR +++)

{\xssc LABEL: {Example WeakTR}}
\label{Example WeakTR-a}

This example shows the expressive weakness of revision based on
distance: not all distance relations can be reconstructed from
the revision operator. Thus, a revision operator does not allow
to ``observe'' all distances relations, so transitivity of $ \xck $ cannot
necessarily be captured in a short condition, requiring
arbitrarily long conditions, see Proposition \ref{Proposition Hamster} (page
\pageref{Proposition Hamster}) .

Note that even when the pseudo-distance is a real distance, the
resulting revision operator $ \xfA_{d}$ does not always permit to
reconstruct the
relations of the distances: revision is a coarse instrument to investigate
distances.

Distances with common start (or end, by symmetry) can always be
compared by looking at the result of revision:

$a \xfA_{d}\{b,b' \}=b$ iff $d(a,b)<d(a,b' ),$

$a \xfA_{d}\{b,b' \}=b' $ iff $d(a,b)>d(a,b' ),$

$a \xfA_{d}\{b,b' \}=\{b,b' \}$ iff $d(a,b)=d(a,b' ).$

This is not the case with arbitrary distances $d(x,y)$ and $d(a,b),$
as this example will show.

\ee

We work in the real plane, with the standard distance, the angles have 120
degrees. $a' $ is closer to $y$ than $x$ is to $y,$ a is closer to $b$
than $x$ is to $y,$
but $a' $ is farther away from $b' $ than $x$ is from $y.$ Similarly for
$b,b'.$
But we cannot distinguish the situation $\{a,b,x,y\}$ and the
situation $\{a',b',x,y\}$ through $ \xfA_{d}.$ (See Diagram \ref{Diagram WeakTR}
(page \pageref{Diagram WeakTR}) ):

Seen from a, the distances are in that order: $y,b,x.$

Seen from $a',$ the distances are in that order: $y,b',x.$

Seen from $b,$ the distances are in that order: $y,a,x.$

Seen from $b',$ the distances are in that order: $y,a',x.$

Seen from $y,$ the distances are in that order: $a/b,x.$

Seen from $y,$ the distances are in that order: $a' /b',x.$

Seen from $x,$ the distances are in that order: $y,a/b.$

Seen from $x,$ the distances are in that order: $y,a' /b'.$

Thus, any $c \xfA_{d}C$ will be the same in both situations (with a
interchanged with
$a',$ $b$ with $b' ).$ The same holds for any $X \xfA_{d}C$ where $X$ has
two elements.

Thus, any $C \xfA_{d}D$ will be the same in both situations, when we
interchange a with
$a',$ and $b$ with $b'.$ So we cannot determine by $ \xfA_{d}$ whether
$d(x,y)>d(a,b)$ or not.
$ \xcz $
\\[3ex]

\vspace{10mm}

\begin{diagram}

{\xssc LABEL: {Diagram WeakTR}}
\label{Diagram WeakTR}
\index{Diagram WeakTR}

\centering
\setlength{\unitlength}{1mm}
{\renewcommand{\dashlinestretch}{30}
\begin{picture}(130,90)(0,0)

\put(5,50){\line(1,0){30}}
\put(35,50){\line(1,2){6}}
\put(35,50){\line(1,-2){6}}

\put(5,50){\circle*{1.5}}
\put(35,50){\circle*{1.5}}

\put(41,62){\circle*{1.5}}
\put(41,38){\circle*{1.5}}

\put(5,47){$x$}
\put(32,47){$y$}
\put(43,61){$a$}
\put(43,37){$b$}

\put(65,50){\line(1,0){35}}
\put(100,50){\line(1,2){12}}
\put(100,50){\line(1,-2){12}}

\put(65,50){\circle*{1.5}}
\put(100,50){\circle*{1.5}}

\put(112,74){\circle*{1.5}}
\put(112,26){\circle*{1.5}}

\put(65,47){$x$}
\put(97,47){$y$}
\put(114,73){$a'$}
\put(114,25){$b'$}

\put(29,10){Indiscernible by revision}

\end{picture}

}

\end{diagram}

\vspace{4mm}

 karl-search= End Example WeakTR
\vspace{7mm}

 *************************************

\vspace{7mm}

\subsubsection{Proposition Hamster}

 {\LARGE karl-search= Start Proposition Hamster }

\index{Proposition Hamster}

\bp

$\hspace{0.01em}$

(+++ Orig. No.:  Proposition Hamster +++)

{\xssc LABEL: {Proposition Hamster}}
\label{Proposition Hamster}

There is no finite characterization of distance based $ \xfA -$operators.

(Attention: this is, of course, false when we fix the left hand side:
the AGM axioms give a finite characterization. So this also shows the
strength of being able to change the left hand side.)

See  \cite{Sch04}.

 karl-search= End Proposition Hamster
\vspace{7mm}

 *************************************

\vspace{7mm}

 karl-search= End ToolBase1-TR-DistRepr
\vspace{7mm}

 *************************************

\vspace{7mm}

 karl-search= End ToolBase1-TR
\vspace{7mm}

 *************************************

\vspace{7mm}

\newpage

\section{
Size
}

\subsection{
}

\subsubsection{ToolBase1-Size}

 {\LARGE karl-search= Start ToolBase1-Size }

{\xssc LABEL: {Section Toolbase1-Size}}
\label{Section Toolbase1-Size}
\index{Section Toolbase1-Size}
\subsubsection{Definition Filter}

 {\LARGE karl-search= Start Definition Filter }

\index{Definition Filter}

\ep

\bd

$\hspace{0.01em}$

(+++ Orig. No.:  Definition Filter +++)

{\xssc LABEL: {Definition Filter}}
\label{Definition Filter}

A filter is an abstract notion of size,
elements of a filter $ \xdf (X)$ on $X$ are called big subsets of $X,$
their complements
are called small, and the rest have medium size. The dual applies to
ideals
$ \xdi (X),$ this is justified by the trivial fact that $\{X-A:A \xbe \xdf
(X)\}$ is an ideal
iff $ \xdf (X)$ is a filter.

In both definitions, the first two conditions (i.e. $ \xCf (FAll),$ $(I
\xCQ ),$ and
$(F \xfB ),$ $(I \xfb ))$ should hold if the notions shall
have anything to do with usual intuition, and there are reasons to
consider
only the weaker, less idealistic, version of the third.

At the same time, we introduce - in rough parallel - coherence conditions
which
describe what might happen when we change the reference or base set $X.$
$(R \xfB )$
is very natural, $(R \xfb )$ is more daring, and $(R \xfb \xfb )$ even
more so.
$(R \xcv disj)$ is a cautious combination of $(R \xfB )$ and $(R \xcv ),$
as we avoid using
the same big set several times in comparison, so $(R \xcv )$ is used more
cautiously
here. See Remark \ref{Remark Ref-Class} (page \pageref{Remark Ref-Class})  for
more details.

Finally, we give a generalized first order quantifier corresponding to a
(weak) filter. The precise connection is formulated in
Definition \ref{Definition Nabla} (page \pageref{Definition Nabla}) , Definition
\ref{Definition N-Model} (page \pageref{Definition N-Model}) ,
Definition \ref{Definition NablaAxioms} (page \pageref{Definition NablaAxioms})
, and
Proposition \ref{Proposition NablaRepr} (page \pageref{Proposition NablaRepr}) ,
respectively their
relativized versions.

Fix now a base set $X \xEd \xCQ.$

A (weak) filter on or over $X$ is a set $ \xdf (X) \xcc \xdp (X),$ s.t.
(FAll), $(F \xfB ),$ $(F \xcs )$
$((FAll),$ $(F \xfB ),$ $(F \xcs ' )$ respectively) hold.

A filter is called a principal filter iff there is $X' \xcc X$ s.t. $ \xdf
=\{A:$ $X' \xcc A \xcc X\}.$

A filter is called an ultrafilter iff for all $X' \xcc X$ $X' \xbe \xdf
(X)$ or $X-X' \xbe \xdf (X).$

A (weak) ideal on or over $X$ is a set $ \xdi (X) \xcc \xdp (X),$ s.t. $(I
\xCQ ),$ $(I \xfb ),$ $(I \xcv )$
$((I \xCQ ),$ $(I \xfb ),$ $(I \xcv ' )$ respectively) hold.

Finally, we set $ \xdm (X):=\{A \xcc X:A \xce \xdi (X),$ $A \xce \xdf
(X)\},$ the ``medium size'' sets, and
$ \xdm^{+}(X):= \xdm (X) \xcv \xdf (X),$
$ \xdm^{+}(X)$ is the set of subsets of $X,$ which are not small, i.e.
have medium or
large size.

For $(R \xfb )$ and $(R \xfb \xfb )$ $ \xCf (-)$ is assumed in the
following
table.

{\footnotesize

\begin{tabular}{|c|c|c|c|}

\hline

\multicolumn{4}{|c|}{Optimum} \xEP

\hline

$(FAll)$
\xEH
$(I \xCQ )$
\xEH
\xEH
\xEP

$X \xbe \xdf (X)$
\xEH
$ \xCQ \xbe \xdi (X)$
\xEH
\xEH
$\xcA x\xbf(x) \xcp \xeA x\xbf(x)$
\xEP

\hline

\multicolumn{4}{|c|}{Improvement} \xEP

\hline

$(F \xfB )$
\xEH
$(I \xfb )$
\xEH
$(R \xfB )$
\xEH
\xEP

$A \xcc B \xcc X,$
\xEH
$A \xcc B \xcc X,$
\xEH
$X \xcc Y$ $ \xch $ $\xdi (X) \xcc \xdi (Y)$
\xEH
$\xeA x\xbf(x) \xcu $
\xEP

$A \xbe \xdf (X)$ $ \xch $
\xEH
$B \xbe \xdi (X)$ $ \xch $
\xEH
\xEH
$ \xcA x(\xbf(x)\xcp\xbq(x)) \xcp$
\xEP

$B \xbe \xdf (X)$
\xEH
$A \xbe \xdi (X)$
\xEH
\xEH
$\xeA x\xbq(x)$
\xEP

\hline

\multicolumn{4}{|c|}{Adding small sets} \xEP

\hline

$(F \xcs )$
\xEH
$(I \xcv )$
\xEH
$(R \xfb)$
\xEH
\xEP

$A,B \xbe \xdf (X)$ $ \xch $
\xEH
$A,B \xbe \xdi (X)$ $ \xch $
\xEH
$A,B \xbe \xdi (X)$ $ \xch $
\xEH
$\xeA x\xbf(x) \xcu \xeA x\xbq(x) \xcp $
\xEP

$A \xcs B \xbe \xdf (X)$
\xEH
$A \xcv B \xbe \xdi (X)$
\xEH
$A-B \xbe \xdi (X-$B)
\xEH
$\xeA x(\xbf(x) \xcu \xbq(x)) $
\xEP

\xEH
\xEH
or:
\xEH
\xEP

\xEH
\xEH
$A \xbe \xdf (X),$ $B \xbe \xdi (X)$ $ \xch $
\xEH
\xEP

\xEH
\xEH
$A-B \xbe \xdf (X-$B)
\xEH
\xEP

\hline

\multicolumn{4}{|c|}{Cautious addition} \xEP

\hline

$(F \xcs ' )$
\xEH
$(I \xcv ' )$
\xEH
$(R \xcv disj)$
\xEH
\xEP

$A,B \xbe \xdf (X)$ $ \xch $
\xEH
$A,B \xbe \xdi (X)$ $ \xch $
\xEH
$A \xbe \xdi (X),$ $B \xbe \xdi (Y),$ $X \xcs Y= \xCQ $ $ \xch $
\xEH
$\xeA x\xbf(x) \xcp \xCN\xeA x\xCN\xbf(x)$
\xEP

$A \xcs B \xEd \xCQ.$
\xEH
$A \xcv B \xEd X.$
\xEH
$A \xcv B \xbe \xdi (X \xcv Y)$
\xEH
and  $\xeA x\xbf(x) \xcp \xcE x\xbf(x)$
\xEP

\hline

\multicolumn{4}{|c|}{Bold addition} \xEP

\hline

Ultrafilter
\xEH
(Dual of) Ultrafilter
\xEH
$(R \xfb \xfb)$
\xEH
\xEP

\xEH
\xEH
$A \xbe \xdi (X),$ $B \xce \xdf (X)$ $ \xch $
\xEH
$\xCN \xeA x\xbf(x) \xcp \xeA x\xCN\xbf(x)$
\xEP

\xEH
\xEH
$A-B \xbe \xdi (X-B)$
\xEH
\xEP

\xEH
\xEH
or:
\xEH
\xEP

\xEH
\xEH
$A \xbe \xdf (X),$ $B \xce \xdf (X)$ $ \xch $
\xEH
\xEP

\xEH
\xEH
$A-B \xbe \xdf (X-$B)
\xEH
\xEP

\xEH
\xEH
or:
\xEH
\xEP

\xEH
\xEH
$A \xbe \xdm^+ (X),$ $X \xbe \xdm^+ (Y)$ $ \xch $
\xEH
\xEP

\xEH
\xEH
$A \xbe \xdm^+ (Y)$ - Transitivity of $\xdm^+$
\xEH
\xEP

\hline

\end{tabular}

}

\ed

These notions are related to nonmonotonic logics as follows:

We can say that, normally, $ \xbf $ implies $ \xbq $ iff in a big subset
of all $ \xbf -$cases, $ \xbq $
holds. In preferential terms, $ \xbf $ implies $ \xbq $ iff $ \xbq $ holds
in all minimal
$ \xbf -$models. If $ \xbm $ is the model choice function of a
preferential structure, i.e.
$ \xbm ( \xbf )$ is the set of minimal $ \xbf -$models, then $ \xbm ( \xbf
)$ will be a (the smallest)
big subset of the set of $ \xbf -$models, and the filter over the $ \xbf
-$models is the
pricipal filter generated by $ \xbm ( \xbf ).$

Due to the finite intersection property, filters and ideals work well with
logics:
If $ \xbf $ holds normally, as it holds in a big subset, and so does $
\xbf ',$ then
$ \xbf \xcu \xbf ' $ will normally hold, too, as the intersection of two
big subsets is
big again. This is a nice property, but not justified in all situations,
consider e.g. simple counting of a finite subset. (The question has a
name,
``lottery paradox'': normally no single participant wins, but someone wins
in the
end.) This motivates the weak versions.

Normality defined by (weak or not) filters is a local concept: the filter
defined on $X$ and the one defined on $X' $ might be totally independent.

Seen more abstractly, set properties like e.g. $(R \xfB )$
allow the transfer of big (or small) subsets from one to
another base set (and the conclusions drawn on this basis), and we call
them
``coherence properties''. They are very important, not only for working with
a logic which respects them, but also for soundness and completeness
questions,
often they are at the core of such problems.

 karl-search= End Definition Filter
\vspace{7mm}

 *************************************

\vspace{7mm}

\subsubsection{Definition Filter2}

 {\LARGE karl-search= Start Definition Filter2 }

\index{Definition Filter2}

\bd

$\hspace{0.01em}$

(+++ Orig. No.:  Definition Filter2 +++)

{\xssc LABEL: {Definition Filter2}}
\label{Definition Filter2}

A filter is an abstract notion of size,
elements of a filter $ \xdf (X)$ on $X$ are called big subsets of $X,$
their complements
are called small, and the rest have medium size. The dual applies to
ideals
$ \xdi (X),$ this is justified by the trivial fact that $\{X-A:A \xbe \xdf
(X)\}$ is an ideal
iff $ \xdf (X)$ is a filter.

In both definitions, the first two conditions (i.e. $ \xCf (FAll),$ $(I
\xCQ ),$ and
$(F \xfB ),$ $(I \xfb ))$ should hold if the notions shall
have anything to do with usual intuition, and there are reasons to
consider
only the weaker, less idealistic, version of the third.

At the same time, we introduce - in rough parallel - coherence conditions
which
describe what might happen when we change the reference or base set $X.$
$(R \xfB )$
is very natural, $(R \xfb )$ is more daring, and $(R \xfb \xfb )$ even
more so.
$(R \xcv disj)$ is a cautious combination of $(R \xfB )$ and $(R \xcv ),$
as we avoid using
the same big set several times in comparison, so $(R \xcv )$ is used more
cautiously
here. See Remark \ref{Remark Ref-Class} (page \pageref{Remark Ref-Class})  for
more details.

Finally, we give a generalized first order quantifier corresponding to a
(weak) filter. The precise connection is formulated in
Definition \ref{Definition Nabla} (page \pageref{Definition Nabla}) , Definition
\ref{Definition N-Model} (page \pageref{Definition N-Model}) ,
Definition \ref{Definition NablaAxioms} (page \pageref{Definition NablaAxioms})
, and
Proposition \ref{Proposition NablaRepr} (page \pageref{Proposition NablaRepr}) ,
respectively their
relativized versions.

Fix now a base set $X \xEd \xCQ.$

A (weak) filter on or over $X$ is a set $ \xdf (X) \xcc \xdp (X),$ s.t.
(FAll), $(F \xfB ),$ $(F \xcs )$
$((FAll),$ $(F \xfB ),$ $(F \xcs ' )$ respectively) hold.

A filter is called a principal filter iff there is $X' \xcc X$ s.t. $ \xdf
=\{A:$ $X' \xcc A \xcc X\}.$

A filter is called an ultrafilter iff for all $X' \xcc X$ $X' \xbe \xdf
(X)$ or $X-X' \xbe \xdf (X).$

A (weak) ideal on or over $X$ is a set $ \xdi (X) \xcc \xdp (X),$ s.t. $(I
\xCQ ),$ $(I \xfb ),$ $(I \xcv )$
$((I \xCQ ),$ $(I \xfb ),$ $(I \xcv ' )$ respectively) hold.

Finally, we set $ \xdm (X):=\{A \xcc X:A \xce \xdi (X),$ $A \xce \xdf
(X)\},$ the ``medium size'' sets, and
$ \xdm^{+}(X):= \xdm (X) \xcv \xdf (X),$
$ \xdm^{+}(X)$ is the set of subsets of $X,$ which are not small, i.e.
have medium or
large size.

For $(R \xfb )$ and $(R \xfb \xfb )$ $ \xCf (-)$ is assumed in the
following
table.

{\footnotesize

\begin{tabular}{|c|c|c|c|}

\hline

Filter
\xEH
Ideal
\xEH
Coherence
\xEH
$\xeA$
\xEP

\hline

\multicolumn{4}{|c|}{Not trivial} \xEP

\hline

$(F \xCQ)$
\xEH
$(I All)$
\xEH
\xEH
$(\xeA \xcE)$
\xEP

$\xCQ \xce \xdf (X)$
\xEH
$ X \xce \xdi (X)$
\xEH
\xEH
$\xeA x \xbf(x) \xcp \xcE x \xbf(x)$
\xEP

\hline

\multicolumn{4}{|c|}{Optimal proportion} \xEP

\hline

$(F All)$
\xEH
$(I \xCQ )$
\xEH
\xEH
$(\xeA All)$
\xEP

$X \xbe \xdf (X)$
\xEH
$ \xCQ \xbe \xdi (X)$
\xEH
\xEH
$\xcA x\xbf(x) \xcp \xeA x\xbf(x)$
\xEP

\hline

\multicolumn{4}{|c|}{Improving proportions} \xEP

\hline

$(F \xfB )$
\xEH
$(I \xfb )$
\xEH
$(R \xfB )$
\xEH
$(\xeA \xfB)$
\xEP

$A \xcc B \xcc X,$
\xEH
$A \xcc B \xcc X,$
\xEH
$X \xcb Y$ $ \xch $ $\xdi (X) \xcc \xdi (Y)$
\xEH
$\xeA x\xbf(x) \xcu $
\xEP

$A \xbe \xdf (X)$ $ \xch $
\xEH
$B \xbe \xdi (X)$ $ \xch $
\xEH
\xEH
$ \xcA x(\xbf(x)\xcp\xbq(x)) \xcp$
\xEP

$B \xbe \xdf (X)$
\xEH
$A \xbe \xdi (X)$
\xEH
\xEH
$\xeA x\xbq(x)$
\xEP

\hline

\multicolumn{4}{|c|}{Improving or keeping proportions} \xEP

\hline

\xEH
\xEH
$(R \xcv disj +)$
\xEH
\xEP

\xEH
\xEH
$A \xbe \xdi (X),$ $B \xbe \xdi (Y),$
\xEH
\xEP

\xEH
\xEH
$(X-A) \xcs (Y-B)= \xCQ $ $ \xch $
\xEH
\xEP

\xEH
\xEH
$A \xcv B \xbe \xdi (X \xcv Y)$
\xEH
\xEP

\hline

\multicolumn{4}{|c|}{Keeping proportions} \xEP

\hline

\xEH
\xEH
$(R \xcv disj)$
\xEH
\xEP

\xEH
\xEH
$A \xbe \xdi (X),$ $B \xbe \xdi (Y),$
\xEH
\xEP

\xEH
\xEH
$X \xcs Y= \xCQ $ $ \xch $
\xEH
\xEP

\xEH
\xEH
$A \xcv B \xbe \xdi (X \xcv Y)$
\xEH
\xEP

\hline

\multicolumn{4}{|c|}{2-Robustness of proportions (small+small $\xEd $ all)}
\xEP

\hline

$(F \xcs 2)$
\xEH
$(I \xcv 2)$
\xEH
$(R \xfb 2)$
\xEH
$(\xeA \xfb 2)$
\xEP

$A,B \xbe \xdf (X)$ $ \xch $
\xEH
$A,B \xbe \xdi (X)$ $ \xch $
\xEH
$A \xbe \xdf (X),$ $X \xbe \xdf (Y)$ $ \xch $
\xEH
$\xeA x\xbf(x) \xcp \xCN\xeA x\xCN\xbf(x)$
\xEP

$A \xcs B \xEd \xCQ.$
\xEH
$A \xcv B \xEd X.$
\xEH
$A \xbe \xdm^+ (Y)$
\xEH
New rules:
\xEP

\xEH
\xEH
\xEH
(1) $\xba\xcn\xbb \xch \xba\xcN\xCN\xbb$
\xEP

\xEH
\xEH
\xEH
(2) $\xba\xcn\xbb, \xba\xcn\xbb' \xch \xba\xcu\xbb\xcN\xCN\xbb'$
\xEP

\hline

\multicolumn{4}{|c|}{$n-$Robustness of proportions ($n*$small $\xEd $ all)}
\xEP

\hline

$(F \xcs n)$
\xEH
$(I \xcv n)$
\xEH
$(R \xfb n)$
\xEH
$(\xeA \xfb n)$
\xEP

$X_1, \ldots, X_n \xbe \xdf (X)$ $ \xch $
\xEH
$X_1, \ldots, X_n \xbe \xdi (X)$ $ \xch $
\xEH
$X_1 \xbe \xdf (X_2), \ldots,$ $X_{n-1} \xbe \xdf (X_n)$ $ \xch $
\xEH
$\xeA x\xbf_1(x) \xcu \ldots \xcu \xeA x\xbf_{n-1}(x)\xcp $
\xEP

$X_1 \xcs \ldots \xcs X_n \xEd \xCQ.$
\xEH
$X_1 \xcv \ldots \xcv X_n \xEd X.$
\xEH
$X_1 \xbe \xdm^+ (X_n)$
\xEH
$\xCN\xeA x(\xCN\xbf_1 \xco \ldots \xco \xCN \xbf_{n-1})(x)$
\xEP

\xEH
\xEH
\xEH
New rules:
\xEP

\xEH
\xEH
\xEH
(1) $\xba\xcn\xbb_1, \ldots, \xba\xcn\xbb_{n-1} \xch $
\xEP

\xEH
\xEH
\xEH
$\xba\xcN(\xCN\xbb_1 \xco \ldots \xco \xCN\xbb_{n-1})$
\xEP

\xEH
\xEH
\xEH
(2) $\xba\xcn\xbb_1, \ldots, \xba\xcn\xbb_n \xch $
\xEP

\xEH
\xEH
\xEH
$\xba \xcu \xbb_1 \xcu \ldots \xcu \xbb_{n-1} \xcN\xCN\xbb_n$ \xEP

\hline

\multicolumn{4}{|c|}{$\xbo-$Robustness of proportions (small+small=small)} \xEP

\hline

$(F \xcs \xbo)$
\xEH
$(I \xcv \xbo)$
\xEH
$(R \xfb \xbo)$
\xEH
$(\xeA \xbo)$
\xEP

$A,B \xbe \xdf (X)$ $ \xch $
\xEH
$A,B \xbe \xdi (X)$ $ \xch $
\xEH
(1) $A,B \xbe \xdi (X)$ $ \xch $
\xEH
$\xeA x\xbf(x) \xcu \xeA x\xbq(x) \xcp $
\xEP

$A \xcs B \xbe \xdf (X)$
\xEH
$A \xcv B \xbe \xdi (X)$
\xEH
$A-B \xbe \xdi (X-$B)
\xEH
$\xeA x(\xbf(x) \xcu \xbq(x)) $
\xEP

\xEH
\xEH
or:
\xEH
\xEP

\xEH
\xEH
(2) $A \xbe \xdf (X),$ $B \xbe \xdi (X)$ $ \xch $
\xEH
\xEP

\xEH
\xEH
$A-B \xbe \xdf (X-$B)
\xEH
\xEP

\xEH
\xEH
or:
\xEH
\xEP

\xEH
\xEH
(3) $A \xbe \xdf (X),$ $X \xbe \xdm^+ (Y)$ $ \xch $
\xEH
\xEP

\xEH
\xEH
$A \xbe \xdm^+ (Y)$
\xEH
\xEP

\xEH
\xEH
or:
\xEH
\xEP

\xEH
\xEH
(4) $A \xbe \xdm^+ (X),$ $X \xbe \xdf (Y)$ $ \xch $
\xEH
\xEP

\xEH
\xEH
$A \xbe \xdm^+ (Y)$
\xEH
\xEP

\hline

\multicolumn{4}{|c|}{Strong robustness of proportions} \xEP

\hline

\xEH
\xEH
$(R \xfb \xfb)$
\xEH
\xEP

\xEH
\xEH
(1) $A \xbe \xdi (X),$ $B \xce \xdf (X)$ $ \xch $
\xEH
\xEP

\xEH
\xEH
$A-B \xbe \xdi (X-B)$
\xEH
\xEP

\xEH
\xEH
or:
\xEH
\xEP

\xEH
\xEH
(2) $A \xbe \xdf (X),$ $B \xce \xdf (X)$ $ \xch $
\xEH
\xEP

\xEH
\xEH
$A-B \xbe \xdf (X-$B)
\xEH
\xEP

\xEH
\xEH
or:
\xEH
\xEP

\xEH
\xEH
(3) $A \xbe \xdm^+ (X),$ $X \xbe \xdm^+ (Y)$ $ \xch $
\xEH
\xEP

\xEH
\xEH
$A \xbe \xdm^+ (Y)$
\xEH
\xEP

\xEH
\xEH
(Transitivity of $\xdm^+$)
\xEH
\xEP

\hline

\multicolumn{4}{|c|}{Ultrafilter} \xEP

\hline

Ultrafilter
\xEH
(Dual of) Ultrafilter
\xEH
\xEH
\xEP

\xEH
\xEH
\xEH
$\xCN \xeA x\xbf(x) \xcp \xeA x\xCN\xbf(x)$
\xEP

\hline

\end{tabular}

}

\subsubsection{
New notations:
}

\ed

$(F \xcs ' ) \xcp (F \xcs 2),$ $(I \xcs ' ) \xcp (I \xcs 2),$ $(R \xfb -)
\xcp (R \xfb 2),$

$(F \xcs ) \xcp (F \xcs \xbo ),$ $(I \xcs ) \xcp (I \xcs \xbo ),$ $(R \xfb
) \xcp (R \xfb \xbo ),$
\subsubsection{
Algebra of size:
}

(1) internal:

(1.1) $n*small=medium$ $ \xch $ $2n*small=all$

(1.2) $n*small=big$ $ \xch $ $(n+1)*small=all$

(1.3) $medium+small \xEd big?$

(1.4) $n*small=all$

So, for (1) the number of ``small'' which can be all seems to be the
decisive measure.

(2) downward:

(2.1) small in $X$ $ \xch $ small in $(X-n*small)$

(2.2) small in $X$ $ \xch $ small in (X-medium)

(3) upward:

(3.1) big in $X,$ $X$ big in $Y$ $ \xch $ big in $Y,$ oder wenigstens $
\xdm^{+}$

(3.2) big in $X,$ $Y=X+n*small$ (in $Y)$ $ \xch $ $big/ \xdm^{+}$ in $Y$

(3.3) big in $X,$ $X$ big in $Y,$ $Y$ big in $Z$  \Xl. $ \xch $ $big/
\xdm^{+}$ in  \Xl.

(3.4) Analogue for starting with medium in $X,$ instead of big in $X.$

Is upward already induction?

What about representativity (true induction), which replaces size?

Is this logical similarity, whereas the others are size similarities?
\subsubsection{
Remarks:
}

Note that $ \xba \xcn \xbb $ $ \xch $ $ \xba \xcN \xCN \xbb $ is less than
$ \xba \xcN \xcT $

From $ \xba \xcn \xbb,$ $ \xba \xcn \xbb ' $ $ \xch $ $ \xba \xcu \xbb
\xcN \xCN \xbb ' $ to $(R \xcp ):$
$A= \xCN \xbb ',$ $X= \xba \xcu \xbb,$ $Y= \xba,$ go backwards.

Is $(R \xfb )$ $X \xbe \xdf (Y),$ $Y \xbe \xdm^{+}(Z)$ $ \xch $ $X \xbe
\xdm^{+}(Z)$ or similar? Yes, still prove!

Other rules:

(1) For $(R \xcv disj)$ $A \xcs B= \xCQ $ not needed, even better without

(2) For $(R( \xfb ))$ 2. also $ \xba \xcn \xbb,$ $ \xba \xcn \xbb ' $ $
\xch $ $ \xba \xcN \xCN \xbb \xco \xCN \xbb ' $

(3) $small+small+small \xEd all$ multitude of such rules,

$A \xbe \xdf (X),$ $X \xbe \xdf (Y),$ $Y \xbe \xdf (Z)$ $ \xch $ $A \xbe
\xdm^{+}(Z)$ etc.

These notions are related to nonmonotonic logics as follows:

We can say that, normally, $ \xbf $ implies $ \xbq $ iff in a big subset
of all $ \xbf -$cases, $ \xbq $
holds. In preferential terms, $ \xbf $ implies $ \xbq $ iff $ \xbq $ holds
in all minimal
$ \xbf -$models. If $ \xbm $ is the model choice function of a
preferential structure, i.e.
$ \xbm ( \xbf )$ is the set of minimal $ \xbf -$models, then $ \xbm ( \xbf
)$ will be a (the smallest)
big subset of the set of $ \xbf -$models, and the filter over the $ \xbf
-$models is the
pricipal filter generated by $ \xbm ( \xbf ).$

Due to the finite intersection property, filters and ideals work well with
logics:
If $ \xbf $ holds normally, as it holds in a big subset, and so does $
\xbf ',$ then
$ \xbf \xcu \xbf ' $ will normally hold, too, as the intersection of two
big subsets is
big again. This is a nice property, but not justified in all situations,
consider e.g. simple counting of a finite subset. (The question has a
name,
``lottery paradox'': normally no single participant wins, but someone wins
in the
end.) This motivates the weak versions.

Normality defined by (weak or not) filters is a local concept: the filter
defined on $X$ and the one defined on $X' $ might be totally independent.

Seen more abstractly, set properties like e.g. $(R \xfB )$
allow the transfer of big (or small) subsets from one to
another base set (and the conclusions drawn on this basis), and we call
them
``coherence properties''. They are very important, not only for working with
a logic which respects them, but also for soundness and completeness
questions,
often they are at the core of such problems.

 karl-search= End Definition Filter2
\vspace{7mm}

 *************************************

\vspace{7mm}

\subsubsection{Remark Ref-Class}

 {\LARGE karl-search= Start Remark Ref-Class }

\index{Remark Ref-Class}

\br

$\hspace{0.01em}$

(+++ Orig. No.:  Remark Ref-Class +++)

{\xssc LABEL: {Remark Ref-Class}}
\label{Remark Ref-Class}

$(R \xfB )$ corresponds to $(I \xfb )$ and $(F \xfB ):$ If $ \xCf A$ is
small in $X \xcc Y,$ then it will a
fortiori be small in the bigger $Y.$

$(R \xfb )$ says that diminishing base sets by a small
amount will keep small subsets small. This goes in the wrong direction, so
we
have to be careful. We cannot diminish arbitrarily, e.g., if $ \xCf A$ is
a small
subset of $B,$ $ \xCf A$ should not be a small subset of $B-(B-A)=A.$ It
still seems quite
safe, if ``small'' is a robust notion, i.e. defined in an abstract way, and
not
anew for each set, and, if ``small'' is sufficiently far from ``big'', as, for
example in a filter.

There is, however, an important conceptual distinction to make here.
Filters
express ``size'' in an abstract way, in the context of
nonmonotonic logics, $ \xba \xcn \xbb $ iff the set of
$ \xba \xcu \xCN \xbb $ is small in $ \xba.$ But here, we were
interested in ``small'' changes in the reference set $X$ (or $ \xba $ in our
example). So
we have two quite different uses of ``size'', one for
nonmonotonic logics, abstractly expressed by
a filter, the other for coherence conditions. It is possible, but not
necessary,
to consider both essentially the same notions. But we should not forget
that we
have two conceptually different uses of size here.

\er

$(R \xfb \xfb )$ is obviously a stronger variant of $(R \xfb ).$

It and its strength is perhaps best understood as transitivity of the
relation xSy $: \xcj $ $x \xbe \xdm^{+}(y).$

Now, (in comparison to $(R \xfb ))$ $A' $ can be a medium size subset of
$B.$
As a matter of fact, $(R \xfb \xfb )$ is a very big strengthening of $(R
\xfb ):$ Consider a
principal filter $ \xdf:=\{X \xcc B:$ $B' \xcc X\},$ $b \xbe B'.$ Then
$\{b\}$ has at least medium size, so
any small set $A \xcc B$ is smaller than $\{b\}$ - and this is, of course,
just
rankedness. If we only
have $(R \xfb ),$ then we need the whole generating set $B' $ to see that
$ \xCf A$ is small.
This is the strong substitution property of rankedness: any $b$ as above
will
show that $ \xCf A$ is small.

The more we see size as an abstract notion, and the more we see
``small'' different from ``big'' (or ``medium'' ), the more we can go from one
base set
to another and find the same sizes - the more we have coherence when we
reason
with small and big subsets.
$(R \xfb )$ works with iterated use of ``small'', just as do filters, but
not
weak filters. So it is not surprising that weak filters and $(R \xfb )$ do
not
cooperate well: Let $A,B,C$ be small subsets of $X$ - pairwise disjoint,
and
$A \xcv B \xcv C=X,$ this is possible. By $(R \xfb )$ $B$ and $C$ will be
small in $X-A$, so again
by $(R \xfb )$ $C$ will be small in $(X-A)-B=C,$ but this is absurd.

If we think that filters are too strong, but we still want some coherence,
i.e.
abstract size, we can consider $(R \xcv disj):$ If $ \xCf A$ is a
small subset of $B,$ and $A' $ of $B',$ and $B$ and $B' $ are disjoint,
then $A \xcv A' $ is a
small subset of $B \xcv B'.$ It
expresses a uniform approach to size, or distributivity, if you like. It
holds,
e.g. when we consider a set to be small iff it is smaller than a certain
fraction. The important point is here that by disjointness, the big
subsets do
not get ``used up''. (This property generalizes in a straightforward way to
the
infinite case.)

 karl-search= End Remark Ref-Class
\vspace{7mm}

 *************************************

\vspace{7mm}

\subsubsection{Fact R-down}

 {\LARGE karl-search= Start Fact R-down }

\index{Fact R-down}

\bfa

$\hspace{0.01em}$

(+++ Orig. No.:  Fact R-down +++)

{\xssc LABEL: {Fact R-down}}
\label{Fact R-down}

The two versions of $(R \xfb )$ and the three versions of $(R \xfb \xfb )$
are each equivalent.
For the third version of $(R \xfb \xfb )$ we use $(I \xfb ).$

 karl-search= End Fact R-down
\vspace{7mm}

 *************************************

\vspace{7mm}

\subsubsection{Fact R-down Proof}

 {\LARGE karl-search= Start Fact R-down Proof }

\index{Fact R-down Proof}

\efa

\subparagraph{
Proof
}

$\hspace{0.01em}$

(+++*** Orig.:  Proof )

For $A,B \xcc X,$ $(X-B)-((X-A)-B)=A-$B.

`` $ \xch $ '': Let $A \xbe \xdf (X),$ $B \xbe \xdi (X),$ so $X-A \xbe \xdi
(X),$ so by prerequisite
$(X-A)-B \xbe \xdi (X-$B), so $A-B=(X-B)-((X-A)-B) \xbe \xdf (X-$B).

`` $ \xci $ '': Let $A,B \xbe \xdi (X),$ so $X-A \xbe \xdf (X),$ so by
prerequisite $(X-A)-B \xbe \xdf (X-$B),
so $A-B=(X-B)-((X-A)-B) \xbe \xdi (X-$B).

The proof for $(R \xfb \xfb )$ is the same for the first two cases.

It remains to show equivalence with the last one. We assume closure under
set difference and union.

$(1) \xch (3):$

Suppose $A \xce \xdm^{+}(Y),$ but $X \xbe \xdm^{+}(Y),$ we show $A \xce
\xdm^{+}(X).$ So $A \xbe \xdi (Y),$ $Y-X \xce \xdf (Y),$
so $A=A-(Y-X) \xbe \xdi (Y-(Y-X))= \xdi (X).$

$(3) \xch (1):$

Suppose $A-B \xce \xdi (X-$B), $B \xce \xdf (X),$ we show $A \xce \xdi
(X).$ By prerequisite $A-B \xbe \xdm^{+}(X-$B),
$X-B \xbe \xdm^{+}(X),$ so $A-B \xbe \xdm^{+}(X),$ so by $(I \xfb )$ $A
\xbe \xdm^{+}(X),$ so $A \xce \xdi (X).$

$ \xcz $
\\[3ex]

 karl-search= End Fact R-down Proof
\vspace{7mm}

 *************************************

\vspace{7mm}

\subsubsection{Proposition Ref-Class-Mu}

 {\LARGE karl-search= Start Proposition Ref-Class-Mu }

\index{Proposition Ref-Class-Mu}

\bp

$\hspace{0.01em}$

(+++ Orig. No.:  Proposition Ref-Class-Mu +++)

{\xssc LABEL: {Proposition Ref-Class-Mu}}
\label{Proposition Ref-Class-Mu}

If $f(X)$ is the smallest $ \xCf A$ s.t. $A \xbe \xdf (X),$ then, given
the property on the
left, the one on the right follows.

Conversely, when we define $ \xdf (X):=\{X':f(X) \xcc X' \xcc X\},$ given
the property on
the right, the one on the left follows. For this direction, we assume
that we can use the full powerset of some base set $U$ - as is the case
for
the model sets of a finite language. This is perhaps not too bold, as
we mainly want to stress here the intuitive connections, without putting
too much weight on definability questions.

{\footnotesize

\begin{tabular}{|c|c|c|c|}

\hline

(1.1)
\xEH
$(R \xfB )$
\xEH
$ \xch $
\xEH
$( \xbm wOR)$
\xEP

\cline{1-1}
\cline{3-3}

(1.2)
\xEH
\xEH
$ \xci $
\xEH
\xEP

\hline

(2.1)
\xEH
$(R \xfB )+(I \xcv )$
\xEH
$ \xch $
\xEH
$( \xbm OR)$
\xEP

\cline{1-1}
\cline{3-3}

(2.2)
\xEH
\xEH
$ \xci $
\xEH
\xEP

\hline

(3.1)
\xEH
$(R \xfB )+(I \xcv )$
\xEH
$ \xch $
\xEH
$( \xbm PR)$
\xEP

\cline{1-1}
\cline{3-3}

(3.2)
\xEH
\xEH
$ \xci $
\xEH
\xEP

\hline

(4.1)
\xEH
$(R \xcv disj )$
\xEH
$ \xch $
\xEH
$( \xbm disjOR)$
\xEP

\cline{1-1}
\cline{3-3}

(4.2)
\xEH
\xEH
$ \xci $
\xEH
\xEP

\hline

(5.1)
\xEH
$(R \xfb)$
\xEH
$ \xch $
\xEH
$( \xbm CM)$
\xEP

\cline{1-1}
\cline{3-3}

(5.2)
\xEH
\xEH
$ \xci $
\xEH
\xEP

\hline

(6.1)
\xEH
$(R \xfb \xfb)$
\xEH
$ \xch $
\xEH
$( \xbm RatM)$
\xEP

\cline{1-1}
\cline{3-3}

(6.2)
\xEH
\xEH
$ \xci $
\xEH
\xEP

\hline

\end{tabular}

}

 karl-search= End Proposition Ref-Class-Mu
\vspace{7mm}

 *************************************

\vspace{7mm}

\subsubsection{Proposition Ref-Class-Mu Proof}

 {\LARGE karl-search= Start Proposition Ref-Class-Mu Proof }

\index{Proposition Ref-Class-Mu Proof}

\ep

\subparagraph{
Proof
}

$\hspace{0.01em}$

(+++*** Orig.:  Proof )

(1.1) $(R \xfB )$ $ \xch $ $( \xbm wOR):$

$X-f(X)$ is small in $X,$ so it is small in $X \xcv Y$ by $(R \xfB ),$ so
$A:=X \xcv Y-(X-f(X)) \xbe \xdf (X \xcv Y),$ but $A \xcc f(X) \xcv Y,$ and
$f(X \xcv Y)$ is the smallest element
of $ \xdf (X \xcv Y),$ so $f(X \xcv Y) \xcc A \xcc f(X) \xcv Y.$

(1.2) $( \xbm wOR)$ $ \xch $ $(R \xfB ):$

Let $X \xcc Y,$ $X':=Y-$X. Let $A \xbe \xdi (X),$ so $X-A \xbe \xdf (X),$
so $f(X) \xcc X-$A, so
$f(X \xcv X' ) \xcc f(X) \xcv X' \xcc (X-A) \xcv X' $ by prerequisite, so
$(X \xcv X' )-((X-A) \xcv X' )=A \xbe \xdi (X \xcv X' ).$

(2.1) $(R \xfB )+(I \xcv )$ $ \xch $ $( \xbm OR):$

$X-f(X)$ is small in $X,$ $Y-f(Y)$ is small in $Y,$ so both are small in
$X \xcv Y$ by
$(R \xfB ),$ so $A:=(X-f(X)) \xcv (Y-f(Y))$ is small in $X \xcv Y$ by $(I
\xcv ),$ but
$X \xcv Y-(f(X) \xcv f(Y)) \xcc A,$ so $f(X) \xcv f(Y) \xbe \xdf (X \xcv
Y),$ so, as $f(X \xcv Y)$ is the smallest
element of $ \xdf (X \xcv Y),$ $f(X \xcv Y) \xcc f(X) \xcv f(Y).$

(2.2) $( \xbm OR)$ $ \xch $ $(R \xfB )+(I \xcv ):$

Let again $X \xcc Y,$ $X':=Y-$X. Let $A \xbe \xdi (X),$ so $X-A \xbe \xdf
(X),$ so $f(X) \xcc X-$A. $f(X' ) \xcc X',$
so $f(X \xcv X' ) \xcc f(X) \xcv f(X' ) \xcc (X-A) \xcv X' $ by
prerequisite, so
$(X \xcv X' )-((X-A) \xcv X' )=A \xbe \xdi (X \xcv X' ).$

$(I \xcv )$ holds by definition.

(3.1) $(R \xfB )+(I \xcv )$ $ \xch $ $( \xbm PR):$

Let $X \xcc Y.$ $Y-f(Y)$ is the largest element of $ \xdi (Y),$ $X-f(X)
\xbe \xdi (X) \xcc \xdi (Y)$ by
$(R \xfB ),$ so $(X-f(X)) \xcv (Y-f(Y)) \xbe \xdi (Y)$ by $(I \xcv ),$ so
by ``largest'' $X-f(X) \xcc Y-f(Y),$
so $f(Y) \xcs X \xcc f(X).$

(3.2) $( \xbm PR)$ $ \xch $ $(R \xfB )+(I \xcv )$

Let again $X \xcc Y,$ $X':=Y-$X. Let $A \xbe \xdi (X),$ so $X-A \xbe \xdf
(X),$ so $f(X) \xcc X-$A, so
by prerequisite $f(Y) \xcs X \xcc X-$A, so $f(Y) \xcc X' \xcv (X-$A), so
$(X \xcv X' )-(X' \xcv (X-A))=A \xbe \xdi (Y).$

Again, $(I \xcv )$ holds by definition.

(4.1) $(R \xcv disj)$ $ \xch $ $( \xbm disjOR):$

If $X \xcs Y= \xCQ,$ then (1) $A \xbe \xdi (X),B \xbe \xdi (Y) \xch A
\xcv B \xbe \xdi (X \xcv Y)$ and
(2) $A \xbe \xdf (X),B \xbe \xdf (Y) \xch A \xcv B \xbe \xdf (X \xcv Y)$
are equivalent. (By $X \xcs Y= \xCQ,$
$(X-A) \xcv (Y-B)=(X \xcv Y)-(A \xcv B).)$
So $f(X) \xbe \xdf (X),$ $f(Y) \xbe \xdf (Y)$ $ \xch $ (by prerequisite)
$f(X) \xcv f(Y) \xbe \xdf (X \xcv Y).$ $f(X \xcv Y)$
is the smallest element of $ \xdf (X \xcv Y),$ so $f(X \xcv Y) \xcc f(X)
\xcv f(Y).$

(4.2) $( \xbm disjOR)$ $ \xch $ $(R \xcv disj):$

Let $X \xcc Y,$ $X':=Y-$X. Let $A \xbe \xdi (X),$ $A' \xbe \xdi (X' ),$
so $X-A \xbe \xdf (X),$ $X' -A' \xbe \xdf (X' ),$
so $f(X) \xcc X-$A, $f(X' ) \xcc X' -A',$ so $f(X \xcv X' ) \xcc f(X)
\xcv f(X' ) \xcc (X-A) \xcv (X' -A' )$ by
prerequisite, so $(X \xcv X' )-((X-A) \xcv (X' -A' ))=A \xcv A' \xbe \xdi
(X \xcv X' ).$

(5.1) $(R \xfb )$ $ \xch $ $( \xbm CM):$

$f(X) \xcc Y \xcc X$ $ \xch $ $X-Y \xbe \xdi (X),$ $X-f(X) \xbe \xdi (X)$
$ \xch_{(R \xfb )}$ $A:=(X-f(X))-(X-Y) \xbe \xdi (Y)$ $ \xch $
$Y-A=f(X)-(X-Y) \xbe \xdf (Y)$ $ \xch $ $f(Y) \xcc f(X)-(X-Y) \xcc f(X).$

(5.2) $( \xbm CM)$ $ \xch $ $(R \xfb )$

Let $A \xbe \xdf (X),$ $B \xbe \xdi (X),$ so $f(X) \xcc X-B \xcc X,$ so by
prerequisite $f(X-B) \xcc f(X).$
As $A \xbe \xdf (X),$ $f(X) \xcc A,$ so $f(X-B) \xcc f(X) \xcc A \xcs
(X-B)=A-$B, and $A-B \xbe \xdf (X-$B).

(6.1) $(R \xfb \xfb )$ $ \xch $ $( \xbm RatM):$

Let $X \xcc Y,$ $X \xcs f(Y) \xEd \xCQ.$ If $Y-X \xbe \xdf (Y),$ then
$A:=(Y-X) \xcs f(Y) \xbe \xdf (Y),$ but by
$X \xcs f(Y) \xEd \xCQ $ $A \xcb f(Y),$ contradicting ``smallest'' of
$f(Y).$ So $Y-X \xce \xdf (Y),$ and
by $(R \xfb \xfb )$ $X-f(Y)=(Y-f(Y))-(Y-X) \xbe \xdi (X),$ so $X \xcs f(Y)
\xbe \xdf (X),$ so $f(X) \xcc f(Y) \xcs X.$

(6.2) $( \xbm RatM)$ $ \xch $ $(R \xfb \xfb )$

Let $A \xbe \xdf (Y),$ $B \xce \xdf (Y).$ $B \xce \xdf (Y)$ $ \xch $ $Y-B
\xce \xdi (Y)$ $ \xch $ $(Y-B) \xcs f(Y) \xEd \xCQ.$
Set $X:=Y-$B, so $X \xcs f(Y) \xEd \xCQ,$ $X \xcc Y,$ so $f(X) \xcc f(Y)
\xcs X$ by prerequisite.
$f(Y) \xcc A$ $ \xch $ $f(X) \xcc f(Y) \xcs X=f(Y)-B \xcc A-$B.

$ \xcz $
\\[3ex]

 karl-search= End Proposition Ref-Class-Mu Proof
\vspace{7mm}

 *************************************

\vspace{7mm}

\subsubsection{Definition Nabla}

 {\LARGE karl-search= Start Definition Nabla }

\index{Definition Nabla}

\bd

$\hspace{0.01em}$

(+++ Orig. No.:  Definition Nabla +++)

{\xssc LABEL: {Definition Nabla}}
\label{Definition Nabla}

Augment the language of first order logic by the new quantifier:
If $ \xbf $ and $ \xbq $ are formulas, then so are $ \xeA x \xbf (x),$ $
\xeA x \xbf (x): \xbq (x),$
for any variable $x.$ The:-versions are the restricted variants.
We call any formula of $ \xdl,$ possibly containing $ \xeA $ a $ \xeA -
\xdl -$formula.

 karl-search= End Definition Nabla
\vspace{7mm}

 *************************************

\vspace{7mm}

\subsubsection{Definition N-Model}

 {\LARGE karl-search= Start Definition N-Model }

\index{Definition N-Model}

\ed

\bd

$\hspace{0.01em}$

(+++ Orig. No.:  Definition N-Model +++)

{\xssc LABEL: {Definition N-Model}}
\label{Definition N-Model}

$( \xdn -$Model)

Let $ \xdl $ be a first order language, and $M$ be a $ \xdl -$structure.
Let $ \xdn (M)$ be
a weak filter, or $ \xdn -$system - $ \xdn $ for normal - over $M.$
Define $<M, \xdn (M)>$ $ \xcm $ $ \xbf $ for any $ \xeA - \xdl -$formula
inductively as usual, with
one additional induction step:

$<M, \xdn (M)>$ $ \xcm $ $ \xeA x \xbf (x)$ iff there is $A \xbe \xdn (M)$
s.t. $ \xcA a \xbe A$ $(<M, \xdn (M)>$ $ \xcm $ $ \xbf [a]).$

 karl-search= End Definition N-Model
\vspace{7mm}

 *************************************

\vspace{7mm}

\subsubsection{Definition NablaAxioms}

 {\LARGE karl-search= Start Definition NablaAxioms }

\index{Definition NablaAxioms}

\ed

\bd

$\hspace{0.01em}$

(+++ Orig. No.:  Definition NablaAxioms +++)

{\xssc LABEL: {Definition NablaAxioms}}
\label{Definition NablaAxioms}

Let any axiomatization of predicate calculus be given. Augment this with
the axiom schemata

(1) $ \xeA x \xbf (x)$ $ \xcu $ $ \xcA x( \xbf (x) \xcp \xbq (x))$ $ \xcp
$ $ \xeA x \xbq (x),$

(2) $ \xeA x \xbf (x)$ $ \xcp $ $ \xCN \xeA x \xCN \xbf (x),$

(3) $ \xcA x \xbf (x)$ $ \xcp $ $ \xeA x \xbf (x)$ and $ \xeA x \xbf (x)$
$ \xcp $ $ \xcE x \xbf (x),$

(4) $ \xeA x \xbf (x)$ $ \xcr $ $ \xeA y \xbf (y)$ if $x$ does not occur
free in $ \xbf (y)$ and $y$ does not
occur free in $ \xbf (x).$

(for all $ \xbf,$ $ \xbq )$.

 karl-search= End Definition NablaAxioms
\vspace{7mm}

 *************************************

\vspace{7mm}

\subsubsection{Proposition NablaRepr}

 {\LARGE karl-search= Start Proposition NablaRepr }

\index{Proposition NablaRepr}

\ed

\bp

$\hspace{0.01em}$

(+++ Orig. No.:  Proposition NablaRepr +++)

{\xssc LABEL: {Proposition NablaRepr}}
\label{Proposition NablaRepr}

The axioms given in Definition \ref{Definition NablaAxioms} (page
\pageref{Definition NablaAxioms})
are sound and complete for the semantics of Definition \ref{Definition N-Model}
(page \pageref{Definition N-Model})

See  \cite{Sch95-1} or  \cite{Sch04}.

 karl-search= End Proposition NablaRepr
\vspace{7mm}

 *************************************

\vspace{7mm}

\ep

\paragraph{
Extension to normal defaults with prerequisites
}

$\hspace{0.01em}$

(+++*** Orig.:  Extension to normal defaults with prerequisites )

{\xssc LABEL: {Section Extension to normal defaults with prerequisites}}
\label{Section Extension to normal defaults with prerequisites}

\subsubsection{Definition Nabla-System}

 {\LARGE karl-search= Start Definition Nabla-System }

\index{Definition Nabla-System}

\bd

$\hspace{0.01em}$

(+++ Orig. No.:  Definition Nabla-System +++)

{\xssc LABEL: {Definition Nabla-System}}
\label{Definition Nabla-System}

Call $ \xdn^{+}(M)=< \xdn (N):N \xcc M>$ a $ \xdn^{+}-system$ or system of
weak filters over $M$ iff
for each $N \xcc M$ $ \xdn (N)$ is a weak filter or $ \xdn -$system over
$N.$
(It suffices to consider the definable subsets of $M.)$

 karl-search= End Definition Nabla-System
\vspace{7mm}

 *************************************

\vspace{7mm}

\subsubsection{Definition N-Model-System}

 {\LARGE karl-search= Start Definition N-Model-System }

\index{Definition N-Model-System}

\ed

\bd

$\hspace{0.01em}$

(+++ Orig. No.:  Definition N-Model-System +++)

{\xssc LABEL: {Definition N-Model-System}}
\label{Definition N-Model-System}

Let $ \xdl $ be a first order language, and $M$ a $ \xdl -$structure. Let
$ \xdn^{+}(M)$ be
a $ \xdn^{+}-system$ over $M.$

Define $<M, \xdn^{+}(M)>$ $ \xcm $ $ \xbf $ for any formula inductively as
usual, with
the additional induction steps:

1. $<M, \xdn^{+}(M)>$ $ \xcm $ $ \xeA x \xbf (x)$ iff there is $A \xbe
\xdn (M)$ s.t. $ \xcA a \xbe A$ $(<M, \xdn^{+}(M)>$ $ \xcm $ $ \xbf [a]),$

2. $<M, \xdn^{+}(M)>$ $ \xcm $ $ \xeA x \xbf (x): \xbq (x)$ iff there is
$A \xbe \xdn (\{x:<M, \xdn^{+}(M)> \xcm \xbf (x)\})$ s.t.
$ \xcA a \xbe A$ $(<M, \xdn^{+}(M)>$ $ \xcm $ $ \xbq [a]).$

 karl-search= End Definition N-Model-System
\vspace{7mm}

 *************************************

\vspace{7mm}

\subsubsection{Definition NablaAxioms-System}

 {\LARGE karl-search= Start Definition NablaAxioms-System }

\index{Definition NablaAxioms-System}

\ed

\bd

$\hspace{0.01em}$

(+++ Orig. No.:  Definition NablaAxioms-System +++)

{\xssc LABEL: {Definition NablaAxioms-System}}
\label{Definition NablaAxioms-System}

Extend the logic of first order predicate calculus by adding the axiom
schemata

(1) a. $ \xeA x \xbf (x)$ $ \xcr $ $ \xeA x(x=x): \xbf (x),$
$b.$ $ \xcA x( \xbs (x) \xcr \xbt (x))$ $ \xcu $ $ \xeA x \xbs (x): \xbf
(x)$ $ \xcp $ $ \xeA x \xbt (x): \xbf (x),$

(2) $ \xeA x \xbf (x): \xbq (x)$ $ \xcu $ $ \xcA x( \xbf (x) \xcu \xbq (x)
\xcp \xbj (x))$ $ \xcp $ $ \xeA x \xbf (x): \xbj (x),$

(3) $ \xcE x \xbf (x)$ $ \xcu $ $ \xeA x \xbf (x): \xbq (x)$ $ \xcp $ $
\xCN \xeA x \xbf (x): \xCN \xbq (x),$

(4) $ \xcA x( \xbf (x) \xcp \xbq (x))$ $ \xcp $ $ \xeA x \xbf (x): \xbq
(x)$
and $ \xeA x \xbf (x): \xbq (x)$ $ \xcp $ $[ \xcE x \xbf (x)$ $ \xcp $ $
\xcE x( \xbf (x) \xcu \xbq (x))],$

(5) $ \xeA x \xbf (x): \xbq (x)$ $ \xcr $ $ \xeA y \xbf (y): \xbq (y)$
(under the usual caveat for substitution).

(for all $ \xbf,$ $ \xbq,$ $ \xbj,$ $ \xbs,$ $ \xbt )$.

 karl-search= End Definition NablaAxioms-System
\vspace{7mm}

 *************************************

\vspace{7mm}

\subsubsection{Proposition NablaRepr-System}

 {\LARGE karl-search= Start Proposition NablaRepr-System }

\index{Proposition NablaRepr-System}

\ed

\bp

$\hspace{0.01em}$

(+++ Orig. No.:  Proposition NablaRepr-System +++)

{\xssc LABEL: {Proposition NablaRepr-System}}
\label{Proposition NablaRepr-System}

The axioms of Definition \ref{Definition NablaAxioms-System} (page
\pageref{Definition NablaAxioms-System})  are
sound and complete for the $ \xdn^{+}-semantics$
of $ \xeA $ as defined in Definition \ref{Definition N-Model-System} (page
\pageref{Definition N-Model-System}) .

See  \cite{Sch95-1} or  \cite{Sch04}.

 karl-search= End Proposition NablaRepr-System
\vspace{7mm}

 *************************************

\vspace{7mm}

\subsubsection{Size-Bib}

 {\LARGE karl-search= Start Size-Bib }

\index{Size-Bib}

\ep

More on different abstract coherent systems based on size,

 \xEI
 \xDH the system of $S.$ Ben-David and $R.$ Ben-Eliyahu (see  \cite{BB94}),
 \xDH the system of the author,
 \xDH the system of $N.$ Friedman and $J.$ Halpern (see  \cite{FH98}).
 \xEJ

can be found in  \cite{Sch04}.

 karl-search= End Size-Bib
\vspace{7mm}

 *************************************

\vspace{7mm}

 karl-search= End ToolBase1-Size
\vspace{7mm}

 *************************************

\vspace{7mm}

\section{
IBRS
}

\subsection{
}

\subsubsection{ToolBase1-IBRS}

 {\LARGE karl-search= Start ToolBase1-IBRS }

{\xssc LABEL: {Section Toolbase1-IBRS}}
\label{Section Toolbase1-IBRS}
\index{Section Toolbase1-IBRS}
\subsubsection{Motivation IBRS}

 {\LARGE karl-search= Start Motivation IBRS }

\index{Motivation IBRS}

The human agent in his daily activity has to deal with many situations
involving change. Chief among them are the following

 \xEh

 \xDH Common sense reasoning from available data. This involves
predication
of what unavailable data is supposed to be (nonmonotonic deduction)
but it is a defeasible prediction, geared towards immediate change.
This is formally known as nonmonotonic reasoning and is studied by
the nonmonotonic community.

 \xDH Belief revision, studied by a very large community. The agent
is unhappy with the totality of his beliefs which he finds internally
unacceptable (usually logically inconsistent but not necessarily so)
and needs to change/revise it.

 \xDH Receiving and updating his data, studied by the update community.

 \xDH Making morally correct decisions, studied by the deontic logic
community.

 \xDH Dealing with hypothetical and counterfactual situations. This
is studied by a large community of philosophers and AI researchers.

 \xDH Considering temporal future possibilities, this is covered by
modal and temporal logic.

 \xDH Dealing with properties that persist through time in the near
future and with reasoning that is constructive. This is covered by
intuitionistic logic.

 \xEj

All the above types of reasoning exist in the human mind and are used
continuously and coherently every hour of the day. The formal modelling
of these types is done by diverse communities which are largely distinct
with no significant communication or cooperation. The formal models
they use are very similar and arise from a more general theory, what
we might call:

``Reasoning with information bearing binary relations''.

 karl-search= End Motivation IBRS
\vspace{7mm}

 *************************************

\vspace{7mm}

\subsubsection{Definition IBRS}

 {\LARGE karl-search= Start Definition IBRS }

\index{Definition IBRS}

\bd

$\hspace{0.01em}$

(+++ Orig. No.:  Definition IBRS +++)

{\xssc LABEL: {Definition IBRS}}
\label{Definition IBRS}

 \xEh
 \xDH An information bearing binary relation frame IBR, has the form
$(S, \xdR ),$ where $S$ is a non-empty set and $ \xdR $ is a subset of
$S,$ where $S$ is defined by induction as follows:

 \xEh

 \xDH $S_{0}=S$

 \xDH $S_{n+1}$ $=$ $S_{n} \xcv (S_{n} \xCK S_{n}).$

 \xDH $S$ $=$ $ \xcV \{S_{n}:n \xbe \xbo \}$

 \xEj

We call elements from $S$ points or nodes, and elements from $ \xdR $
arrows.
Given $(S, \xdR ),$ we also set $ \xdP ((S, \xdR )):=S,$ and $ \xdA ((S,
\xdR )):= \xdR.$

If $ \xba $ is an arrow, the origin and destination of $ \xba $ are
defined
as usual, and we write $ \xba:x \xcp y$ when $x$ is the origin, and $y$
the destination of the arrow $ \xba.$ We also write $o( \xba )$ and $d(
\xba )$ for
the origin and destination of $ \xba.$

 \xDH Let $Q$ be a set of atoms, and $ \xdL $ be a set of labels (usually
$\{0,1\}$ or $[0,1]).$ An information assignment $h$ on $(S, \xdR )$
is a function $h:Q \xCK \xdR \xcp \xdL.$

 \xDH An information bearing system IBRS, has the form
$(S, \xdR,h,Q, \xdL ),$ where $S,$ $ \xdR,$ $h,$ $Q,$ $ \xdL $ are as
above.

 \xEj

See Diagram \ref{Diagram IBRS} (page \pageref{Diagram IBRS})  for an
illustration.

\vspace{10mm}

\begin{diagram}

\centering
\setlength{\unitlength}{0.00083333in}
{\renewcommand{\dashlinestretch}{30}
\begin{picture}(4961,5004)(0,0)
\path(1511,1583)(611,3683)
\blacken\path(685.845,3584.520)(611.000,3683.000)(630.696,3560.885)(672.451,3539.613)(685.845,3584.520)
\path(1511,1583)(2411,3683)
\blacken\path(2391.304,3560.885)(2411.000,3683.000)(2336.155,3584.520)(2349.549,3539.613)(2391.304,3560.885)
\path(3311,1583)(4361,4133)
\blacken\path(4343.050,4010.616)(4361.000,4133.000)(4287.570,4033.461)(4301.603,3988.750)(4343.050,4010.616)
\path(3316,1574)(2416,3674)
\blacken\path(2490.845,3575.520)(2416.000,3674.000)(2435.696,3551.885)(2477.451,3530.613)(2490.845,3575.520)
\path(986,2783)(2621,2783)
\blacken\path(2501.000,2753.000)(2621.000,2783.000)(2501.000,2813.000)(2465.000,2783.000)(2501.000,2753.000)
\path(2486,2783)(2786,2783)
\blacken\path(2666.000,2753.000)(2786.000,2783.000)(2666.000,2813.000)(2630.000,2783.000)(2666.000,2753.000)
\path(3311,1583)(2051,2368)
\blacken\path(2168.714,2330.008)(2051.000,2368.000)(2136.987,2279.083)(2183.406,2285.509)(2168.714,2330.008)
\path(2166,2288)(1906,2458)
\blacken\path(2022.854,2417.439)(1906.000,2458.000)(1990.019,2367.221)(2036.567,2372.629)(2022.854,2417.439)

\put(1511,1358) {{\xssc $a$}}
\put(3311,1358) {{\xssc $d$}}
\put(3311,1058)  {{\xssc $(p,q)=(1,0)$}}
\put(1511,1058)  {{\xssc $(p,q)=(0,0)$}}
\put(2411,3758){{\xssc $c$}}
\put(4361,4433){{\xssc $(p,q)=(1,1)$}}
\put(4361,4208){{\xssc $e$}}
\put(2411,3983){{\xssc $(p,q)=(0,1)$}}
\put(611,3983) {{\xssc $(p,q)=(0,1)$}}
\put(611,3758) {{\xssc $b$}}
\put(1211,2883){{\xssc $(p,q)=(1,1)$}}
\put(260,2333) {{\xssc $(p,q)=(1,1)$}}
\put(2261,1583) {{\xssc $(p,q)=(1,1)$}}
\put(1286,3233){{\xssc $(p,q)=(1,1)$}}
\put(2711,3083){{\xssc $(p,q)=(1,1)$}}
\put(3836,2633){{\xssc $(p,q)=(1,1)$}}

\put(300,700)
{{\rm\bf
A simple example of an information bearing system.
}}

\end{picture}
}
{\xssc LABEL: {Diagram IBRS-a}}
\label{Diagram IBRS-a}
\index{Diagram IBRS}

\end{diagram}

We have here:
\[\begin{array}{l}
S =\{a,b,c,d,e\}.\\
\xdR = S \cup \{(a,b), (a,c), (d,c), (d,e)\} \cup \{((a,b), (d,c)),
(d,(a,c))\}.\\
Q = \{p,q\}
\end{array}
\]
The values of $h$ for $p$ and $q$ are as indicated in the figure. For
example $h(p,(d,(a,c))) =1$.

\vspace{4mm}

 karl-search= End Definition IBRS
\vspace{7mm}

 *************************************

\vspace{7mm}

\subsubsection{Comment IBRS}

 {\LARGE karl-search= Start Comment IBRS }

\index{Comment IBRS}

\ed

\bcom

$\hspace{0.01em}$

(+++ Orig. No.:  Comment +++)

{\xssc LABEL: {Comment}}
\label{Comment}

{\xssc LABEL: {Comment IBRS}}
\label{Comment IBRS}

The elements in Figure Diagram \ref{Diagram IBRS} (page \pageref{Diagram IBRS}) 
can be interpreted in
many ways,
depending on the area of application.

 \xEh

 \xDH The points in $S$ can be interpreted as possible worlds, or
as nodes in an argumentation network or nodes in a neural net or
states, etc.

 \xDH The direct arrows from nodes to nodes can be interpreted as
accessibility relation, attack or support arrows in an argumentation
networks, connection in a neural nets, a preferential ordering in
a nonmonotonic model, etc.

 \xDH The labels on the nodes and arrows can be interpreted as fuzzy
values in the accessibility relation or weights in the neural net
or strength of arguments and their attack in argumentation nets, or
distances in a counterfactual model, etc.

 \xDH The double arrows can be interpreted as feedback loops to nodes
or to connections, or as reactive links changing the system which are
activated
as we pass between the nodes.

 \xEj

 karl-search= End Comment IBRS
\vspace{7mm}

 *************************************

\vspace{7mm}

\subsection{IBRS as abstraction}
\subsubsection{IBRS as abstraction}

 {\LARGE karl-search= Start IBRS as abstraction }

\index{IBRS as abstraction}
{\xssc LABEL: {IBRS as abstraction}}
\label{IBRS as abstraction}

\ecom

Thus, IBRS can be used as a source of information for various logics
based on the atoms in $Q.$ We now illustrate by listing several such
logics.
\subsubsection*{Modal Logic}

One can consider the figure as giving rise to two modal logic models.
One with actual world a and one with $d,$ these being the two minimal
points of the relation. Consider a language with $ \xcX q.$ how
do we evaluate $a \xcm \xcX q?$

The modal logic will have to give an algorithm for calculating the
values.

Say we choose algorithm $ \xda_{1}$ for $a \xcm \xcX q,$ namely:

[
$ \xda_{1}(a, \xcX q)=1$
]
iff for all $x \xbe S$ such that $a=x$ or $(a,x) \xbe \xdR $ we have
$h(q,x)=1.$

According to $ \xda_{1}$ we get that $ \xcX q$ is false at a. $ \xda_{1}$
gives rise to a $T-$modal logic. Note that the reflexivity is not
anchored at the relation $ \xdR $ of the network but in the algorithm
$ \xda_{1}$ in the way we evaluate. We say $(S, \xdR, \Xl.) \xcm \xcX $
$q$ iff $ \xcX q$ holds in all minimal points of $(S, \xdR ).$

For orderings without minimal points we may choose a subset of
distinguished
points.
\subsubsection*{Nonmonotonic Deduction}

We can ask whether $p \xcn q$ according to algorithm $ \xda_{2}$ defined
below. $ \xda_{2}$ says that $p \xcn q$ holds iff $q$ holds in all
minimal models of $p.$ Let us check the value of $ \xda_{2}$ in this
case:

Let $S_{p}=\{s \xbe S \xfA h(p,s)=1\}.$ Thus $S_{p}=\{d,e\}.$

The minimal points of $S_{p}$ are $\{d\}.$ Since $h(q,d)=0,$ we have that
$p \xcN q.$

Note that in the cases of modal logic and
nonmonotonic logic we ignored the arrows $(d,(a,c))$ (i.e. the double
arrow from $d$ to the arc $(a,c))$ and the $h$ values to arcs. These
values do not play a part in the traditional modal or nonmonotonic
logic. They do play a part in other logics. The attentive reader may
already suspect that we have her an opportunity for generalisation
of say nonmonotonic logic, by giving a role to arc annotations.
\subsubsection*{Argumentation Nets}

Here the nodes of $S$ are interpreted as arguments. The atoms $\{p,q\}$
can be interpreted as types of arguments and the arrows e.g. $(a,b) \xbe
\xdR $
as indicating that the argument a is attacking the argument $b.$

So, for example, let
\begin{quote}

a $=$ we must win votes.

$b$ $=$ death sentence for murderers.

$c$ $=$ We must allow abortion for teenagers

$d$ $=$ Bible forbids taking of life.

$q$ $=$ the argument is a social argument

$p$ $=$ the argument is a religious argument.

$(d,(a,c))$ $=$ there should be no connection between winning votes
and
abortion.

$((a,b),(d,c))$ $=$ If we attack the death sentence in order to win
votes then we must stress (attack) that there should be no connection
between religion (Bible) and social issues.
\end{quote}

Thus we have according to this model that supporting abortion can
lose votes. The argument for abortion is a social one and the argument
from the Bible against it is a religious one.

We can extract information from this IBRS using two algorithms. The
modal logic one can check whether for example every social argument
is attacked by a religious argument. The answer is no, since the social
argument $b$ is attacked only by a which is not a religious argument.

We can also use algorithm $ \xda_{3}$ (following Dung) to extract the
winning arguments of this system. The arguments a and $d$ are winning
since they are not attacked. $d$ attacks the connection between a
and $c$ (i.e. stops a attacking $c).$

The attack of a on $b$ is successful and so $b$ is out. However the
arc $(a,b)$ attacks the arc $(d,c).$ So $c$ is not attacked at all
as both arcs leading into it are successfully eliminated. So $c$
is in. $e$ is out because it is attacked by $d.$

So the winning arguments are $\{a,c,d\}$

In this model we ignore the annotations on arcs. To be consistent
in our mathematics we need to say that $h$ is a partial function on
$ \xdR.$ The best way is to give more specific definition on IBRS to
make it suitable for each logic.

See also  \cite{Gab08b} and  \cite{BGW05}.
\subsubsection*{Counterfactuals}

The traditional semantics for counterfactuals involves closeness of
worlds. The clauses $y \xcm p \xfo q,$ where $ \xfo $
is a counterfactual implication is that $q$ holds in all worlds $y' $
``near enough'' to $y$ in which $p$ holds. So if we interpret the
annotation on arcs as distances then we can define ``near'' as distance
$ \xck $ 2, we get: $a \xcm p \xfo q$ iff in all worlds
of $p-$distance $ \xck 2$ if $p$ holds so does $q.$ Note that the distance
depends on $p.$

In this case we get that $a \xcm p \xfo q$ holds. The
distance function can also use the arrows from arcs to arcs, etc.
There are many opportunities for generalisation in our IBRS set up.
\subsubsection*{Intuitionistic Persistence}

We can get an intuitionistic Kripke model out of this IBRS by letting,
for $t,s \xbe S,t \xbr_{0}s$ iff $t=s$ or $[tRs \xcu \xcA q \xbe Q(h(q,t)
\xck h(q,s))].$
We get that

[
$r_{0}=\{(y,y) \xfA y \xbe S\} \xcv \{(a,b),(a,c),(d,e)\}.$
]

Let $ \xbr $ be the transitive closure of $ \xbr_{0}.$ Algorithm $
\xda_{4}$
evaluates $p \xch q$ in this model, where $ \xch $ is intuitionistic
implication.

$ \xda_{4}:$ $p \xch q$ holds at the IBRS iff $p \xch q$ holds
intuitionistically
at every $ \xbr -$minimal point $of(S, \xbr ).$

 karl-search= End IBRS as abstraction
\vspace{7mm}

 *************************************

\vspace{7mm}

\subsubsection{Reac-Sem}

 {\LARGE karl-search= Start Reac-Sem }

{\xssc LABEL: {Section Reac-Sem}}
\label{Section Reac-Sem}
\index{Section Reac-Sem}

\subsection{
Introduction
}

\subsubsection{Reac-Sem-Intro}

 {\LARGE karl-search= Start Reac-Sem-Intro }

{\xssc LABEL: {Section Reac-Sem-Intro}}
\label{Section Reac-Sem-Intro}
\index{Section Reac-Sem-Intro}

(1) Nodes and arrows

As we may have counterarguments not only against nodes, but also against
arrows, they must be treated basically the same way, i.e. in some way
there has
to be a positive, but also a negative influence on both. So arrows cannot
just
be concatenation between the contents of nodes, or so.

We will differentiate between nodes and arrows by labelling arrows in
addition
with a time delay. We see nodes as situations, where the output is
computed
instantenously from the input, whereas arrows describe some ``force'' or
``mechanism'' which may need some time to ``compute'' the result from the
input.

Consequently, if $ \xba $ is an arrow, and $ \xbb $ an arrow pointing to $
\xba,$ then it
should point to the input of $ \xba,$ i.e. before the time lapse.
Conversely,
any arrow originating in $ \xba $ should originate after the time lapse.

Apart this distinction, we will treat nodes and arrows the same way, so
the
following discussion will apply to both - which we call just ``objects''.

(2) Defeasibility

The general idea is to code each object, say $X,$ by $I(X):U(X) \xcp
C(X):$ If $I(X)$
holds then, unless $U(X)$ holds, consequence $C(X)$ will hold. (We adopted
Reiter's
notation for defaults, as IBRS have common points with the former.)

The situation is slightly more complicated, as there can be several
counterarguments, so $U(X)$ really is an ``or''. Likewise, there can be
several
supporting arguments, so $I(X)$ also is an ``or''.

A counterargument must not always be an argument against a specific
supporting
argument, but it can be. Thus, we should admit both possibilties. As we
can use
arrows to arrows, the second case is easy to treat (as is the dual, a
supporting
argument can be against a specific counterargument). How do we treat the
case
of unspecific pro- and counterarguments? Probably the easiest way is to
adopt
Dung's idea: an object is in, if it has at least one support, and no
counterargument - see  \cite{Dun95}.
Of course, other possibilities may be adopted, counting, use
of labels, etc., but we just consider the simple case here.

(3) Labels

In the general case, objects stand for some kind of defeasible
transmission.
We may in some cases see labels as restricting this transmission to
certain
values. For instance, if the label is $p=1$ and $q=0,$ then the $p-$part
may be
transmitted and the $q-$part not.

Thus, a transmission with a label can sometimes be considered as a family
of
transmissions, which ones are active is indicated by the label.

\be

$\hspace{0.01em}$

(+++ Orig. No.:  Example 2.1 +++)

{\xssc LABEL: {Example 2.1}}
\label{Example 2.1}

In fuzzy Kripke models, labels are elements of $[0,1].$ $p=0.5$ as label
for a node
$m' $ which stands for a fuzzy model means that the value of $p$ is 0.5.
$p=0.5$ as
label for an arrow from $m$ to $m' $ means that $p$ is transmitted with
value 0.5.
Thus, when we look from $m$ to $m',$ we see $p$ with value
$0.5*0.5=0.25.$ So, we have
$ \xcx p$ with value 0.25 at $m$ - if, e.g., $m,m' $ are the only models.

\ee

(4) Putting things together

If an arrow leaves an object, the object's output will be connected to the
(only) positive input of the arrow. (An arrow has no negative inputs from
objects it leaves.) If a positive arrow enters an object, it is connected
to
one of the positive inputs of the object, analogously for negative arrows
and
inputs.

When labels are present, they are transmitted through some operation.

 karl-search= End Reac-Sem-Intro
\vspace{7mm}

 *************************************

\vspace{7mm}


\subsection{
Formal definition
}

\subsubsection{Reac-Sem-Def}

 {\LARGE karl-search= Start Reac-Sem-Def }

{\xssc LABEL: {Section Reac-Sem-Def}}
\label{Section Reac-Sem-Def}
\index{Section Reac-Sem-Def}

\bd

$\hspace{0.01em}$

(+++ Orig. No.:  Definition 2.1 +++)

{\xssc LABEL: {Definition 2.1}}
\label{Definition 2.1}

In the most general case, objects of IBRS have the form:
$(<I_{1},L_{1}>, \Xl,<I_{n},L_{n}>):(<U_{1},L'_{1}>, \Xl
,<U_{n},L'_{n}>),$ where the $L_{i},L'_{i}$ are labels
and the $I_{i},U_{i}$ might be just truth values, but can also be more
complicated,
a (possibly infinite) sequence of some values. Connected objects have,
of course, to have corresponding such sequences. In addition, the object
$X$ has a criterion for each input, whether it is valid or not (in the
simple
case,
this will just be the truth value $'' true'' ).$ If there is at least one
positive
valid input $I_{i},$ and no valid negative input $U_{i},$ then the output
$C(X)$ and its
label are calculated on the basis of the valid inputs and their labels. If
the
object is an arrow, this will take some time, $t,$ otherwise, this is
instantaneous.

\ed

Evaluating a diagram

An evaluation is relative to a fixed input, i.e. some objects will be
given
certain values, and the diagram is left to calculate the others. It may
well
be that it oscillates, i.e. shows a cyclic behaviour. This may be true for
a
subset of the diagram, or the whole diagram. If it is restricted to an
unimportant part, we might neglect this. Whether it oscillates or not can
also depend on the time delays of the arrows (see Example \ref{Example 2.2}
(page \pageref{Example 2.2}) ).

We therefore define for a diagram $ \xbD $

$ \xba \xcn_{ \xbD } \xbb $ iff

(a) $ \xba $ is a (perhaps partial) input - where the other values are set
``not valid''

(b) $ \xbb $ is a (perhaps partial) output

(c) after some time, $ \xbb $ is stable, i.e. all still possible
oscillations
do not affect $ \xbb $

(d) the other possible input values do not matter, i.e. whatever the
input,
the result is the same.

In the cases examined here more closely, all input values will be defined.

 karl-search= End Reac-Sem-Def
\vspace{7mm}

 *************************************

\vspace{7mm}


\subsection{
A circuit semantics for simple IBRS without labels
}

\subsubsection{Reac-Sem-Circuit}

 {\LARGE karl-search= Start Reac-Sem-Circuit }

{\xssc LABEL: {Section Reac-Sem-Circuit}}
\label{Section Reac-Sem-Circuit}
\index{Section Reac-Sem-Circuit}

It is standard to implement the usual logical connectives by electronic
circuits. These components are called gates. Circuits with feedback
sometimes
show undesirable behaviour when the initial conditions are not specified.
(When
we switch a circuit on, the outputs of the individual gates can have
arbitrary
values.) The technical
realization of these initial values shows the way to treat defaults. The
initial
values are set via resistors (in the order of 1 $k \xbO )$ between the
point in the
circuit we want to intialize and the desired tension (say 0 Volt for
false,
5 Volt for true). They are called pull-down or pull-up resistors (for
default
0 or 5 Volt). When a ``real'' result comes in, it will override the tension
applied via the resistor.

Closer inspection reveals that we have here a 3 level default situation:
The initial value will be the weakest, which can be overridden by any
``real''
signal, but a positive argument can be overridden by a negative one. Thus,
the biggest resistor will be for the initialization, the smaller one for
the
supporting arguments, and the negative arguments have full power.

Technical details will be left to the experts.

We give now an example which shows that the delays of the arrows can
matter.
In one situation, a stable state is reached, in another, the circuit
begins to
oscillate.

\be

$\hspace{0.01em}$

(+++ Orig. No.:  Example 2.2 +++)

{\xssc LABEL: {Example 2.2}}
\label{Example 2.2}

(In engineering terms, this is a variant of a JK flip-flop with $R*S=0,$ a
circuit
with feedback.)

We have 8 measuring points.

$In1,In2$ are the overall input, $Out1,Out2$ the overall output,
$A1,A2,A3,A4$ are auxiliary internal points. All points can be true or
false.

The logical structure is as follows:

A1 $=$ $In1 \xcu Out1,$ A2 $=$ $In2 \xcu Out2,$

A3 $=$ $A1 \xco Out2,$ A4 $=$ $A2 \xco Out1,$

Out1 $=$ $ \xCN A3,$ Out2 $=$ $ \xCN A4.$

Thus, the circuit is symmetrical, with In1 corresponding to In2, A1 to A2,
A3 to A4, Out1 to Out2.

The input is held constant. See Diagram \ref{Diagram Gate-Sem} (page
\pageref{Diagram Gate-Sem}) .

\vspace{10mm}

\begin{diagram}

{\xssc LABEL: {Diagram Gate-Sem}}
\label{Diagram Gate-Sem-a}
\index{Diagram Gate-Sem}

\centering
\setlength{\unitlength}{1mm}
{\renewcommand{\dashlinestretch}{30}
\begin{picture}(150,170)(0,0)


\put(15,130){\arc{10}{-1.57}{1.57}}
\path(15,125)(15,135)
\put(16.3,129.3){\xssc{$\xcu$}}
\put(22,132){\xssc{$A1$}}
\path(13,132)(15,133)(13,134)
\path(13,126)(15,127)(13,128)

\put(45,133){\arc{10}{-1.57}{1.57}}
\path(45,128)(45,138)
\put(46.3,132.3){\xssc{$\xco$}}
\put(52,135){\xssc{$A3$}}
\path(43,135)(45,136)(43,137)
\path(43,129)(45,130)(43,131)

\path(75,128)(75,138)(83,133)(75,128)
\put(76.3,132.3){\xssc{$\xCN$}}
\path(73,132)(75,133)(73,134)

\path(0,127)(15,127)
\path(20,130)(45,130)
\path(50,133)(75,133)
\path(83,133)(108,133)
\path(106,132)(108,133)(106,134)
\put(93,133){\circle*{1}}
\put(101,133){\circle*{1}}

\put(0,124){\xssc{$In1$}}
\put(110,132){\xssc{$Out1$}}

\path(15,133)(5,133)(5,150)(101,150)(101,133)
\path(45,136)(35,136)(35,143)(85,143)(93,77)


\put(15,80){\arc{10}{-1.57}{1.57}}
\path(15,75)(15,85)
\put(16.3,79.3){\xssc{$\xcu$}}
\put(22,76){\xssc{$A2$}}
\path(13,82)(15,83)(13,84)
\path(13,76)(15,77)(13,78)

\put(45,77){\arc{10}{-1.57}{1.57}}
\path(45,72)(45,82)
\put(46.3,76.3){\xssc{$\xco$}}
\put(52,73){\xssc{$A4$}}
\path(43,79)(45,80)(43,81)
\path(43,73)(45,74)(43,75)

\path(75,72)(75,82)(83,77)(75,72)
\put(76.3,76.3){\xssc{$\xCN$}}
\path(73,78)(75,77)(73,76)

\path(0,83)(15,83)
\path(20,80)(45,80)
\path(50,77)(75,77)
\path(83,77)(108,77)
\path(106,76)(108,77)(106,78)
\put(93,77){\circle*{1}}
\put(101,77){\circle*{1}}

\put(0,85.5){\xssc{$In2$}}
\put(110,76){\xssc{$Out2$}}

\path(15,77)(5,77)(5,60)(101,60)(101,77)
\path(45,74)(35,74)(35,67)(85,67)(93,133)

\put(30,20) {{\rm\bf Gate Semantics}}

\end{picture}
}

\end{diagram}

\vspace{4mm}

\ee

We suppose that the output of the individual gates is present $n$ time
slices
after the input was present. $n$ will in the first circuit be equal to 1
for all
gates, in the second circuit equal to 1 for all but the AND gates, which
will
take 2 time slices. Thus, in both cases, e.g. Out1 at time $t$ will be the
negation of A3 at time $t-1.$ In the first case, A1 at time $t$ will be
the
conjunction of In1 and Out1 at time $t-1,$ and in the second case the
conjunction
of In1 and Out1 at time $t-2.$

We initialize In1 as true, all others as false. (The initial value of A3
and A4
does not matter, the behaviour is essentially the same for all such
values.)

The first circuit will oscillate with a period of 4, the second circuit
will go
to a stable state.

We have the following transition tables (time slice shown at left):

Circuit 1, $delay=1$ everywhere:

\begin{tabular}{lccccccccl}
   & In1 & In2 & A1 & A2 & A3 & A4 & Out1 & Out2 & \\
1: &   T &   F &  F &  F &  F &  F &    F &    F & \\
2: &   T &   F &  F &  F &  F &  F &    T &    T & \\
3: &   T &   F &  T &  F &  T &  T &    T &    T & \\
4: &   T &   F &  T &  F &  T &  T &    F &    F & \\
5: &   T &   F &  F &  F &  T &  F &    F &    F & oscillation starts \\
6: &   T &   F &  F &  F &  F &  F &    F &    T & \\
7: &   T &   F &  F &  F &  T &  F &    T &    T & \\
8: &   T &   F &  T &  F &  T &  T &    F &    T & \\
9: &   T &   F &  F &  F &  T &  F &    F &    F &
back to start of oscillation \\
\end{tabular}

Circuit 2, $delay=1$ everywhere, except for AND with $delay=2:$

(Thus, A1 and A2 are held at their intial value up to time 2, then they
are
calculated using the values of time $t-2.)$

\begin{tabular}{lccccccccl}
   & In1 & In2 & A1 & A2 & A3 & A4 & Out1 & Out2 & \\
1: &   T &   F &  F &  F &  F &  F &    F &    F & \\
2: &   T &   F &  F &  F &  F &  F &    T &    T & \\
3: &   T &   F &  F &  F &  T &  T &    T &    T & \\
4: &   T &   F &  T &  F &  T &  T &    F &    F & \\
5: &   T &   F &  T &  F &  T &  F &    F &    F & \\
6: &   T &   F &  F &  F &  T &  F &    F &    T & stable state reached \\
7: &   T &   F &  F &  F &  T &  F &    F &    T & \\
\end{tabular}

Note that state 6 of circuit 2 is also stable in circuit 1, but it is
never
reached in that circuit.

 karl-search= End Reac-Sem-Circuit
\vspace{7mm}

 *************************************

\vspace{7mm}

 karl-search= End Reac-Sem
\vspace{7mm}

 *************************************

\vspace{7mm}

 karl-search= End ToolBase1-IBRS
\vspace{7mm}

 *************************************

\vspace{7mm}

\newpage

\section{
Generalized preferential structures
}

\subsection{
}

\subsubsection{ToolBase1-HigherPref}

 {\LARGE karl-search= Start ToolBase1-HigherPref }

{\xssc LABEL: {Section Toolbase1-HigherPref}}
\label{Section Toolbase1-HigherPref}
\index{Section Toolbase1-HigherPref}
\subsubsection{Comment Gen-Pref}

 {\LARGE karl-search= Start Comment Gen-Pref }

\index{Comment Gen-Pref}

\bcom

$\hspace{0.01em}$

(+++ Orig. No.:  Comment Gen-Pref +++)

{\xssc LABEL: {Comment Gen-Pref}}
\label{Comment Gen-Pref}

A counterargument to $ \xba $ is NOT an argument for $ \xCN \xba $ (this
is asking for too
much), but just showing one case where $ \xCN \xba $ holds. In
preferential structures,
an argument for $ \xba $ is a set of level 1 arrows, eliminating $ \xCN
\xba -$models. A
counterargument is one level 2 arrow, attacking one such level 1 arrow.

Of course, when we have copies, we may need many successful attacks, on
all
copies, to achieve the goal. As we may have copies of level 1 arrows, we
may need many level 2 arrows to destroy them all.

 karl-search= End Comment Gen-Pref
\vspace{7mm}

 *************************************

\vspace{7mm}

\subsubsection{Definition Generalized preferential structure}

 {\LARGE karl-search= Start Definition Generalized preferential structure }

\index{Definition Generalized preferential structure}

\ecom

\bd

$\hspace{0.01em}$

(+++ Orig. No.:  Definition Generalized preferential structure +++)

{\xssc LABEL: {Definition Generalized preferential structure}}
\label{Definition Generalized preferential structure}

An IBR is called a generalized preferential structure iff the origins of
all
arrows are points. We will usually write $x,y$ etc. for points, $ \xba,$
$ \xbb $ etc. for
arrows.

 karl-search= End Definition Generalized preferential structure
\vspace{7mm}

 *************************************

\vspace{7mm}

\subsubsection{Definition Level-n-Arrow}

 {\LARGE karl-search= Start Definition Level-n-Arrow }

\index{Definition Level-n-Arrow}

\ed

\bd

$\hspace{0.01em}$

(+++ Orig. No.:  Definition Level-n-Arrow +++)

{\xssc LABEL: {Definition Level-n-Arrow}}
\label{Definition Level-n-Arrow}

Consider a generalized preferential structure $ \xdx.$

(1) Level $n$ arrow:

Definition by upward induction.

If $ \xba:x \xcp y,$ $x,y$ are points, then $ \xba $ is a level 1 arrow.

If $ \xba:x \xcp \xbb,$ $x$ is a point, $ \xbb $ a level $n$ arrow, then
$ \xba $ is a level $n+1$ arrow.
$(o( \xba )$ is the origin, $d( \xba )$ is the destination of $ \xba.)$

$ \xbl ( \xba )$ will denote the level of $ \xba.$

(2) Level $n$ structure:

$ \xdx $ is a level $n$ structure iff all arrows in $ \xdx $ are at most
level $n$ arrows.

We consider here only structures of some arbitrary but finite level $n.$

(3) We define for an arrow $ \xba $ by induction $O( \xba )$ and $D( \xba
).$

If $ \xbl ( \xba )=1,$ then $O( \xba ):=\{o( \xba )\},$ $D( \xba ):=\{d(
\xba )\}.$

If $ \xba:x \xcp \xbb,$ then $D( \xba ):=D( \xbb ),$ and $O( \xba
):=\{x\} \xcv O( \xbb ).$

Thus, for example, if $ \xba:x \xcp y,$ $ \xbb:z \xcp \xba,$ then $O(
\xbb ):=\{x,z\},$ $D( \xbb )=\{y\}.$

 karl-search= End Definition Level-n-Arrow
\vspace{7mm}

 *************************************

\vspace{7mm}

\subsubsection{Example Inf-Level}

 {\LARGE karl-search= Start Example Inf-Level }

\index{Example Inf-Level}

\ed

We will not consider here diagrams with arbitrarily high levels. One
reason is
that diagrams like the following will have an unclear meaning:

\be

$\hspace{0.01em}$

(+++ Orig. No.:  Example Inf-Level +++)

{\xssc LABEL: {Example Inf-Level}}
\label{Example Inf-Level}

$< \xba,1>:x \xcp y,$

$< \xba,n+1>:x \xcp < \xba,n>$ $(n \xbe \xbo ).$

Is $y \xbe \xbm (X)?$

 karl-search= End Example Inf-Level
\vspace{7mm}

 *************************************

\vspace{7mm}

\subsubsection{Definition Valid-Arrow}

 {\LARGE karl-search= Start Definition Valid-Arrow }

\index{Definition Valid-Arrow}

\ee

\bd

$\hspace{0.01em}$

(+++ Orig. No.:  Definition Valid-Arrow +++)

{\xssc LABEL: {Definition Valid-Arrow}}
\label{Definition Valid-Arrow}

Let $ \xdx $ be a generalized preferential structure of (finite) level
$n.$

We define (by downward induction):

(1) Valid $X-to-Y$ arrow:

Let $X,Y \xcc \xdP ( \xdx ).$

$ \xba \xbe \xdA ( \xdx )$ is a valid $X-to-Y$ arrow iff

(1.1) $O( \xba ) \xcc X,$ $D( \xba ) \xcc Y,$

(1.2) $ \xcA \xbb:x' \xcp \xba.(x' \xbe X$ $ \xch $ $ \xcE \xbg:x''
\xcp \xbb.( \xbg $ is a valid $X-to-Y$ arrow)).

We will also say that $ \xba $ is a valid arrow in $X,$ or just valid in
$X,$ iff $ \xba $ is a
valid $X-to-X$ arrow.

(2) Valid $X \xch Y$ arrow:

Let $X \xcc Y \xcc \xdP ( \xdx ).$

$ \xba \xbe \xdA ( \xdx )$ is a valid $X \xch Y$ arrow iff

(2.1) $o( \xba ) \xbe X,$ $O( \xba ) \xcc Y,$ $D( \xba ) \xcc Y,$

(2.2) $ \xcA \xbb:x' \xcp \xba.(x' \xbe Y$ $ \xch $ $ \xcE \xbg:x''
\xcp \xbb.( \xbg $ is a valid $X \xch Y$ arrow)).

(Note that in particular $o( \xbg ) \xbe X,$ and that $o( \xbb )$ need not
be in $X,$ but
can be in the bigger $Y.)$

 karl-search= End Definition Valid-Arrow
\vspace{7mm}

 *************************************

\vspace{7mm}

\subsubsection{Fact Higher-Validity}

 {\LARGE karl-search= Start Fact Higher-Validity }

\index{Fact Higher-Validity}

\ed

\bfa

$\hspace{0.01em}$

(+++ Orig. No.:  Fact Higher-Validity +++)

{\xssc LABEL: {Fact Higher-Validity}}
\label{Fact Higher-Validity}

(1) If $ \xba $ is a valid $X \xch Y$ arrow, then $ \xba $ is a valid
$Y-to-Y$ arrow.

(2) If $X \xcc X' \xcc Y' \xcc Y \xcc \xdP ( \xdx )$ and $ \xba \xbe \xdA
( \xdx )$ is a valid $X \xch Y$ arrow, and
$O( \xba ) \xcc Y',$ $D( \xba ) \xcc Y',$ then $ \xba $ is a valid $X'
\xch Y' $ arrow.

 karl-search= End Fact Higher-Validity
\vspace{7mm}

 *************************************

\vspace{7mm}

\subsubsection{Fact Higher-Validity Proof}

 {\LARGE karl-search= Start Fact Higher-Validity Proof }

\index{Fact Higher-Validity Proof}

\efa

\paragraph{
Proof Fact Higher-Validity
}

$\hspace{0.01em}$

(+++*** Orig.:  Proof Fact Higher-Validity )

{\xssc LABEL: {Section Proof Fact Higher-Validity}}
\label{Section Proof Fact Higher-Validity}

Let $ \xba $ be a valid $X \xch Y$ arrow. We show (1) and (2) together by
downward induction
(both are trivial).

By prerequisite $o( \xba ) \xbe X \xcc X',$ $O( \xba ) \xcc Y' \xcc Y,$
$D( \xba ) \xcc Y' \xcc Y.$

Case 1: $ \xbl ( \xba )=n.$ So $ \xba $ is a valid $X' \xch Y' $ arrow,
and a valid $Y-to-Y$ arrow.

Case 2: $ \xbl ( \xba )=n-1.$ So there is no $ \xbb:x' \xcp \xba,$ $y
\xbe Y,$ so $ \xba $ is a valid
$Y-to-Y$ arrow. By $Y' \xcc Y$ $ \xba $ is a valid $X' \xch Y' $ arrow.

Case 3: Let the result be shown down to $m,$ $n>m>1,$ let $ \xbl ( \xba
)=m-1.$
So $ \xcA \xbb:x' \xcp \xba (x' \xbe Y$ $ \xch $ $ \xcE \xbg:x'' \xcp
\xbb (x'' \xbe X$ and $ \xbg $ is a valid $X \xch Y$ arrow)).
By induction hypothesis $ \xbg $ is a valid $Y-to-Y$ arrow, and a valid
$X' \xch Y' $ arrow. So $ \xba $ is a valid $Y-to-Y$ arrow, and by $Y'
\xcc Y,$ $ \xba $ is a valid
$X' \xch Y' $ arrow.

$ \xcz $
\\[3ex]

 karl-search= End Fact Higher-Validity Proof
\vspace{7mm}

 *************************************

\vspace{7mm}

\subsubsection{Definition Higher-Mu}

 {\LARGE karl-search= Start Definition Higher-Mu }

\index{Definition Higher-Mu}

\bd

$\hspace{0.01em}$

(+++ Orig. No.:  Definition Higher-Mu +++)

{\xssc LABEL: {Definition Higher-Mu}}
\label{Definition Higher-Mu}

Let $ \xdx $ be a generalized preferential structure of level $n,$ $X \xcc
\xdP ( \xdx ).$

$ \xbm (X):=\{x \xbe X:$ $ \xcE <x,i>. \xCN \xcE $ valid $X-to-X$ arrow $
\xba:x' \xcp <x,i>\}.$

 karl-search= End Definition Higher-Mu
\vspace{7mm}

 *************************************

\vspace{7mm}

\subsubsection{Comment Smooth-Gen}

 {\LARGE karl-search= Start Comment Smooth-Gen }

\index{Comment Smooth-Gen}

\ed

\bcom

$\hspace{0.01em}$

(+++ Orig. No.:  Comment Smooth-Gen +++)

{\xssc LABEL: {Comment Smooth-Gen}}
\label{Comment Smooth-Gen}

The purpose of smoothness is to guarantee cumulativity. Smoothness
achieves Cumulativity by mirroring all information present in $X$ also in
$ \xbm (X).$
Closer inspection shows that smoothness does more than necessary. This is
visible when there are copies (or, equivalently, non-injective labelling
functions). Suppose we have two copies of $x \xbe X,$ $<x,i>$ and $<x,i'
>,$ and there
is $y \xbe X,$ $ \xba:<y,j> \xcp <x,i>,$ but there is no $ \xba ':<y'
,j' > \xcp <x,i' >,$ $y' \xbe X.$ Then
$ \xba:<y,j> \xcp <x,i>$
is irrelevant, as $x \xbe \xbm (X)$ anyhow. So mirroring $ \xba:<y,j>
\xcp <x,i>$ in $ \xbm (X)$ is not
necessary, i.e. it is not necessary to have some $ \xba ':<y',j' > \xcp
<x,i>,$ $y' \xbe \xbm (X).$

On the other hand, Example \ref{Example Need-Smooth} (page \pageref{Example
Need-Smooth})  shows that,
if we want smooth structures to correspond to the property $( \xbm CUM),$
we
need at least some valid arrows from $ \xbm (X)$ also for higher level
arrows.
This ``some'' is made precise (essentially) in Definition \ref{Definition
X-Sub-X'} (page \pageref{Definition X-Sub-X'}) .

From a more philosophical point of view,
when we see the (inverted) arrows of preferential structures as attacks on
non-minimal elements, then we should see smooth structures as always
having
attacks also from valid (minimal) elements. So, in general structures,
also
attacks from non-valid elements are valid, in smooth structures we always
also have attacks from valid elements.

The analogon to usual smooth structures, on level 2, is then that any
successfully attacked level 1 arrow is also attacked from a minimal point.

 karl-search= End Comment Smooth-Gen
\vspace{7mm}

 *************************************

\vspace{7mm}

\subsubsection{Definition X-Sub-X'}

 {\LARGE karl-search= Start Definition X-Sub-X' }

\index{Definition X-Sub-X'}

\ecom

\bd

$\hspace{0.01em}$

(+++ Orig. No.:  Definition X-Sub-X' +++)

{\xssc LABEL: {Definition X-Sub-X'}}
\label{Definition X-Sub-X'}

Let $ \xdx $ be a generalized preferential structure.

$X \xes X' $ iff

(1) $X \xcc X' \xcc \xdP ( \xdx ),$

(2) $ \xcA x \xbe X' -X$ $ \xcA <x,i>$ $ \xcE \xba:x' \xcp <x,i>( \xba $
is a valid $X \xch X' $ arrow),

(3) $ \xcA x \xbe X$ $ \xcE <x,i>$

$ \xDC $ $( \xcA \xba:x' \xcp <x,i>(x' \xbe X' $ $ \xch $ $ \xcE \xbb
:x'' \xcp \xba.( \xbb $ is a valid $X \xch X' $ arrow))).

Note that (3) is not simply the negation of (2):

Consider a level 1 structure. Thus all level 1 arrows are valid, but the
source
of the arrows must not be neglected.

(2) reads now: $ \xcA x \xbe X' -X$ $ \xcA <x,i>$ $ \xcE \xba:x' \xcp
<x,i>.x' \xbe X$

(3) reads: $ \xcA x \xbe X$ $ \xcE <x,i>$ $ \xCN \xcE \xba:x' \xcp
<x,i>.x' \xbe X' $

This is intended: intuitively, $X= \xbm (X' ),$ and minimal elements must
not be
attacked at all, but non-minimals must be attacked from $X$ - which is a
modified
version of smoothness.

 karl-search= End Definition X-Sub-X'
\vspace{7mm}

 *************************************

\vspace{7mm}

\subsubsection{Remark X-Sub-X'}

 {\LARGE karl-search= Start Remark X-Sub-X' }

\index{Remark X-Sub-X'}

\ed

\br

$\hspace{0.01em}$

(+++ Orig. No.:  Remark X-Sub-X' +++)

{\xssc LABEL: {Remark X-Sub-X'}}
\label{Remark X-Sub-X'}

We note the special case of Definition \ref{Definition X-Sub-X'} (page
\pageref{Definition X-Sub-X'})  for level
3 structures,
as it will be used later. We also write it immediately for the intended
case
$ \xbm (X) \xes X,$ and explicitly with copies.

$x \xbe \xbm (X)$ iff

(1) $ \xcE <x,i> \xcA < \xba,k>:<y,j> \xcp <x,i>$

$ \xDC (y \xbe X$ $ \xcp $ $ \xcE < \xbb ',l' >:<z',m' > \xcp < \xba
,k>.$

$ \xDC \xDC (z' \xbe \xbm (X)$ $ \xcu $ $ \xCN \xcE < \xbg ',n' >:<u',p'
> \xcp < \xbb ',l' >.u' \xbe X))$

See Diagram \ref{Diagram Essential-Smooth-3-1-2} (page \pageref{Diagram
Essential-Smooth-3-1-2}) .

$x \xbe X- \xbm (X)$ iff

(2) $ \xcA <x,i> \xcE < \xba ',k' >:<y',j' > \xcp <x,i>$

$ \xDC (y' \xbe \xbm (X)$ $ \xcu $

$ \xDC \xDC (a)$ $ \xCN \xcE < \xbb ',l' >:<z',m' > \xcp < \xba ',k'
>.z' \xbe X$

$ \xDC \xDC or$

$ \xDC \xDC (b)$ $ \xcA < \xbb ',l' >:<z',m' > \xcp < \xba ',k' >$

$ \xDC \xDC \xDC (z' \xbe X$ $ \xcp $ $ \xcE < \xbg ',n' >:<u',p' > \xcp
< \xbb ',l' >.u' \xbe \xbm (X))$ )

See Diagram \ref{Diagram Essential-Smooth-3-2} (page \pageref{Diagram
Essential-Smooth-3-2}) .

 karl-search= End Remark X-Sub-X'
\vspace{7mm}

 *************************************

\vspace{7mm}

\vspace{15mm}

\begin{diagram}

{\xssc LABEL: {Diagram Essential-Smooth-3-1-2}}
\label{Diagram Essential-Smooth-3-1-2-a}
\index{Diagram Essential-Smooth-3-1-2}

\centering
\setlength{\unitlength}{1mm}
{\renewcommand{\dashlinestretch}{30}
\begin{picture}(100,130)(0,0)
\put(50,80){\circle{80}}
\path(10,80)(90,80)

\put(100,100){{\xssc $X$}}
\put(100,60){{\xssc $\xbm (X)$}}

\path(50,100)(50,60)
\path(48.5,63)(50,60)(51.5,63)
\put(50,101){\circle*{0.3}}
\put(50,59){\circle*{0.3}}
\put(50,57){{\xssc $<x,i>$}}
\put(50,102){{\xssc $<y,j>$}}
\put(52,90){{\xssc $<\xba,k>$}}

\path(20,70)(48,70)
\path(46,71)(48,70)(46,69)
\put(19,70){\circle*{0.3}}
\put(13,67){{\xssc $<z',m'>$}}
\put(26,71){{\xssc $<\xbb',l'>$}}

\path(30,30)(30,68)
\path(29,66)(30,68)(31,66)
\put(30,29){\circle*{0.3}}
\put(31,53){{\xssc $<\xbg ',n'>$}}
\put(30,26.5){{\xssc $<u',p'>$}}

\put(30,10) {{\rm\bf Case 3-1-2}}

\end{picture}
}
\end{diagram}

\vspace{4mm}

\vspace{10mm}

\begin{diagram}

{\xssc LABEL: {Diagram Essential-Smooth-3-2}}
\label{Diagram Essential-Smooth-3-2-a}
\index{Diagram Essential-Smooth-3-2}

\centering
\setlength{\unitlength}{1mm}
{\renewcommand{\dashlinestretch}{30}
\begin{picture}(100,110)(0,0)
\put(50,60){\circle{80}}
\path(10,60)(90,60)

\put(100,80){{\xssc $X$}}
\put(100,40){{\xssc $\xbm (X)$}}

\path(50,80)(50,40)
\path(48.5,77)(50,80)(51.5,77)
\put(50,81){\circle*{0.3}}
\put(50,39){\circle*{0.3}}
\put(50,37){{\xssc $<y',j'>$}}
\put(50,82){{\xssc $<x,i>$}}
\put(52,70){{\xssc $<\xba',k'>$}}

\path(20,70)(48,70)
\path(46,71)(48,70)(46,69)
\put(19,70){\circle*{0.3}}
\put(13,67){{\xssc $<z',m'>$}}
\put(26,71){{\xssc $<\xbb',l'>$}}

\path(30,30)(30,68)
\path(29,66)(30,68)(31,66)
\put(30,29){\circle*{0.3}}
\put(31,53){{\xssc $<\xbg ',n'>$}}
\put(31,26.5){{\xssc $<u',p'>$}}

\put(30,0) {{\rm\bf Case 3-2}}

\end{picture}
}
\end{diagram}

\vspace{4mm}

\subsubsection{Fact X-Sub-X'}

 {\LARGE karl-search= Start Fact X-Sub-X' }

\index{Fact X-Sub-X'}

\er

\bfa

$\hspace{0.01em}$

(+++ Orig. No.:  Fact X-Sub-X' +++)

{\xssc LABEL: {Fact X-Sub-X'}}
\label{Fact X-Sub-X'}

(1) If $X \xes X',$ then $X= \xbm (X' ),$

(2) $X \xes X',$ $X \xcc X'' \xcc X' $ $ \xch $ $X \xes X''.$ (This
corresponds to $( \xbm CUM).)$

(3) $X \xes X',$ $X \xcc Y',$ $Y \xes Y',$ $Y \xcc X' $ $ \xch $ $X=Y.$
(This corresponds to $( \xbm \xcc \xcd ).)$

 karl-search= End Fact X-Sub-X'
\vspace{7mm}

 *************************************

\vspace{7mm}

\subsubsection{Fact X-Sub-X' Proof}

 {\LARGE karl-search= Start Fact X-Sub-X' Proof }

\index{Fact X-Sub-X' Proof}

\efa

\paragraph{
Proof Fact $X-Sub-X' $
}

$\hspace{0.01em}$

(+++*** Orig.:  Proof Fact X-Sub-X' )

{\xssc LABEL: {Section Proof Fact X-Sub-X'}}
\label{Section Proof Fact X-Sub-X'}

\subparagraph{
Proof
}

$\hspace{0.01em}$

(+++*** Orig.:  Proof )

(1) Trivial by Fact \ref{Fact Higher-Validity} (page \pageref{Fact
Higher-Validity})  (1).

(2)

We have to show

(a) $ \xcA x \xbe X'' -X$ $ \xcA <x,i>$ $ \xcE \xba:x' \xcp <x,i>( \xba $
is a valid $X \xch X'' $ arrow), and

(b) $ \xcA x \xbe X$ $ \xcE <x,i>$
$( \xcA \xba:x' \xcp <x,i>(x' \xbe X'' $ $ \xch $ $ \xcE \xbb:x'' \xcp
\xba.( \xbb $ is a valid $X \xch X'' $ arrow))).

Both follow from the corresponding condition for $X \xch X',$ the
restriction of the
universal quantifier, and Fact \ref{Fact Higher-Validity} (page \pageref{Fact
Higher-Validity})  (2).

(3)

Let $x \xbe X-$Y.

(a) By $x \xbe X \xes X',$ $ \xcE <x,i>$ s.t.
$( \xcA \xba:x' \xcp <x,i>(x' \xbe X' $ $ \xch $ $ \xcE \xbb:x'' \xcp
\xba.( \xbb $ is a valid $X \xch X' $ arrow))).

(b) By $x \xce Y \xes \xcE \xba_{1}:x' \xcp <x,i>$ $ \xba_{1}$ is a valid
$Y \xch Y' $ arrow, in
particular $x' \xbe Y \xcc X'.$ Moreover, $ \xbl ( \xba_{1})=1.$

So by (a) $ \xcE \xbb_{2}:x'' \xcp \xba_{1}.( \xbb_{2}$ is a valid $X \xch
X' $ arrow), in particular $x'' \xbe X \xcc Y',$
moreover $ \xbl ( \xbb_{2})=2.$

It follows by induction from the definition of valid $A \xch B$ arrows
that

$ \xcA n \xcE \xba_{2m+1},$ $ \xbl ( \xba_{2m+1})=2m+1,$ $ \xba_{2m+1}$ a
valid $Y \xch Y' $ arrow and

$ \xcA n \xcE \xbb_{2m+2},$ $ \xbl ( \xbb_{2m+2})=2m+2,$ $ \xbb_{2m+2}$ a
valid $X \xch X' $ arrow,

which is impossible, as $ \xdx $ is a structure of finite level.

$ \xcz $
\\[3ex]

 karl-search= End Fact X-Sub-X' Proof
\vspace{7mm}

 *************************************

\vspace{7mm}

\subsubsection{Definition Totally-Smooth}

 {\LARGE karl-search= Start Definition Totally-Smooth }

\index{Definition Totally-Smooth}

\bd

$\hspace{0.01em}$

(+++ Orig. No.:  Definition Totally-Smooth +++)

{\xssc LABEL: {Definition Totally-Smooth}}
\label{Definition Totally-Smooth}

Let $ \xdx $ be a generalized preferential structure, $X \xcc \xdP ( \xdx
).$

$ \xdx $ is called totally smooth for $X$ iff

(1) $ \xcA \xba:x \xcp y \xbe \xdA ( \xdx )(O( \xba ) \xcv D( \xba ) \xcc
X$ $ \xch $ $ \xcE \xba ':x' \xcp y.x' \xbe \xbm (X))$

(2) if $ \xba $ is valid, then there must also exist such $ \xba ' $ which
is valid.

(y a point or an arrow).

If $ \xdy \xcc \xdP ( \xdx ),$ then $ \xdx $ is called $ \xdy -$totally
smooth iff for all $X \xbe \xdy $
$ \xdx $ is totally smooth for $X.$

 karl-search= End Definition Totally-Smooth
\vspace{7mm}

 *************************************

\vspace{7mm}

\subsubsection{Example Totally-Smooth}

 {\LARGE karl-search= Start Example Totally-Smooth }

\index{Example Totally-Smooth}

\ed

\be

$\hspace{0.01em}$

(+++ Orig. No.:  Example Totally-Smooth +++)

{\xssc LABEL: {Example Totally-Smooth}}
\label{Example Totally-Smooth}

$X:=\{ \xba:a \xcp b,$ $ \xba ':b \xcp c,$ $ \xba '':a \xcp c,$ $ \xbb
:b \xcp \xba ' \}$ is not totally smooth,

$X:=\{ \xba:a \xcp b,$ $ \xba ':b \xcp c,$ $ \xba '':a \xcp c,$ $ \xbb
:b \xcp \xba ',$ $ \xbb ':a \xcp \xba ' \}$ is totally smooth.

 karl-search= End Example Totally-Smooth
\vspace{7mm}

 *************************************

\vspace{7mm}

\subsubsection{Example Need-Smooth}

 {\LARGE karl-search= Start Example Need-Smooth }

\index{Example Need-Smooth}

\ee

\be

$\hspace{0.01em}$

(+++ Orig. No.:  Example Need-Smooth +++)

{\xssc LABEL: {Example Need-Smooth}}
\label{Example Need-Smooth}

Consider $ \xba ':a \xcp b,$ $ \xba '':b \xcp c,$ $ \xba:a \xcp c,$ $
\xbb:a \xcp \xba.$

Then $ \xbm (\{a,b,c\})=\{a\},$ $ \xbm (\{a,c\})=\{a,c\}.$
Thus, $( \xbm CUM)$ does not hold in this structure.
Note that there is no valid arrow from $ \xbm (\{a,b,c\})$ to $c.$

 karl-search= End Example Need-Smooth
\vspace{7mm}

 *************************************

\vspace{7mm}

\subsubsection{Definition Essentially-Smooth}

 {\LARGE karl-search= Start Definition Essentially-Smooth }

\index{Definition Essentially-Smooth}

\ee

\bd

$\hspace{0.01em}$

(+++ Orig. No.:  Definition Essentially-Smooth +++)

{\xssc LABEL: {Definition Essentially-Smooth}}
\label{Definition Essentially-Smooth}

Let $ \xdx $ be a generalized preferential structure, $X \xcc \xdP ( \xdx
).$

$ \xdx $ is called essentially smooth for $X$ iff $ \xbm (X) \xes X.$

If $ \xdy \xcc \xdP ( \xdx ),$ then $ \xdx $ is called $ \xdy
-$essentially smooth iff for all $X \xbe \xdy $
$ \xbm (X) \xes X.$

 karl-search= End Definition Essentially-Smooth
\vspace{7mm}

 *************************************

\vspace{7mm}

\subsubsection{Example Total-vs-Essential}

 {\LARGE karl-search= Start Example Total-vs-Essential }

\index{Example Total-vs-Essential}

\ed

\be

$\hspace{0.01em}$

(+++ Orig. No.:  Example Total-vs-Essential +++)

{\xssc LABEL: {Example Total-vs-Essential}}
\label{Example Total-vs-Essential}

It is easy to see that we can distinguish total and essential smoothness
in richer structures, as the following Example shows:

We add an accessibility relation $R,$ and consider only those models which
are accessible.

Let e.g. $a \xcp b \xcp <c,0>,$ $<c,1>,$ without transitivity. Thus, only
$c$ has two
copies. This structure is essentially smooth, but of course not totally
so.

Let now mRa, mRb, $mR<c,0>,$ $mR<c,1>,$ $m' Ra,$ $m' Rb,$ $m' R<c,0>.$

Thus, seen from $m,$ $ \xbm (\{a,b,c\})=\{a,c\},$ but seen from $m',$ $
\xbm (\{a,b,c\})=\{a\},$
but $ \xbm (\{a,c\})=\{a,c\},$ contradicting (CUM).

$ \xcz $
\\[3ex]

 karl-search= End Example Total-vs-Essential
\vspace{7mm}

 *************************************

\vspace{7mm}

 karl-search= End ToolBase1-HigherPref
\vspace{7mm}

 *************************************

\vspace{7mm}

\section{
New remarks
}

\ee

 \xEh

 \xDH $eM \xdf_{1}$ entails: $ \xbm (Y) \xcc X \xcc Y \xch \xbm (X) \xcc
\xbm (Y)$

 \xDH Let $X \xcc Y \xcc Z.$

$eM \xdi_{1}:$ $X \xbe \xdi (Y) \xch X \xbe \xdi (Z)$

$eM \xdi_{2}:$ $X \xbe \xdm^{-}(Y) \xch X \xbe \xdm^{-}(Z)$

$eM \xdf_{1}:$ $X \xbe \xdf (Z) \xch X \xbe \xdf (Y)$

$eM \xdf_{2}:$ $X \xbe \xdm^{+}(Z) \xch X \xbe \xdm^{+}(Y)$

 \xDH We have: $eM \xdi_{1} \xcj eM \xdf_{2},$ $eM \xdi_{2} \xcj eM
\xdf_{1}.$

 \xDH We can represent the semantics for $n*s$ by reactive structures: the
choice of
one big subset disables the other choices.

 \xDH Some such structures can be represented by permutations (sometimes
not all)
of elements chosen.

 \xDH (iM) is done automatically, we have basically a preferential
structure,
but one which is switched on/off. The basic preferential idea is in
the fact that small sets are upward small. And Big sets are downward big,
this corresponds to Cautious Monotony.

 \xDH Was genau entspricht $(eM \xdi ),$ $(eM \xdf )$?

As the versions (1) suffice, we work with them only.

$(eM \xdi )$ without any domain prerequisites for $ \xcn:$

$(CUT' )$ $ \xba \xcn \xbb,$ $ \xba \xcl \xba ',$ $ \xba ' \xcu \xCN
\xbb \xcl \xba $ $ \xch $ $ \xba ' \xcn \xbb $

$(eM \xdi )$ with $( \xcv ):$

(wOR) $ \xba \xcn \xbb,$ $ \xba ' \xcl \xbb $ $ \xch $ $ \xba \xco \xba '
\xcn \xbb $

$(eM \xdf )$ without any domain prerequisites for $ \xcn:$

$(CM' )$ $ \xba \xcn \xbb,$ $ \xba ' \xcl \xba,$ $ \xba \xcu \xbb \xcl
\xba ' $ $ \xch $ $ \xba ' \xcn \xbb.$

 \xDH $(LLE)+(SC)+(RW)+(CUT' )+(CM' )$ characterize basic systems.

(iM) holds by (RW).

$(eM \xdi ):$ Let $A \xcc X \xcc Y,$ A be small in $X,$ $A:=M( \xba \xcu
\xCN \xbb )=M( \xba ' \xcu \xCN \xbb )$ $ \xch $
$ \xba ' \xcn \xbb $ $ \xch $ $M( \xba ' \xcu \xCN \xbb )$ is small in
$Y.$

$(eM \xdf )$ analogously.

 \xDH subideal situations and size: Optimum: smallest big subset,
subideal:
bigger big set, least ideal: all. This is ordered by logical strength.

 \xEj
\section{
Main table
}

$I$ should translate all logical rules into algebraic rules.

E.g.

$ \xba \xcn \xbb $ corresponds to $A-B \xbe \xdi (A)$

$ \xba \xcN \xbb $ to $A-B \xbe \xdm^{+}(A)$

$ \xba \xcl \xbb $ to $A-B= \xCQ $

And then give only the algebraic versions of the rules $(AND_{x}),$
$(OR_{x}),$ $(CM_{x}).$

Define $ \xdm^{+},$ $ \xdm^{-}.$
\newpage

\begin{turn}{90}

{\tiny

\begin{tabular}{|c|c@{:}c|c|c|c|c|c|c|}

\hline

\xEH
Ideal
\xEH
Filter
\xEH
$ \xdm^+ $
\xEH
$ \xeA $
\xEH
div. rules
\xEH
AND
\xEH
OR
\xEH
CM/Rat.Mon.
\xEP

\hline
\hline

\multicolumn{9}{|c|}{Optimal proportion} \xEP

\hline

(Opt)
\xEH
$ \xCQ \xbe \xdi (X)$
\xEH
$X \xbe \xdf (X)$
\xEH
\xEH
$ \xcA x \xbf \xcp \xeA x \xbf$
\xEH
$(SC)$
\xEH
\xEH
\xEH
\xEP

\xEH
\xEH
\xEH
\xEH
\xEH
$ \xba \xcl \xbb \xch \xba \xcn \xbb $
\xEH
\xEH
\xEH
\xEP

\hline
\hline

\multicolumn{9}{|c|}{Monotony: Improving proportions} \xEP

\hline

(iM)
\xEH
$A \xcc B \xbe \xdi (X)$ $ \xch $
\xEH
$A \xbe \xdf (X)$, $A \xcc B $
\xEH
\xEH
$\xeA x \xbf \xcu \xcA(\xbf \xcp \xbf')$
\xEH
$(RW)$
\xEH
\xEH
\xEH
\xEP

\xEH
$A \xbe \xdi (X)$
\xEH
$ \xch $ $B \xbe \xdf (X)$
\xEH
\xEH
$ \xcp $ $ \xeA x \xbf'$
\xEH
$ \xba \xcn \xbb, \xbb \xcl \xbb' \xch $
\xEH
\xEH
\xEH
\xEP

\xEH
\xEH
\xEH
\xEH
\xEH
$ \xba \xcn \xbb' $
\xEH
\xEH
\xEH
\xEP

\hline

$(eM \xdi )$
\xEH
Let $X \xcc Y$
\xEH
\xEH
\xEH
$\xeA x (\xbf: \xbq) \xcu$
\xEH
\xEH
\xEH
$(wOR)$
\xEH
\xEP

\xEH
(1) $\xdi (X) \xcc \xdi (Y)$
\xEH
\xEH
\xEH
$\xcA x (\xbf' \xcp \xbq) \xcp$
\xEH
\xEH
\xEH
$ \xba \xcn \xbb,$ $ \xba ' \xcl \xbb $ $ \xch $
\xEH
\xEP

\xEH
(2) $ \xdm^- (X) \xcc \xdm^- (Y)$
\xEH
\xEH
\xEH
$\xeA x (\xbf \xco \xbf': \xbq)$
\xEH
\xEH
\xEH
$ \xba \xco \xba ' \xcn \xbb $
\xEH
\xEP

\hline

$(eM \xdf )$
\xEH
\xEH
Let $X \xcc Y$
\xEH
\xEH
$\xeA x (\xbf: \xbq) \xcu$
\xEH
\xEH
\xEH
\xEH
$(wCM)$
\xEP

\xEH
\xEH
(1) $\xdf (Y) \xcs \xdp (X) \xcc \xdf (X)$
\xEH
\xEH
$\xcA x (\xbq \xcp \xbq') \xcp$
\xEH
\xEH
\xEH
\xEH
$ \xba \xcn \xbb,$ $ \xbb \xcl \xbb ' $ $ \xch $
\xEP

\xEH
\xEH
(2) $\xdm^+ (Y) \xcs \xdp (X) \xcc $
\xEH
\xEH
$\xeA x (\xbf \xcu \xbq': \xbq)$
\xEH
\xEH
\xEH
\xEH
$ \xba \xcu \xbb ' \xcn \xbb $
\xEP

\xEH
\xEH
$\xdm^+ (X)$
\xEH
\xEH
\xEH
\xEH
\xEH
\xEH
\xEP

\hline
\hline

\multicolumn{9}{|c|}{Keeping proportions} \xEP

\hline

$(\xCd)$
\xEH
$(\xdi \xcv disj)$
\xEH
$(\xdf \xcv disj)$
\xEH
\xEH
\xEH
\xEH
\xEH
$(disjOR)$
\xEH
\xEP

\xEH
$A \xbe \xdi (X),$ $B \xbe \xdi (Y),$
\xEH
$A \xbe \xdf (X),$ $B \xbe \xdf (Y),$
\xEH
\xEH
\xEH
\xEH
\xEH
$ \xbf \xcn \xbq,$ $ \xbf ' \xcn \xbq $
\xEH
\xEP

\xEH
$X \xcs Y= \xCQ $ $ \xch $
\xEH
$X \xcs Y= \xCQ $ $ \xch $
\xEH
\xEH
\xEH
\xEH
\xEH
$ \xbf \xcl \xCN \xbf ',$ $ \xch $
\xEH
\xEP

\xEH
$A \xcv B \xbe \xdi (X \xcv Y)$
\xEH
$A \xcv B \xbe \xdf (X \xcv Y)$
\xEH
\xEH
\xEH
\xEH
\xEH
$ \xbf \xco \xbf ' \xcn \xbq $
\xEH
\xEP

\hline
\hline

\multicolumn{9}{|c|}{Robustness of proportions: $n*small \xEd All$} \xEP

\hline

$(1*s)$
\xEH
$(\xdi_1)$
\xEH
$(\xdf_1)$
\xEH
\xEH
$(\xeA_1)$
\xEH
$(CP)$
\xEH
$(AND_{1})$
\xEH
\xEH
\xEP

\xEH
$X \xce \xdi (X)$
\xEH
$ \xCQ \xce \xdf (X)$
\xEH
\xEH
$ \xeA x \xbf \xcp \xcE x \xbf $
\xEH
$\xbf\xcn\xcT \xch \xbf\xcl\xcT$
\xEH
$ \xba \xcn \xbb $ $ \xch $ $ \xba \xcL \xCN \xbb $
\xEH
\xEH
\xEP

\hline

$(2*s)$
\xEH
$(\xdi_2)$
\xEH
$(\xdf_2)$
\xEH
\xEH
$(\xeA_2)$
\xEH
\xEH
$(AND_{2})$
\xEH
\xEH
$(CM_{2})$
\xEP

\xEH
$A,B \xbe \xdi (X) \xch $
\xEH
$A,B \xbe \xdf (X) \xch $
\xEH
\xEH
(1)
\xEH
\xEH
(1) $ \xba \xcn \xbb,$ $ \xba \xcn \xbb ' $ $ \xch $
\xEH
\xEH
$ \xba \xcn \xbb,$ $ \xba \xcn \xbb ' $ $ \xch $
\xEP

\xEH
$ A \xcv B \xEd X$
\xEH
$A \xcs B \xEd \xCQ $
\xEH
\xEH
$ \xeA x \xbf \xcu \xeA x \xbq $
\xEH
\xEH
$ \xba \xcL \xCN \xbb \xco \xCN \xbb ' $
\xEH
\xEH
$ \xba \xcu \xbb \xcL \xCN \xbb ' $
\xEP

\xEH
\xEH
\xEH
\xEH
$ \xcp $ $ \xcE x( \xbf \xcu \xbq )$
\xEH
\xEH
(2) $ \xba \xcn \xbb $ $ \xch $ $ \xba \xcN \xCN \xbb $
\xEH
\xEH
\xEP

\xEH
\xEH
\xEH
\xEH
(2)
\xEH
\xEH
\xEH
\xEH
\xEP

\xEH
\xEH
\xEH
\xEH
$\xeA x\xbf(x) \xcp \xCN\xeA x\xCN\xbf(x)$
\xEH
\xEH
\xEH
\xEH
\xEP

\hline

$(3*s)$
\xEH
$(\xdi_3)$
\xEH
$(\xdf_3)$
\xEH
$( \xdm^{+}_{3})$
\xEH
$(\xeA_3)$
\xEH
$( \xdl^{+}_{3})$
\xEH
$(AND_{3})$
\xEH
$(OR_{3})$
\xEH
$(CM_{3})$
\xEP

\xEH
$A,B,C \xbe \xdi (X) \xch $
\xEH
$A,B,C \xbe \xdf (X) \xch $
\xEH
$A \xbe \xdf (X),$ $X \xbe \xdf (Y)$
\xEH
$ \xeA x \xbf \xcu \xeA x \xbq \xcu \xeA x \xbs $
\xEH
$ \xbg \xcn \xbb,$ $ \xbg \xcu \xbb \xcn \xba $
\xEH
(1)
\xEH
$ \xba \xcn \xbb,$ $ \xba ' \xcn \xbb $ $ \xch $
\xEH
(1)
\xEP

\xEH
$ A \xcv B \xcv C \xEd X$
\xEH
$ A \xcs B \xcs C \xEd \xCQ $
\xEH
$Y \xbe \xdf (Z)$
\xEH
$ \xcp $
\xEH
$ \xch $ $ \xbg \xcN \xCN \xba $
\xEH
$ \xba \xcn \xbb,$ $ \xba \xcn \xbb ',$ $ \xba \xcn \xbb '' $
\xEH
$ \xba \xco \xba ' \xcN \xCN \xbb $
\xEH
$ \xba \xcn \xbb,$ $ \xba \xcn \xbb ',$ $ \xba \xcn \xbb '' $
\xEP

\xEH
\xEH
\xEH
 $ \xch $ $A \xbe \xdm^{+}(Z)$
\xEH
$ \xcE x( \xbf \xcu \xbq \xcu \xbs )$
\xEH
\xEH
 $ \xch $ $ \xba \xcL \xCN \xbb \xco \xCN \xbb ' \xco \xCN \xbb '' $
\xEH
\xEH
 $ \xch $ $ \xba \xcu \xbb \xcu \xbb ' \xcL \xCN \xbb '' $
\xEP

\xEH
\xEH
\xEH
\xEH
\xEH
\xEH
(2)
\xEH
\xEH
(2)
\xEP

\xEH
\xEH
\xEH
\xEH
\xEH
\xEH
$ \xba \xcn \xbb,$ $ \xba \xcn \xbb ' $ $ \xch $
\xEH
\xEH
$ \xba \xcn \xbb,$ $ \xba \xcn \xbb ' $ $ \xch $
\xEP

\xEH
\xEH
\xEH
\xEH
\xEH
\xEH
$ \xba \xcN \xCN \xbb \xco \xCN \xbb ' $
\xEH
\xEH
$ \xba \xcu \xbb \xcN \xCN \xbb ' $
\xEP

\hline

$(n*s)$
\xEH
$(\xdi_n)$
\xEH
$(\xdf_n)$
\xEH
$( \xdm^{+}_{n})$
\xEH
$(\xeA_n)$
\xEH
$( \xdl^{+}_{n})$
\xEH
$(AND_{n})$
\xEH
$(OR_{n})$
\xEH
$(CM_{n})$
\xEP

\xEH
$A_{1},.,A_{n} \xbe \xdi (X) $
\xEH
$A_{1},.,A_{n} \xbe \xdi (X) $
\xEH
$X_{1} \xbe \xdf (X_{2}),., $
\xEH
$ \xeA x \xbf_{1} \xcu.  \xcu \xeA x \xbf_{n} $
\xEH
$ \xba_{n} \xcn \xba_{n-1},$
\xEH
(1)
\xEH
$ \xba_{1} \xcn \xbb,., \xba_{n-1} \xcn \xbb $
\xEH
(1)
\xEP

\xEH
$ \xch $
\xEH
$ \xch $
\xEH
$ X_{n-1} \xbe \xdf (X_{n})$ $ \xch $
\xEH
$ \xcp $
\xEH
$ \xba_{n} \xcu \xba_{n-1} \xcn \xba_{n-2}$
\xEH
$ \xba \xcn \xbb_{1},., \xba \xcn \xbb_{n}$ $ \xch $
\xEH
$ \xch $
\xEH
$ \xba \xcn \xbb_{1},., \xba \xcn \xbb_{n}$ $ \xch $
\xEP

\xEH
$ A_{1} \xcv.  \xcv A_{n} \xEd X $
\xEH
$A_{1} \xcs.  \xcs A_{n} \xEd \xCQ$
\xEH
$X_{1} \xbe \xdm^{+}(X_{n})$
\xEH
$ \xcE x (\xbf_{1} \xcu.  \xcu \xbf_{n}) $
\xEH
$ \ldots $
\xEH
$ \xba \xcL \xCN \xbb_{1} \xco.  \xco \xCN \xbb_{n}$
\xEH
$ \xba_{1} \xco.  \xco \xba_{n-1} \xcN $
\xEH
$ \xba \xcu \xbb_{1} \xcu.  \xcu \xbb_{n-1} \xcL \xCN \xbb_{n}$
\xEP

\xEH
\xEH
\xEH
\xEH
\xEH
$ \xba_{n} \xcu.  \xcu \xba_{2} \xcn \xba_{1}$
\xEH
(2)
\xEH
$ \xCN \xbb $
\xEH
(2)
\xEP

\xEH
\xEH
\xEH
\xEH
\xEH
$ \xch $ $ \xba_{n} \xcN \xCN \xba_{1}$
\xEH
$ \xba \xcn \xbb_{1},., \xba \xcn \xbb_{n-1}$ $ \xch $
\xEH
\xEH
$ \xba \xcn \xbb_{1},., \xba \xcn \xbb_{n-1}$ $ \xch $
\xEP

\xEH
\xEH
\xEH
\xEH
\xEH
\xEH
$ \xba \xcN \xCN \xbb_{1} \xco.  \xco \xCN \xbb_{n-1}$
\xEH
\xEH
$ \xba \xcu \xbb_{1} \xcu.  \xcu \xbb_{n-2} \xcN $
\xEP

\xEH
\xEH
\xEH
\xEH
\xEH
\xEH
\xEH
\xEH
$ \xCN \xbb_{n-1}$
\xEP

\hline

$(< \xbo*s)$
\xEH
$(\xdi_\xbo)$
\xEH
$(\xdf_\xbo)$
\xEH
$( \xdm^{+}_{ \xbo })$
\xEH
$(\xeA_{\xbo})$
\xEH
$( \xdl^{+}_{ \xbo })$
\xEH
$(AND_{ \xbo })$
\xEH
$(OR_{ \xbo })$
\xEH
$(CM_{ \xbo })$
\xEP

\xEH
$A,B \xbe \xdi (X) \xch $
\xEH
$A,B \xbe \xdf (X) \xch $
\xEH
(1)
\xEH
$ \xeA x \xbf \xcu \xeA x \xbq \xcp $
\xEH
(1)
\xEH
$ \xba \xcn \xbb,$ $ \xba \xcn \xbb ' $ $ \xch $
\xEH
$ \xba \xcn \xbb,$ $ \xba ' \xcn \xbb $ $ \xch $
\xEH
$ \xba \xcn \xbb,$ $ \xba \xcn \xbb ' $ $ \xch $
\xEP

\xEH
$ A \xcv B \xbe \xdi (X)$
\xEH
$ A \xcs B \xbe \xdf (X)$
\xEH
$A \xbe \xdf (X),$ $X \xbe \xdm^{+}(Y)$
\xEH
$ \xeA x( \xbf \xcu \xbq )$
\xEH
$ \xbg \xcN \xCN \xbb,$ $ \xbg \xcu \xbb \xcn \xba $
\xEH
$ \xba \xcn \xbb \xcu \xbb ' $
\xEH
$ \xba \xco \xba ' \xcn \xbb $
\xEH
$ \xba \xcu \xbb \xcn \xbb ' $
\xEP

\xEH
\xEH
\xEH
$ \xch $ $A \xbe \xdm^{+}(Y)$
\xEH
\xEH
$ \xch $ $ \xbg \xcN \xCN \xba $
\xEH
\xEH
\xEH
\xEP

\xEH
\xEH
\xEH
(2)
\xEH
\xEH
(2)
\xEH
\xEH
\xEH
\xEP

\xEH
\xEH
\xEH
$A \xbe \xdm^{+}(X),$ $X \xbe \xdf (Y)$
\xEH
\xEH
$ \xbg \xcn \xbb,$ $ \xbg \xcu \xbb \xcN \xCN \xba $
\xEH
\xEH
\xEH
\xEP

\xEH
\xEH
\xEH
$ \xch $ $A \xbe \xdm^{+}(Y)$
\xEH
\xEH
$ \xch $ $ \xbg \xcN \xCN \xba $
\xEH
\xEH
\xEH
\xEP

\xEH
\xEH
\xEH
(3)
\xEH
\xEH
(3)
\xEH
\xEH
\xEH
\xEP

\xEH
\xEH
\xEH
$A \xbe \xdf (X),$ $X \xbe \xdf (Y)$
\xEH
\xEH
$ \xbg \xcu \xbb \xcn \xba,$ $ \xbg \xcn \xbb $
\xEH
\xEH
\xEH
\xEP

\xEH
\xEH
\xEH
$ \xch $ $A \xbe \xdf (Y)$
\xEH
\xEH
$ \xch $ $ \xbg \xcn \xba $
\xEH
\xEH
\xEH
\xEP

\xEH
\xEH
\xEH
(4)
\xEH
\xEH
\xEH
\xEH
\xEH
\xEP

\xEH
\xEH
\xEH
$A,B \xbe \xdi (X)$ $ \xch $
\xEH
\xEH
\xEH
\xEH
\xEH
\xEP

\xEH
\xEH
\xEH
$A-B \xbe \xdi (X-$B)
\xEH
\xEH
\xEH
\xEH
\xEH
\xEP

\xEH
\xEH
\xEH
(5)
\xEH
\xEH
\xEH
\xEH
\xEH
\xEP

\xEH
\xEH
\xEH
$A \xbe \xdf (X),$ $B \xbe \xdi (X)$
\xEH
\xEH
\xEH
\xEH
\xEH
\xEP

\xEH
\xEH
\xEH
$ \xch $ $A-B \xbe \xdf (X-B)$
\xEH
\xEH
\xEH
\xEH
\xEH
\xEP

\hline
\hline

\multicolumn{9}{|c|}{Robustness of $\xdm^+$} \xEP

\hline

$(\xdm^{++})$
\xEH
\xEH
\xEH
$(\xdm^{++})$
\xEH
\xEH
analogue $( \xdl^{+}_{ \xbo })$
\xEH
\xEH
\xEH
$(RatM)$
\xEP

\xEH
\xEH
\xEH
(1)
\xEH
\xEH
\xEH
\xEH
\xEH
$ \xbf \xcn \xbq,  \xbf \xcN \xCN \xbq '   \xch $
\xEP

\xEH
\xEH
\xEH
$A \xbe \xdi (X),$ $B \xce \xdf (X)$
\xEH
\xEH
\xEH
\xEH
\xEH
$ \xbf \xcu \xbq ' \xcn \xbq $
\xEP

\xEH
\xEH
\xEH
$ \xch $ $A-B \xbe \xdi (X-B)$
\xEH
\xEH
\xEH
\xEH
\xEH
\xEP

\xEH
\xEH
\xEH
(2)
\xEH
\xEH
\xEH
\xEH
\xEH
\xEP

\xEH
\xEH
\xEH
$A \xbe \xdf (X),$ $B \xce \xdf (X)$
\xEH
\xEH
\xEH
\xEH
\xEH
\xEP

\xEH
\xEH
\xEH
$ \xch $ $A-B \xbe \xdf (X-$B)
\xEH
\xEH
\xEH
\xEH
\xEH
\xEP

\xEH
\xEH
\xEH
(3)
\xEH
\xEH
\xEH
\xEH
\xEH
\xEP

\xEH
\xEH
\xEH
$A \xbe \xdm^+ (X),$
\xEH
\xEH
\xEH
\xEH
\xEH
\xEP

\xEH
\xEH
\xEH
$X \xbe \xdm^+ (Y)$
\xEH
\xEH
\xEH
\xEH
\xEH
\xEP

\xEH
\xEH
\xEH
$ \xch $ $A \xbe \xdm^+ (Y)$
\xEH
\xEH
\xEH
\xEH
\xEH
\xEP

\hline

\hline

\end{tabular}

}

\end{turn}

\section{
Comments
}

The usual rules $ \xCf (AND)$ etc. are named here $(AND_{ \xbo }),$ as
they are in a
natural ascending line of similar rules, based on strengthening of the
filter/ideal properties.

$ \xba \xcn \xbb $ $: \xcj $ $M( \xba \xcu \xbb ) \xbe \xdf (M( \xba ))$ $
\xcj $ $M( \xba \xcu \xCN \xbb ) \xbe \xdi (M( \xba )).$

Thus $ \xba \xcN \xbb $ $ \xcj $ $M( \xba \xcu \xCN \xbb ) \xbe
\xdm^{+}(M( \xba )).$

$ \xba \xcl \xbb $ $ \xcj $ $M( \xba \xcu \xbb )=M( \xba ).$
\subsection{
Regularities
}

The rules are divided into 5 groups:

 \xEh

 \xDH $ \xCf (Opt),$ which says that All is optimal - i.e. when there are
no
exceptions, then a rule holds.

 \xDH 3 monotony rules:

 \xEh
 \xDH $ \xCf (iM)$ is inner monotony, a subset of a small set is small
 \xDH $(eM \xdi )$ external monotony for ideals: enlarging the base set
keeps small
sets small
 \xDH $(em \xdf )$ external monotony for filters: a big subset stays big
when the base
set shrinks.
 \xEj

These three rules are very natural if ``size'' is anything coherent over
change
of base sets. In particular, they can be seen as weakening.

 \xDH $( \xCd )$ keeps proportions, it is here mainly to point the
possibility out.

 \xDH a group of rules $x*s,$ which say how many small sets will not yet
add to
the base set.
 \xDH Rational monotony, which can best be understood as robustness of $
\xdm^{+}$ -
where $ \xdm^{+}$ is the set of subsets, which are not small, i.e. big or
medium
size.

 \xEj

There are more regularities in the table:

Starting at $3*s,$ the properties can
be expressed nicely by ever stronger conditions $ \xdm^{+}.$

The conditions $ \xdi_{x}$ (or, equivalently, $ \xdf_{x})$ correspond
directly to the
conditions $(AND)_{x}.$ The other logical and algebraic conditions in the
same line
can be obtained using the weakening rules (monotony). Thus, $(AND)_{x},$
somewhat
surpisingly, reveals itself as, in this sense, the strongest rule of the
line $x.$

See  \Xl. below.

Thus, we can summarize:

We can obtain all rules except $ \xCf (RatM)$ from $ \xCf (Opt),$ the
monotony rules, and
$x*s.$
\subsection{
The position of RatMon:
}

RatM does not fit intoadding small sets. We have exhausted the combination
of small sets by $(< \xbo *s),$ unless we go to languages with infinitary
formulas.

The next idea would be to add medium size sets. But, by definition,
$2*medium$
can be all. Adding small and medium sets would not help either: Suppose we
have a rule $medium+n*small \xEd all.$ Taking the complement of the first
medium
set, which is again medium, we have the rule $2*n*small \xEd all.$ So we
do not
see any meaningful new internal rule. i.e. without changing the base set.
\section{
Coherent systems
}

\bd

$\hspace{0.01em}$

(+++ Orig. No.:  Definition CoherentSystem +++)

{\xssc LABEL: {Definition CoherentSystem}}
\label{Definition CoherentSystem}

A coherent system of sizes $ \xdc \xds $ consists of a universe $U,$ $
\xCQ \xce \xdy \xcc \xdp (U),$ and
for all $X \xbe \xdy $ $ \xdi (X)$ (dually $ \xdf (X)).$

We say that $ \xdc \xds $ satisfies a certain property iff all $X,Y \xbe
\xdy $ satisfy this
property.

$ \xdc \xds $ is called basic or level 1 iff it satisfies (iM), $(eM \xdi
),$ $(eM \xdf ),$ $(1*s).$

$ \xdc \xds $ is level $n$ iff it satisfies (iM), $(eM \xdi ),$ $(eM \xdf
),$ $(n*s).$

\ed

\bfa

$\hspace{0.01em}$

(+++ Orig. No.:  Fact Not-2*s +++)

{\xssc LABEL: {Fact Not-2*s}}
\label{Fact Not-2*s}

Let a $ \xdc \xds $ be given s.t. $ \xdy = \xdp (U).$ If $X \xbe \xdy $
satisfies $ \xdm^{++},$ but not
$(< \xbo *s),$ then there is $Y \xbe \xdy $ which does not satisfy
$(2*s).$

\efa

\subparagraph{
Proof
}

$\hspace{0.01em}$

(+++*** Orig.:  Proof )

As $X$ does not satisfy $(< \xbo *s),$ there are small $A,B \xcc X$ s.t.
$A \xcv B \xbe \xdm^{+}.$
Consider now $A \xcv B$ as base set $Y.$ By $( \xdm^{++})$ for $X,$ $A,B
\xce \xdm^{+}(A \xcv B),$ so
$A,B \xbe \xdi (A \xcv B),$ so $(2*s)$ does not hold for $A \xcv B.$ $
\xcz $
\\[3ex]

\bfa

$\hspace{0.01em}$

(+++ Orig. No.:  Fact Independence-eM +++)

{\xssc LABEL: {Fact Independence-eM}}
\label{Fact Independence-eM}

(1) $(eM \xdi )$ and $(eM \xdf )$ (1) are formally independent, though
intuitively
equivalent.

(2) $(eM \xdf )$ $(1)+(2)$ $ \xch $ $(eM \xdi ).$

\efa

\subparagraph{
Proof
}

$\hspace{0.01em}$

(+++*** Orig.:  Proof )

(1) Let $U:=\{x,y,z\},$ $X:=\{x,z\},$ $ \xdy:=\{U,X\}.$

(1.1) Let $ \xdf (U):=\{A \xcc U:z \xbe A\},$ $ \xdf (X):=\{X\}.$
$(eM \xdi )$ holds for $X$ and $U,$ but $\{z\} \xbe \xdf (U),$ $\{z\} \xcc
X,$ $\{z\} \xce \xdf (X),$ so $(eM \xdf )$ fails.

(1.2) Let $ \xdf (U):=\{U\},$ $ \xdf (X):=\{A \xcc X:z \xbe A\}.$
$(eM \xdf )$ holds trivially, but $(eM \xdi )$ fails, as $\{x\} \xbe \xdi
(X),$ but $\{x\} \xce \xdi (U).$

(2) Let $A \xcc X \xcc Y$ $A \xbe \xdi (X).$ If $A \xbe \xdm (Y),$ then $A
\xbe \xdm (X),$ likewise if $A \xbe \xdf (Y),$
so $A \xbe \xdi (Y).$

$ \xcz $
\\[3ex]

\bfa

$\hspace{0.01em}$

(+++ Orig. No.:  Fact Level-n-n+1 +++)

{\xssc LABEL: {Fact Level-n-n+1}}
\label{Fact Level-n-n+1}

A level $n$ system is strictly weaker than a level $n+1$ system.

\efa

\subparagraph{
Proof
}

$\hspace{0.01em}$

(+++*** Orig.:  Proof )

Consider $U:=\{1, \Xl,n+1\},$ $ \xdy:= \xdp (U)-\{ \xCQ \}.$ Let $ \xdi
(U):=\{ \xCQ \} \xcv \{\{x\}:x \xbe U\},$
$ \xdi (X):=\{ \xCQ \}$ for $X \xEd U.$
(iM), $(eM \xdi ),$ $(eM \xdf )$ hold trivially, so does $(1*s).$
$(n*s)$ holds trivially for $X \xEd U,$ but also for $U.$ $((n*1)*s)$ does
not hold
for $U.$ $ \xcz $
\\[3ex]
\section{
Ideals, filters, and logical rules
}

$( \xdi_{n})$ says $A_{1}, \Xl,A_{n} \xbe \xdi (X)$ $ \xch $ $A_{1} \xcv
\Xl  \xcv A_{n} \xEd X.$

The proofs for $(< \xbo *s)$ are analogous.
\subsection{
There are infinitely many new rules
}

Note that our schemata allow us to generate infintely many new rules, here
is
an example:

Start with A, add $s_{1,1},$ $s_{1,2}$ two sets small in $A \xcv s_{1,1}$
$(A \xcv s_{1,2}$ respectively).
Consider now $A \xcv s_{1,1} \xcv s_{1,2}$ and $s_{2}$ s.t. $s_{2}$ is
small in $A \xcv s_{1,1} \xcv s_{1,2} \xcv s_{2}.$
Continue with $s_{3,1},$ $s_{3,2}$ small in $A \xcv s_{1,1} \xcv s_{1,2}
\xcv s_{2} \xcv s_{3,1}$ etc.

Without additional properties, this system creates a new rule, which is
not equivalent to any usual rules.
\newpage
\section{
Facts about $\xdm$
}
\subsubsection{Fact R-down-neu}

 {\LARGE karl-search= Start Fact R-down-neu }

\index{Fact R-down-neu}

\bfa

$\hspace{0.01em}$

(+++ Orig. No.:  Fact R-down-neu +++)

{\xssc LABEL: {Fact R-down-neu}}
\label{Fact R-down-neu}

$( \xdm^{+}_{ \xbo })$ (4) and (5) and the three versions of $(
\xdm^{++})$ are each equivalent.

For the third version of $( \xdm^{++})$ we use $(eM \xdi )$ and $(eM \xdf
).$

 karl-search= End Fact R-down-neu
\vspace{7mm}

 *************************************

\vspace{7mm}

\subsubsection{Fact R-down-neu Proof}

 {\LARGE karl-search= Start Fact R-down-neu Proof }

\index{Fact R-down-neu Proof}

\efa

\subparagraph{
Proof
}

$\hspace{0.01em}$

(+++*** Orig.:  Proof )

For $A,B \xcc X,$ $(X-B)-((X-A)-B)=A-$B.

`` $ \xch $ '': Let $A \xbe \xdf (X),$ $B \xbe \xdi (X),$ so $X-A \xbe \xdi
(X),$ so by prerequisite
$(X-A)-B \xbe \xdi (X-$B), so $A-B=(X-B)-((X-A)-B) \xbe \xdf (X-$B).

`` $ \xci $ '': Let $A,B \xbe \xdi (X),$ so $X-A \xbe \xdf (X),$ so by
prerequisite $(X-A)-B \xbe \xdf (X-$B),
so $A-B=(X-B)-((X-A)-B) \xbe \xdi (X-$B).

The proof for $( \xdm^{++})$ is the same for the first two cases.

It remains to show equivalence with the last one. We assume closure under
set difference and union.

$(1) \xch (3):$

Suppose $A \xce \xdm^{+}(Y),$ but $X \xbe \xdm^{+}(Y),$ we show $A \xce
\xdm^{+}(X).$ So $A \xbe \xdi (Y),$ $Y-X \xce \xdf (Y),$
so $A=A-(Y-X) \xbe \xdi (Y-(Y-X))= \xdi (X).$

$(3) \xch (1):$

Suppose $A-B \xce \xdi (X-$B), $B \xce \xdf (X),$ we show $A \xce \xdi
(X).$ By prerequisite $A-B \xbe \xdm^{+}(X-$B),
$X-B \xbe \xdm^{+}(X),$ so $A-B \xbe \xdm^{+}(X),$ so by $(eM \xdi )$ and
$(eM \xdf )$ $A \xbe \xdm^{+}(X),$ so $A \xce \xdi (X).$

$ \xcz $
\\[3ex]

 karl-search= End Fact R-down-neu Proof
\vspace{7mm}

 *************************************

\vspace{7mm}

\bfa

$\hspace{0.01em}$

(+++ Orig. No.:  Fact 3->1-3->2 +++)

{\xssc LABEL: {Fact 3->1-3->2}}
\label{Fact 3->1-3->2}

(3) $ \xch $ (1) of $( \xdm^{+}_{ \xbo }),$ (3) $ \xch $ (2) of $(
\xdm^{+}_{ \xbo }).$

\efa

\subparagraph{
Proof
}

$\hspace{0.01em}$

(+++*** Orig.:  Proof )

$A \xbe \xdi (Y)$ $ \xch $ $X=(X-A) \xcv A \xbe \xdi (Y).$ The other
implication is analogous. $ \xcz $
\\[3ex]
\section{
Equivalences between size and logic
}
\subsubsection{Proposition Ref-Class-Mu-neu}

 {\LARGE karl-search= Start Proposition Ref-Class-Mu-neu }

\index{Proposition Ref-Class-Mu-neu}

\bp

$\hspace{0.01em}$

(+++ Orig. No.:  Proposition Ref-Class-Mu-neu +++)

{\xssc LABEL: {Proposition Ref-Class-Mu-neu}}
\label{Proposition Ref-Class-Mu-neu}

If $f(X)$ is the smallest $ \xCf A$ s.t. $A \xbe \xdf (X),$ then, given
the property on the
left, the one on the right follows.

Conversely, when we define $ \xdf (X):=\{X':f(X) \xcc X' \xcc X\},$ given
the property on
the right, the one on the left follows. For this direction, we assume
that we can use the full powerset of some base set $U$ - as is the case
for
the model sets of a finite language. This is perhaps not too bold, as
we mainly want to stress here the intuitive connections, without putting
too much weight on definability questions.

We assume $ \xCf (iM)$ to hold.

{\footnotesize

\begin{tabular}{|c|c|c|c|}

\hline

(1.1)
\xEH
$(eM\xdi )$
\xEH
$ \xch $
\xEH
$( \xbm wOR)$
\xEP

\cline{1-1}
\cline{3-3}

(1.2)
\xEH
\xEH
$ \xci $
\xEH
\xEP

\hline

(2.1)
\xEH
$(eM\xdi )+(I_\xbo )$
\xEH
$ \xch $
\xEH
$( \xbm OR)$
\xEP

\cline{1-1}
\cline{3-3}

(2.2)
\xEH
\xEH
$ \xci $
\xEH
\xEP

\hline

(3.1)
\xEH
$(eM\xdi )+(I_\xbo )$
\xEH
$ \xch $
\xEH
$( \xbm PR)$
\xEP

\cline{1-1}
\cline{3-3}

(3.2)
\xEH
\xEH
$ \xci $
\xEH
\xEP

\hline

(4.1)
\xEH
$(I \xcv disj )$
\xEH
$ \xch $
\xEH
$( \xbm disjOR)$
\xEP

\cline{1-1}
\cline{3-3}

(4.2)
\xEH
\xEH
$ \xci $
\xEH
\xEP

\hline

(5.1)
\xEH
$(\xdm^+_\xbo)$
\xEH
$ \xch $
\xEH
$( \xbm CM)$
\xEP

\cline{1-1}
\cline{3-3}

(5.2)
\xEH
\xEH
$ \xci $
\xEH
\xEP

\hline

(6.1)
\xEH
$(\xdm^{++})$
\xEH
$ \xch $
\xEH
$( \xbm RatM)$
\xEP

\cline{1-1}
\cline{3-3}

(6.2)
\xEH
\xEH
$ \xci $
\xEH
\xEP

\hline

(7.1)
\xEH
$(I_\xbo )$
\xEH
$ \xch $
\xEH
$( \xbm AND)$
\xEP

\cline{1-1}
\cline{3-3}

(7.2)
\xEH
\xEH
$ \xci $
\xEH
\xEP

\hline

\end{tabular}

}

 karl-search= End Proposition Ref-Class-Mu-neu
\vspace{7mm}

 *************************************

\vspace{7mm}

\subsubsection{Proposition Ref-Class-Mu-neu Proof}

 {\LARGE karl-search= Start Proposition Ref-Class-Mu-neu Proof }

\index{Proposition Ref-Class-Mu-neu Proof}

\ep

\subparagraph{
Proof
}

$\hspace{0.01em}$

(+++*** Orig.:  Proof )

(1.1) $(eM \xdi )$ $ \xch $ $( \xbm wOR):$

$X-f(X)$ is small in $X,$ so it is small in $X \xcv Y$ by $(eM \xdi ),$ so
$A:=X \xcv Y-(X-f(X)) \xbe \xdf (X \xcv Y),$ but $A \xcc f(X) \xcv Y,$ and
$f(X \xcv Y)$ is the smallest element
of $ \xdf (X \xcv Y),$ so $f(X \xcv Y) \xcc A \xcc f(X) \xcv Y.$

(1.2) $( \xbm wOR)$ $ \xch $ $(eM \xdi ):$

Let $X \xcc Y,$ $X':=Y-$X. Let $A \xbe \xdi (X),$ so $X-A \xbe \xdf (X),$
so $f(X) \xcc X-$A, so
$f(X \xcv X' ) \xcc f(X) \xcv X' \xcc (X-A) \xcv X' $ by prerequisite, so
$(X \xcv X' )-((X-A) \xcv X' )=A \xbe \xdi (X \xcv X' ).$

(2.1) $(eM \xdi )+(I_{ \xbo })$ $ \xch $ $( \xbm OR):$

$X-f(X)$ is small in $X,$ $Y-f(Y)$ is small in $Y,$ so both are small in
$X \xcv Y$ by
$(eM \xdi ),$ so $A:=(X-f(X)) \xcv (Y-f(Y))$ is small in $X \xcv Y$ by
$(I_{ \xbo }),$ but
$X \xcv Y-(f(X) \xcv f(Y)) \xcc A,$ so $f(X) \xcv f(Y) \xbe \xdf (X \xcv
Y),$ so, as $f(X \xcv Y)$ is the smallest
element of $ \xdf (X \xcv Y),$ $f(X \xcv Y) \xcc f(X) \xcv f(Y).$

(2.2) $( \xbm OR)$ $ \xch $ $(eM \xdi )+(I_{ \xbo }):$

Let again $X \xcc Y,$ $X':=Y-$X. Let $A \xbe \xdi (X),$ so $X-A \xbe \xdf
(X),$ so $f(X) \xcc X-$A. $f(X' ) \xcc X',$
so $f(X \xcv X' ) \xcc f(X) \xcv f(X' ) \xcc (X-A) \xcv X' $ by
prerequisite, so
$(X \xcv X' )-((X-A) \xcv X' )=A \xbe \xdi (X \xcv X' ).$

$(I_{ \xbo })$ holds by definition.

(3.1) $(eM \xdi )+(I_{ \xbo })$ $ \xch $ $( \xbm PR):$

Let $X \xcc Y.$ $Y-f(Y)$ is the largest element of $ \xdi (Y),$ $X-f(X)
\xbe \xdi (X) \xcc \xdi (Y)$ by
$(eM \xdi ),$ so $(X-f(X)) \xcv (Y-f(Y)) \xbe \xdi (Y)$ by $(I_{ \xbo }),$
so by ``largest'' $X-f(X) \xcc Y-f(Y),$
so $f(Y) \xcs X \xcc f(X).$

(3.2) $( \xbm PR)$ $ \xch $ $(eM \xdi )+(I_{ \xbo })$

Let again $X \xcc Y,$ $X':=Y-$X. Let $A \xbe \xdi (X),$ so $X-A \xbe \xdf
(X),$ so $f(X) \xcc X-$A, so
by prerequisite $f(Y) \xcs X \xcc X-$A, so $f(Y) \xcc X' \xcv (X-$A), so
$(X \xcv X' )-(X' \xcv (X-A))=A \xbe \xdi (Y).$

Again, $(I_{ \xbo })$ holds by definition.

(4.1) $(I \xcv disj)$ $ \xch $ $( \xbm disjOR):$

If $X \xcs Y= \xCQ,$ then (1) $A \xbe \xdi (X),B \xbe \xdi (Y) \xch A
\xcv B \xbe \xdi (X \xcv Y)$ and
(2) $A \xbe \xdf (X),B \xbe \xdf (Y) \xch A \xcv B \xbe \xdf (X \xcv Y)$
are equivalent. (By $X \xcs Y= \xCQ,$
$(X-A) \xcv (Y-B)=(X \xcv Y)-(A \xcv B).)$
So $f(X) \xbe \xdf (X),$ $f(Y) \xbe \xdf (Y)$ $ \xch $ (by prerequisite)
$f(X) \xcv f(Y) \xbe \xdf (X \xcv Y).$ $f(X \xcv Y)$
is the smallest element of $ \xdf (X \xcv Y),$ so $f(X \xcv Y) \xcc f(X)
\xcv f(Y).$

(4.2) $( \xbm disjOR)$ $ \xch $ $(I \xcv disj):$

Let $X \xcc Y,$ $X':=Y-$X. Let $A \xbe \xdi (X),$ $A' \xbe \xdi (X' ),$
so $X-A \xbe \xdf (X),$ $X' -A' \xbe \xdf (X' ),$
so $f(X) \xcc X-$A, $f(X' ) \xcc X' -A',$ so $f(X \xcv X' ) \xcc f(X)
\xcv f(X' ) \xcc (X-A) \xcv (X' -A' )$ by
prerequisite, so $(X \xcv X' )-((X-A) \xcv (X' -A' ))=A \xcv A' \xbe \xdi
(X \xcv X' ).$

(5.1) $( \xdm^{+}_{ \xbo })$ $ \xch $ $( \xbm CM):$

$f(X) \xcc Y \xcc X$ $ \xch $ $X-Y \xbe \xdi (X),$ $X-f(X) \xbe \xdi (X)$
$ \xch $ (by $( \xdm^{+}_{ \xbo }),$ (4))
$A:=(X-f(X))-(X-Y) \xbe \xdi (Y)$ $ \xch $
$Y-A=f(X)-(X-Y) \xbe \xdf (Y)$ $ \xch $ $f(Y) \xcc f(X)-(X-Y) \xcc f(X).$

(5.2) $( \xbm CM)$ $ \xch $ $( \xdm^{+}_{ \xbo })$

Let $A \xbe \xdf (X),$ $B \xbe \xdi (X),$ so $f(X) \xcc X-B \xcc X,$ so by
prerequisite $f(X-B) \xcc f(X).$
As $A \xbe \xdf (X),$ $f(X) \xcc A,$ so $f(X-B) \xcc f(X) \xcc A \xcs
(X-B)=A-$B, and $A-B \xbe \xdf (X-$B), so
$( \xdm^{+}_{ \xbo }),$ (5) holds.

(6.1) $( \xdm^{++})$ $ \xch $ $( \xbm RatM):$

Let $X \xcc Y,$ $X \xcs f(Y) \xEd \xCQ.$ If $Y-X \xbe \xdf (Y),$ then
$A:=(Y-X) \xcs f(Y) \xbe \xdf (Y),$ but by
$X \xcs f(Y) \xEd \xCQ $ $A \xcb f(Y),$ contradicting ``smallest'' of
$f(Y).$ So $Y-X \xce \xdf (Y),$ and
by $( \xdm^{++})$ $X-f(Y)=(Y-f(Y))-(Y-X) \xbe \xdi (X),$ so $X \xcs f(Y)
\xbe \xdf (X),$ so $f(X) \xcc f(Y) \xcs X.$

(6.2) $( \xbm RatM)$ $ \xch $ $( \xdm^{++})$

Let $A \xbe \xdf (Y),$ $B \xce \xdf (Y).$ $B \xce \xdf (Y)$ $ \xch $ $Y-B
\xce \xdi (Y)$ $ \xch $ $(Y-B) \xcs f(Y) \xEd \xCQ.$
Set $X:=Y-$B, so $X \xcs f(Y) \xEd \xCQ,$ $X \xcc Y,$ so $f(X) \xcc f(Y)
\xcs X$ by prerequisite.
$f(Y) \xcc A$ $ \xch $ $f(X) \xcc f(Y) \xcs X=f(Y)-B \xcc A-$B.

(7.1) $( \xdi_{ \xbo })$ $ \xch $ $( \xbm AND)$

Trivial.

(7.2) $( \xbm AND)$ $ \xch $ $( \xdi_{ \xbo })$

Trivial.

$ \xcz $
\\[3ex]

 karl-search= End Proposition Ref-Class-Mu-neu Proof
\vspace{7mm}

 *************************************

\vspace{7mm}

\bfa

$\hspace{0.01em}$

(+++ Orig. No.:  Fact i-Rule +++)

{\xssc LABEL: {Fact i-Rule}}
\label{Fact i-Rule}

So $( \xdi_{n})$ is equivalent to the rule:

$ \xba \xcn \xbb_{1}, \Xl, \xba \xcn \xbb_{n}$ $ \xch $ $ \xba \xcL \xCN
\xbb_{1} \xco  \Xl  \xco \xCN \xbb_{n}.$

\efa

\subparagraph{
Proof
}

$\hspace{0.01em}$

(+++*** Orig.:  Proof )

Let $ \xba \xcn \xbb_{1}, \Xl, \xba \xcn \xbb_{n},$ so
$M( \xba \xcu ( \xCN \xbb_{1} \xco  \Xl  \xco \xCN \xbb_{n}))$ $=$ $M(
\xba \xcu \xCN \xbb_{1}) \xcv  \Xl  \xcv M( \xba \xcu \xCN \xbb_{n})$ $
\xEd $ $M( \xba ),$
or $ \xba \xcL \xCN \xbb_{1} \xco  \Xl  \xco \xCN \xbb_{n}.$

The converse is analogue.

$ \xcz $
\\[3ex]

\bfa

$\hspace{0.01em}$

(+++ Orig. No.:  Fact i-Reformulation +++)

{\xssc LABEL: {Fact i-Reformulation}}
\label{Fact i-Reformulation}

$( \xdi_{n})$ can be reformulated to $A_{1}, \Xl,A_{n-1} \xbe \xdi (X)$ $
\xch $ $X-(A_{1} \xcv  \Xl  \xcv A_{n-1}) \xce \xdi (X).$

This translates then to $ \xba \xcn \xbb_{1}, \Xl, \xba \xcn \xbb_{n-1}$
$ \xch $ $ \xba \xcN \xCN \xbb_{1} \xco  \Xl  \xco \xCN \xbb_{n-1}.$

$ \xcz $
\\[3ex]
\section{
Strength of $(AND)$
}

\efa

\bfa

$\hspace{0.01em}$

(+++ Orig. No.:  Fact i+eM->m +++)

{\xssc LABEL: {Fact i+eM->m}}
\label{Fact i+eM->m-a}

$( \xdi_{n})+(eM \xdi )$ entail $ \xdm^{+}_{n}$

\efa

\subparagraph{
Proof
}

$\hspace{0.01em}$

(+++*** Orig.:  Proof )

By prerequisite, $X_{1} \xcc  \Xl  \xcc X_{n},$ and $X_{2}-X_{1} \xbe \xdi
(X_{2}), \Xl,X_{n}-X_{n-1} \xbe \xdi (X_{n}),$ so by
$(eM \xdi )$ $X_{2}-X_{1}, \Xl,X_{n}-X_{n-1} \xbe \xdi (X_{n}).$ By $(
\xdi_{n})$ $X_{n}-((X_{n}-X_{n-1}) \xcv  \Xl  \xcv (X_{2}-X_{1}))$ $=$
$X_{1} \xbe \xdm^{+}(X_{n}).$ $ \xcz $
\\[3ex]

\bfa

$\hspace{0.01em}$

(+++ Orig. No.:  Fact i+eM->m(3) +++)

{\xssc LABEL: {Fact i+eM->m}}
\label{Fact i+eM->m}

$( \xdi_{ \xbo })+(eM \xdi )$ entail $ \xdm^{+}_{ \xbo }$ (3)

\efa

\subparagraph{
Proof
}

$\hspace{0.01em}$

(+++*** Orig.:  Proof )

$A \xbe \xdf (X)$ $ \xch $ $X-A \xbe \xdi (X) \xcc \xdi (Y),$ $Y-X \xbe
\xdi (Y)$ $ \xch $ $Y-A=(Y-X) \xcv (X-A) \xbe \xdi (Y)$ $ \xch $
$A \xbe \xdf (Y).$ $ \xcz $
\\[3ex]

\bfa

$\hspace{0.01em}$

(+++ Orig. No.:  Fact i+em->m(4) +++)

{\xssc LABEL: {Fact i+em->m}}
\label{Fact i+em->m-b}

$( \xdi_{ \xbo })+(eM \xdf )$ entail $ \xdm^{+}_{ \xbo }$ (4)

\efa

\subparagraph{
Proof
}

$\hspace{0.01em}$

(+++*** Orig.:  Proof )

$A,B \xbe \xdi (X)$ $ \xch $ $A \xcv B \xbe \xdi (X)$ $ \xch $ $X-(A \xcv
B) \xbe \xdf (X)$ $ \xch $ (by $(eM \xdf ))$ $X-(A \xcv B) \xbe \xdf
(X-$B) $ \xch $
$(X-B)-(X-(A \xcv B))=A-B \xbe \xdi (X-$B). $ \xcz $
\\[3ex]

\bfa

$\hspace{0.01em}$

(+++ Orig. No.:  Fact on omega +++)

{\xssc LABEL: {Fact on omega}}
\label{Fact on omega}

First, all versions $(._{n})$ for all $n \xbe \xbo $ hold.

Note that in the following conditions, transitivity is ``built in'', so
repetition is implicit.

\efa

Prove $ \xdm^{+}_{ \xbo }$ from $(CUM_{ \xbo })$ and $(AND_{ \xbo }):$

(a) (6.1) is equivalent to: $A \xbe \xdf (X)$ $ \xch $ $(A \xbe \xdi (Y)
\xch X \xbe \xdi (Y)),$ follows from:
$X-A$ is small in $X,$ so in $Y,$ A small in $Y,$ so $X=(X-A) \xcv A$
small in $Y.$

(6.2) $B \xbe \xdi (X)$ $ \xch $ $(A \xbe \xdi (X) \xch A \xbe \xdi
(X-$B)) or
$B \xbe \xdi (X)$ $ \xch $ $(A \xbe \xdm^{+}(X-B) \xch A \xbe
\xdm^{+}(X)),$ but $B \xbe \xdi (X)$ $ \xcj $ $X-B \xbe \xdf (X).$

\bfa

$\hspace{0.01em}$

(+++ Orig. No.:  Fact em+i +++)

{\xssc LABEL: {Fact em+i}}
\label{Fact em+i}

Using $(eM \xdi ),$ we conclude from $( \xdi_{n}):$

$X_{1}-B \xbe \xdi (X_{1}), \Xl,X_{n-1}-B \xbe \xdi (X_{n-1})$ $ \xch $
$(X_{1}-B) \xcv  \Xl  \xcv (X_{n-1}-B) \xce \xdi (X_{1} \xcv  \Xl  \xcv
X_{n-1}).$

This is equivalent to $ \xba_{1} \xcn \xbb, \Xl, \xba_{n-1} \xcn \xbb $
$ \xch $ $ \xba_{1} \xco  \Xl  \xco \xba_{n-1} \xcN \xCN \xbb.$

$ \xcz $
\\[3ex]

\efa

\bfa

$\hspace{0.01em}$

(+++ Orig. No.:  Fact Reformulation +++)

{\xssc LABEL: {Fact Reformulation}}
\label{Fact Reformulation}

We reformulate

$( \xdf_{n}):$ $B_{1} \xbe \xdf (X), \Xl,B_{n} \xbe \xdf (X)$ $ \xch $
$B_{1} \xcs  \Xl  \xcs B_{n} \xEd \xCQ $ or $B_{1} \xcs  \Xl  \xcs B_{n-1}
\xcC X-B_{n},$ thus

$ \xba \xcn \xbb_{1}, \Xl, \xba \xcn \xbb_{n}$ $ \xch $ $ \xba \xcu
\xbb_{1} \xcu  \Xl  \xcu \xbb_{n-1} \xcL \xCN \xbb_{n}.$

Or: $A_{1} \xbe \xdi (X), \Xl,A_{n-1} \xbe \xdi (X)$ $ \xch $ $X-(A_{1}
\xcv  \Xl  \xcv A_{n-1}) \xce \xdi (X),$ so by $(eM \xdf )$ (2)
$X-(A_{1} \xcv  \Xl  \xcv A_{n-1}) \xce \xdi (X-(A_{1} \xcv  \Xl  \xcv
A_{n-2})),$ thus

$ \xba \xcn \xbb_{1}, \Xl, \xba \xcn \xbb_{n-1}$ $ \xch $ $ \xba \xcu
\xbb_{1} \xcu  \Xl  \xcu \xbb_{n-2} \xcN \xCN \xbb_{n-1}.$

(This seems to be the only time we use (2) of $(eM \xdf ).$ Check all
proofs
of $(eM \xdf )$ if (2) holds, too.)

$ \xcz $
\\[3ex]

\efa

\bfa

$\hspace{0.01em}$

(+++ Orig. No.:  Fact i+eM->Rules +++)

{\xssc LABEL: {Fact i+eM->Rules}}
\label{Fact i+eM->Rules}

Let $n \xcg 3.$

(1) In the presence of (iM), $(eM \xdi ),$ $(eM \xdf ),$ $(AND_{n})$
implies $(OR_{n})$ and $(CM_{n}).$

(2) (iM), $(eM \xdi ),$ $(eM \xdf ),$ $(OR_{n})$ do not imply $(AND_{n}).$

(3) (iM), $(eM \xdi ),$ $(eM \xdf ),$ $(CM_{n})$ do not imply $(AND_{n}).$

\efa

\subparagraph{
Proof
}

$\hspace{0.01em}$

(+++*** Orig.:  Proof )

(1)

(2)

Consider

$U:=\{1, \Xl,n+1\},$ $ \xdi (U):=\{ \xCQ \} \xcv \{\{i\}:1 \xck i \xck
n\},$

$X:=\{1, \Xl,n\},$ $ \xdi (X):=\{ \xCQ \} \xcv \{\{i\}:1 \xck i \xck
n\},$

$ \xdi (Y):=\{ \xCQ \}$ for all other $Y \xcc U.$

$(AND_{n})$ fails for $X,$ (iM), $(eM \xdi ),$ $(eM \xdf ),$ $(OR_{n})$
hold:

$(AND_{n})$ fails: trivial.

(iM) holds: trivial.

$(eM \xdi )$ holds: trivial

$(eM \xdf )$ (1) holds: trivial as all big subsets of any $Y$ are either
$Y$ or $Y-\{x\}$
for some $x \xbe Y,$ so if $Y-\{x\} \xcc Y' \xcc Y,$ then $Y-\{x\}=Y'.$

$(eM \xdf )$ (2) holds: If $Y=X$ or $Y=U,$ then for $A \xcc Y$ $A \xbe
\xdm^{+}(Y)$ iff A contains at least
2 elements. If $Y \xEd X,$ $Y \xEd U,$ then for $A \xcc Y$ $A \xbe
\xdm^{+}(Y)$ iff A is not empty. These
properties are inherited downward.

$(OR_{n})$ holds: Let $ \xba_{1} \xcn \xbb, \Xl, \xba_{n-1} \xcn \xbb,$
we have to show $ \xba_{1} \xco  \Xl  \xco \xba_{n-1} \xcN \xCN \xbb.$
If $ \xba_{i}$ is neither $X$ nor $U,$ $ \xba_{i} \xcn \xbb $ is $
\xba_{i} \xcl \xbb.$ So the only exceptions
to $ \xba_{i} \xcl \xbb $ can be some $\{i\}:1 \xck i \xck n,$ but there
can be only one such $i,$
as otherwise $ \xba_{i} \xcu \xCN \xbb $ would be $\{i,i' \},$ which is
not small. But $ \xba_{i}$ contains
at least 2 $j,j' $ s.t. $j,j' \xbe \xba_{i} \xcu \xbb,$ and 2-element
sets are not small, so
$ \xba_{1} \xco  \Xl  \xco \xba_{n-1} \xcN \xCN \xbb.$

It seems that even $(OR_{ \xbo })$ holds.

(3)

Consider

$U:=\{1, \Xl,n+1\},$ $ \xdf (U):=\{A \xcc U:n+1 \xbe A\},$

$X:=\{1, \Xl,n\},$ $ \xdi (X):=\{ \xCQ \} \xcv \{\{i\}:1 \xck i \xck
n\},$

for all other $Y \xcc U$ let

$ \xdf (Y):=\{A \xcc Y:n+1 \xbe A\}$ if $n+1 \xbe Y,$ and

$ \xdf (Y):=\{Y\}$ if $n+1 \xce Y.$

Thus, $(AND_{n})$ fails in $X,$ for all $Y \xcc U,$ $Y \xEd X$ $(< \xbo
*s)$ holds, and (iM), $(eM \xdi ),$
$(eM \xdf )$ hold. $(CM_{n})$ (2) holds, too:

$(AND_{n})$ fails: trivial.

(im) holds: trivial.

$(eM \xdi )$ holds: trivial

$(eM \xdf )$ (1) holds: If $A \xcc Y' \xcc Y$ is big in $Y,$ $n+1 \xbe Y,$
then $n+1 \xbe A,$ so $n+1 \xbe Y',$
and $A \xcc Y' $ is big. If $A \xcc Y' \xcc Y$ is big in $Y,$ $n+1 \xce
Y,$ then $A=Y'.$

$(eM \xdf )$ (2) holds: If $n+1 \xbe Y,$ then there are no medium size
sets. If $Y=X,$
then for $A \xcc Y$ $A \xbe \xdm^{+}(Y)$ iff A contains at least 2
elements. Otherwise,
$A \xbe \xdm^{+}(Y)$ iff A is not empty. These properties inherit
downwards.

$(CM_{n})$ (2) holds:

$((CM_{n})$ version (1) fails here.)

Let $ \xba \xcn \xbb_{1}, \Xl, \xba \xcn \xbb_{n-1},$ we have to show $
\xba \xcu \xbb_{1} \xcu  \Xl  \xcu \xbb_{n-2} \xcN \xCN \xbb_{n-1}.$

Let $ \xba $ correspond to $Y$ with $n+1 \xbe Y.$ Then $n+1 \xbe \xba \xcu
\xbb_{i}$ for all $i,$ so
$ \xba \xcu \xbb_{1} \xcu  \Xl  \xcu \xbb_{n-2} \xcn \xbb_{n-1}.$

Let $ \xba $ correspond to $Y \xEd X$ with $n+1 \xce Y.$ Then $ \xba \xcn
\xbb_{i}$ is $ \xba \xcl \xbb_{i}.$

Let $ \xba $ correspond to $X.$

If $ \xba \xcu \xbb_{1} \xcu  \Xl  \xcu \xbb_{n-2}$ still corresponds to
$X= \xba,$ then $ \xba \xcn \xbb_{n-1}.$

If not, $ \xba \xcu \xbb_{1} \xcu  \Xl  \xcu \xbb_{n-2} \xcn \xCN
\xbb_{n-1}$ is $ \xba \xcu \xbb_{1} \xcu  \Xl  \xcu \xbb_{n-2} \xcl \xCN
\xbb_{n-1}.$
$ \xba \xcu \xbb_{1} \xcu  \Xl  \xcu \xbb_{n-2}$ contains at least 2
elements, and $ \xba \xcu \xCN \xbb_{n-1}$ at most
one element, so $ \xba \xcu \xbb_{1} \xcu  \Xl  \xcu \xbb_{n-2} \xcl \xCN
\xbb_{n-1}$ cannot be.

Does $(CM_{ \xbo })$ hold here?

$ \xcz $
\\[3ex]

\bfa

$\hspace{0.01em}$

(+++ Orig. No.:  Fact More-Rules +++)

{\xssc LABEL: {Fact More-Rules}}
\label{Fact More-Rules}

(1) In the presence of (iM), $(eM \xdi ),$ $(eM \xdf )$ $(AND_{ \xbo })$
imply $(OR_{ \xbo }),$ $(CM_{ \xbo }),$ $( \xdm^{+}_{ \xbo }).$

(2) (iM), $(eM \xdi ),$ $(eM \xdf ),$ $(OR_{ \xbo })$ do not imply $(AND_{
\xbo }).$

(3) (iM), $(eM \xdi ),$ $(eM \xdf ),$ $(CM_{ \xbo })$ do not imply $(AND_{
\xbo }).$

(4) (iM), $(eM \xdi ),$ $(eM \xdf ),$ $( \xdm^{+}_{ \xbo })$ do not imply
$(AND_{ \xbo }).$

(5) (iM), $(eM \xdi ),$ $(eM \xdf ),$ $(OR_{ \xbo }),$ $(CM_{ \xbo })$
imply $(AND_{ \xbo }).$

\efa

\subparagraph{
Proof
}

$\hspace{0.01em}$

(+++*** Orig.:  Proof )

(5): Let $A,B \xcc X$ small, then A small in $X-$B, $B$ small in $X-$A, so
$A \xcv B$ small in
$(X-B) \xcv (X-$A).

$(CUM_{ \xbo }):$ By $(AND_{ \xbo })$ is $ \xCN \xbb \xco \xCN \xbb ' $
small, so $X-( \xCN \xbb \xco \xCN \xbb ' )$ is big in $X,$
thus a fortiori big in $X-( \xCN \xbb ).$
(The argument does not work properly with small sets!)

$(OR_{ \xbo }):$ $ \xba \xcu \xCN \xbb $ is small in $ \xba,$ $ \xba '
\xcu \xCN \xbb $ is small in $ \xba ',$ so a fortiori
small in $ \xba \xco \xba '.$
(The argument does not work properly with big sets!)

$A \xbe \xdf (X),$ $X \xbe \xdf (Y)$ $ \xch $ $A \xbe \xdm^{+}(Y)$ is also
weaker, as first small is ``diluted''

$( \xdm^{+}_{ \xbo })$ (6.3): $X-A$ is small in $X,$ so a fortiori in $Y.$

\section{
Div.
}

$++++++++++++++++++++++++++++++++++++++++++++++$
\subsubsection{Plausibility Logic}

 {\LARGE karl-search= Start Plausibility Logic }

\index{Plausibility Logic}
\subsection{Plausibility Logic}
{\xssc LABEL: {Section Plausibility Logic}}
\label{Section Plausibility Logic}

\paragraph{
Discussion of plausibility logic
}

$\hspace{0.01em}$

(+++*** Orig.:  Discussion of plausibility logic )

{\xssc LABEL: {Section Discussion of plausibility logic}}
\label{Section Discussion of plausibility logic}

Plausibility logic was introduced by $D.$ Lehmann  \cite{Leh92a},
 \cite{Leh92b}
as a sequent
calculus in a propositional language without connectives. Thus, a
plausibility
logic language $ \xdl $ is just a set, whose elements correspond to
propositional
variables, and a sequent has the form $X \xcn Y,$ where $X,$ $Y$ are $
\ul{finite}$ subsets of
$ \xdl,$ thus, in the intuitive reading, $ \xcU X \xcn \xcO Y.$ (We use $
\xcn $ instead of the $ \xcl $ used
in  \cite{Leh92a},  \cite{Leh92b} and continue to reserve $
\xcl $ for classical logic.)

\paragraph{
The details:
}

$\hspace{0.01em}$

(+++*** Orig.:  The details: )

{\xssc LABEL: {Section The details:}}
\label{Section The details:}

\bn

$\hspace{0.01em}$

(+++ Orig. No.:  Notation Plausi-1 +++)

{\xssc LABEL: {Notation Plausi-1}}
\label{Notation Plausi-1}

We abuse notation, and write $X \xcn a$ for $X \xcn \{a\},$ $X,a \xcn Y$
for $X \xcv \{a\} \xcn Y,$ $ab \xcn Y$ for
$\{a,b\} \xcn Y,$ etc. When discussing plausibility logic, $X,Y,$ etc.
will denote finite
subsets of $ \xdl,$ $a,b,$ etc. elements of $ \xdl.$

We first define the logical properties we will examine.

\en

\bd

$\hspace{0.01em}$

(+++ Orig. No.:  Definition Plausi-1 +++)

{\xssc LABEL: {Definition Plausi-1}}
\label{Definition Plausi-1}

$X$ and $Y$ will be finite subsets of $ \xdl,$ a, etc. elements of $ \xdl
.$
The base axiom and rules of plausibility logic are
(we use the prefix ``Pl'' to differentiate them from the usual ones):

(PlI) (Inclusion): $X \xcn a$ for all $a \xbe X,$

(PlRM) (Right Monotony): $X \xcn Y$ $ \xch $ $X \xcn a,Y,$

(PlCLM) (Cautious Left Monotony): $X \xcn a,$ $X \xcn Y$ $ \xch $ $X,a
\xcn Y,$

(PlCC) (Cautious Cut): $X,a_{1} \Xl a_{n} \xcn Y,$ and for all $1 \xck i
\xck n$ $X \xcn a_{i},Y$ $ \xch $ $X \xcn Y,$

and as a special case of (PlCC):

(PlUCC) (Unit Cautious Cut): $X,a \xcn Y$, $X \xcn a,Y$ $ \xch $ $X \xcn
Y.$

and we denote by PL, for plausibility logic, the full system, i.e.
$(PlI)+(PlRM)+(PlCLM)+(PlCC).$

\ed

We now adapt the definition of a preferential model to plausibility logic.
This is the central definition on the semantic side.

\bd

$\hspace{0.01em}$

(+++ Orig. No.:  Definition Plausi-2 +++)

{\xssc LABEL: {Definition Plausi-2}}
\label{Definition Plausi-2}

Fix a plausibility logic language $ \xdl.$ A model for $ \xdl $ is then
just an arbitrary
subset of $ \xdl.$

If $ \xdm:=<M, \xeb >$ is a preferential model s.t. $M$ is a set of
(indexed) $ \xdl -$models,
then for a finite set $X \xcc \xdl $ (to be imagined on the left hand side
of $ \xcn $!), we
define

(a) $m \xcm X$ iff $X \xcc m$

(b) $M(X)$ $:=$ $\{m$: $<m,i> \xbe M$ for some $i$ and $m \xcm X\}$

(c) $ \xbm (X)$ $:=$ $\{m \xbe M(X)$: $ \xcE <m,i> \xbe M. \xCN \xcE <m'
,i' > \xbe M$ $(m' \xbe M(X)$ $ \xcu $ $<m',i' > \xeb <m,i>)\}$

(d) $X \xcm_{ \xdm }Y$ iff $ \xcA m \xbe \xbm (X).m \xcs Y \xEd \xCQ.$

(a) reflects the intuitive reading of $X$ as $ \xcU X,$ and (d) that of
$Y$ as $ \xcO Y$ in
$X \xcn Y.$ Note that $X$ is a set of ``formulas'', and $ \xbm (X)= \xbm_{
\xdm }(M(X)).$

\ed

We note as trivial consequences of the definition.

\bfa

$\hspace{0.01em}$

(+++ Orig. No.:  Fact Plausi-1 +++)

{\xssc LABEL: {Fact Plausi-1}}
\label{Fact Plausi-1}

(a) $a \xcm_{ \xdm }b$ iff for all $m \xbe \xbm (a).b \xbe m$

(b) $X \xcm_{ \xdm }Y$ iff $ \xbm (X) \xcc \xcV \{M(b):b \xbe Y\}$

(c) $m \xbe \xbm (X)$ $ \xcu $ $X \xcc X' $ $ \xcu $ $m \xbe M(X' )$ $
\xcp $ $m \xbe \xbm (X' )$.

\efa

We note without proof: $(PlI)+(PlRM)+(PlCC)$ is complete (and sound) for
preferential models

We note the following fact for smooth preferential models:

\bfa

$\hspace{0.01em}$

(+++ Orig. No.:  Fact Plausi-2 +++)

{\xssc LABEL: {Fact Plausi-2}}
\label{Fact Plausi-2}

Let $ \xCf U,X,Y$ be any sets, $ \xdm $ be smooth for at least $\{Y,X\}$
and
let $ \xbm (Y) \xcc U \xcv X,$ $ \xbm (X) \xcc U,$ then $X \xcs Y \xcs
\xbm (U) \xcc \xbm (Y).$ (This is, of course,
a special case of $( \xbm Cum1),$ see Definition \ref{Definition Cum-Alpha}
(page \pageref{Definition Cum-Alpha}) .

\efa

\be

$\hspace{0.01em}$

(+++ Orig. No.:  Example Plausi-1 +++)

{\xssc LABEL: {Example Plausi-1}}
\label{Example Plausi-1}

Let $ \xdl:=\{a,b,c,d,e,f\},$ and
$ \xdx $ $:=$ $\{a \xcn b$, $b \xcn a$, $a \xcn c$, $a \xcn fd$, $dc
\xcn ba$, $dc \xcn e$, $fcba \xcn e\}.$
(fd stands for $f,d$ etc.) Note that, intuitively, left of $ \xcn $ stands
a
conjunction, right of $ \xcn $ a disjunction - in the tradition of sequent
calculus notation.
We show that $ \xdx $ does not have a smooth representation.

\ee

\bfa

$\hspace{0.01em}$

(+++ Orig. No.:  Fact Plausi-3 +++)

{\xssc LABEL: {Fact Plausi-3}}
\label{Fact Plausi-3}

$ \xdx $ does not entail $a \xcn e.$

\efa

See  \cite{Sch96-3} for a proof.

Suppose now that there is a smooth preferential model $ \xdm =<M, \xeb >$
for plausibility
logic which represents $ \xcn,$ i.e. for all $ \xCf X,Y$ finite subsets
of $ \xdl $ $X \xcn Y$ iff
$X \xcm_{ \xdm }Y.$ (See Definition \ref{Definition Plausi-2} (page
\pageref{Definition Plausi-2})  and Fact \ref{Fact Plausi-1} (page \pageref{Fact
Plausi-1}) .)

$a \xcn a,$ $a \xcn b,$ $a \xcn c$ implies for $m \xbe \xbm (a)$ $a,b,c
\xbe m.$ Moreover, as $a \xcn df,$ then also
$d \xbe m$ or $f \xbe m.$ As $a \xcN e,$ there must be $m \xbe \xbm (a)$
s.t. $e \xce m.$ Suppose now $m \xbe \xbm (a)$
with $f \xbe m.$ So $a,b,c,f \xbe m,$ thus by $m \xbe \xbm (a)$ and Fact
\ref{Fact Plausi-1} (page \pageref{Fact Plausi-1}) ,
$m \xbe \xbm (a,b,c,f).$ But
$fcba \xcn e,$ so $e \xbe m.$ We thus have shown that $m \xbe \xbm (a)$
and $f \xbe m$ implies $e \xbe m.$
Consequently, there must be $m \xbe \xbm (a)$ s.t. $d \xbe m,$ $e \xce m.$
Thus, in particular, as $cd \xcn e,$ there is $m \xbe \xbm (a),$ $a,b,c,d
\xbe m,$ $m \xce \xbm (cd).$
But by $cd \xcn ab,$ and $b \xcn a,$ $ \xbm (cd) \xcc M(a) \xcv M(b)$ and
$ \xbm (b) \xcc M(a)$ by
Fact \ref{Fact Plausi-1} (page \pageref{Fact Plausi-1}) .
Let now $T:=M(cd),$ $R:=M(a),$ $S:=M(b),$ and $ \xbm_{ \xdm }$ be the
choice function of the
minimal elements in the structure $ \xdm,$ we then have by $ \xbm (S)=
\xbm_{ \xdm }(M(S))$:

1. $ \xbm_{ \xdm }(T) \xcc R \xcv S,$

2. $ \xbm_{ \xdm }(S) \xcc R,$

3. there is $m \xbe S \xcs T \xcs \xbm_{ \xdm }(R),$ but $m \xce \xbm_{
\xdm }(T),$

but this contradicts above Fact \ref{Fact Plausi-2} (page \pageref{Fact
Plausi-2}) .

 karl-search= End Plausibility Logic
\vspace{7mm}

 *************************************

\vspace{7mm}

\subsubsection{Arieli-Avron}

 {\LARGE karl-search= Start Arieli-Avron }

\index{Arieli-Avron}
\subsection{A comment on work by Arieli and Avron}
{\xssc LABEL: {Section Arieli-Avron}}
\label{Section Arieli-Avron}

We turn to a similar case, published in  \cite{AA00}.
Definitions are due to  \cite{AA00}, for motivation the reader is
referred there.

We follow here the convention of Arieli and Avron and use upper-case
Greek letters for sets of formulae. At the same time this different
notation should remind the
reader that sets of formulae are read as conjunctions on the left of $
\xcn,$
and as disjunctions on the right of $ \xcn.$

\bd

$\hspace{0.01em}$

(+++ Orig. No.:  Definition Arieli-Avron-1 +++)

{\xssc LABEL: {Definition Arieli-Avron-1}}
\label{Definition Arieli-Avron-1}

(1) A Scott consequence relation, abbreviated scr, is a binary relation $
\xcl $
between sets of formulae, that satisfies the following conditions:

(s-R) if $ \xbG \xcs \xbD \xEd \xCQ,$ the $ \xbG \xcl \xbD $
(M) if $ \xbG \xcl \xbD $ and $ \xbG \xcc \xbG ',$ $ \xbD \xcc \xbD ',$
then $ \xbG ' \xcl \xbD ' $
(C) if $ \xbG \xcl \xbq, \xbD $ and $ \xbG ', \xbq \xcl \xbD ',$ then $
\xbG, \xbG ' \xcl \xbD, \xbD ' $

(2) A Scott cautious consequence relation, abbreviated sccr, is a binary
relation $ \xcn $ between nonempty sets of formulae, that satisfies the
following
conditions:

(s-R) if $ \xbG \xcs \xbD \xEd \xCQ,$ the $ \xbG \xcn \xbD $
(CM) if $ \xbG \xcn \xbD $ and $ \xbG \xcn \xbq,$ then $ \xbG, \xbq \xcn
\xbD $
(CC) if $ \xbG \xcn \xbq $ and $ \xbG, \xbq \xcn \xbD,$ then $ \xbG \xcn
\xbD.$

\ed

\be

$\hspace{0.01em}$

(+++ Orig. No.:  Example Arieli-Avron-1 +++)

{\xssc LABEL: {Example Arieli-Avron-1}}
\label{Example Arieli-Avron-1}

We have two consequence relations, $ \xcl $ and $ \xcn.$

The rules to consider are

$LCC^{n}$ $ \frac{ \xbG \xcn \xbq_{1}, \xbD  \Xl  \xbG \xcn \xbq_{n}, \xbD
\xbG, \xbq_{1}, \Xl, \xbq_{n} \xcn \xbD }{ \xbG \xcn \xbD }$

$RW^{n}$ $ \frac{ \xbG \xcn \xbq_{i}, \xbD i=1 \Xl n \xbG, \xbq_{1}, \Xl
, \xbq_{n} \xcl \xbf }{ \xbG \xcn \xbf, \xbD }$

Cum $ \xbG, \xbD \xEd \xCQ,$ $ \xbG \xcl \xbD $ $ \xcp $ $ \xbG \xcn
\xbD $

RM $ \xbG \xcn \xbD $ $ \xcp $ $ \xbG \xcn \xbq, \xbD $

CM $ \frac{ \xbG \xcn \xbq \xbG \xcn \xbD }{ \xbG, \xbq \xcn \xbD }$

$s-R$ $ \xbG \xcs \xbD \xEd \xCQ $ $ \xcp $ $ \xbG \xcn \xbD $

$M$ $ \xbG \xcl \xbD,$ $ \xbG \xcc \xbG ',$ $ \xbD \xcc \xbD ' $ $ \xcp
$ $ \xbG ' \xcl \xbD ' $

$C$ $ \frac{ \xbG_{1} \xcl \xbq, \xbD_{1} \xbG_{2}, \xbq \xcl \xbD_{2}}{
\xbG_{1}, \xbG_{2} \xcl \xbD_{1}, \xbD_{2}}$

Let $ \xdl $ be any set.
Define now $ \xbG \xcl \xbD $ iff $ \xbG \xcs \xbD \xEd \xCQ.$
Then $s-R$ and $M$ for $ \xcl $ are trivial. For $C:$ If $ \xbG_{1} \xcs
\xbD_{1} \xEd \xCQ $ or $ \xbG_{1} \xcs \xbD_{1} \xEd \xCQ,$ the
result is trivial. If not, $ \xbq \xbe \xbG_{1}$ and $ \xbq \xbe
\xbD_{2},$ which implies the result.
So $ \xcl $ is a scr.

Consider now the rules for a sccr which is $ \xcl -$plausible for this $
\xcl.$
Cum is equivalent to $s-$R, which is essentially (PlI) of Plausibility
Logic.
Consider $RW^{n}.$ If $ \xbf $ is one of the $ \xbq_{i},$ then the
consequence $ \xbG \xcn \xbf, \xbD $ is a case
of one of the other hypotheses. If not, $ \xbf \xbe \xbG,$ so $ \xbG \xcn
\xbf $ by $s-$R, so $ \xbG \xcn \xbf, \xbD $
by RM (if $ \xbD $ is finite). So, for this $ \xcl,$ $RW^{n}$ is a
consequence of $s-R$ $+$ RM.

We are left with $LCC^{n},$ RM, CM, $s-$R, it was shown in  \cite{Sch04}
and  \cite{Sch96-3}
that this does not suffice to guarantee smooth representability, by
failure of
$( \xbm Cum1)$ - see Definition \ref{Definition Cum-Alpha} (page
\pageref{Definition Cum-Alpha}) .

 karl-search= End Arieli-Avron
\vspace{7mm}

 *************************************

\vspace{7mm}

\ee

$++++++++++++++++++++++++++++++++++++++++++++++$
\subsubsection{Comment Cum-Union}

 {\LARGE karl-search= Start Comment Cum-Union }

\index{Comment Cum-Union}

\bcom

$\hspace{0.01em}$

(+++ Orig. No.:  Comment Cum-Union +++)

{\xssc LABEL: {Comment Cum-Union}}
\label{Comment Cum-Union}

We show here that, without sufficient closure
properties, there is an infinity of versions of cumulativity, which
collapse
to usual cumulativity when the domain is closed under finite unions.
Closure properties thus reveal themselves as a powerful tool to show
independence of properties.

We work in some fixed arbitrary set $Z,$ all sets considered will be
subsets of $Z.$

Unless said otherwise, we use without further
mentioning $( \xbm PR)$ and $( \xbm \xcc ).$

 karl-search= End Comment Cum-Union
\vspace{7mm}

 *************************************

\vspace{7mm}

\subsubsection{Definition Cum-Alpha}

 {\LARGE karl-search= Start Definition Cum-Alpha }

\index{Definition Cum-Alpha}

\ecom

\bd

$\hspace{0.01em}$

(+++ Orig. No.:  Definition Cum-Alpha +++)

{\xssc LABEL: {Definition Cum-Alpha}}
\label{Definition Cum-Alpha}

For any ordinal $ \xba,$ we define

$( \xbm Cum \xba ):$

If for all $ \xbb \xck \xba $ $ \xbm (X_{ \xbb }) \xcc U \xcv \xcV \{X_{
\xbg }: \xbg < \xbb \}$ hold, then so does
$ \xcS \{X_{ \xbg }: \xbg \xck \xba \} \xcs \xbm (U) \xcc \xbm (X_{ \xba
}).$

$( \xbm Cumt \xba ):$

If for all $ \xbb \xck \xba $ $ \xbm (X_{ \xbb }) \xcc U \xcv \xcV \{X_{
\xbg }: \xbg < \xbb \}$ hold, then so does
$X_{ \xba } \xcs \xbm (U) \xcc \xbm (X_{ \xba }).$

( `` $t$ '' stands for transitive, see Fact \ref{Fact Cum-Alpha} (page
\pageref{Fact Cum-Alpha}) , (2.2)
below.)

$( \xbm Cum \xca )$ and $( \xbm Cumt \xca )$ will be the class of all $(
\xbm Cum \xba )$ or $( \xbm Cumt \xba )$ -
read their ``conjunction'', i.e. if we say that $( \xbm Cum \xca )$ holds,
we mean that
all $( \xbm Cum \xba )$ hold.

 karl-search= End Definition Cum-Alpha
\vspace{7mm}

 *************************************

\vspace{7mm}

\subsubsection{Note Cum-Alpha}

 {\LARGE karl-search= Start Note Cum-Alpha }

\index{Note Cum-Alpha}

\ed

\paragraph{
Note
}

$\hspace{0.01em}$

(+++*** Orig.:  Note )

{\xssc LABEL: {Section Note}}
\label{Section Note}

The first conditions thus have the form:

$( \xbm Cum0)$ $ \xbm (X_{0}) \xcc U$ $ \xcp $ $X_{0} \xcs \xbm (U) \xcc
\xbm (X_{0}),$

$( \xbm Cum1)$ $ \xbm (X_{0}) \xcc U,$ $ \xbm (X_{1}) \xcc U \xcv X_{0}$ $
\xcp $ $X_{0} \xcs X_{1} \xcs \xbm (U) \xcc \xbm (X_{1}),$

$( \xbm Cum2)$ $ \xbm (X_{0}) \xcc U,$ $ \xbm (X_{1}) \xcc U \xcv X_{0},$
$ \xbm (X_{2}) \xcc U \xcv X_{0} \xcv X_{1}$ $ \xcp $ $X_{0} \xcs X_{1}
\xcs X_{2} \xcs \xbm (U) \xcc \xbm (X_{2}).$

$( \xbm Cumt \xba )$ differs from $( \xbm Cum \xba )$ only in the
consequence, the intersection contains
only the last $X_{ \xba }$ - in particular, $( \xbm Cum0)$ and $( \xbm
Cumt0)$ coincide.

Recall that condition $( \xbm Cum1)$ is the crucial condition in  \cite{Leh92a},
which
failed, despite $( \xbm CUM),$ but which has to hold in all smooth models.
This
condition $( \xbm Cum1)$ was the starting point of the investigation.

We briefly mention some major results on above conditions, taken from
Fact \ref{Fact Cum-Alpha} (page \pageref{Fact Cum-Alpha})  and
shown there - we use the same numbering:

(1.1) $( \xbm Cum \xba )$ $ \xcp $ $( \xbm Cum \xbb )$ for all $ \xbb \xck
\xba $

(1.2) $( \xbm Cumt \xba )$ $ \xcp $ $( \xbm Cumt \xbb )$ for all $ \xbb
\xck \xba $

(2.1) All $( \xbm Cum \xba )$ hold in smooth preferential structures

(2.2) All $( \xbm Cumt \xba )$ hold in transitive smooth preferential
structures

(3.1) $( \xbm Cum \xbb )$ $+$ $( \xcv )$ $ \xcp $ $( \xbm Cum \xba )$ for
all $ \xbb \xck \xba $

(3.2) $( \xbm Cumt \xbb )$ $+$ $( \xcv )$ $ \xcp $ $( \xbm Cumt \xba )$
for all $ \xbb \xck \xba $

(5.2) $( \xbm Cum \xba )$ $ \xcp $ $( \xbm CUM)$ for all $ \xba $

(5.3) $( \xbm CUM)$ $+$ $( \xcv )$ $ \xcp $ $( \xbm Cum \xba )$ for all $
\xba $

 karl-search= End Note Cum-Alpha
\vspace{7mm}

 *************************************

\vspace{7mm}

\subsubsection{Definition HU-All}

 {\LARGE karl-search= Start Definition HU-All }

\index{Definition HU-All}

The following inductive definition of $H(U,u)$ and of the property $ \xCf
(HU,u)$
concerns closure under $( \xbm Cum \xca ),$ its main
property is formulated in Fact \ref{Fact HU} (page \pageref{Fact HU}) , its main
interest
is its use in the proof of Proposition D-4.4.6.

\bd

$\hspace{0.01em}$

(+++ Orig. No.:  Definition HU-All +++)

{\xssc LABEL: {Definition HU-All}}
\label{Definition HU-All}

$(H(U,u)_{ \xba },$ $H(U)_{ \xba },$ $ \xCf (HU,u),$ $ \xCf (HU).)$

$H(U,u)_{0}$ $:=$ $U,$

$H(U,u)_{ \xba +1}$ $:=$ $H(U,u)_{ \xba }$ $ \xcv $ $ \xcV \{X:$ $u \xbe
X$ $ \xcu $ $ \xbm (X) \xcc H(U,u)_{ \xba }\},$

$H(U,u)_{ \xbl }$ $:=$ $ \xcV \{H(U,u)_{ \xba }: \xba < \xbl \}$ for
$limit( \xbl ),$

$H(U,u)$ $:=$ $ \xcV \{H(U,u)_{ \xba }: \xba < \xbk \}$ for $ \xbk $
sufficiently big $(card(Z)$ suffices, as

the procedure trivializes, when we cannot add any new elements).

$ \xCf (HU,u)$ is the property:

$u \xbe \xbm (U),$ $u \xbe Y- \xbm (Y)$ $ \xcp $ $ \xbm (Y) \xcC H(U,u)$ -
of course for all $u$ and $U.$
$(U,Y \xbe \xdy ).$

For the case with $( \xcv ),$ we further define, independent of $u,$

$H(U)_{0}$ $:=$ $U,$

$H(U)_{ \xba +1}$ $:=$ $H(U)_{ \xba }$ $ \xcv $ $ \xcV \{X:$ $ \xbm (X)
\xcc H(U)_{ \xba }\},$

$H(U)_{ \xbl }$ $:=$ $ \xcV \{H(U)_{ \xba }: \xba < \xbl \}$ for $limit(
\xbl ),$

$H(U)$ $:=$ $ \xcV \{H(U)_{ \xba }: \xba < \xbk \}$ again for $ \xbk $
sufficiently big

$ \xCf (HU)$ is the property:

$u \xbe \xbm (U),$ $u \xbe Y- \xbm (Y)$ $ \xcp $ $ \xbm (Y) \xcC H(U)$ -
of course for all $U.$
$(U,Y \xbe \xdy ).$

\ed

Obviously, $H(U,u) \xcc H(U),$ so $(HU) \xcp (HU,u).$

 karl-search= End Definition HU-All
\vspace{7mm}

 *************************************

\vspace{7mm}

\subsubsection{Definition St-Tree}

 {\LARGE karl-search= Start Definition St-Tree }

\index{Definition St-Tree}

\bd

$\hspace{0.01em}$

(+++ Orig. No.:  Definition St-Tree +++)

{\xssc LABEL: {Definition St-Tree}}
\label{Definition St-Tree}

(ST-trees and $( \xbm ST))$

Let $u \xbe \xbm (U).$

A tree $t$ (of height $ \xck \xbo )$ is an ST-tree for $<u,U>$ iff

(1) the nodes are pairs $<x,X>$ s.t. $x \xbe \xbm (X)$

(2) Level 0:

$<u,U>$ is the root.

(3) Level $n \xcp n+1$

Let $<x_{n},X_{n}>$ be at level $n$ (so $x_{n} \xbe \xbm (X_{n})).$

(3.1) For all $X_{n+1} \xbe \xdy $ s.t. $x_{n} \xbe X_{n+1}- \xbm
(X_{n+1})$
there is a successor $<x_{n+1},X_{n+1}>$ of $<x_{n},X_{n}>$ in $t$ with
$x_{n+1} \xbe \xbm (X_{n+1})$ and
$x_{n+1} \xce H(X_{n},x_{n})$ and for all predessors $<x',X' >$ of
$<x_{n},X_{n}>$ also $x_{n+1} \xce H(X',x' ).$

(3.2) For all $X_{n+1} \xbe \xdy $ s.t. $x_{n} \xbe \xbm (X_{n+1})$ and
s.t. there is a successor
$<x'_{n+1},X'_{n+1}>$ of $<x_{n},X_{n}>$ in $t$ with $x'_{n+1} \xbe
X_{n+1}- \xbm (X_{n+1})$
there is a successor $<x_{n+1},X_{n+1}>$ of $<x_{n},X_{n}>$ in $t$ with
$x_{n+1} \xbe \xbm (X_{n+1})$ and
$x_{n+1} \xce H(X_{n},x_{n})$ and for all predessors $<x',X' >$ of
$<x_{n},X_{n}>$ also $x_{n+1} \xce H(X',x' ).$

Finally $( \xbm ST)$ is the condition:

For all $U \xbe \xdy,$ $u \xbe \xbm (U)$ there is an ST-tree for $<u,U>.$

 karl-search= End Definition St-Tree
\vspace{7mm}

 *************************************

\vspace{7mm}

\subsubsection{Example Inf-Cum-Alpha}

 {\LARGE karl-search= Start Example Inf-Cum-Alpha }

\index{Example Inf-Cum-Alpha}

\ed

\be

$\hspace{0.01em}$

(+++ Orig. No.:  Example Inf-Cum-Alpha +++)

{\xssc LABEL: {Example Inf-Cum-Alpha}}
\label{Example Inf-Cum-Alpha}

This important example shows that the conditions $( \xbm Cum \xba )$ and
$( \xbm Cumt \xba )$ defined in Definition \ref{Definition Cum-Alpha} (page
\pageref{Definition Cum-Alpha})
are all different in the absence of $( \xcv ),$ in its
presence they all collapse (see Fact \ref{Fact Cum-Alpha} (page \pageref{Fact
Cum-Alpha})  below).
More precisely,
the following (class of) examples shows that the $( \xbm Cum \xba )$
increase in
strength. For any finite or infinite ordinal $ \xbk >0$ we construct an
example s.t.

(a) $( \xbm PR)$ and $( \xbm \xcc )$ hold

(b) $( \xbm CUM)$ holds

(c) $( \xcS )$ holds

(d) $( \xbm Cumt \xba )$ holds for $ \xba < \xbk $

(e) $( \xbm Cum \xbk )$ fails.

\ee

\subparagraph{
Proof
}

$\hspace{0.01em}$

(+++*** Orig.:  Proof )

We define a suitable base set and a non-transitive binary relation $ \xeb
$
on this set, as well as a suitable set $ \xdx $ of subsets, closed under
arbitrary
intersections, but not under finite unions, and define $ \xbm $ on these
subsets
as usual in preferential structures by $ \xeb.$ Thus, $( \xbm PR)$ and $(
\xbm \xcc )$ will hold.
It will be immediate that $( \xbm Cum \xbk )$ fails, and we will show that
$( \xbm CUM)$ and
$( \xbm Cumt \xba )$ for $ \xba < \xbk $ hold by examining the cases.

For simplicity, we first define a set of generators for $ \xdx,$ and
close under
$( \xcS )$ afterwards. The set $U$ will have a special position, it is the
``useful''
starting point to construct chains corresponding to above definitions of
$( \xbm Cum \xba )$ and $( \xbm Cumt \xba ).$

In the sequel,
$i,j$ will be successor ordinals, $ \xbl $ etc. limit ordinals, $ \xba,$
$ \xbb,$ $ \xbk $ any ordinals,
thus e.g. $ \xbl \xck \xbk $ will imply that $ \xbl $ is a limit ordinal $
\xck \xbk,$ etc.

\paragraph{
The base set and the relation $ \xeb $:
}

$\hspace{0.01em}$

(+++*** Orig.:  The base set and the relation  b: )

{\xssc LABEL: {Section The base set and the relation  b:}}
\label{Section The base set and the relation  b:}

$ \xbk >0$ is fixed, but arbitrary. We go up to $ \xbk >0.$

The base set is $\{a,b,c\}$ $ \xcv $ $\{d_{ \xbl }: \xbl \xck \xbk \}$ $
\xcv $ $\{x_{ \xba }: \xba \xck \xbk +1\}$ $ \xcv $ $\{x'_{ \xba }: \xba
\xck \xbk \}.$
$a \xeb b \xeb c,$ $x_{ \xba } \xeb x_{ \xba +1},$ $x_{ \xba } \xeb x'_{
\xba },$ $x'_{0} \xeb x_{ \xbl }$ (for any $ \xbl )$ - $ \xeb $ is NOT
transitive.

\paragraph{
The generators:
}

$\hspace{0.01em}$

(+++*** Orig.:  The generators: )

{\xssc LABEL: {Section The generators:}}
\label{Section The generators:}

$U:=\{a,c,x_{0}\} \xcv \{d_{ \xbl }: \xbl \xck \xbk \}$ - i.e. $ \Xl
.\{d_{ \xbl }:lim( \xbl ) \xcu \xbl \xck \xbk \},$

$X_{i}:=\{c,x_{i},x'_{i},x_{i+1}\}$ $(i< \xbk ),$

$X_{ \xbl }:=\{c,d_{ \xbl },x_{ \xbl },x'_{ \xbl },x_{ \xbl +1}\} \xcv
\{x'_{ \xba }: \xba < \xbl \}$ $( \xbl < \xbk ),$

$X'_{ \xbk }:=\{a,b,c,x_{ \xbk },x'_{ \xbk },x_{ \xbk +1}\}$ if $ \xbk $
is a successor,

$X'_{ \xbk }:=\{a,b,c,d_{ \xbk },x_{ \xbk },x'_{ \xbk },x_{ \xbk +1}\}
\xcv \{x'_{ \xba }: \xba < \xbk \}$ if $ \xbk $ is a limit.

Thus, $X'_{ \xbk }=X_{ \xbk } \xcv \{a,b\}$ if $X_{ \xbk }$ were defined.

Note that there is only one $X'_{ \xbk },$ and $X_{ \xba }$ is defined
only for $ \xba < \xbk,$ so we will
not have $X_{ \xba }$ and $X'_{ \xba }$ at the same time.

Thus, the values of the generators under $ \xbm $ are:

$ \xbm (U)=U,$

$ \xbm (X_{i})=\{c,x_{i}\},$

$ \xbm (X_{ \xbl })=\{c,d_{ \xbl }\} \xcv \{x'_{ \xba }: \xba < \xbl \},$

$ \xbm (X'_{i})=\{a,x_{i}\}$ $(i>0,$ $i$ has to be a successor),

$ \xbm (X'_{ \xbl })=\{a,d_{ \xbl }\} \xcv \{x'_{ \xba }: \xba < \xbl \}.$

(We do not assume that the domain is closed under $ \xbm.)$

\paragraph{
Intersections:
}

$\hspace{0.01em}$

(+++*** Orig.:  Intersections: )

{\xssc LABEL: {Section Intersections:}}
\label{Section Intersections:}

We consider first pairwise intersections:

(1) $U \xcs X_{0}=\{c,x_{0}\},$

(2) $U \xcs X_{i}=\{c\},$ $i>0,$

(3) $U \xcs X_{ \xbl }=\{c,d_{ \xbl }\},$

(4) $U \xcs X'_{i}=\{a,c\}$ $(i>0),$

(5) $U \xcs X'_{ \xbl }=\{a,c,d_{ \xbl }\},$

(6) $X_{i} \xcs X_{j}:$

(6.1) $j=i+1$ $\{c,x_{i+1}\},$

(6.2) else $\{c\},$

(7) $X_{i} \xcs X_{ \xbl }:$

(7.1) $i< \xbl $ $\{c,x'_{i}\},$

(7.2) $i= \xbl +1$ $\{c,x_{ \xbl +1}\},$

(7.3) $i> \xbl +1$ $\{c\},$

(8) $X_{ \xbl } \xcs X_{ \xbl ' }:$ $\{c\} \xcv \{x'_{ \xba }: \xba \xck
min( \xbl, \xbl ' )\}.$

As $X'_{ \xbk }$ occurs only once, $X_{ \xba } \xcs X'_{ \xbk }$ etc. give
no new results.

Note that $ \xbm $ is constant on all these pairwise intersections.

Iterated intersections:

As $c$ is an element of all sets, sets of the type $\{c,z\}$ do not give
any
new results. The possible subsets of $\{a,c,d_{ \xbl }\}:$ $\{c\},$
$\{a,c\},$ $\{c,d_{ \xbl }\}$ exist
already. Thus, the only source of new sets via iterated intersections is
$X_{ \xbl } \xcs X_{ \xbl ' }=\{c\} \xcv \{x'_{ \xba }: \xba \xck min(
\xbl, \xbl ' )\}.$ But, to intersect them, or with some
old sets, will not generate any new sets either. Consequently, the example
satisfies $( \xcS )$ for $ \xdx $ defined by $U,$ $X_{i}$ $(i< \xbk ),$
$X_{ \xbl }$ $( \xbl < \xbk ),$ $X'_{ \xbk },$ and above
paiwise intersections.

We will now verify the positive properties. This is tedious, but
straightforward, we have to check the different cases.

\paragraph{
Validity of $( \xbm CUM)$:
}

$\hspace{0.01em}$

(+++*** Orig.:  Validity of ( mCUM): )

{\xssc LABEL: {Section Validity of ( mCUM):}}
\label{Section Validity of ( mCUM):}

Consider the prerequisite $ \xbm (X) \xcc Y \xcc X.$ If $ \xbm (X)=X$ or
if $X- \xbm (X)$ is a singleton,
$X$ cannot give a violation of $( \xbm CUM).$ So we are left with the
following
candidates for $X:$

(1) $X_{i}:=\{c,x_{i},x'_{i},x_{i+1}\},$ $ \xbm (X_{i})=\{c,x_{i}\}$

Interesting candidates for $Y$ will have 3 elements, but they will all
contain
a. (If $ \xbk < \xbo:$ $U=\{a,c,x_{0}\}.)$

(2) $X_{ \xbl }:=\{c,d_{ \xbl },x_{ \xbl },x'_{ \xbl },x_{ \xbl +1}\} \xcv
\{x'_{ \xba }: \xba < \xbl \},$ $ \xbm (X_{ \xbl })=\{c,d_{ \xbl }\} \xcv
\{x'_{ \xba }: \xba < \xbl \}$

The only sets to contain $d_{ \xbl }$ are $X_{ \xbl },$ $U,$ $U \xcs X_{
\xbl }.$ But $a \xbe U,$ and
$U \xcs X_{ \xbl }$ ist finite. $(X_{ \xbl }$ and $X'_{ \xbl }$ cannot be
present at the same time.)

(3) $X'_{i}:=\{a,b,c,x_{i},x'_{i},x_{i+1}\},$ $ \xbm (X'_{i})=\{a,x_{i}\}$

a is only in $U,$ $X'_{i},$ $U \xcs X'_{i}=\{a,c\},$ but $x_{i} \xce U,$
as $i>0.$

(4) $X'_{ \xbl }:=\{a,b,c,d_{ \xbl },x_{ \xbl },x'_{ \xbl },x_{ \xbl +1}\}
\xcv \{x'_{ \xba }: \xba < \xbl \},$ $ \xbm (X'_{ \xbl })=\{a,d_{ \xbl }\}
\xcv \{x'_{ \xba }: \xba < \xbl \}$

$d_{ \xbl }$ is only in $X'_{ \xbl }$ and $U,$ but $U$ contains no $x'_{
\xba }.$

Thus, $( \xbm CUM)$ holds trivially.

\paragraph{
$( \xbm Cumt \xba )$ hold for $ \xba < \xbk $:
}

$\hspace{0.01em}$

(+++*** Orig.:  ( mCumt a) hold for  a< k: )

{\xssc LABEL: {Section ( mCumt a) hold for  a< k:}}
\label{Section ( mCumt a) hold for  a< k:}

To simplify language, we say that we reach $Y$ from $X$ iff $X \xEd Y$ and
there is a
sequence $X_{ \xbb },$ $ \xbb \xck \xba $ and $ \xbm (X_{ \xbb }) \xcc X
\xcv \xcV \{X_{ \xbg }: \xbg < \xbb \},$ and $X_{ \xba }=Y,$ $X_{0}=X.$
Failure of $( \xbm Cumt \xba )$ would then mean that there are $X$ and
$Y,$ we can reach
$Y$ from $X,$ and $x \xbe ( \xbm (X) \xcs Y)- \xbm (Y).$ Thus, in a
counterexample, $Y= \xbm (Y)$ is
impossible, so none of the intersections can be such $Y.$

To reach $Y$ from $X,$ we have to get started from $X,$ i.e. there must be
$Z$ s.t.
$ \xbm (Z) \xcc X,$ $Z \xcC X$ (so $ \xbm (Z) \xEd Z).$ Inspection of the
different cases shows that
we cannot reach any set $Y$ from any case of the intersections, except
from
(1), (6.1), (7.2).

If $Y$ contains a globally minimal element (i.e. there is no smaller
element in
any set), it can only be reached from any $X$ which already contains this
element. The globally minimal elements are a, $x_{0},$ and the $d_{ \xbl
},$ $ \xbl \xck \xbk.$

By these observations, we see that $X_{ \xbl }$ and $X'_{ \xbk }$ can only
be reached from $U.$
From no $X_{ \xba }$ $U$ can be reached, as the globally minimal a is
missing. But
$U$ cannot be reached from $X'_{ \xbk }$ either, as the globally minimal
$x_{0}$ is missing.

When we look at the relation $ \xeb $ defining $ \xbm,$ we see that we
can reach $Y$ from $X$
only by going upwards, adding bigger elements. Thus, from $X_{ \xba },$ we
cannot reach
any $X_{ \xbb },$ $ \xbb < \xba,$ the same holds for $X'_{ \xbk }$ and
$X_{ \xbb },$ $ \xbb < \xbk.$ Thus, from $X'_{ \xbk },$ we
cannot go anywhere interesting (recall that the intersections are not
candidates
for a $Y$ giving a contradiction).

Consider now $X_{ \xba }.$ We can go up to any $X_{ \xba +n},$ but not to
any $X_{ \xbl },$ $ \xba < \xbl,$ as
$d_{ \xbl }$ is missing, neither to $X'_{ \xbk },$ as a is missing. And we
will be stopped by
the first $ \xbl > \xba,$ as $x_{ \xbl }$ will be missing to go beyond
$X_{ \xbl }.$ Analogous observations
hold for the remaining intersections (1), (6.1), (7.2). But in all these
sets we
can reach, we will not destroy minimality of any element of $X_{ \xba }$
(or of the
intersections).

Consequently, the only candidates for failure will all start with $U.$ As
the only
element of $U$ not globally minimal is $c,$ such failure has to have $c
\xbe Y- \xbm (Y),$ so
$Y$ has to be $X'_{ \xbk }.$ Suppose we omit one of the $X_{ \xba }$ in
the sequence going up to
$X'_{ \xbk }.$ If $ \xbk \xcg \xbl > \xba,$ we cannot reach $X_{ \xbl }$
and beyond, as $x'_{ \xba }$ will be missing.
But we cannot go to $X_{ \xba +n}$ either, as $x_{ \xba +1}$ is missing.
So we will be stopped
at $X_{ \xba }.$ Thus, to see failure, we need the full sequence
$U=X_{0},$ $X'_{ \xbk }=Y_{ \xbk },$
$Y_{ \xba }=X_{ \xba }$ for $0< \xba < \xbk.$

\paragraph{
$( \xbm Cum \xbk )$ fails:
}

$\hspace{0.01em}$

(+++*** Orig.:  ( mCum k) fails: )

{\xssc LABEL: {Section ( mCum k) fails:}}
\label{Section ( mCum k) fails:}

The full sequence $U=X_{0},$ $X'_{ \xbk }=Y_{ \xbk },$ $Y_{ \xba }=X_{
\xba }$ for $0< \xba < \xbk $ shows this, as
$c \xbe \xbm (U) \xcs X'_{ \xbk },$ but $c \xce \xbm (X'_{ \xbk }).$

Consequently, the example satisfies $( \xcS ),$ $( \xbm CUM),$ $( \xbm
Cumt \xba )$ for $ \xba < \xbk,$ and
$( \xbm Cum \xbk )$ fails.

$ \xcz $
\\[3ex]

 karl-search= End Example Inf-Cum-Alpha
\vspace{7mm}

 *************************************

\vspace{7mm}

\subsubsection{Fact Cum-Alpha}

 {\LARGE karl-search= Start Fact Cum-Alpha }

\index{Fact Cum-Alpha}

\bfa

$\hspace{0.01em}$

(+++ Orig. No.:  Fact Cum-Alpha +++)

{\xssc LABEL: {Fact Cum-Alpha}}
\label{Fact Cum-Alpha}

We summarize some properties of $( \xbm Cum \xba )$ and $( \xbm Cumt \xba
)$ - sometimes with some
redundancy. Unless said otherwise, $ \xba,$ $ \xbb $ etc. will be
arbitrary ordinals.

For (1) to (6) $( \xbm PR)$ and $( \xbm \xcc )$ are assumed to hold, for
(7) only
$( \xbm \xcc ).$

(1) Downward:

(1.1) $( \xbm Cum \xba )$ $ \xcp $ $( \xbm Cum \xbb )$ for all $ \xbb \xck
\xba $

(1.2) $( \xbm Cumt \xba )$ $ \xcp $ $( \xbm Cumt \xbb )$ for all $ \xbb
\xck \xba $

(2) Validity of $( \xbm Cum \xba )$ and $( \xbm Cumt \xba )$:

(2.1) All $( \xbm Cum \xba )$ hold in smooth preferential structures

(2.2) All $( \xbm Cumt \xba )$ hold in transitive smooth preferential
structures

(2.3) $( \xbm Cumt \xba )$ for $0< \xba $ do not necessarily hold in
smooth structures without
transitivity, even in the presence of $( \xcS )$

(3) Upward:

(3.1) $( \xbm Cum \xbb )$ $+$ $( \xcv )$ $ \xcp $ $( \xbm Cum \xba )$ for
all $ \xbb \xck \xba $

(3.2) $( \xbm Cumt \xbb )$ $+$ $( \xcv )$ $ \xcp $ $( \xbm Cumt \xba )$
for all $ \xbb \xck \xba $

(3.3) $\{( \xbm Cumt \xbb ): \xbb < \xba \}$ $+$ $( \xbm CUM)$ $+$ $( \xcS
)$ $ \xcP $ $( \xbm Cum \xba )$ for $ \xba >0.$

(4) Connection $( \xbm Cum \xba )/( \xbm Cumt \xba )$:

(4.1) $( \xbm Cumt \xba )$ $ \xcp $ $( \xbm Cum \xba )$

(4.2) $( \xbm Cum \xba )$ $+$ $( \xcS )$ $ \xcP $ $( \xbm Cumt \xba )$

(4.3) $( \xbm Cum \xba )$ $+$ $( \xcv )$ $ \xcp $ $( \xbm Cumt \xba )$

(5) $( \xbm CUM)$ and $( \xbm Cumi)$:

(5.1) $( \xbm CUM)$ $+$ $( \xcv )$ entail:

(5.1.1) $ \xbm (A) \xcc B$ $ \xcp $ $ \xbm (A \xcv B)= \xbm (B)$

(5.1.2) $ \xbm (X) \xcc U,$ $U \xcc Y$ $ \xcp $ $ \xbm (Y \xcv X)= \xbm
(Y)$

(5.1.3) $ \xbm (X) \xcc U,$ $U \xcc Y$ $ \xcp $ $ \xbm (Y) \xcs X \xcc
\xbm (U)$

(5.2) $( \xbm Cum \xba )$ $ \xcp $ $( \xbm CUM)$ for all $ \xba $

(5.3) $( \xbm CUM)$ $+$ $( \xcv )$ $ \xcp $ $( \xbm Cum \xba )$ for all $
\xba $

(5.4) $( \xbm CUM)$ $+$ $( \xcs )$ $ \xcp $ $( \xbm Cum0)$

(6) $( \xbm CUM)$ and $( \xbm Cumt \xba )$:

(6.1) $( \xbm Cumt \xba )$ $ \xcp $ $( \xbm CUM)$ for all $ \xba $

(6.2) $( \xbm CUM)$ $+$ $( \xcv )$ $ \xcp $ $( \xbm Cumt \xba )$ for all $
\xba $

(6.3) $( \xbm CUM)$ $ \xcP $ $( \xbm Cumt \xba )$ for all $ \xba >0$

(7) $( \xbm Cum0)$ $ \xcp $ $( \xbm PR)$

 karl-search= End Fact Cum-Alpha
\vspace{7mm}

 *************************************

\vspace{7mm}

\subsubsection{Fact Cum-Alpha Proof}

 {\LARGE karl-search= Start Fact Cum-Alpha Proof }

\index{Fact Cum-Alpha Proof}

\efa

\subparagraph{
Proof
}

$\hspace{0.01em}$

(+++*** Orig.:  Proof )

We prove these facts in a different order: (1), (2), (5.1), (5.2), (4.1),
(6.1),
(6.2), (5.3), (3.1), (3.2), (4.2), (4.3), (5.4), (3.3), (6.3), (7).

(1.1)

For $ \xbb < \xbg \xck \xba $ set $X_{ \xbg }:=X_{ \xbb }.$ Let the
prerequisites of $( \xbm Cum \xbb )$ hold. Then for
$ \xbg $ with $ \xbb < \xbg \xck \xba $ $ \xbm (X_{ \xbg }) \xcc X_{ \xbb
}$ by $( \xbm \xcc ),$ so the prerequisites
of $( \xbm Cum \xba )$ hold, too, so by $( \xbm Cum \xba )$ $ \xcS \{X_{
\xbd }: \xbd \xck \xbb \} \xcs \xbm (U)$ $=$
$ \xcS \{X_{ \xbd }: \xbd \xck \xba \} \xcs \xbm (U)$ $ \xcc $ $ \xbm (X_{
\xba })$ $=$ $ \xbm (X_{ \xbb }).$

(1.2)

Analogous.

(2.1)

Proof by induction.

$( \xbm Cum0)$ Let $ \xbm (X_{0}) \xcc U,$ suppose there is $x \xbe \xbm
(U) \xcs (X_{0}- \xbm (X_{0})).$ By smoothness,
there is $y \xeb x,$ $y \xbe \xbm (X_{0}) \xcc U,$ $contradiction$ (The
same arguments works for copies: all
copies of $x$ must be minimized by some $y \xbe \xbm (X_{0}),$ but at
least one copy of $x$
has to be minimal in $U.)$

Suppose $( \xbm Cum \xbb )$ hold for all $ \xbb < \xba.$ We show $( \xbm
Cum \xba ).$ Let the prerequisites
of $( \xbm Cum \xba )$ hold, then those for $( \xbm Cum \xbb ),$ $ \xbb <
\xba $ hold, too. Suppose there is
$x \xbe \xbm (U) \xcs \xcS \{X_{ \xbg }: \xbg \xck \xba \}- \xbm (X_{ \xba
}).$ So by $( \xbm Cum \xbb )$ for $ \xbb < \xba $ $x \xbe \xbm (X_{ \xbb
})$
moreover $x \xbe \xbm (U).$ By smoothness, there is $y \xbe \xbm (X_{ \xba
}) \xcc U \xcv \xcV \{X_{ \xbb }: \xbb < \xba \},$ $y \xeb x,$
but this is a contradiction. The same argument works again for copies.

(2.2)

We use the following Fact:
Let, in a smooth transitive structure, $ \xbm (X_{ \xbb })$ $ \xcc $ $U
\xcv \xcV \{X_{ \xbg }: \xbg < \xbb \}$ for all $ \xbb \xck \xba,$
and let $x \xbe \xbm (U).$ Then there is no $y \xeb x,$ $y \xbe U \xcv
\xcV \{X_{ \xbg }: \xbg \xck \xba \}.$

Proof of the Fact by induction:
$ \xba =0:$ $y \xbe U$ is impossible: if $y \xbe X_{0},$ then if $y \xbe
\xbm (X_{0}) \xcc U,$ which is impossible,
or there is $z \xbe \xbm (X_{0}),$ $z \xeb y,$ so $z \xeb x$ by
transitivity, but $ \xbm (X_{0}) \xcc U.$
Let the result hold for all $ \xbb < \xba,$ but fail for $ \xba,$
so $ \xCN \xcE y \xeb x.y \xbe U \xcv \xcV \{X_{ \xbg }: \xbg < \xba \},$
but $ \xcE y \xeb x.y \xbe U \xcv \xcV \{X_{ \xbg }: \xbg \xck \xba \},$
so $y \xbe X_{ \xba }.$
If $y \xbe \xbm (X_{ \xba }),$ then $y \xbe U \xcv \xcV \{X_{ \xbg }: \xbg
< \xba \},$ but this is impossible, so
$y \xbe X_{ \xba }- \xbm (X_{ \xba }),$ let by smoothness $z \xeb y,$ $z
\xbe \xbm (X_{ \xba }),$ so by transitivity $z \xeb x,$ $contradiction.$
The result is easily modified for the case with copies.

Let the prerequisites of $( \xbm Cumt \xba )$ hold, then those of the Fact
will hold,
too. Let now $x \xbe \xbm (U) \xcs (X_{ \xba }- \xbm (X_{ \xba })),$ by
smoothness, there must be $y \xeb x,$
$y \xbe \xbm (X_{ \xba }) \xcc U \xcv \xcV \{X_{ \xbg }: \xbg < \xba \},$
contradicting the Fact.

(2.3)

Let $ \xba >0,$ and consider the following structure over $\{a,b,c\}:$
$U:=\{a,c\},$
$X_{0}:=\{b,c\},$ $X_{ \xba }:= \Xl:=X_{1}:=\{a,b\},$ and their
intersections, $\{a\},$ $\{b\},$ $\{c\},$ $ \xCQ $ with
the order $c \xeb b \xeb a$ (without transitivity). This is preferential,
so $( \xbm PR)$ and
$( \xbm \xcc )$ hold.
The structure is smooth for $U,$ all $X_{ \xbb },$ and their
intersections.
We have $ \xbm (X_{0}) \xcc U,$ $ \xbm (X_{ \xbb }) \xcc U \xcv X_{0}$ for
all $ \xbb \xck \xba,$ so $ \xbm (X_{ \xbb }) \xcc U \xcv \xcV \{X_{ \xbg
}: \xbg < \xbb \}$
for all $ \xbb \xck \xba $ but $X_{ \xba } \xcs \xbm (U)=\{a\} \xcC \{b\}=
\xbm (X_{ \xba })$ for $ \xba >0.$

(5.1)

(5.1.1) $ \xbm (A) \xcc B$ $ \xcp $ $ \xbm (A \xcv B) \xcc \xbm (A) \xcv
\xbm (B) \xcc B$ $ \xcp_{( \xbm CUM)}$ $ \xbm (B)= \xbm (A \xcv B).$

(5.1.2) $ \xbm (X) \xcc U \xcc Y$ $ \xcp $ (by (1)) $ \xbm (Y \xcv X)=
\xbm (Y).$

(5.1.3) $ \xbm (Y) \xcs X$ $=$ (by (2)) $ \xbm (Y \xcv X) \xcs X$ $ \xcc $
$ \xbm (Y \xcv X) \xcs (X \xcv U)$ $ \xcc $ (by $( \xbm PR))$
$ \xbm (X \xcv U)$ $=$ (by (1)) $ \xbm (U).$

(5.2)

Using (1.1), it suffices to show $( \xbm Cum0)$ $ \xcp $ $( \xbm CUM).$
Let $ \xbm (X) \xcc U \xcc X.$ By $( \xbm Cum0)$ $X \xcs \xbm (U) \xcc
\xbm (X),$ so by $ \xbm (U) \xcc U \xcc X$ $ \xcp $ $ \xbm (U) \xcc \xbm
(X).$
$U \xcc X$ $ \xcp $ $ \xbm (X) \xcs U \xcc \xbm (U),$ but also $ \xbm (X)
\xcc U,$ so $ \xbm (X) \xcc \xbm (U).$

(4.1)

Trivial.

(6.1)

Follows from (4.1) and (5.2).

(6.2)

Let the prerequisites of $( \xbm Cumt \xba )$ hold.

We first show by induction $ \xbm (X_{ \xba } \xcv U) \xcc \xbm (U).$

Proof:

$ \xba =0:$ $ \xbm (X_{0}) \xcc U$ $ \xcp $ $ \xbm (X_{0} \xcv U)= \xbm
(U)$ by (5.1.1).
Let for all $ \xbb < \xba $ $ \xbm (X_{ \xbb } \xcv U) \xcc \xbm (U) \xcc
U.$ By prerequisite,
$ \xbm (X_{ \xba }) \xcc U \xcv \xcV \{X_{ \xbb }: \xbb < \xba \},$ thus $
\xbm (X_{ \xba } \xcv U)$ $ \xcc $ $ \xbm (X_{ \xba }) \xcv \xbm (U)$ $
\xcc $ $ \xcV \{U \xcv X_{ \xbb }: \xbb < \xba \},$

so $ \xcA \xbb < \xba $ $ \xbm (X_{ \xba } \xcv U) \xcs (U \xcv X_{ \xbb
})$ $ \xcc $ $ \xbm (U)$ by (5.1.3), thus $ \xbm (X_{ \xba } \xcv U) \xcc
\xbm (U).$

Consequently, under the above prerequisites, we have $ \xbm (X_{ \xba }
\xcv U)$ $ \xcc $ $ \xbm (U)$ $ \xcc $
$U$ $ \xcc $ $U \xcv X_{ \xba },$ so by $( \xbm CUM)$ $ \xbm (U)= \xbm
(X_{ \xba } \xcv U),$ and, finally,
$ \xbm (U) \xcs X_{ \xba }= \xbm (X_{ \xba } \xcv U) \xcs X_{ \xba } \xcc
\xbm (X_{ \xba })$ by $( \xbm PR).$

Note that finite unions take us over the limit step, essentially, as all
steps collapse, and $ \xbm (X_{ \xba } \xcv U)$ will always be $ \xbm
(U),$ so there are no real
changes.

(5.3)

Follows from (6.2) and (4.1).

(3.1)

Follows from (5.2) and (5.3).

(3.2)

Follows from (6.1) and (6.2).

(4.2)

Follows from (2.3) and (2.1).

(4.3)

Follows from (5.2) and (6.2).

(5.4)

$ \xbm (X) \xcc U$ $ \xcp $ $ \xbm (X) \xcc U \xcs X \xcc X$ $ \xcp $ $
\xbm (X \xcs U)= \xbm (X)$ $ \xcp $
$X \xcs \xbm (U)=(X \xcs U) \xcs \xbm (U) \xcc \xbm (X \xcs U)= \xbm (X)$

(3.3)

See Example \ref{Example Inf-Cum-Alpha} (page \pageref{Example Inf-Cum-Alpha}) .

(6.3)

This is a consequence of (3.3).

(7)

Trivial. Let $X \xcc Y,$ so by $( \xbm \xcc )$ $ \xbm (X) \xcc X \xcc Y,$
so by $( \xbm Cum0)$ $X \xcs \xbm (Y) \xcc \xbm (X).$

$ \xcz $
\\[3ex]

 karl-search= End Fact Cum-Alpha Proof
\vspace{7mm}

 *************************************

\vspace{7mm}

\subsubsection{Fact Cum-Alpha-HU}

 {\LARGE karl-search= Start Fact Cum-Alpha-HU }

\index{Fact Cum-Alpha-HU}

\bfa

$\hspace{0.01em}$

(+++ Orig. No.:  Fact Cum-Alpha-HU +++)

{\xssc LABEL: {Fact Cum-Alpha-HU}}
\label{Fact Cum-Alpha-HU}

Assume $( \xbm \xcc ).$

We have for $( \xbm Cum \xca )$ and $ \xCf (HU,u)$:

(1) $x \xbe \xbm (Y),$ $ \xbm (Y) \xcc H(U,x)$ $ \xcp $ $Y \xcc H(U,x)$

(2) $( \xbm Cum \xca )$ $ \xcp $ $ \xCf (HU,u)$

(3) $ \xCf (HU,u)$ $ \xcp $ $( \xbm Cum \xca )$

 karl-search= End Fact Cum-Alpha-HU
\vspace{7mm}

 *************************************

\vspace{7mm}

\subsubsection{Fact Cum-Alpha-HU Proof}

 {\LARGE karl-search= Start Fact Cum-Alpha-HU Proof }

\index{Fact Cum-Alpha-HU Proof}

\efa

\bfa

$\hspace{0.01em}$

(+++ Orig. No.:  Fact Cum-Alpha-HU Proof +++)

{\xssc LABEL: {Fact Cum-Alpha-HU Proof}}
\label{Fact Cum-Alpha-HU Proof}

(1)

Trivial by definition of $H(U,x).$

(2)

Let $x \xbe \xbm (U),$ $x \xbe Y,$ $ \xbm (Y) \xcc H(U,x)$ (and thus $Y
\xcc H(U,x)$ by definition).
Thus, we have a sequence $X_{0}:=U,$ $ \xbm (X_{ \xbb }) \xcc U \xcv \xcV
\{X_{ \xbg }: \xbg < \xbb \}$ and $Y=X_{ \xba }$ for some $ \xba $
(after $X_{0},$ enumerate arbitrarily $H(U,x)_{1},$ then $H(U,x)_{2},$
etc., do nothing at
limits). So $x \xbe \xcS \{X_{ \xbg }: \xbg \xck \xba \} \xcs \xbm (U),$
and $x \xbe \xbm (X_{ \xba })= \xbm (Y)$ by $( \xbm Cum \xca ).$
Remark: The same argument shows that we can replace `` $x \xbe X$ ''
equivalently by
`` $x \xbe \xbm (X)$ '' in the definition of $H(U,x)_{ \xba +1},$ as was
done in Definition 3.7.5
in  \cite{Sch04}.

(3)

Suppose $( \xbm Cum \xba )$ fails, we show that then so does $ \xCf
(HUx).$ As $( \xbm Cum \xba )$ fails, for
all $ \xbb \xck \xba $ $ \xbm (X_{ \xbb }) \xcc U \xcv \xcV \{X_{ \xbg }:
\xbg < \xbb \},$ but there is $x \xbe \xcS \{X_{ \xbg }: \xbg \xck \xba \}
\xcs \xbm (U),$
$x \xce \xbm (X_{ \xba }).$ Thus for all $ \xbb \xck \xba $ $ \xbm (X_{
\xbb }) \xcc X_{ \xbb } \xcc H(U,x),$ moreover $x \xbe \xbm (U),$
$x \xbe X_{ \xba }- \xbm (X_{ \xba }),$ but $ \xbm (X_{ \xba }) \xcc
H(U,x),$ so $ \xCf (HUx)$ fails.

$ \xcz $
\\[3ex]

 karl-search= End Fact Cum-Alpha-HU Proof
\vspace{7mm}

 *************************************

\vspace{7mm}

\subsubsection{Fact HU}

 {\LARGE karl-search= Start Fact HU }

\index{Fact HU}

\efa

\bfa

$\hspace{0.01em}$

(+++ Orig. No.:  Fact HU +++)

{\xssc LABEL: {Fact HU}}
\label{Fact HU}

We continue to show results for $H(U)$ and $H(U,u).$

Let A, $X,$ $U,$ $U',$ $Y$ and all $A_{i}$ be in $ \xdy.$

(0) $H(U)$ and $H(U,u)$

(0.1) $H(U,u) \xcc H(U)$

(0.2) $(HU) \xcp (HU,u)$

(0.3) $( \xcv )$ $+$ $( \xbm PR)$ entail $H(U) \xcc H(U,u)$

(0.4) $( \xcv )$ $+$ $( \xbm PR)$ entail $(HU,u) \xcp (HU)$

(1) $( \xbm \xcc )$ and $ \xCf (HU)$ entail:

(1.1) $( \xbm PR)$

(1.2) $( \xbm CUM)$

(2) $ \xCf (HU)$ $+$ $( \xcv )$ $ \xcp $ $ \xCf (HU,u)$

(3) $( \xbm \xcc )$ and $( \xbm PR)$ entail:

(3.1) $A= \xcV \{A_{i}:i \xbe I\}$ $ \xcp $ $ \xbm (A) \xcc \xcV \{ \xbm
(A_{i}):i \xbe I\},$

(3.2) $U \xcc H(U),$ and $U \xcc U' \xcp H(U) \xcc H(U' ),$

(3.3) $ \xbm (U \xcv Y)-H(U) \xcc \xbm (Y)$ - if $ \xbm (U \xcv Y)$ is
defined, in particular, if $( \xcv )$
holds.

(4) $( \xcv ),$ $( \xbm \xcc ),$ $( \xbm PR),$ $( \xbm CUM)$ entail:

(4.1) $H(U)=H_{1}(U)$

(4.2) $U \xcc A,$ $ \xbm (A) \xcc H(U)$ $ \xcp $ $ \xbm (A) \xcc U,$

(4.3) $ \xbm (Y) \xcc H(U)$ $ \xcp $ $Y \xcc H(U)$ and $ \xbm (U \xcv Y)=
\xbm (U),$

(4.4) $x \xbe \xbm (U),$ $x \xbe Y- \xbm (Y)$ $ \xcp $ $Y \xcC H(U)$ (and
thus $ \xCf (HU)),$

(4.5) $Y \xcC H(U)$ $ \xcp $ $ \xbm (U \xcv Y) \xcC H(U).$

(5) $( \xcv ),$ $( \xbm \xcc ),$ $ \xCf (HU)$ entail

(5.1) $H(U)=H_{1}(U)$

(5.2) $U \xcc A,$ $ \xbm (A) \xcc H(U)$ $ \xcp $ $ \xbm (A) \xcc U,$

(5.3) $ \xbm (Y) \xcc H(U)$ $ \xcp $ $Y \xcc H(U)$ and $ \xbm (U \xcv Y)=
\xbm (U),$

(5.4) $x \xbe \xbm (U),$ $x \xbe Y- \xbm (Y)$ $ \xcp $ $Y \xcC H(U),$

(5.5) $Y \xcC H(U)$ $ \xcp $ $ \xbm (U \xcv Y) \xcC H(U).$

 karl-search= End Fact HU
\vspace{7mm}

 *************************************

\vspace{7mm}

\subsubsection{Fact HU Proof}

 {\LARGE karl-search= Start Fact HU Proof }

\index{Fact HU Proof}

\efa

\bfa

$\hspace{0.01em}$

(+++ Orig. No.:  Fact HU Proof +++)

{\xssc LABEL: {Fact HU Proof}}
\label{Fact HU Proof}

(0.1) and (0.2) trivial by definition.

(0.3) Proof by induction. Let $X \xbe \xdy,$ $ \xbm (X) \xcc H(U)_{ \xba
},$ then $U \xcv X \xbe \xdy,$
$ \xbm (U \xcv X) \xcc_{( \xbm PR)} \xbm (U) \xcv \xbm (X) \xcc_{( \xbm
\xcc )}H(U)_{ \xba }=H(U,u)_{ \xba }$ by induction hypothesis,
and $u \xbe U \xcv X.$

(0.4) Immediate by (0.3).

(1.1) By $ \xCf (HU),$ if $ \xbm (Y) \xcc H(U),$ then $ \xbm (U) \xcs Y
\xcc \xbm (Y).$ But, if $Y \xcc U,$ then
$ \xbm (Y) \xcc H(U)$ by $( \xbm \xcc ).$

(1.2) Let $ \xbm (U) \xcc X \xcc U.$ Then by (1.1) $ \xbm (U)= \xbm (U)
\xcs X \xcc \xbm (X).$ By prerequisite,
$ \xbm (U) \xcc U \xcc H(X),$ so $ \xbm (X)= \xbm (X) \xcs U \xcc \xbm
(U)$ by $( \xbm \xcc ).$

(2) By (1.2), $ \xCf (HU)$ entails $( \xbm CUM),$ so by $( \xcv )$ and
Fact \ref{Fact Cum-Alpha} (page \pageref{Fact Cum-Alpha}) , (5.2) $( \xbm Cum
\xca )$
holds, so by Fact \ref{Fact HU} (page \pageref{Fact HU}) , (2) $ \xCf (HUx)$
holds.

(3.1) $ \xbm (A) \xcs A_{j} \xcc \xbm (A_{j}) \xcc \xcV \xbm (A_{i}),$ so
by $ \xbm (A) \xcc A= \xcV A_{i}$ $ \xbm (A) \xcc \xcV \xbm (A_{i}).$

(3.2) trivial.

(3.3) $ \xbm (U \xcv Y)-H(U)$ $ \xcc_{(3.2)}$ $ \xbm (U \xcv Y)-U$ $ \xcc
$ (by $( \xbm \xcc )$ and (3.1))
$ \xbm (U \xcv Y) \xcs Y$ $ \xcc_{( \xbm PR)}$ $ \xbm (Y).$

(4.1) We show that, if $X \xcc H_{2}(U),$ then $X \xcc H_{1}(U),$ more
precisely, if $ \xbm (X) \xcc H_{1}(U),$
then already $X \xcc H_{1}(U),$ so the construction stops already at
$H_{1}(U).$
Suppose then $ \xbm (X) \xcc \xcV \{Y: \xbm (Y) \xcc U\},$ and let $A:=X
\xcv U.$ We show that $ \xbm (A) \xcc U,$ so
$X \xcc A \xcc H_{1}(U).$ Let $a \xbe \xbm (A).$ By (3.1), $ \xbm (A) \xcc
\xbm (X) \xcv \xbm (U).$ If $a \xbe \xbm (U) \xcc U,$ we
are done. If $a \xbe \xbm (X),$ there is $Y$ s.t. $ \xbm (Y) \xcc U$ and
$a \xbe Y,$ so $a \xbe \xbm (A) \xcs Y.$
By Fact \ref{Fact Cum-Alpha} (page \pageref{Fact Cum-Alpha}) , (5.1.3), we have
for $Y$ s.t.
$ \xbm (Y) \xcc U$ and $U \xcc A$ $ \xbm (A) \xcs Y \xcc \xbm (U).$ Thus
$a \xbe \xbm (U),$ and we are done again.

(4.2) Let $U \xcc A,$ $ \xbm (A) \xcc H(U)=H_{1}(U)$ by (4.1). So $ \xbm
(A)$ $=$ $ \xcV \{ \xbm (A) \xcs Y: \xbm (Y) \xcc U\}$ $ \xcc $
$ \xbm (U)$ $ \xcc $ $U,$ again by Fact \ref{Fact Cum-Alpha} (page \pageref{Fact
Cum-Alpha}) ,
(5.1.3).

(4.3) Let $ \xbm (Y) \xcc H(U),$ then by $ \xbm (U) \xcc H(U)$ and (3.1)
$ \xbm (U \xcv Y) \xcc \xbm (U) \xcv \xbm (Y) \xcc H(U),$ so by (4.2) $
\xbm (U \xcv Y) \xcc U$ and $U \xcv Y \xcc H(U).$
Moreover, $ \xbm (U \xcv Y) \xcc U \xcc U \xcv Y$ $ \xcp_{( \xbm CUM)}$ $
\xbm (U \xcv Y)= \xbm (U).$

(4.4) If not, $Y \xcc H(U),$ so $ \xbm (Y) \xcc H(U),$ so $ \xbm (U \xcv
Y)= \xbm (U)$ by (4.3),
but $x \xbe Y- \xbm (Y)$ $ \xcp_{( \xbm PR)}$ $x \xce \xbm (U \xcv Y)=
\xbm (U),$ $contradiction.$

(4.5) $ \xbm (U \xcv Y) \xcc H(U)$ $ \xcp_{(4.3)}$ $U \xcv Y \xcc H(U).$

(5) Trivial by (1) and (4).

$ \xcz $
\\[3ex]

 karl-search= End Fact HU Proof
\vspace{7mm}

 *************************************

\vspace{7mm}

\section{
Validity
}

\efa

$++++++++++++++++++++++++++++++++++++++++++++++$
\subsubsection{Introduction Path-Validity}

 {\LARGE karl-search= Start Introduction Path-Validity }

\index{Introduction Path-Validity}

\paragraph{
Introduction to Path-Validity
}

$\hspace{0.01em}$

(+++*** Orig.:  Introduction to Path-Validity )

{\xssc LABEL: {Section Introduction to Path-Validity}}
\label{Section Introduction to Path-Validity}

All definitions are relative to a fixed diagram $ \xbG.$

For simplicity, we consider $ \xbG $ to be just a set of points and
arrows, thus e.g. $x \xcp y \xbe \xbG $ and $x \xbe \xbG $ are defined,
when $x$ is a point in $ \xbG,$ and
$x \xcp y$ an arrow in $ \xbG.$

Recall that we have two types of arrows, positive and
negative ones.

We first define generalized and potential paths, and
finally validity of paths, written $ \xbG \xcm \xbs,$ if $ \xbs $ is a
path, as well as $ \xbG \xcm xy,$
if $ \xbG \xcm \xbs $ and $ \xbs:x \Xl. \xcp y.$

 karl-search= End Introduction Path-Validity
\vspace{7mm}

 *************************************

\vspace{7mm}

$++++++++++++++++++++++++++++++++++++++++++++++$
\subsubsection{Definition Gen-Path}

 {\LARGE karl-search= Start Definition Gen-Path }

\index{Definition Gen-Path}

\bd

$\hspace{0.01em}$

(+++ Orig. No.:  Definition Gen-Path +++)

{\xssc LABEL: {Definition Gen-Path}}
\label{Definition Gen-Path}

 \xEh

 \xDH Generalized paths:

A generalized path is an uninterrupted chain of positive or negative
arrows
pointing in the same direction, more precisely:

 \xEh
 \xDH The empty path is a generalized path.
 \xDH If $x \xcp p \xbe \xbG,$ then $x \xcp p$ is a generalized path,
 \xDH if $x \xcP p \xbe \xbG,$ then $x \xcP p$ is a generalized path.
 \xDH If $x \xFB \xcp p$ is a generalized path, and $p \xcp q \xbe \xbG $
, then
$x \xFB \xcp p \xcp q$ is a generalized path,
 \xDH if $x \xFB \xcp p$ is a generalized path, and $p \xcP q \xbe \xbG $
, then
$x \xFB \xcp p \xcP q$ is a generalized path.
 \xDH If the starting point of a generalized path $ \xbs $ is $x,$ and $y$
its endpoint,
we say that $ \xbs $ is a generalized path from $x$ to $y,$ and write $
\xbs:x \xfq y$
 \xEj

 \xDH Concatenation:

If $ \xbs $ and $ \xbt $ are two generalized paths, and the end point of $
\xbs $ is the same
as the starting point of $ \xbt,$ then $ \xbs \xDM \xbt $ is the
concatenation of $ \xbs $ and $ \xbt.$

 \xDH Subpath:

If $ \xbs = \xbt \xDM \xbt ' \xDM \xbt '' $ is a generalized path, $ \xbt
$ and $ \xbt '' $ are generalized paths
(possibly empty), then $ \xbt ' $ is a subpath of $ \xbs.$

 \xDH $[x,y]:$

If $x,y$ are nodes in $ \xbG,$ then $[x,y]$ is the set of all subpaths of
all generalized
paths from $x$ to $y.$ Note that $ \xcc $ is a well-founded relation on
the set of
$[x,y]$ of $ \xbG,$ so we can do induction on $[x,y]$ and $ \xcc.$

 \xDH Potential paths (pp.):

A generalized path, which contains at most one negative arrow, and then at
the
end, is a potential path. If the last link is positive, it is a positive
potential path, if not, a negative one.

 \xEj

 karl-search= End Definition Gen-Path
\vspace{7mm}

 *************************************

\vspace{7mm}

\ed

$++++++++++++++++++++++++++++++++++++++++++++++$
\subsubsection{Definition Arrow-Origin}

 {\LARGE karl-search= Start Definition Arrow-Origin }

\index{Definition Arrow-Origin}

\bd

$\hspace{0.01em}$

(+++ Orig. No.:  Definition Arrow-Origin +++)

{\xssc LABEL: {Definition Arrow-Origin}}
\label{Definition Arrow-Origin}

This definition is for IBRS - otherwise it is trivial.

The definition is by recursion. Intuitively, we go back until we find a
node.

If $ \xba:x \xcp y$ is an arrow from node to node or arrow, then $or(
\xba ):=x.$

If $ \xbb: \xba \xcp y$ is an arrow from arrow to node or arrow, then
$or( \xbb ):=or( \xba ).$

 karl-search= End Definition Arrow-Origin
\vspace{7mm}

 *************************************

\vspace{7mm}

\ed

$++++++++++++++++++++++++++++++++++++++++++++++$
\subsubsection{Definition Path-Validity}

 {\LARGE karl-search= Start Definition Path-Validity }

\index{Definition Path-Validity}

\bd

$\hspace{0.01em}$

(+++ Orig. No.:  Definition Path-Validity +++)

{\xssc LABEL: {Definition Path-Validity}}
\label{Definition Path-Validity}

Inductive definition of $ \xbG \xcm \xbs $ or validity of path. At the
same time,
we construct dynamically an IBRS - which is just a reformulation of the
same mechanism.

Let $ \xbs:x \xfq y$ be a potential path, and let validity,
as well as the construction of new arrows in the IBRS,
be determined by induction for all $ \xbs ':x' \xfq y' $ with $[x',y' ]
\xcb [x,y].$

 \xEh

 \xDH Case $I,$ $ \xbs:x \xcp y$ (or $x \xcP y)$ is a direct link in $
\xbG:$

Then $ \xbG \xcm \xbs,$ and we add the arrow $ \xba:x \xcp y,$ with two
labels: $v$ for validity,
and $+/-$ if $ \xbs:x \xcp y$ or $ \xbs:x \xcP y$ - we denote this arrow
$ \xba_{ \xbs }.$

(Recall that we have no hard contradictions in $ \xbG.)$

 \xDH Case II, $ \xbs $ is a compound potential path:

 \xEh

 \xDH Case II.1, $ \xbs $ is a positive pp. $x \xfq u \xcp y:$

Let $ \xbs ':=x \xfq u,$ so $ \xbs = \xbs ' \xDM u \xcp y.$

Then, intuitively, $ \xbG \xcm \xbs $ iff

 \xEh
 \xDH $ \xbs $ is a candidate by upward chaining,
 \xDH $ \xbs $ is not precluded by more specific contradicting
information,
 \xDH all potential contradictions are themselves precluded by information
contradicting them.

 \xEj

Formally,

$ \xbG \xcm \xbs $ and we add an arrow $ \xba_{ \xbs }: \xba_{ \xbs ' }
\xcp y$ with labels $v$ and $+$

iff

 \xEh

 \xDH $ \xbG \xcm \xbs ' $ and $u \xcp y \xbe \xbG.$

For IBRS, the prerequisite is that there is an
arrow $ \xba_{ \xbs ' }$ s.t. $or( \xba )=x,$ the destination of $ \xba $
is $y,$ the labels of $ \xba $ are
$v$ and $+.$

(The initial segment must be a path, as we have an upward chaining
approach.
This is decided by the induction hypothesis.)

 \xDH There are no $v,$ $ \xbt,$ $ \xbt ' $ such that $v \xcP y \xbe \xbG
$ and $ \xbG \xcm \xbt:=x \xfq v$ and
$ \xbG \xcm \xbt ':=v \xfq u$ - there are arrows $ \xba_{ \xbt }$ with
$or( \xba_{ \xbt })=x,$ destination $v$
and $ \xba_{ \xbt ' }$ with $or( \xba_{ \xbt ' })=$ $v,$ destination $u$
for the IBRS
$( \xbt $ may be the empty path, i.e. $x=v.)$

$( \xbs $ itself is not precluded by split validity preclusion and a
contradictory
link. Note that $ \xbt \xDM v \xcP y$ need not be valid, it suffices
that it is a better candidate (by $ \xbt ' ).)$

 \xDH all potentially conflicting paths are precluded by information
contradicting them:

For all $v$ and $ \xbt $ such that $v \xcP y \xbe \xbG $ and $ \xbG \xcm
\xbt:=x \xfq v$ (i.e. for all potentially
conflicting paths $ \xbt \xDM v \xcP y)$ - $ \xba_{ \xbt }$ with $or(
\xba_{ \xbt })=x$ and destination $v$ -
there is $z$ such that $z \xcp y \xbe \xbG $ and either

$z=x$

(the potentially conflicting pp. is itself precluded by a direct link,
which is
thus valid)

or

there are $ \xbG \xcm \xbr:=x \xfq z$ and $ \xbG \xcm \xbr ':=z \xfq v$
for suitable $ \xbr $ and $ \xbr ' $ -
$ \xba_{ \xbr }$ and $ \xba_{ \xbr ' }$ with suitable origin and
destination for IBRS.

 \xEj

 \xEh

 \xDH Case II.2, the negative case, i.e.
$ \xbs $ a negative pp. $x \xfq u \xcP y,$ $ \xbs ':=x \xfq u,$ $ \xbs =
\xbs ' \xDM u \xcP y$
is entirely symmetrical.

 \xEj

 \xEj

 \xEj

\ed

Note that the new arrows $ \xba $ allow us to reconstruct the whole path,
if needed.

 karl-search= End Definition Path-Validity
\vspace{7mm}

 *************************************

\vspace{7mm}

$++++++++++++++++++++++++++++++++++++++++++++++$
\subsubsection{Remark Path-Validity}

 {\LARGE karl-search= Start Remark Path-Validity }

\index{Remark Path-Validity}

\br

$\hspace{0.01em}$

(+++ Orig. No.:  Remark Path-Validity +++)

{\xssc LABEL: {Remark Path-Validity}}
\label{Remark Path-Validity}

The following remarks all concern preclusion.

(1) Thus, in the case of preclusion, there is a valid path from $x$ to
$z,$
and $z$ is more specific than $v,$ so
$ \xbt \xDM v \xcP y$ is precluded. Again, $ \xbr \xDM z \xcp y$ need not
be a valid path, but it is
a better candidate than $ \xbt \xDM v \xcP y$ is, and as $ \xbt \xDM v
\xcP y$ is in simple contradiction,
this suffices.

(2) Our definition is stricter than many usual ones, in the
following sense: We require - according to our general picture to treat
only
direct links as information - that the preclusion ``hits'' the precluded
path
at the end, i.e. $v \xcP y \xbe \xbG,$ and $ \xbr ' $ hits $ \xbt \xDM v
\xcP y$ at $v.$ In other definitions,
it is possible that the preclusion hits at some $v',$ which is somewhere
on the
path $ \xbt,$ and not necessarily at its end. For instance, in the Tweety
Diagram,

\begin{verbatim}

  *****************************

  **  Index unter Hauptteil  **

  *****************************

\end{verbatim}

$see$ Diagram \ref{Diagram Tweety} (page \pageref{Diagram Tweety}) , if there
were a node $b' $ between
$b$ and $d,$ we will need
path $c \xcp b \xcp b' $ to be valid, (obvious) validity of the arrow $c
\xcp b$ will not
suffice.

(3) If we allow $ \xbr $ to be the empty path, then the case $z=x$ is a
subcase of the
present one.

(4) Our conceptual analysis has led to a very important simplification of
the
definition of validity. If we adopt on-path preclusion, we have to
remember
all paths which led to the information source to be considered: In the
Tweety
diagram, we have to remember that there is an arrow $a \xcp b,$ it is not
sufficient
to note that we somehow came from a to $b$ by a valid path, as the path $a
\xcp c \xcp b \xcp d$
is precluded, but not the path $a \xcp b \xcp d.$ If we adopt total path
preclusion, we
have to remember the valid path $a \xcp c \xcp b$ to see that it precludes
$a \xcp c \xcp d.$ If
we allow preclusion to ``hit'' below the last node, we also have to remember
the entire path which is precluded. Thus, in all those cases, whole paths
(which can be very long) have to be remembered, but NOT in our definition.

We only need to remember (consider the Tweety diagram):

(a) we want to know if $a \xcp b \xcp d$ is valid, so we have to remember
a, $b,$ $d.$
Note that the (valid) path from a to $b$ can be composed and very long.

(b) we look at possible preclusions, so we have to remember $a \xcp c \xcP
d,$ again
the (valid) path from a to $c$ can be very long.

(c) we have to remember that the path from $c$ to $b$ is valid (this was
decided
by induction before).

So in all cases (the last one is even simpler), we need only remember the
starting node, a (or $c),$ the last node of the valid paths, $b$ (or $c),$
and the
information $b \xcp d$ or $c \xcP d$ - i.e. the size of what has to be
recalled is
$ \xck 3.$ (Of course, there may be many possible preclusions, but in all
cases we have to look at a very limited situation, and not arbitrarily
long
paths.)

$ \xcz $
\\[3ex]

 karl-search= End Remark Path-Validity
\vspace{7mm}

 *************************************

\vspace{7mm}

\newpage
\subsubsection{Diagram I-U-reac-x}
\subsubsection{Diagram I-U-reac-x}

 {\LARGE karl-search= Start Diagram I-U-reac-x }

\vspace{10mm}

\begin{diagram}

{\xssc LABEL: {Diagram I-U-reac-x}}
\label{Diagram I-U-reac-x}
\index{Diagram I-U-reac-x}

\centering
\setlength{\unitlength}{1mm}
{\renewcommand{\dashlinestretch}{30}
\begin{picture}(150,150)(0,0)

\put(10,60){\circle*{1}}
\put(70,60){\circle*{1}}
\put(130,60){\circle*{1}}

\put(70,10){\circle*{1}}
\put(70,110){\circle*{1}}

\put(40,50){\circle*{1}}
\put(40,70){\circle*{1}}
\put(100,50){\circle*{1}}
\put(100,70){\circle*{1}}

\path(69,11)(12,58)
\path(13.5,55.4)(12,58)(14.8,57)
\path(71,11)(128,58)
\path(125.2,57)(128,58)(126.5,55.4)
\path(12,62)(69,109)
\path(66.2,108)(69,109)(67.5,106.4)
\path(128,62)(71,109)
\path(72.5,106.4)(71,109)(73.8,108)

\path(70,61)(70,108)
\path(68,86)(72,86)
\path(69,105.2)(70,108)(71,105.2)

\path(11,60.5)(39,69.5)
\path(14,60.4)(11,60.5)(13.4,62.3)
\path(11,59.5)(39,50.5)
\path(13.4,57.7)(11,59.5)(14,59.6)
\path(41,69.5)(69,60.5)
\path(43.4,67.7)(41,69.5)(44,69.6)
\path(41,50.5)(69,59.5)
\path(44,50.4)(41,50.5)(43.4,52.3)

\path(24,67)(26,63)

\path(71,60.5)(99,69.5)
\path(71,59.5)(99,50.5)
\path(101,69.5)(129,60.5)
\path(101,50.5)(129,59.5)
\path(84,67)(86,63)

\path(74,60.4)(71,60.5)(73.4,62.3)
\path(73.4,57.7)(71,59.5)(74,59.6)
\path(103.4,67.7)(101,69.5)(104,69.6)
\path(104,50.4)(101,50.5)(103.4,52.3)

\path(100,52)(100,68)
\path(99.2,65.9)(100,68)(100.8,65.9)

\put(69,5){{\xssc $x$}}
\put(69,115){{\xssc $y$}}

\put(7,60){{\xssc $a$}}
\put(69,57){{\xssc $b$}}
\put(132,60){{\xssc $c$}}

\put(39,47){{\xssc $d$}}
\put(99,47){{\xssc $e$}}
\put(39,72){{\xssc $f$}}
\put(99,72){{\xssc $g$}}

\path(92,27)(140,27)(140,75)(89,75)(89,68)
\put(89,75){\vector(0,-1){7}}
\put(89,75){\vector(0,-1){6}}

\path(90,25)(142,25)(142,77)(72,77)
\put(116,77){\vector(-1,0){44}}
\put(116,77){\vector(-1,0){43}}

\path(88,23)(144,23)(144,120)(22,120)(22,66)
\put(22,120){\vector(0,-1){54}}
\put(22,120){\vector(0,-1){53}}

\path(80,17)(80,3)(22,3)(22,54)
\put(22,3){\vector(0,1){51}}
\put(22,3){\vector(0,1){50}}

\end{picture}
}

\end{diagram}

\vspace{4mm}

 karl-search= End Diagram I-U-reac-x
\vspace{7mm}

 *************************************

\vspace{7mm}

\clearpage

\er

$++++++++++++++++++++++++++++++++++++++++++++++++++++++$
\subsubsection{Diagram I-U-reac-c}
\subsubsection{Diagram I-U-reac-c}

 {\LARGE karl-search= Start Diagram I-U-reac-c }

\vspace{10mm}

\begin{diagram}

{\xssc LABEL: {Diagram I-U-reac-c}}
\label{Diagram I-U-reac-c}
\index{Diagram I-U-reac-c}

\centering
\setlength{\unitlength}{1mm}
{\renewcommand{\dashlinestretch}{30}
\begin{picture}(150,150)(0,0)

\put(10,60){\circle*{1}}
\put(70,60){\circle*{1}}
\put(130,60){\circle*{1}}

\put(70,10){\circle*{1}}
\put(70,110){\circle*{1}}

\put(40,50){\circle*{1}}
\put(40,70){\circle*{1}}
\put(100,50){\circle*{1}}
\put(100,70){\circle*{1}}

\path(69,11)(12,58)
\path(13.5,55.4)(12,58)(14.8,57)
\path(71,11)(128,58)
\path(125.2,57)(128,58)(126.5,55.4)
\path(12,62)(69,109)
\path(66.2,108)(69,109)(67.5,106.4)
\path(128,62)(71,109)
\path(72.5,106.4)(71,109)(73.8,108)

\path(70,61)(70,108)
\path(68,86)(72,86)
\path(69,105.2)(70,108)(71,105.2)

\path(11,60.5)(39,69.5)
\path(14,60.4)(11,60.5)(13.4,62.3)
\path(11,59.5)(39,50.5)
\path(13.4,57.7)(11,59.5)(14,59.6)
\path(41,69.5)(69,60.5)
\path(43.4,67.7)(41,69.5)(44,69.6)
\path(41,50.5)(69,59.5)
\path(44,50.4)(41,50.5)(43.4,52.3)

\path(24,67)(26,63)

\path(71,60.5)(99,69.5)
\path(71,59.5)(99,50.5)
\path(101,69.5)(129,60.5)
\path(101,50.5)(129,59.5)
\path(84,67)(86,63)

\path(74,60.4)(71,60.5)(73.4,62.3)
\path(73.4,57.7)(71,59.5)(74,59.6)
\path(103.4,67.7)(101,69.5)(104,69.6)
\path(104,50.4)(101,50.5)(103.4,52.3)

\path(100,52)(100,68)
\path(99.2,65.9)(100,68)(100.8,65.9)

\put(69,5){{\xssc $x$}}
\put(69,115){{\xssc $y$}}

\put(7,60){{\xssc $a$}}
\put(69,57){{\xssc $b$}}
\put(132,60){{\xssc $c$}}

\put(39,47){{\xssc $d$}}
\put(99,47){{\xssc $e$}}
\put(39,72){{\xssc $f$}}
\put(99,72){{\xssc $g$}}

\path(109,68)(109,75)(89,75)(89,68)
\put(89,75){\vector(0,-1){7}}
\put(89,75){\vector(0,-1){6}}
\path(115,56)(115,76)(86,76)(86,67)
\put(86,76){\vector(0,-1){9}}
\put(86,76){\vector(0,-1){8}}
\path(116,56.5)(116,77)(72,77)
\put(116,77){\vector(-1,0){44}}
\put(116,77){\vector(-1,0){43}}
\path(117,57)(117,120)(22,120)(22,66)
\put(22,120){\vector(0,-1){54}}
\put(22,120){\vector(0,-1){53}}
\path(117,54)(117,3)(22,3)(22,54)
\put(22,3){\vector(0,1){51}}
\put(22,3){\vector(0,1){50}}

\end{picture}
}

\end{diagram}

\vspace{4mm}

 karl-search= End Diagram I-U-reac-c
\vspace{7mm}

 *************************************

\vspace{7mm}

\clearpage

$++++++++++++++++++++++++++++++++++++++++++++++++++++++$
\subsubsection{Diagram U-D-reactive}
\subsubsection{Diagram U-D-reactive}

 {\LARGE karl-search= Start Diagram U-D-reactive }

\vspace{10mm}

\begin{diagram}

{\xssc LABEL: {Diagram U-D-reactive}}
\label{Diagram U-D-reactive}
\index{Diagram U-D-reactive}

\unitlength1.0mm
\begin{picture}(130,100)

\newsavebox{\ZWEIvierreac}
\savebox{\ZWEIvierreac}(140,110)[bl]
{

\put(0,95){{\rm\bf The problem of downward chaining - reactive}}

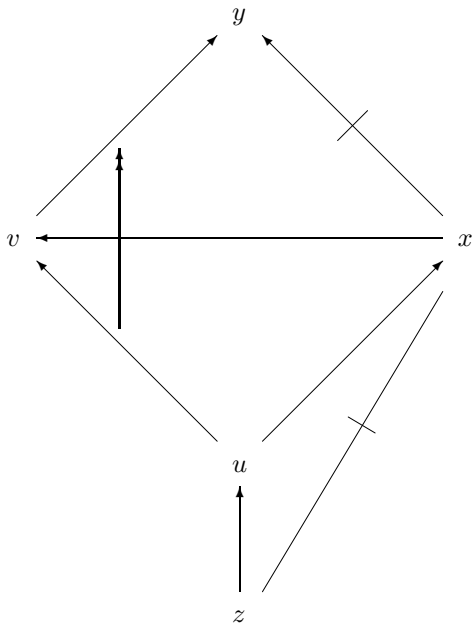
\put(43,27){\vector(1,1){24}}
\put(37,27){\vector(-1,1){24}}
\put(13,57){\vector(1,1){24}}
\put(67,57){\vector(-1,1){24}}

\put(53,67){\line(1,1){4}}

\put(67,54){\vector(-1,0){54}}

\put(40,7){\vector(0,1){14}}
\put(43,7){\line(3,5){24}}
\put(58,28.1){\line(-5,3){3.6}}

\put(24,42){\vector(0,1){24}}
\put(24,42){\vector(0,1){22.5}}

\put(39,3){$z$}
\put(39,23){$u$}
\put(9,53){$v$}
\put(69,53){$x$}
\put(39,83){$y$}

}

\put(0,0){\usebox{\ZWEIvierreac}}
\end{picture}

\end{diagram}

\vspace{4mm}

 karl-search= End Diagram U-D-reactive
\vspace{7mm}

 *************************************

\vspace{7mm}

\clearpage

$++++++++++++++++++++++++++++++++++++++++++++++++++++++$
\subsubsection{Diagram Dov-Is-2}
\subsubsection{Diagram Dov-Is-2}

 {\LARGE karl-search= Start Diagram Dov-Is-2 }

\vspace{10mm}

\begin{diagram}

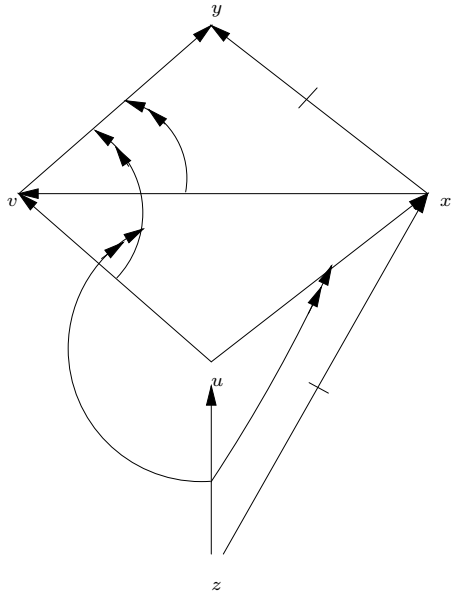
{\xssc LABEL: {Diagram Dov-Is-2}}
\label{Diagram Dov-Is-2}
\index{Diagram Dov-Is-2}

\centering
\setlength{\unitlength}{0.00083333in}
{\renewcommand{\dashlinestretch}{30}

\begin{picture}(3253,3820)(0,-10)
\put(-6848.000,6274.000){\arc{20135.789}{0.4303}{0.5839}}
\blacken\path(2278.077,1952.642)(2302.000,2074.000)(2223.782,1978.177)(2278.077,1952.642)
\put(903.476,2620.958){\arc{992.406}{4.9432}{6.4593}}
\blacken\path(1139.941,3090.381)(1017.000,3104.000)(1119.068,3034.128)(1139.941,3090.381)
\put(531.951,2405.800){\arc{1183.940}{5.2342}{7.0408}}
\blacken\path(941.882,2873.151)(827.000,2919.000)(906.791,2824.483)(941.882,2873.151)
\put(1490.621,1557.466){\arc{1660.763}{1.5028}{4.2591}}
\blacken\path(1038.516,2217.568)(1127.000,2304.000)(1008.252,2269.376)(1038.516,2217.568)
\path(2902,2524)(352,2524)
\blacken\path(472.000,2554.000)(352.000,2524.000)(472.000,2494.000)(472.000,2554.000)
\path(1552,1474)(352,2524)
\blacken\path(462.064,2467.557)(352.000,2524.000)(422.554,2422.402)(462.064,2467.557)
\path(1552,1474)(2902,2524)
\blacken\path(2825.696,2426.647)(2902.000,2524.000)(2788.860,2474.008)(2825.696,2426.647)
\path(2902,2524)(1552,3574)
\blacken\path(1665.140,3524.008)(1552.000,3574.000)(1628.304,3476.647)(1665.140,3524.008)
\path(352,2524)(1552,3574)
\blacken\path(1481.446,3472.402)(1552.000,3574.000)(1441.936,3517.557)(1481.446,3472.402)
\path(1552,274)(1552,1324)
\blacken\path(1582.000,1204.000)(1552.000,1324.000)(1522.000,1204.000)(1582.000,1204.000)
\path(1627,274)(2902,2524)
\blacken\path(2868.939,2404.807)(2902.000,2524.000)(2816.738,2434.388)(2868.939,2404.807)
\path(2157,1774)(2237,1939)
\blacken\path(2211.642,1817.934)(2237.000,1939.000)(2157.653,1844.111)(2211.642,1817.934)
\path(1262,2959)(1162,3049)
\blacken\path(1271.264,2991.023)(1162.000,3049.000)(1231.126,2946.425)(1271.264,2991.023)
\path(1037,2694)(957,2814)
\blacken\path(1048.526,2730.795)(957.000,2814.000)(998.603,2697.513)(1048.526,2730.795)
\path(867,2119)(1002,2219)
\blacken\path(923.430,2123.466)(1002.000,2219.000)(887.716,2171.679)(923.430,2123.466)
\path(2212,3184)(2097,3059)
\path(2162,1344)(2282,1279)

\put(1552,1324){\xssc $u$}
\put(1552,49)  {\xssc $z$}
\put(1552,3649){\xssc $y$}
\put(2977,2449){\xssc $x$}
\put(277,2449) {\xssc $v$}

\end{picture}

}

\end{diagram}

\vspace{4mm}

 karl-search= End Diagram Dov-Is-2
\vspace{7mm}

 *************************************

\vspace{7mm}

\clearpage

$++++++++++++++++++++++++++++++++++++++++++++++++++++++$
\subsubsection{Diagram Dov-Is-1}
\subsubsection{Diagram Dov-Is-1}

 {\LARGE karl-search= Start Diagram Dov-Is-1 }

\vspace{10mm}

\begin{diagram}

{\xssc LABEL: {Diagram Dov-Is-1}}
\label{Diagram Dov-Is-1}
\index{Diagram Dov-Is-1}

\centering
\setlength{\unitlength}{0.00083333in}
{\renewcommand{\dashlinestretch}{30}

\begin{picture}(2630,4046)(0,-10)
\put(1289.500,1811.500){\arc{2059.430}{4.5293}{7.6709}}
\blacken\path(1217.372,2868.600)(1102.000,2824.000)(1224.788,2809.060)(1217.372,2868.600)
\path(1102,274)(1852,1324)
\blacken\path(1806.663,1208.915)(1852.000,1324.000)(1757.839,1243.789)(1806.663,1208.915)
\path(1102,274)(352,1324)
\blacken\path(446.161,1243.789)(352.000,1324.000)(397.337,1208.915)(446.161,1243.789)
\path(352,1324)(1102,2374)
\blacken\path(1056.663,2258.915)(1102.000,2374.000)(1007.839,2293.789)(1056.663,2258.915)
\path(1852,1324)(1102,2374)
\blacken\path(1196.161,2293.789)(1102.000,2374.000)(1147.337,2258.915)(1196.161,2293.789)
\path(1102,2374)(1102,3274)
\blacken\path(1132.000,3154.000)(1102.000,3274.000)(1072.000,3154.000)(1132.000,3154.000)
\path(1412,2854)(1257,2844)
\blacken\path(1374.820,2881.664)(1257.000,2844.000)(1378.683,2821.788)(1374.820,2881.664)

\path(2452,274)(1852,1324)
\blacken\path(1937.584,1234.695)(1852.000,1324.000)(1885.489,1204.927)(1937.584,1234.695)

\put(1950,1300){\xssc $c$}
\put(1252,2299){\xssc $d$}
\put(1177,3199){\xssc $e$}
\put(277,1249) {\xssc $b$}
\put(1102,49)  {\xssc $a$}
\put(2477,59)  {\xssc $a'$}

\put(20,3500) {{\rm\bf Reactive graph}}

\end{picture}
}

\end{diagram}

\vspace{4mm}

 karl-search= End Diagram Dov-Is-1
\vspace{7mm}

 *************************************

\vspace{7mm}

\clearpage

$++++++++++++++++++++++++++++++++++++++++++++++++++++++$
\subsubsection{Diagram CJ-O2}
\subsubsection{Diagram CJ-O2}

 {\LARGE karl-search= Start Diagram CJ-O2 }

\vspace{10mm}

\begin{diagram}

{\xssc LABEL: {Diagram CJ-O2}}
\label{Diagram CJ-O2}
\index{Diagram CJ-O2}

\centering
\setlength{\unitlength}{1mm}
{\renewcommand{\dashlinestretch}{30}
\begin{picture}(150,100)(0,0)

\path(10,10)(90,10)(90,50)(10,10)

\put(60,40){\circle{40}}

\path(50,22.7)(70,22.7)

\put(8,6){{\xssc $m$}}

\put(58,40){{\xssc $M(\xbq)$}}

\put(55,15){{\xssc $ob(a(m)\xcs M(\xbq))$}}

\put(20,3) {{\rm\bf $ob(a(m)\xcs M(\xbq))$}}

\end{picture}
}

\end{diagram}

\vspace{4mm}

 karl-search= End Diagram CJ-O2
\vspace{7mm}

 *************************************

\vspace{7mm}

\clearpage

$++++++++++++++++++++++++++++++++++++++++++++++++++++++$
\subsubsection{Diagram CJ-O1}
\subsubsection{Diagram CJ-O1}

 {\LARGE karl-search= Start Diagram CJ-O1 }

\vspace{10mm}

\begin{diagram}

{\xssc LABEL: {Diagram CJ-O1}}
\label{Diagram CJ-O1}
\index{Diagram CJ-O1}

\centering
\setlength{\unitlength}{1mm}
{\renewcommand{\dashlinestretch}{30}
\begin{picture}(150,100)(0,0)

\path(10,10)(90,10)(90,50)(10,10)

\path(60,10)(60,30)(80,30)(80,10)

\put(8,6){{\xssc $m$}}

\put(62,15){{\xssc $ob(a(m))$}}

\put(20,3) {{\rm\bf $ob(a(m))$}}

\end{picture}
}

\end{diagram}

\vspace{4mm}

 karl-search= End Diagram CJ-O1
\vspace{7mm}

 *************************************

\vspace{7mm}

\clearpage

$++++++++++++++++++++++++++++++++++++++++++++++++++++++$
\subsubsection{Diagram FunnyRule}
\subsubsection{Diagram FunnyRule}

 {\LARGE karl-search= Start Diagram FunnyRule }

\vspace{10mm}

\begin{diagram}

{\xssc LABEL: {Diagram FunnyRule}}
\label{Diagram FunnyRule}
\index{Diagram FunnyRule}

\centering
\setlength{\unitlength}{1mm}
{\renewcommand{\dashlinestretch}{30}
\begin{picture}(150,100)(0,0)

\path(10,10)(10,70)(70,70)(70,10)(10,10)

\path(70,50)(80,50)(80,30)(70,30)
\path(70,40)(80,40)
\path(80,45)(90,45)(90,35)(80,35)
\path(90,50)(100,50)(100,30)(90,30)(90,50)
\path(90,40)(100,40)
\path(100,45)(110,45)(110,35)(100,35)

\put(40,40){{\xssc $\xba$}}

\put(20,3) {{\rm\bf FunnyRule for $()$}}

\end{picture}
}

\end{diagram}

\vspace{4mm}

 karl-search= End Diagram FunnyRule
\vspace{7mm}

 *************************************

\vspace{7mm}

\clearpage

$++++++++++++++++++++++++++++++++++++++++++++++++++++++$
\subsubsection{Diagram External2}
\subsubsection{Diagram External2}

 {\LARGE karl-search= Start Diagram External2 }

\vspace{10mm}

\begin{diagram}

{\xssc LABEL: {Diagram External2}}
\label{Diagram External2}
\index{Diagram External2}

\centering
\setlength{\unitlength}{1mm}
{\renewcommand{\dashlinestretch}{30}
\begin{picture}(150,100)(0,0)

\path(10,10)(10,90)(90,90)(90,10)(10,10)

\path(90,60)(110,60)(110,40)(90,40)
\path(110,50)(125,50)(125,30)(90,30)

\put(40,40){{\xssc $\xba$}}
\put(95,50){{\xssc $\xba \xcu \xCN \xbb$}}
\put(95,35){{\xssc $\xba \xcu \xCN \xbb'$}}

\put(20,3) {{\rm\bf External2 for $()$}}

\end{picture}
}

\end{diagram}

\vspace{4mm}

 karl-search= End Diagram External2
\vspace{7mm}

 *************************************

\vspace{7mm}

\clearpage

$++++++++++++++++++++++++++++++++++++++++++++++++++++++$
\subsubsection{Diagram External}
\subsubsection{Diagram External}

 {\LARGE karl-search= Start Diagram External }

\vspace{10mm}

\begin{diagram}

{\xssc LABEL: {Diagram External}}
\label{Diagram External}
\index{Diagram External}

\centering
\setlength{\unitlength}{1mm}
{\renewcommand{\dashlinestretch}{30}
\begin{picture}(150,120)(0,0)

\path(10,10)(10,90)(90,90)(90,10)(10,10)

\path(40,90)(40,110)(60,110)(60,90)
\path(90,60)(110,60)(110,40)(90,40)

\put(40,40){{\xssc $\xba$}}
\put(45,100){{\xssc $\xba \xcu \xCN \xbb$}}
\put(95,50){{\xssc $\xba \xcu \xCN \xbb'$}}

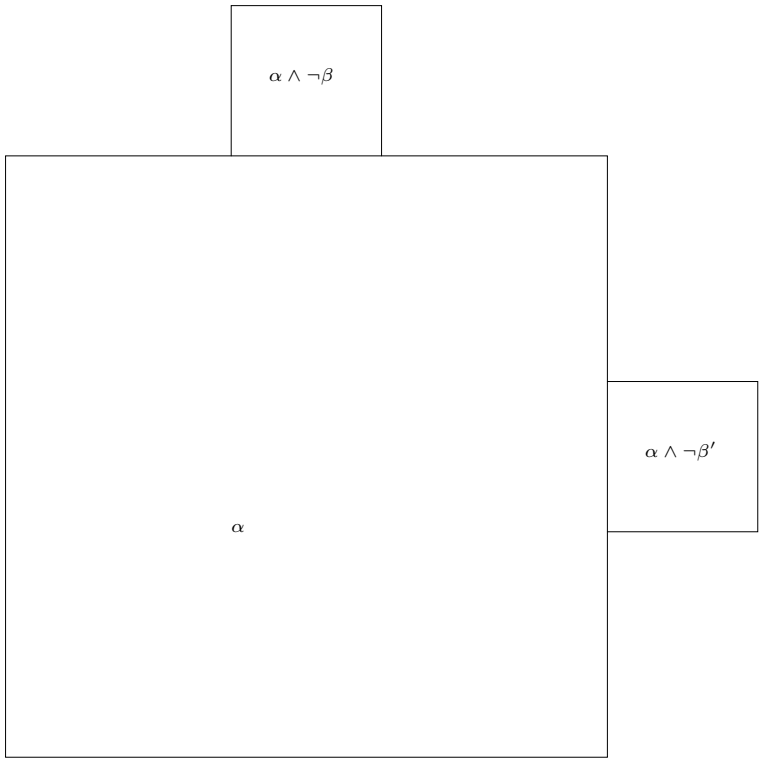
\put(20,3) {{\rm\bf External for $(OR)$}}

\end{picture}
}

\end{diagram}

\vspace{4mm}

 karl-search= End Diagram External
\vspace{7mm}

 *************************************

\vspace{7mm}

\clearpage

$++++++++++++++++++++++++++++++++++++++++++++++++++++++$
\subsubsection{Diagram Internal}
\subsubsection{Diagram Internal}

 {\LARGE karl-search= Start Diagram Internal }

\vspace{10mm}

\begin{diagram}

{\xssc LABEL: {Diagram Internal}}
\label{Diagram Internal}
\index{Diagram Internal}

\centering
\setlength{\unitlength}{1mm}
{\renewcommand{\dashlinestretch}{30}
\begin{picture}(150,100)(0,0)

\path(10,10)(10,90)(90,90)(90,10)(10,10)

\path(40,90)(40,70)(60,70)(60,90)
\path(90,60)(70,60)(70,40)(90,40)

\put(40,40){{\xssc $\xba$}}
\put(45,80){{\xssc $\xba \xcu \xCN \xbb$}}
\put(75,50){{\xssc $\xba \xcu \xCN \xbb'$}}

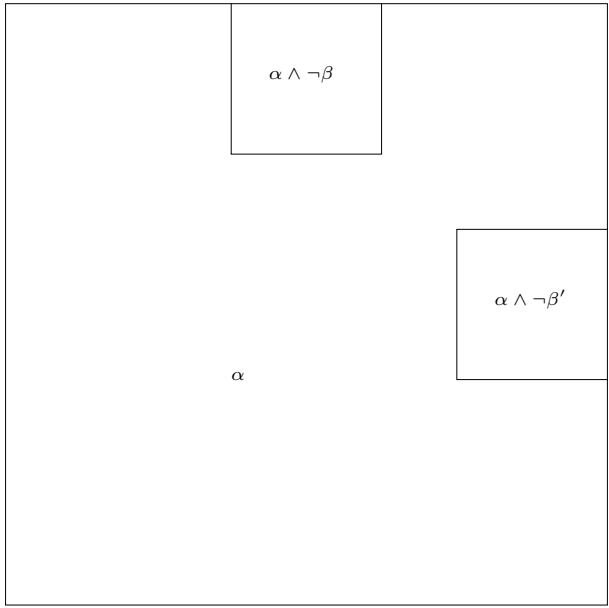
\put(20,3) {{\rm\bf Internal for $(CUM)$ and $(AND)$}}

\end{picture}
}

\end{diagram}

\vspace{4mm}

 karl-search= End Diagram Internal
\vspace{7mm}

 *************************************

\vspace{7mm}

\clearpage

$++++++++++++++++++++++++++++++++++++++++++++++++++++++$
\subsubsection{Diagram InherUniv}
\subsubsection{Diagram InherUniv}

 {\LARGE karl-search= Start Diagram InherUniv }

\vspace{10mm}

\begin{diagram}

{\xssc LABEL: {Diagram InherUniv}}
\label{Diagram InherUniv}
\index{Diagram InherUniv}

\centering
\setlength{\unitlength}{1mm}
{\renewcommand{\dashlinestretch}{30}
\begin{picture}(150,150)(0,0)

\put(10,60){\circle*{1}}
\put(70,60){\circle*{1}}
\put(130,60){\circle*{1}}

\put(70,10){\circle*{1}}
\put(70,110){\circle*{1}}

\put(40,50){\circle*{1}}
\put(40,70){\circle*{1}}
\put(100,50){\circle*{1}}
\put(100,70){\circle*{1}}

\path(69,11)(12,58)
\path(13.5,55.4)(12,58)(14.8,57)
\path(71,11)(128,58)
\path(125.2,57)(128,58)(126.5,55.4)
\path(12,62)(69,109)
\path(66.2,108)(69,109)(67.5,106.4)
\path(128,62)(71,109)
\path(72.5,106.4)(71,109)(73.8,108)

\path(70,61)(70,108)
\path(68,86)(72,86)
\path(69,105.2)(70,108)(71,105.2)

\path(11,60.5)(39,69.5)
\path(14,60.4)(11,60.5)(13.4,62.3)
\path(11,59.5)(39,50.5)
\path(13.4,57.7)(11,59.5)(14,59.6)
\path(41,69.5)(69,60.5)
\path(43.4,67.7)(41,69.5)(44,69.6)
\path(41,50.5)(69,59.5)
\path(44,50.4)(41,50.5)(43.4,52.3)

\path(24,67)(26,63)

\path(71,60.5)(99,69.5)
\path(71,59.5)(99,50.5)
\path(101,69.5)(129,60.5)
\path(101,50.5)(129,59.5)
\path(84,67)(86,63)

\path(74,60.4)(71,60.5)(73.4,62.3)
\path(73.4,57.7)(71,59.5)(74,59.6)
\path(103.4,67.7)(101,69.5)(104,69.6)
\path(104,50.4)(101,50.5)(103.4,52.3)

\path(100,52)(100,68)
\path(99.2,65.9)(100,68)(100.8,65.9)

\put(69,5){{\xssc $x$}}
\put(69,115){{\xssc $y$}}

\put(7,60){{\xssc $a$}}
\put(69,57){{\xssc $b$}}
\put(132,60){{\xssc $c$}}

\put(39,47){{\xssc $d$}}
\put(99,47){{\xssc $e$}}
\put(39,72){{\xssc $f$}}
\put(99,72){{\xssc $g$}}

\put(110,65){{\xssc $STOP$}}
\put(110,54){{\xssc $STOP$}}


\end{picture}
}

\end{diagram}

\vspace{4mm}

 karl-search= End Diagram InherUniv
\vspace{7mm}

 *************************************

\vspace{7mm}

\clearpage

$++++++++++++++++++++++++++++++++++++++++++++++++++++++$
\subsubsection{Diagram Now-1}
\subsubsection{Diagram Now-1}

 {\LARGE karl-search= Start Diagram Now-1 }

\vspace{10mm}

\begin{minipage}{14cm}

\begin{diagram}

{\xssc LABEL: {Diagram Now-1}}
\label{Diagram Now-1}
\index{Diagram Now-1}

\centering
\setlength{\unitlength}{0.00083333in}
{\renewcommand{\dashlinestretch}{30}
\begin{picture}(3286,2735)(0,0)
\put(1323,1848){\ellipse{2556}{1730}}
\path(1324,0525)(874,1575)
\path(948.845,1476.520)(874.000,1575.000)(893.696,1452.885)
\path(1324,0525)(1924,1575)
\path(1890.511,1455.927)(1924.000,1575.000)(1838.416,1485.695)
\path(1324,0525)(3274,1575)
\path(3182.566,1491.694)(3274.000,1575.000)(3154.120,1544.522)
\put(1324,2075){{\xssc minimal (preferred) worlds of $A$}}
\put(874,1650) {{\xssc $t_1$}}
\put(1924,1650){{\xssc $t_2$}}
\put(3274,1650){{\xssc $t_3\ldots$}}
\put(1324,0300)   {{\xssc now}}
\put(100,100) {{\rm\bf Diagram 1}}
\end{picture}
}

\end{diagram}

\end{minipage}

\vspace{4mm}

 karl-search= End Diagram Now-1
\vspace{7mm}

 *************************************

\vspace{7mm}

\clearpage

$++++++++++++++++++++++++++++++++++++++++++++++++++++++$
\subsubsection{Diagram Now-2}
\subsubsection{Diagram Now-2}

 {\LARGE karl-search= Start Diagram Now-2 }

\vspace{10mm}

\begin{diagram}

{\xssc LABEL: {Diagram Now-2}}
\label{Diagram Now-2}
\index{Diagram Now-2}

\centering
\setlength{\unitlength}{0.00083333in}
{\renewcommand{\dashlinestretch}{30}
\begin{picture}(3677,4235)(0,0)
\put(1286,3348){\ellipse{2556}{1730}}
\dashline{60.000}(837,3075)(537,2175)
\dashline{60.000}(1812,3075)(1587,2175)
\path(1062,1425)(1587,2175)
\path(1542.761,2059.488)(1587.000,2175.000)(1493.608,2093.896)
\path(1062,1425)(537,2175)
\path(630.392,2093.896)(537.000,2175.000)(581.239,2059.488)
\path(1887,0525)(2337,1125)
\path(2289.000,1011.000)(2337.000,1125.000)(2241.000,1047.000)
\path(1887,0525)(1062,1125)
\path(1176.693,1078.681)(1062.000,1125.000)(1141.403,1030.157)
\path(1887,0525)(3462,1125)
\path(3360.541,1054.246)(3462.000,1125.000)(3339.182,1110.315)
\put(837,3150) {{\xssc $t_1$}}
\put(1887,3150){{\xssc $t_2$}}
\put(1287,3675){{\xssc minimal $A$ worlds}}
\put(1062,1200) {{\xssc $r_1$}}
\put(2337,1200) {{\xssc $r_2$}}
\put(3462,1200) {{\xssc $r_3$}}
\put(1887,300)   {{\xssc now}}
\put(100,100) {{\rm\bf Diagram 2}}
\end{picture}
}
\end{diagram}

\vspace{4mm}

 karl-search= End Diagram Now-2
\vspace{7mm}

 *************************************

\vspace{7mm}

\clearpage

$++++++++++++++++++++++++++++++++++++++++++++++++++++++$
\subsubsection{Diagram IBRS}
\subsubsection{Diagram IBRS}

 {\LARGE karl-search= Start Diagram IBRS }

\vspace{10mm}

\begin{diagram}

\centering
\setlength{\unitlength}{0.00083333in}
{\renewcommand{\dashlinestretch}{30}
\begin{picture}(4961,5004)(0,0)
\path(1511,1583)(611,3683)
\blacken\path(685.845,3584.520)(611.000,3683.000)(630.696,3560.885)(672.451,3539.613)(685.845,3584.520)
\path(1511,1583)(2411,3683)
\blacken\path(2391.304,3560.885)(2411.000,3683.000)(2336.155,3584.520)(2349.549,3539.613)(2391.304,3560.885)
\path(3311,1583)(4361,4133)
\blacken\path(4343.050,4010.616)(4361.000,4133.000)(4287.570,4033.461)(4301.603,3988.750)(4343.050,4010.616)
\path(3316,1574)(2416,3674)
\blacken\path(2490.845,3575.520)(2416.000,3674.000)(2435.696,3551.885)(2477.451,3530.613)(2490.845,3575.520)
\path(986,2783)(2621,2783)
\blacken\path(2501.000,2753.000)(2621.000,2783.000)(2501.000,2813.000)(2465.000,2783.000)(2501.000,2753.000)
\path(2486,2783)(2786,2783)
\blacken\path(2666.000,2753.000)(2786.000,2783.000)(2666.000,2813.000)(2630.000,2783.000)(2666.000,2753.000)
\path(3311,1583)(2051,2368)
\blacken\path(2168.714,2330.008)(2051.000,2368.000)(2136.987,2279.083)(2183.406,2285.509)(2168.714,2330.008)
\path(2166,2288)(1906,2458)
\blacken\path(2022.854,2417.439)(1906.000,2458.000)(1990.019,2367.221)(2036.567,2372.629)(2022.854,2417.439)

\put(1511,1358) {{\xssc $a$}}
\put(3311,1358) {{\xssc $d$}}
\put(3311,1058)  {{\xssc $(p,q)=(1,0)$}}
\put(1511,1058)  {{\xssc $(p,q)=(0,0)$}}
\put(2411,3758){{\xssc $c$}}
\put(4361,4433){{\xssc $(p,q)=(1,1)$}}
\put(4361,4208){{\xssc $e$}}
\put(2411,3983){{\xssc $(p,q)=(0,1)$}}
\put(611,3983) {{\xssc $(p,q)=(0,1)$}}
\put(611,3758) {{\xssc $b$}}
\put(1211,2883){{\xssc $(p,q)=(1,1)$}}
\put(260,2333) {{\xssc $(p,q)=(1,1)$}}
\put(2261,1583) {{\xssc $(p,q)=(1,1)$}}
\put(1286,3233){{\xssc $(p,q)=(1,1)$}}
\put(2711,3083){{\xssc $(p,q)=(1,1)$}}
\put(3836,2633){{\xssc $(p,q)=(1,1)$}}

\put(300,700)
{{\rm\bf
A simple example of an information bearing system.
}}

\end{picture}
}
{\xssc LABEL: {Diagram IBRS}}
\label{Diagram IBRS}
\index{Diagram IBRS}

\end{diagram}

We have here:
\[\begin{array}{l}
S =\{a,b,c,d,e\}.\\
\xdR = S \cup \{(a,b), (a,c), (d,c), (d,e)\} \cup \{((a,b), (d,c)),
(d,(a,c))\}.\\
Q = \{p,q\}
\end{array}
\]
The values of $h$ for $p$ and $q$ are as indicated in the figure. For
example $h(p,(d,(a,c))) =1$.

\vspace{4mm}

 karl-search= End Diagram IBRS
\vspace{7mm}

 *************************************

\vspace{7mm}

\clearpage

$++++++++++++++++++++++++++++++++++++++++++++++++++++++$
\subsubsection{Diagram Pischinger}
\subsubsection{Diagram Pischinger}

 {\LARGE karl-search= Start Diagram Pischinger }

\vspace{10mm}

\begin{diagram}

{\xssc LABEL: {Diagram Pischinger}}
\label{Diagram Pischinger}
\index{Diagram Pischinger}

\centering
\setlength{\unitlength}{1mm}
{\renewcommand{\dashlinestretch}{30}
\begin{picture}(150,150)(0,0)

\path(30,10)(30,110)(110,110)(110,10)(30,10)
\path(30,30)(110,30)
\path(30,50)(110,50)
\path(30,70)(110,70)
\path(30,90)(110,90)

\put(70,10){\arc{20}{-3.14}{0}}
\put(70,30){\arc{20}{-3.14}{0}}
\put(70,50){\arc{20}{-3.14}{0}}
\put(70,70){\arc{20}{-3.14}{0}}
\put(70,90){\arc{20}{-3.14}{0}}

\put(60,80){\circle{30}}

\path(5,80)(25,80)
\path(22.3,81)(25,80)(22.3,79)
\put(5,80){\circle*{1}}
\put(25,80){\circle*{1}}

\put(15,75){\xssc{$R$}}
\put(5,75){\xssc{$t$}}
\put(25,75){\xssc{$s$}}

\put(10,130){\xssc{The overall structure is visible from $t$}}
\put(10,123){\xssc{Only the inside of the circle is visible from $s$}}
\put(10,116){\xssc{Half-circles are the sets of minimal elements of layers}}

\put(20,3) {{\rm\bf $\xda-$ ranked structure and accessibility}}

\end{picture}
}

\end{diagram}

\vspace{4mm}

 karl-search= End Diagram Pischinger
\vspace{7mm}

 *************************************

\vspace{7mm}

\clearpage

$++++++++++++++++++++++++++++++++++++++++++++++++++++++$
\subsubsection{Diagram ReacA}
\subsubsection{Diagram ReacA}

 {\LARGE karl-search= Start Diagram ReacA }

\vspace{10mm}

\begin{diagram}

{\xssc LABEL: {Diagram ReacA}}
\label{Diagram ReacA}
\index{Diagram ReacA}

\centering
\setlength{\unitlength}{1mm}
{\renewcommand{\dashlinestretch}{30}
\begin{picture}(150,150)(0,0)

\put(50,15){\circle*{1}}
\put(20,50){\circle*{1}}
\put(80,50){\circle*{1}}
\put(50,85){\circle*{1}}
\put(50,115){\circle*{1}}

\put(50,10){\xssc{$a$}}
\put(15,50){\xssc{$b$}}
\put(82,50){\xssc{$c$}}
\put(52,85){\xssc{$d$}}
\put(52,115){\xssc{$e$}}

\path(51,16)(79,49)
\path(76.4,47.5)(79,49)(78,46.2)
\path(49,16)(21,49)
\path(22,46.2)(21,49)(23.6,47.5)
\path(21,51)(49,84)
\path(46.4,82.5)(49,84)(48,81.2)
\path(79,51)(51,84)
\path(52,81.2)(51,84)(53.6,82.5)
\path(50,86)(50,114)
\path(49,111.2)(50,114)(51,111.2)

\path(65,32.5)(100,32.5)(100,100)(51,100)
\path(53.8,99)(51,100)(53.8,101)
\path(55.8,99)(53,100)(55.8,101)


\end{picture}
}

\end{diagram}

\vspace{4mm}

 karl-search= End Diagram ReacA
\vspace{7mm}

 *************************************

\vspace{7mm}

\clearpage

$++++++++++++++++++++++++++++++++++++++++++++++++++++++$
\subsubsection{Diagram ReacB}
\subsubsection{Diagram ReacB}

 {\LARGE karl-search= Start Diagram ReacB }

\vspace{10mm}

\begin{diagram}

{\xssc LABEL: {Diagram ReacB}}
\label{Diagram ReacB}
\index{Diagram ReacB}

\centering
\setlength{\unitlength}{1mm}
{\renewcommand{\dashlinestretch}{30}
\begin{picture}(150,150)(0,0)

\put(50,15){\circle*{1}}
\put(20,50){\circle*{1}}
\put(80,50){\circle*{1}}
\put(20,80){\circle*{1}}
\put(20,110){\circle*{1}}
\put(80,80){\circle*{1}}

\put(50,10){\xssc{$a$}}
\put(15,50){\xssc{$ab$}}
\put(82,50){\xssc{$ac$}}
\put(13,80){\xssc{$abd$}}
\put(82,80){\xssc{$acd$}}
\put(11,110){\xssc{$abde$}}

\path(51,16)(79,49)
\path(76.4,47.5)(79,49)(78,46.2)
\path(49,16)(21,49)
\path(22,46.2)(21,49)(23.6,47.5)
\path(20,51)(20,79)
\path(19,76.2)(20,79)(21,76.2)
\path(20,81)(20,109)
\path(19,106.2)(20,109)(21,106.2)
\path(80,51)(80,79)
\path(79,76.2)(80,79)(81,76.2)

\path(5,70)(95,70)(95,90)(5,90)(5,70)


\end{picture}
}

\end{diagram}

\vspace{4mm}

 karl-search= End Diagram ReacB
\vspace{7mm}

 *************************************

\vspace{7mm}

\clearpage

$++++++++++++++++++++++++++++++++++++++++++++++++++++++$
\subsubsection{Diagram CumSem}
\subsubsection{Diagram CumSem}

 {\LARGE karl-search= Start Diagram CumSem }

\vspace{10mm}

\begin{diagram}

{\xssc LABEL: {Diagram CumSem}}
\label{Diagram CumSem}
\index{Diagram CumSem}

\centering
\setlength{\unitlength}{1mm}
{\renewcommand{\dashlinestretch}{30}
\begin{picture}(150,150)(0,0)

\put(60,60){\circle{50}}
\path(35,60)(85,60)
\path(40,45)(80,45)

\put(90,75){\xssc{$M(T)$}}
\put(90,55){\xssc{$M(T')$}}
\put(56,40){\xssc{$\xbm(M(T))$}}


\end{picture}
}

\end{diagram}

\vspace{4mm}

 karl-search= End Diagram CumSem
\vspace{7mm}

 *************************************

\vspace{7mm}

\clearpage

$++++++++++++++++++++++++++++++++++++++++++++++++++++++$
\subsubsection{Diagram Eta-Rho-1}
\subsubsection{Diagram Eta-Rho-1}

 {\LARGE karl-search= Start Diagram Eta-Rho-1 }

\vspace{10mm}

\begin{diagram}

{\xssc LABEL: {Diagram Eta-Rho-1}}
\label{Diagram Eta-Rho-1}
\index{Diagram Eta-Rho-1}

\centering
\setlength{\unitlength}{1mm}
{\renewcommand{\dashlinestretch}{30}
\begin{picture}(150,100)(0,0)

\put(30,50){\circle{50}}
\put(100,50){\circle{40}}
\put(100,50){\circle{70}}

\put(30,50){\circle*{1}}
\put(70,50){\circle*{1}}

\path(31,50)(69,50)
\path(66.6,50.9)(69,50)(66.6,49.1)

\put(30,60){\xssc{$X$}}
\put(100,80){\xssc{$\xbh(X)$}}
\put(100,60){\xssc{$\xbr(X)$}}

\put(30,10) {{\rm\bf Attacking structure}}

\end{picture}
}

\end{diagram}

\vspace{4mm}

 karl-search= End Diagram Eta-Rho-1
\vspace{7mm}

 *************************************

\vspace{7mm}

\clearpage

$++++++++++++++++++++++++++++++++++++++++++++++++++++++$
\subsubsection{Diagram Structure-rho-eta}
\subsubsection{Diagram Structure-rho-eta}

 {\LARGE karl-search= Start Diagram Structure-rho-eta }

\vspace{10mm}

\begin{diagram}

{\xssc LABEL: {Diagram Structure-rho-eta}}
\label{Diagram Structure-rho-eta}
\index{Diagram Structure-rho-eta}

\centering
\setlength{\unitlength}{1mm}
{\renewcommand{\dashlinestretch}{30}
\begin{picture}(150,100)(0,0)

\put(30,50){\circle{50}}
\put(100,50){\circle{40}}
\put(100,50){\circle{70}}

\put(30,50){\circle*{1}}
\put(90,50){\circle*{1}}
\put(92,50){{\xssc $< x, f, X>$}}
\put(32,50) {{\xssc $x'$}}

\put(60,20){\arc{84.5}{-2.34}{-0.8}}
\put(70,62){{\xssc $< a, x_0''>$}}
\put(60,80){\arc{84.5}{0.8}{2.34}}
\put(70,38){{\xssc $< a, x_1''>$}}

\path(88.9,51.7)(90,50)(88,50.5)
\path(88,49.5)(90,50)(88.9,48.3)

\put(35,65){\circle*{1}}
\put(35,35){\circle*{1}}

\put(36,65) {{\xssc $x_0''$}}
\put(36,35){{\xssc $x''_1$}}

\path(35,65)(40,58)
\path(39.7,59.2)(40,58)(39.0,58.7)
\put(25,62){{\xssc $< \beta, x_1''>$}}

\path(35,35)(40,42)
\path(39,41.3)(40,42)(39.7,40.8)
\put(25,38) {{\xssc $< \beta, x_o''>$}}

\put(15,50){\xssc{$X$}}
\put(100,80){\xssc{$\xbh(X)$}}
\put(100,60){\xssc{$\xbr(X)$}}

\put(30,10) {{\rm\bf Attacking structure}}

\end{picture}
}

\end{diagram}

\vspace{4mm}

 karl-search= End Diagram Structure-rho-eta
\vspace{7mm}

 *************************************

\vspace{7mm}

\clearpage

$++++++++++++++++++++++++++++++++++++++++++++++++++++++$
\subsubsection{Diagram Gate-Sem}
\subsubsection{Diagram Gate Semantics}

 {\LARGE karl-search= Start Diagram Gate Semantics }

\vspace{10mm}

\begin{diagram}

{\xssc LABEL: {Diagram Gate-Sem}}
\label{Diagram Gate-Sem}
\index{Diagram Gate-Sem}

\centering
\setlength{\unitlength}{1mm}
{\renewcommand{\dashlinestretch}{30}
\begin{picture}(150,170)(0,0)


\put(15,130){\arc{10}{-1.57}{1.57}}
\path(15,125)(15,135)
\put(16.3,129.3){\xssc{$\xcu$}}
\put(22,132){\xssc{$A1$}}
\path(13,132)(15,133)(13,134)
\path(13,126)(15,127)(13,128)

\put(45,133){\arc{10}{-1.57}{1.57}}
\path(45,128)(45,138)
\put(46.3,132.3){\xssc{$\xco$}}
\put(52,135){\xssc{$A3$}}
\path(43,135)(45,136)(43,137)
\path(43,129)(45,130)(43,131)

\path(75,128)(75,138)(83,133)(75,128)
\put(76.3,132.3){\xssc{$\xCN$}}
\path(73,132)(75,133)(73,134)

\path(0,127)(15,127)
\path(20,130)(45,130)
\path(50,133)(75,133)
\path(83,133)(108,133)
\path(106,132)(108,133)(106,134)
\put(93,133){\circle*{1}}
\put(101,133){\circle*{1}}

\put(0,124){\xssc{$In1$}}
\put(110,132){\xssc{$Out1$}}

\path(15,133)(5,133)(5,150)(101,150)(101,133)
\path(45,136)(35,136)(35,143)(85,143)(93,77)


\put(15,80){\arc{10}{-1.57}{1.57}}
\path(15,75)(15,85)
\put(16.3,79.3){\xssc{$\xcu$}}
\put(22,76){\xssc{$A2$}}
\path(13,82)(15,83)(13,84)
\path(13,76)(15,77)(13,78)

\put(45,77){\arc{10}{-1.57}{1.57}}
\path(45,72)(45,82)
\put(46.3,76.3){\xssc{$\xco$}}
\put(52,73){\xssc{$A4$}}
\path(43,79)(45,80)(43,81)
\path(43,73)(45,74)(43,75)

\path(75,72)(75,82)(83,77)(75,72)
\put(76.3,76.3){\xssc{$\xCN$}}
\path(73,78)(75,77)(73,76)

\path(0,83)(15,83)
\path(20,80)(45,80)
\path(50,77)(75,77)
\path(83,77)(108,77)
\path(106,76)(108,77)(106,78)
\put(93,77){\circle*{1}}
\put(101,77){\circle*{1}}

\put(0,85.5){\xssc{$In2$}}
\put(110,76){\xssc{$Out2$}}

\path(15,77)(5,77)(5,60)(101,60)(101,77)
\path(45,74)(35,74)(35,67)(85,67)(93,133)

\put(30,20) {{\rm\bf Gate Semantics}}

\end{picture}
}

\end{diagram}

\vspace{4mm}

 karl-search= End Diagram Gate Semantics
\vspace{7mm}

 *************************************

\vspace{7mm}

\clearpage

$++++++++++++++++++++++++++++++++++++++++++++++++++++++$
\subsubsection{Diagram A-Ranked}
\subsubsection{Diagram A-Ranked}

 {\LARGE karl-search= Start Diagram A-Ranked }

\vspace{20mm}

\begin{diagram}

{\xssc LABEL: {Diagram A-Ranked}}
\label{Diagram A-Ranked}
\index{Diagram A-Ranked}

\centering
\setlength{\unitlength}{0.00083333in}
{\renewcommand{\dashlinestretch}{30}
\begin{picture}(2390,3581)(0,0)
\put(1212.000,2028.000){\arc{1110.000}{3.4719}{5.9529}}
\put(979.651,835.201){\arc{1095.700}{3.7717}{5.6071}}
\put(949,2283){\ellipse{1874}{2550}}
\path(12,2208)(1887,2208)
\path(1137,3108)(1137,2358)
\path(1107.000,2478.000)(1137.000,2358.000)(1167.000,2478.000)
\path(1137,3108)(312,1758)
\path(348.976,1876.037)(312.000,1758.000)(400.172,1844.750)
\path(1437,2358)(1662,1758)
\path(1591.775,1859.826)(1662.000,1758.000)(1647.955,1880.893)
\path(1437,2358)(1137,1758)
\path(1163.833,1878.748)(1137.000,1758.000)(1217.498,1851.915)
\path(1137,1758)(1137,1158)
\path(1107.000,1278.000)(1137.000,1158.000)(1167.000,1278.000)
\put(2037,2658){{\xssc $A'$, layer of lesser quality}}
\put(2037,1458){{\xssc $A$, best layer}}
\put(-700,800){{\xssc
Each layer behaves inside like any preferential structure.}}
\put(-700,600){{\xssc
Amongst each other, layers behave like ranked structures.}}

\put(100,200) {{\rm\bf $\xda-$ ranked structure}}

\end{picture}
}
\end{diagram}

\vspace{4mm}

 karl-search= End Diagram A-Ranked
\vspace{7mm}

 *************************************

\vspace{7mm}

\clearpage

$++++++++++++++++++++++++++++++++++++++++++++++++++++++$
\subsubsection{Diagram C-Validity}
\subsubsection{Diagram C-Validity}

 {\LARGE karl-search= Start Diagram C-Validity }

\vspace{20mm}

\begin{diagram}

{\xssc LABEL: {Diagram C-Validity}}
\label{Diagram C-Validity}
\index{Diagram C-Validity}

\centering
\setlength{\unitlength}{0.00083333in}
{\renewcommand{\dashlinestretch}{30}
\begin{picture}(3841,3665)(0,200)
\put(2407.000,2628.000){\arc{1110.000}{3.4719}{5.9529}}
\put(2407.000,1578.000){\arc{1110.000}{3.4719}{5.9529}}
\put(2398.149,661.355){\arc{1503.312}{3.9463}{5.4292}}
\put(-24.763,2227.012){\arc{4826.591}{5.7407}{6.7841}}
\put(2369,2283){\ellipse{1874}{2550}}
\path(1507,2808)(3232,2808)
\path(1507,1758)(3232,1758)
\path(907,3263)(2562,1963)
\path(2449.100,2013.534)(2562.000,1963.000)(2486.163,2060.718)
\path(157,3248)(2072,1213)
\path(1967.915,1279.831)(2072.000,1213.000)(2011.611,1320.950)
\path(272,3258)(847,3258)
\put(210,3258){\circle*{20}}
\put(880,3258){\circle*{20}}
\path(727.000,3228.000)(847.000,3258.000)(727.000,3288.000)
\put(3457,3108){{\xssc $A''$}}
\put(3457,2208){{\xssc $A'$}}
\put(3457,1308){{\xssc $A$}}
\put(4000,1093){{\xssc $m \xcm \xdc$}}
\put(4000,8033){{\xssc $m' \xcM \xdc$}}
\put(2407,2883){{\xssc $\mu(A'')$}}
\put(2407,1833){{\xssc $\mu(A')$}}
\put(2407,1093){{\xssc $\mu(A)$}}
\put(132,3318){{\xssc $m$}}
\put(882,3313){{\xssc $m'$}}
\put(1972,3518){{\xssc $B$}}
\put(10,800){{\xssc Here, the ``best'' element $m$ sees is in
$B$, so $\xdc$ holds in $m$.}}
\put(10,600){{\xssc The ``best'' element $m'$ sees is not in
$B$, so $\xdc$ does not hold in $m'$.}}

\put(100,250) {{\rm\bf Validity of $\xdc$ from $m$ and $m'$}}

\end{picture}
}
\end{diagram}

\vspace{4mm}

 karl-search= End Diagram C-Validity
\vspace{7mm}

 *************************************

\vspace{7mm}

\clearpage

$++++++++++++++++++++++++++++++++++++++++++++++++++++++$
\subsubsection{Diagram IBRS-Base}
\subsubsection{Diagram IBRS-Base}

 {\LARGE karl-search= Start Diagram IBRS-Base }

\vspace{30mm}

\begin{diagram}

{\xssc LABEL: {Diagram IBRS-Base}}
\label{Diagram IBRS-Base}
\index{Diagram IBRS-Base}

\centering
\setlength{\unitlength}{1mm}
{\renewcommand{\dashlinestretch}{30}
\begin{picture}(110,150)(0,0)

\put(50,50){\circle{80}}
\put(50,50){\circle{40}}

\path(50,40)(50,60)
\path(49,57.3)(50,60)(51,57.3)
\put(50,39.2){\circle*{0.3}}
\put(50,60.8){\circle*{0.3}}
\put(50,38){\xssc{$a$}}
\put(50,62){\xssc{$b$}}

\path(20,50)(49,50)
\path(46.2,51)(49,50)(46.2,49)
\put(19.2,50){\circle*{0.3}}
\put(17,50){\xssc{$c$}}

\put(60,50){\xssc{$X$}}
\put(60,80){\xssc{$Y$}}

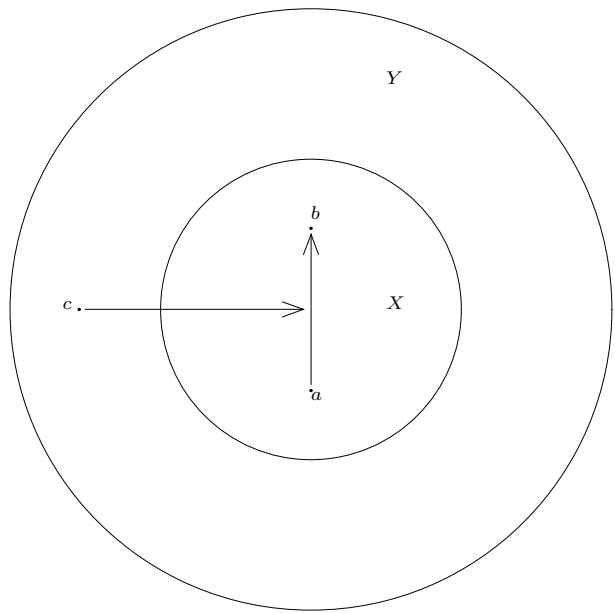
\put(30,2) {{\rm\bf Attacking an arrow}}

\end{picture}
}

\end{diagram}

\vspace{4mm}

 karl-search= End Diagram IBRS-Base
\vspace{7mm}

 *************************************

\vspace{7mm}

\clearpage

$++++++++++++++++++++++++++++++++++++++++++++++++++++++$
\subsubsection{Diagram Essential-Smooth-Repr}
\subsubsection{Diagram Essential Smooth Repr}

 {\LARGE karl-search= Start Diagram Essential Smooth Repr }

\vspace{30mm}

\begin{diagram}

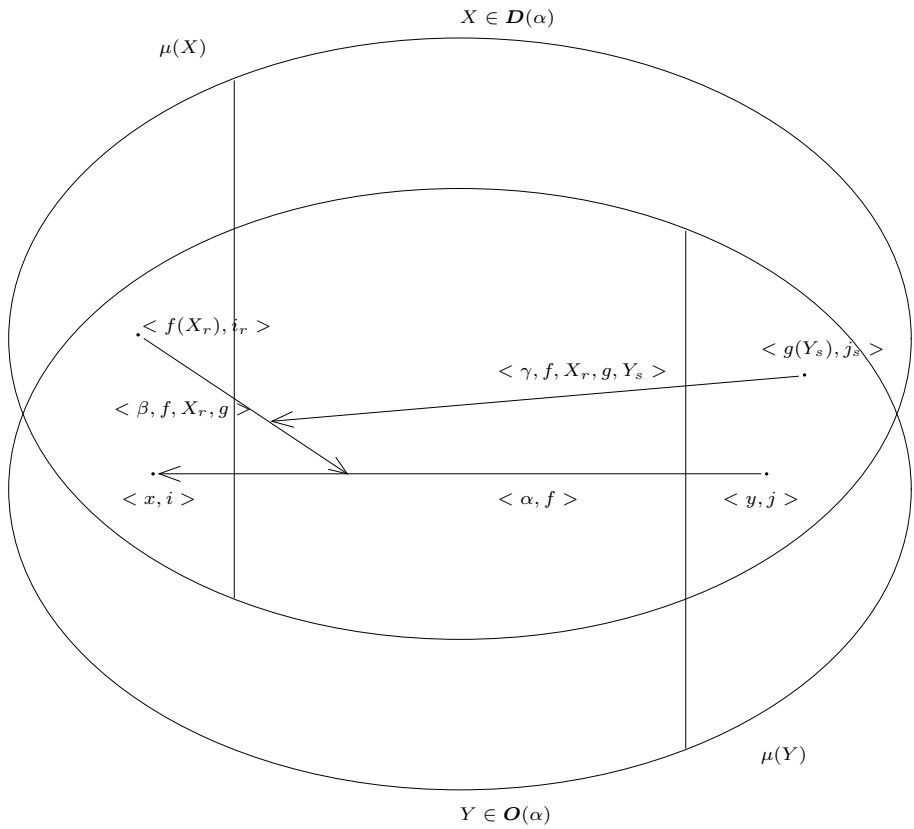
{\xssc LABEL: {Diagram Essential-Smooth-Repr}}
\label{Diagram Essential-Smooth-Repr}
\index{Diagram Essential-Smooth-Repr}

\centering
\setlength{\unitlength}{1mm}
{\renewcommand{\dashlinestretch}{30}
\begin{picture}(110,150)(0,0)
\put(50,80){\ellipse{120}{80}}
\put(50,100){\ellipse{120}{80}}
\path(20,134.3)(20,65.5)
\path(80,114.3)(80,45.5)
\put(50,142){\xssc{$X \xbe \xdD(\xba)$}}
\put(50,36){\xssc{$Y \xbe \xdO(\xba)$}}
\put(10,138){\xssc{$\xbm(X)$}}
\put(90,44){\xssc{$\xbm(Y)$}}

\path(10,82)(90,82)
\put(9.2,82){\circle*{0.3}}
\put(90.8,82){\circle*{0.3}}
\path(12.8,81)(10,82)(12.8,83)
\put(5,78){{\xssc $<x,i>$}}
\put(85,78){{\xssc $<y,j>$}}
\put(55,78){{\xssc $<\xba,f>$}}

\path(8,100)(35,82)
\path(32.5,82.6)(35,82)(33.5,84.2)
\put(7.2,100.5){\circle*{0.3}}
\put(7.7,101){{\xssc $<f(X_r),i_r>$}}
\put(4,90){{\xssc $<\xbb,f,X_r,g>$}}

\path(25,89)(95,95)
\path(27.9,88.2)(25,89)(27.7,90.26)
\put(95.8,95.2){\circle*{0.3}}
\put(90,98){{\xssc $<g(Y_s),j_s>$}}
\put(55,95){{\xssc $<\xbg,f,X_r,g,Y_s>$}}

\put(10,15) {{\rm\bf The construction}}

\end{picture}
}

\end{diagram}

\vspace{4mm}

 karl-search= End Diagram Essential Smooth Repr
\vspace{7mm}

 *************************************

\vspace{7mm}

\clearpage

$++++++++++++++++++++++++++++++++++++++++++++++++++++++$
\subsubsection{Diagram Essential-Smooth-2-1-2}
\subsubsection{Diagram Essential Smooth 2-1-2}

 {\LARGE karl-search= Start Diagram Essential Smooth 2-1-2 }

\vspace{10mm}

\begin{diagram}

{\xssc LABEL: {Diagram Essential-Smooth-2-1-2}}
\label{Diagram Essential-Smooth-2-1-2}
\index{Diagram Essential-Smooth-2-1-2}

\centering
\setlength{\unitlength}{1mm}
{\renewcommand{\dashlinestretch}{30}
\begin{picture}(100,100)(0,0)
\put(50,60){\circle{80}}
\path(10,60)(90,60)

\put(100,80){{\xssc $X$}}
\put(100,40){{\xssc $\xbm (X)$}}

\path(50,80)(50,40)
\path(48.5,43)(50,40)(51.5,43)
\put(50,81){\circle*{0.3}}
\put(50,39){\circle*{0.3}}
\put(50,37){{\xssc $<x,i>$}}
\put(50,82){{\xssc $<y,j>$}}
\put(52,70){{\xssc $<\xba,k>$}}

\path(20,50)(48,50)
\path(46,51)(48,50)(46,49)
\put(19,50){\circle*{0.3}}
\put(13,47){{\xssc $<z',m'>$}}

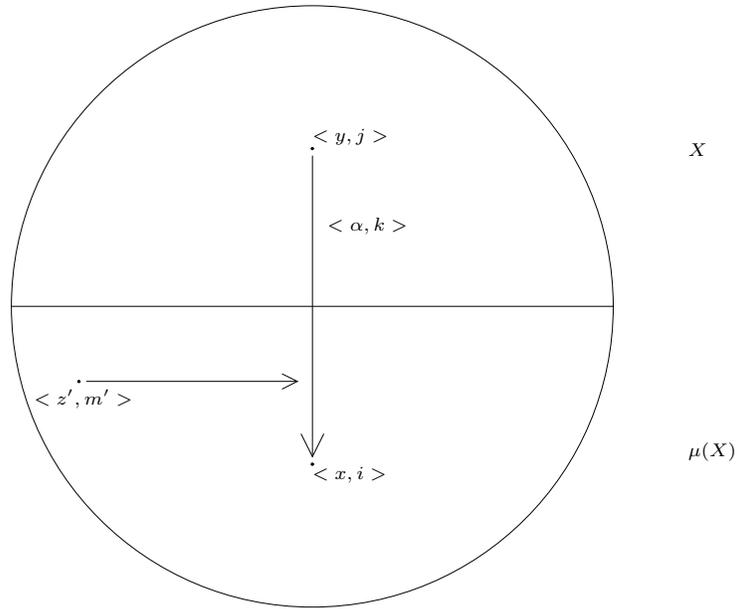
\put(30,10) {{\rm\bf Case 2-1-2}}

\end{picture}
}
\end{diagram}

\vspace{4mm}

 karl-search= End Diagram Essential Smooth 2-1-2
\vspace{7mm}

 *************************************

\vspace{7mm}

\clearpage

$++++++++++++++++++++++++++++++++++++++++++++++++++++++$
\subsubsection{Diagram Essential-Smooth-3-1-2}
\subsubsection{Diagram Essential Smooth 3-1-2}

 {\LARGE karl-search= Start Diagram Essential Smooth 3-1-2 }

\vspace{15mm}

\begin{diagram}

{\xssc LABEL: {Diagram Essential-Smooth-3-1-2}}
\label{Diagram Essential-Smooth-3-1-2}
\index{Diagram Essential-Smooth-3-1-2}

\centering
\setlength{\unitlength}{1mm}
{\renewcommand{\dashlinestretch}{30}
\begin{picture}(100,130)(0,0)
\put(50,80){\circle{80}}
\path(10,80)(90,80)

\put(100,100){{\xssc $X$}}
\put(100,60){{\xssc $\xbm (X)$}}

\path(50,100)(50,60)
\path(48.5,63)(50,60)(51.5,63)
\put(50,101){\circle*{0.3}}
\put(50,59){\circle*{0.3}}
\put(50,57){{\xssc $<x,i>$}}
\put(50,102){{\xssc $<y,j>$}}
\put(52,90){{\xssc $<\xba,k>$}}

\path(20,70)(48,70)
\path(46,71)(48,70)(46,69)
\put(19,70){\circle*{0.3}}
\put(13,67){{\xssc $<z',m'>$}}
\put(26,71){{\xssc $<\xbb',l'>$}}

\path(30,30)(30,68)
\path(29,66)(30,68)(31,66)
\put(30,29){\circle*{0.3}}
\put(31,53){{\xssc $<\xbg ',n'>$}}
\put(30,26.5){{\xssc $<u',p'>$}}

\put(30,10) {{\rm\bf Case 3-1-2}}

\end{picture}
}
\end{diagram}

\vspace{4mm}

 karl-search= End Diagram Essential Smooth 3-1-2
\vspace{7mm}

 *************************************

\vspace{7mm}

\clearpage

$++++++++++++++++++++++++++++++++++++++++++++++++++++++$
\subsubsection{Diagram Essential-Smooth-3-2}
\subsubsection{Diagram Essential Smooth 3-2}

 {\LARGE karl-search= Start Diagram Essential Smooth 3-2 }

\vspace{10mm}

\begin{diagram}

{\xssc LABEL: {Diagram Essential-Smooth-3-2}}
\label{Diagram Essential-Smooth-3-2}
\index{Diagram Essential-Smooth-3-2}

\centering
\setlength{\unitlength}{1mm}
{\renewcommand{\dashlinestretch}{30}
\begin{picture}(100,110)(0,0)
\put(50,60){\circle{80}}
\path(10,60)(90,60)

\put(100,80){{\xssc $X$}}
\put(100,40){{\xssc $\xbm (X)$}}

\path(50,80)(50,40)
\path(48.5,77)(50,80)(51.5,77)
\put(50,81){\circle*{0.3}}
\put(50,39){\circle*{0.3}}
\put(50,37){{\xssc $<y',j'>$}}
\put(50,82){{\xssc $<x,i>$}}
\put(52,70){{\xssc $<\xba',k'>$}}

\path(20,70)(48,70)
\path(46,71)(48,70)(46,69)
\put(19,70){\circle*{0.3}}
\put(13,67){{\xssc $<z',m'>$}}
\put(26,71){{\xssc $<\xbb',l'>$}}

\path(30,30)(30,68)
\path(29,66)(30,68)(31,66)
\put(30,29){\circle*{0.3}}
\put(31,53){{\xssc $<\xbg ',n'>$}}
\put(31,26.5){{\xssc $<u',p'>$}}

\put(30,0) {{\rm\bf Case 3-2}}

\end{picture}
}
\end{diagram}

\vspace{4mm}

 karl-search= End Diagram Essential Smooth 3-2
\vspace{7mm}

 *************************************

\vspace{7mm}

\clearpage

$++++++++++++++++++++++++++++++++++++++++++++++++++++++$
\subsubsection{Diagram Essential-Smooth-2-2}
\subsubsection{Diagram Essential Smooth 2-2}

 {\LARGE karl-search= Start Diagram Essential Smooth 2-2 }

\vspace{10mm}

\begin{diagram}

{\xssc LABEL: {Diagram Essential-Smooth-2-2}}
\label{Diagram Essential-Smooth-2-2}
\index{Diagram Essential-Smooth-2-2}

\centering
\setlength{\unitlength}{1mm}
{\renewcommand{\dashlinestretch}{30}
\begin{picture}(100,100)(0,-10)
\put(50,50){\circle{80}}
\path(10,50)(90,50)

\put(100,70){{\xssc $X$}}
\put(100,30){{\xssc $\xbm (X)$}}

\path(50,70)(50,30)
\path(48.5,67)(50,70)(51.5,67)
\put(50,71){\circle*{0.3}}
\put(50,29){\circle*{0.3}}
\put(50,27){{\xssc $<y',j'>$}}
\put(50,72){{\xssc $<x,i>$}}
\put(52,60){{\xssc $<\xba',k'>$}}

\path(0,60)(48,60)
\path(46,61)(48,60)(46,59)
\put(-1,60){\circle*{0.3}}
\put(-1,57){{\xssc $<z,m>$}}
\put(26,61){{\xssc $<\xbb,l>$}}

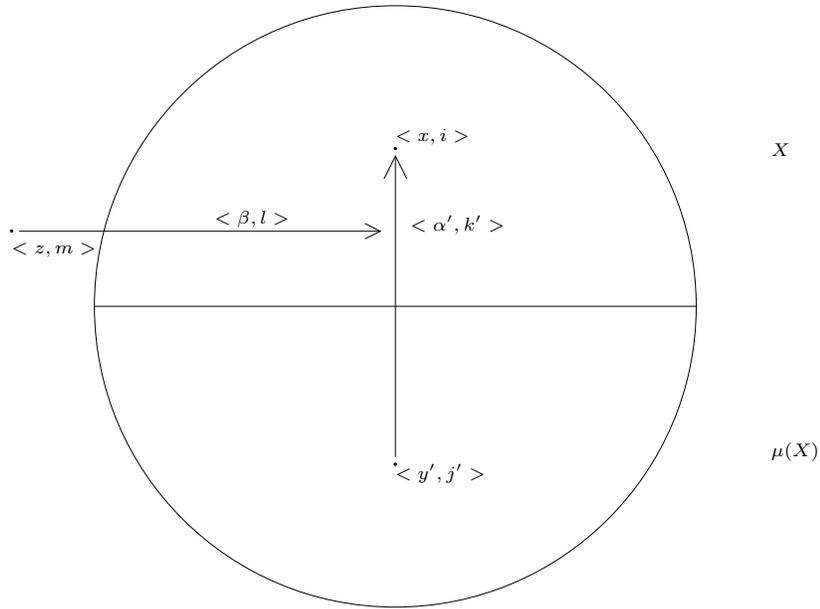
\put(30,0) {{\rm\bf Case 2-2}}

\end{picture}
}
\end{diagram}

\vspace{4mm}

 karl-search= End Diagram Essential Smooth 2-2
\vspace{7mm}

 *************************************

\vspace{7mm}

\clearpage

$++++++++++++++++++++++++++++++++++++++++++++++++++++++$
\subsubsection{Diagram Condition-rho-eta}
\subsubsection{Diagram Condition rho-eta}

 {\LARGE karl-search= Start Diagram Condition rho-eta }

\vspace{10mm}

\begin{diagram}

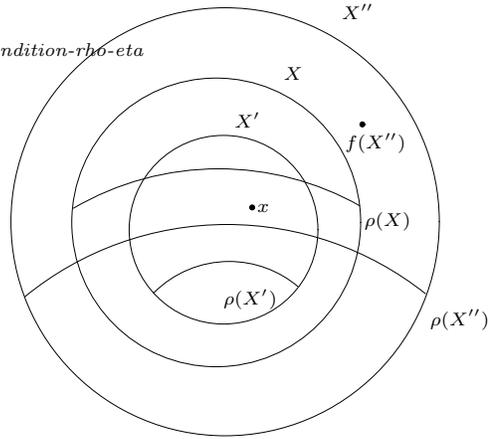
{\xssc LABEL: {Diagram Condition-rho-eta}}
\label{Diagram Condition-rho-eta}
\index{Diagram Condition-rho-eta}

\centering
\setlength{\unitlength}{0.00083333in}
{\renewcommand{\dashlinestretch}{30}
\begin{picture}(2727,2755)(0,-500)
\put(1304.562,-118.818){\arc{3584.794}{4.1827}{5.2252}}
\put(1356.377,-664.180){\arc{3978.201}{4.0248}{5.3839}}
\put(1370.659,430.666){\arc{1327.952}{3.9224}{5.4136}}
\put(1344,1343){\ellipse{2672}{2672}}
\put(1334,1293){\ellipse{1178}{1178}}
\put(1289,1338){\ellipse{1802}{1802}}
\put(2200,1950){\circle*{30}}
\put(2094,1800){{\xssc $f(X'')$}}
\put(2214,1318){{\xssc $\rho(X)$}}
\put(2074,2608){{\xssc $X''$}}
\put(1709,2228){{\xssc $X$}}
\put(2624,698) {{\xssc $\rho(X'')$}}
\put(1334,828) {{\xssc $\rho(X')$}}
\put(1404,1928){{\xssc $X'$}}
\put(1514,1433){\circle*{30}}
\put(1544,1403){{\xssc $x$}}
\put(100,2800){{\xssc For simplicity, $\xbh(X)=X$ here}}

\put(150,-400) {{\rm\bf The complicated case}}

\end{picture}
}
\end{diagram}

\vspace{4mm}

 karl-search= End Diagram Condition rho-eta
\vspace{7mm}

 *************************************

\vspace{7mm}

\clearpage

$++++++++++++++++++++++++++++++++++++++++++++++++++++++$
\subsubsection{Diagram Struc-old-rho-eta}
\subsubsection{Diagram Struc-old-rho-eta}

 {\LARGE karl-search= Start Diagram Struc-old-rho-eta }

\vspace{10mm}

\begin{diagram}

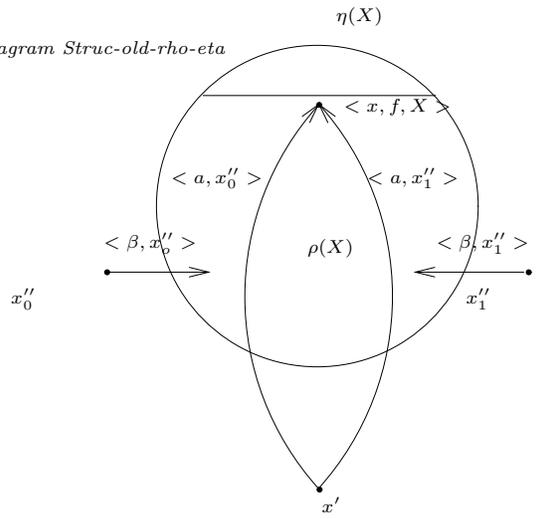
{\xssc LABEL: {Diagram Struc-old-rho-eta}}
\label{Diagram Struc-old-rho-eta}
\index{Diagram Struc-old-rho-eta}

\centering
\setlength{\unitlength}{0.00083333in}
{\renewcommand{\dashlinestretch}{30}
\begin{picture}(3422,3437)(0,-500)
\path(1000,2800)(2450,2800)
\put(386.337,1547.287){\arc{3587.021}{5.5542}{7.0164}}
\path(1823.829,2668.965)(1724.000,2742.000)(1777.715,2630.579)
\put(1724.000,2742.000){\circle*{30}}
\put(3048.125,1544.214){\arc{3578.416}{2.4123}{3.8751}}
\path(1664.387,2631.016)(1719.000,2742.000)(1618.584,2669.774)
\put(1725,340){\circle*{30}}

\put(1711,2111){\ellipse{2006}{2006}}
\path(3000,1697)(2324,1697)
\put(3030,1697){\circle*{30}}
\path(2444.000,1727.000)(2324.000,1697.000)(2444.000,1667.000)
\path(414,1697)(1034,1697)
\put(400,1697){\circle*{30}}
\path(914.000,1667.000)(1034.000,1697.000)(914.000,1727.000)

\put(1879,2700){{\xssc $< x, f, X>$}}
\put(1829,3262){{\xssc $\xbh(X)$}}
\put(1650,1817){{\xssc $\rho(X)$}}
\put(800,2252){{\xssc $< a, x_0''>$}}
\put(2024,2252){{\xssc $< a, x_1''>$}}
\put(2454,1847){{\xssc $< \beta, x_1''>$}}
\put(380,1847) {{\xssc $< \beta, x_o''>$}}
\put(-200, 1500) {{\xssc $x_0''$}}
\put(2639,1500){{\xssc $x''_1$}}
\put(1739,200) {{\xssc $x'$}}

\put(100,-400) {{\rm\bf The full structure}}

\end{picture}
}
\end{diagram}

\vspace{4mm}

 karl-search= End Diagram Struc-old-rho-eta
\vspace{7mm}

 *************************************

\vspace{7mm}

\clearpage

$++++++++++++++++++++++++++++++++++++++++++++++++++++++$
\subsubsection{Diagram Smooth-Level-3}
\subsubsection{Diagram Smooth Level-3}

 {\LARGE karl-search= Start Diagram Smooth Level-3 }

\vspace{10mm}

\begin{diagram}

{\xssc LABEL: {Diagram Smooth-Level-3}}
\label{Diagram Smooth-Level-3}
\index{Diagram Smooth-Level-3}

\centering
\setlength{\unitlength}{0.00083333in}
{\renewcommand{\dashlinestretch}{30}
\begin{picture}(3584,3681)(0,-500)
\put(1980.858,1684.979){\arc{2896.137}{3.0863}{5.1276}}
\path(503.037,1724.492)(535.000,1605.000)(563.029,1725.476)
\put(613.871,2084.678){\arc{1045.415}{1.8998}{5.0662}}
\path(671.312,2573.903)(795.000,2575.000)(685.347,2632.239)
\put(1755.902,1748.510){\arc{3480.502}{1.0969}{3.3803}}
\path(69.787,2036.400)(65.000,2160.000)(11.059,2048.688)
\put(1915.859,2762.471){\arc{1397.908}{3.4506}{5.9440}}
\path(1268.561,3097.293)(1250.000,2975.000)(1323.927,3074.171)
\put(1892.516,1773.709){\arc{3365.589}{4.8901}{7.4299}}
\path(2313.635,3433.788)(2190.000,3430.000)(2300.872,3375.161)
\path(605,1595)(2535,2930)
\path(2453.376,2837.062)(2535.000,2930.000)(2419.243,2886.407)
\path(600,1540)(2545,245)
\put(550,1560){\circle*{30}}
\path(2428.488,286.533)(2545.000,245.000)(2461.741,336.476)
\path(2560,2870)(1820,765)
\path(1831.496,888.158)(1820.000,765.000)(1888.100,868.259)
\path(2530,270)(1745,2345)
\put(2570,220){\circle*{30}}
\put(2600,2950){\circle*{30}}
\path(1815.520,2243.378)(1745.000,2345.000)(1759.401,2222.148)
\put(3315,2745){{\xssc $\gamma_1$}}
\put(2570,2845){{\xssc $y$}}
\put(2575,60)  {{\xssc $y'$}}
\put(1960,950) {{\xssc $\beta_1$}}
\put(1915,1985){{\xssc $\beta_2$}}
\put(2030,2445){{\xssc $\alpha_1$}}
\put(520,1435) {{\xssc $x$}}
\put(1960,3510){{\xssc $\beta_3$}}
\put(750,185)  {{\xssc $\gamma_2$}}
\put(170,2530) {{\xssc $\beta_4$}}
\put(630,2110) {{\xssc $\alpha_3$}}
\put(1920,500) {{\xssc $\alpha_2$}}

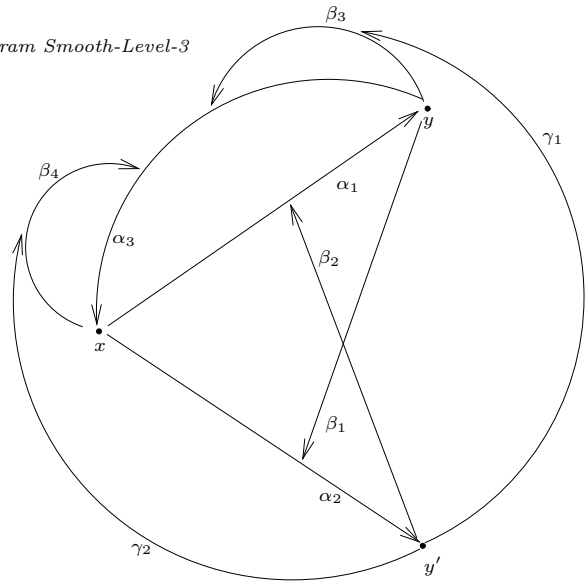
\put(100,-400) {{\rm\bf Solution by smooth level 3 structure}}

\end{picture}
}
\end{diagram}

\vspace{4mm}

 karl-search= End Diagram Smooth Level-3
\vspace{7mm}

 *************************************

\vspace{7mm}

\clearpage

$++++++++++++++++++++++++++++++++++++++++++++++++++++++$
\subsubsection{Diagram Tweety}
\subsubsection{Diagram Tweety}

 {\LARGE karl-search= Start Diagram Tweety }

\vspace{10mm}

\begin{diagram}

{\xssc LABEL: {Diagram Tweety}}
\label{Diagram Tweety}
\index{Diagram Tweety}

\unitlength1.0mm
\begin{picture}(130,100)(0,0)

\newsavebox{\SECHSacht}
\savebox{\SECHSacht}(140,110)[bl]
{

\put(0,95){{\rm\bf The Tweety diagram}}

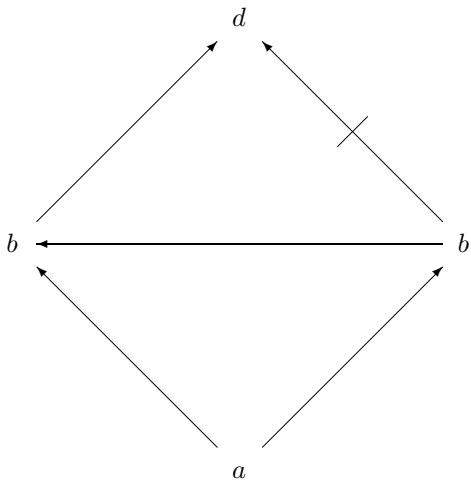
\put(43,27){\vector(1,1){24}}
\put(37,27){\vector(-1,1){24}}
\put(13,57){\vector(1,1){24}}
\put(67,57){\vector(-1,1){24}}

\put(53,67){\line(1,1){4}}

\put(67,54){\vector(-1,0){54}}

\put(39,23){$a$}
\put(9,53){$b$}
\put(69,53){$b$}
\put(39,83){$d$}

}

\put(0,0){\usebox{\SECHSacht}}
\end{picture}

\end{diagram}

\vspace{4mm}

 karl-search= End Diagram Tweety
\vspace{7mm}

 *************************************

\vspace{7mm}

\clearpage

$++++++++++++++++++++++++++++++++++++++++++++++++++++++$
\subsubsection{Diagram Nixon Diamond}
\subsubsection{Diagram Nixon Diamond}

 {\LARGE karl-search= Start Diagram Nixon Diamond }

\vspace{10mm}

\begin{diagram}

{\xssc LABEL: {Diagram Nixon-Diamond}}
\label{Diagram Nixon-Diamond}
\index{Diagram Nixon-Diamond}

\unitlength1.0mm
\begin{picture}(130,100)

\newsavebox{\ZWEIzwei}
\savebox{\ZWEIzwei}(140,110)[bl]
{

\put(0,95){{\rm\bf The Nixon Diamond}}

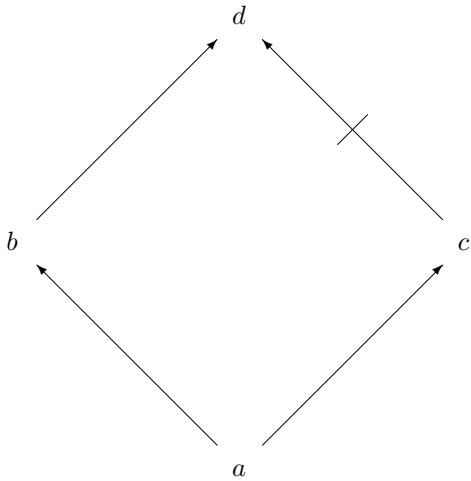
\put(43,27){\vector(1,1){24}}
\put(37,27){\vector(-1,1){24}}
\put(13,57){\vector(1,1){24}}
\put(67,57){\vector(-1,1){24}}

\put(53,67){\line(1,1){4}}

\put(39,23){$a$}
\put(9,53){$b$}
\put(69,53){$c$}
\put(39,83){$d$}

}

\put(0,0){\usebox{\ZWEIzwei}}
\end{picture}

\end{diagram}

\vspace{4mm}

 karl-search= End Diagram Nixon Diamond
\vspace{7mm}

 *************************************

\vspace{7mm}

\clearpage

$++++++++++++++++++++++++++++++++++++++++++++++++++++++$
\subsubsection{Diagram The complicated case}
\subsubsection{Diagram The complicated case}

 {\LARGE karl-search= Start Diagram The complicated case }

\vspace{10mm}

\begin{diagram}

{\xssc LABEL: {Diagram Complicated-Case}}
\label{Diagram Complicated-Case}
\index{Diagram Complicated-Case}

\unitlength1.0mm
\begin{picture}(130,100)

\newsavebox{\Preclusion}
\savebox{\Preclusion}(140,90)[bl]
{

\multiput(43,8)(1,1){5}{\circle*{.3}}
\put(48,13){\vector(1,1){17}}

\multiput(37,8)(-1,1){5}{\circle*{.3}}
\put(32,13){\vector(-1,1){17}}

\put(13,38){\vector(1,1){24}}
\put(67,38){\vector(-1,1){24}}

\put(40,37){\vector(0,1){23}}

\multiput(40,8)(0,1){5}{\circle*{.3}}
\put(40,13){\vector(0,1){17}}

\multiput(66,34)(-1,0){5}{\circle*{.3}}
\put(61,34){\vector(-1,0){17}}

\put(39,3){$x$}
\put(9,33){$u$}
\put(39,33){$v$}
\put(69,33){$z$}
\put(39,63){$y$}

\put(38,50){\line(1,0){3.7}}

\put(10,0) {{\rm\bf The complicated case}}

}

\put(0,0){\usebox{\Preclusion}}
\end{picture}

\end{diagram}

\vspace{4mm}

 karl-search= End Diagram The complicated case
\vspace{7mm}

 *************************************

\vspace{7mm}

\clearpage

$++++++++++++++++++++++++++++++++++++++++++++++++++++++$
\subsubsection{Diagram Upward vs. downward chaining}
\subsubsection{Diagram Upward vs. downward chaining}

 {\LARGE karl-search= Start Diagram Upward vs. downward chaining }

\vspace{10mm}

\begin{diagram}

{\xssc LABEL: {Diagram Up-Down-Chaining}}
\label{Diagram Up-Down-Chaining}
\index{Diagram Up-Down-Chaining}

\unitlength1.0mm
\begin{picture}(130,100)

\newsavebox{\ZWEIvier}
\savebox{\ZWEIvier}(140,110)[bl]
{

\put(0,95){{\rm\bf The problem of downward chaining}}

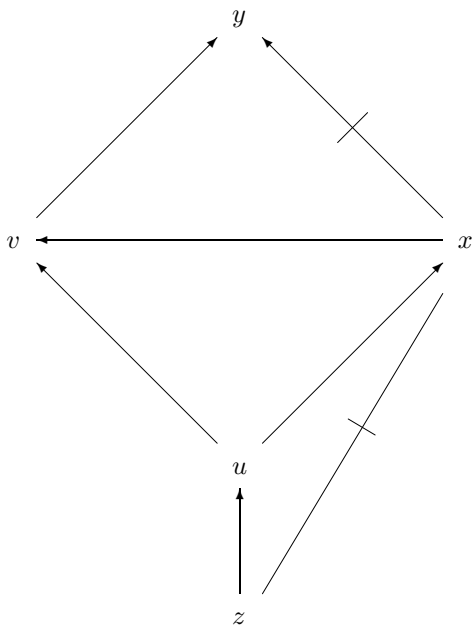
\put(43,27){\vector(1,1){24}}
\put(37,27){\vector(-1,1){24}}
\put(13,57){\vector(1,1){24}}
\put(67,57){\vector(-1,1){24}}

\put(53,67){\line(1,1){4}}

\put(67,54){\vector(-1,0){54}}

\put(40,7){\vector(0,1){14}}
\put(43,7){\line(3,5){24}}
\put(58,28.1){\line(-5,3){3.6}}

\put(39,3){$z$}
\put(39,23){$u$}
\put(9,53){$v$}
\put(69,53){$x$}
\put(39,83){$y$}

}

\put(0,0){\usebox{\ZWEIvier}}
\end{picture}

\end{diagram}

\vspace{4mm}

 karl-search= End Diagram Upward vs. downward chaining
\vspace{7mm}

 *************************************

\vspace{7mm}

\clearpage

$++++++++++++++++++++++++++++++++++++++++++++++++++++++$
\subsubsection{Diagram Split vs. total validity preclusion}
\subsubsection{Diagram Split vs. total validity preclusion}

 {\LARGE karl-search= Start Diagram Split vs. total validity preclusion }

\vspace{10mm}

\begin{diagram}

{\xssc LABEL: {Diagram Split-Total-Preclusion}}
\label{Diagram Split-Total-Preclusion}
\index{Diagram Split-Total-Preclusion}

\unitlength1.0mm
\begin{picture}(130,100)

\newsavebox{\SECHSneun}
\savebox{\SECHSneun}(140,90)[bl]
{

\put(0,75){{\rm\bf Split vs. total validity preclusion}}

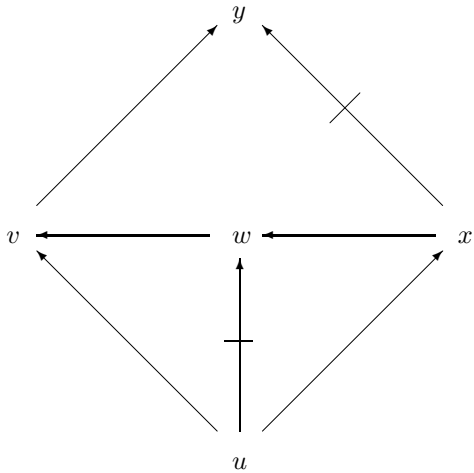
\put(43,8){\vector(1,1){24}}
\put(37,8){\vector(-1,1){24}}
\put(13,38){\vector(1,1){24}}
\put(67,38){\vector(-1,1){24}}

\put(40,8){\vector(0,1){23}}
\put(66,34){\vector(-1,0){23}}
\put(36,34){\vector(-1,0){23}}

\put(39,3){$u$}
\put(9,33){$v$}
\put(39,33){$w$}
\put(69,33){$x$}
\put(39,63){$y$}

\put(38,20){\line(1,0){3.7}}
\put(52,49){\line(1,1){4}}

}

\put(0,0){\usebox{\SECHSneun}}
\end{picture}

\end{diagram}

\vspace{4mm}

 karl-search= End Diagram Split vs. total validity preclusion
\vspace{7mm}

 *************************************

\vspace{7mm}

\clearpage

$++++++++++++++++++++++++++++++++++++++++++++++++++++++$
\subsubsection{Diagram Information transfer}
\subsubsection{Diagram Information transfer}

 {\LARGE karl-search= Start Diagram Information transfer }

\vspace{10mm}

\begin{diagram}

{\xssc LABEL: {Diagram Information-Transfer}}
\label{Diagram Information-Transfer}
\index{Diagram Information-Transfer}

\unitlength1.0mm
\begin{picture}(130,110)

\newsavebox{\Tweety}
\savebox{\Tweety}(130,110)[bl]
{

\put(57,18){\vector(1,1){24}}
\put(51,18){\vector(-1,1){24}}
\put(27,51){\vector(1,1){24}}
\put(81,51){\vector(-1,1){24}}

\put(67,61){\line(1,1){4}}

\put(81,47){\vector(-1,0){54}}

\put(24,51){\vector(0,1){22}}
\put(54,81){\vector(0,1){22}}

\put(53,16){$a$}
\put(23,46){$b$}
\put(83,46){$c$}
\put(53,76){$d$}
\put(23,76){$e$}
\put(53,106){$f$}

\put(10,0) {{\rm\bf Information transfer}}

}

\put(0,0){\usebox{\Tweety}}
\end{picture}

\end{diagram}

\vspace{4mm}

 karl-search= End Diagram Information transfer
\vspace{7mm}

 *************************************

\vspace{7mm}

\clearpage

$++++++++++++++++++++++++++++++++++++++++++++++++++++++$
\subsubsection{Diagram Multiple and conflicting information}
\subsubsection{Diagram Multiple and conflicting information}

 {\LARGE karl-search= Start Diagram Multiple and conflicting information }

\vspace{10mm}

\begin{diagram}

{\xssc LABEL: {Diagram Multiple}}
\label{Diagram Multiple}
\index{Diagram Multiple}

\unitlength1.0mm
\begin{picture}(130,100)

\newsavebox{\sets}
\savebox{\sets}(140,90)[bl]
{

\put(37,8){\vector(-1,1){22}}

\put(13,38){\vector(1,1){24}}
\put(67,38){\vector(-1,1){24}}

\put(40,37){\vector(0,1){23}}

\put(40,8){\vector(0,1){22}}

\put(36,34){\vector(-1,0){22}}

\put(36,4){\vector(-1,0){22}}

\put(39,3){$X$}
\put(9,33){$Y$}
\put(39,33){$Y'$}
\put(69,33){$Y''$}
\put(39,63){$Z$}
\put(9,3){$U$}

\put(38,50){\line(1,0){3.7}}

\put(10,90) {{\rm\bf Multiple and conflicting information}}

}

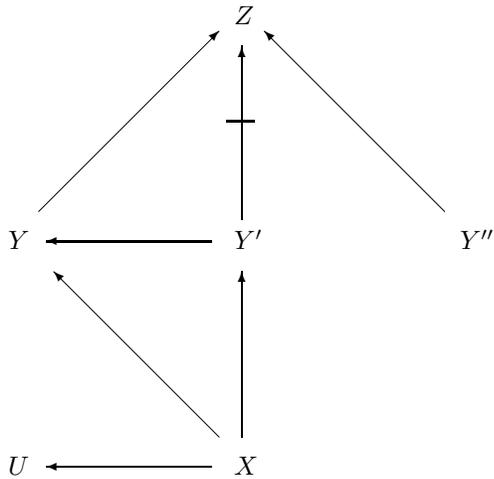
\put(0,0){\usebox{\sets}}
\end{picture}

\end{diagram}

\vspace{4mm}

 karl-search= End Diagram Multiple and conflicting information
\vspace{7mm}

 *************************************

\vspace{7mm}

\clearpage

$++++++++++++++++++++++++++++++++++++++++++++++++++++++$
\subsubsection{Diagram Valid paths vs. valid conclusions}
\subsubsection{Diagram Valid paths vs. valid conclusions}

 {\LARGE karl-search= Start Diagram Valid paths vs. valid conclusions }

\vspace{10mm}

\begin{diagram}

{\xssc LABEL: {Diagram Paths-Conclusions}}
\label{Diagram Paths-Conclusions}
\index{Diagram Paths-Conclusions}

\unitlength1.0mm
\begin{picture}(130,100)

\newsavebox{\setsx}
\savebox{\setsx}(140,90)[bl]
{

\put(37,8){\vector(-1,1){22}}
\put(41,8){\vector(1,1){22}}

\put(13,38){\vector(1,1){24}}
\put(67,38){\vector(-1,1){24}}

\put(40,37){\vector(0,1){23}}

\put(40,8){\vector(0,1){22}}

\put(36,34){\vector(-1,0){22}}
\put(66,34){\vector(-1,0){22}}

\put(39,3){$X$}
\put(9,33){$Y$}
\put(39,33){$Y'$}
\put(69,33){$Y''$}
\put(39,63){$Z$}

\put(38,50){\line(1,0){3.7}}

\put(10,90) {{\rm\bf Valid paths vs. valid conclusions}}

}

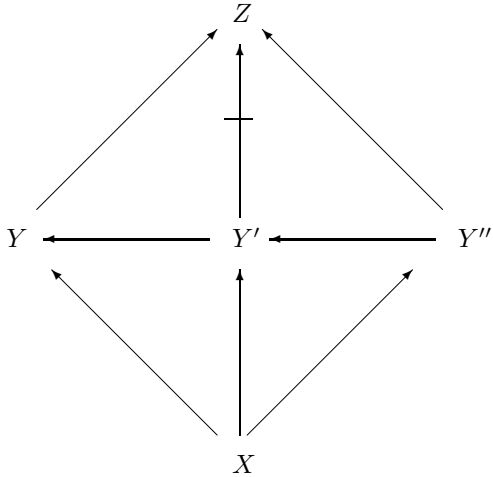
\put(0,0){\usebox{\setsx}}
\end{picture}

\end{diagram}

\vspace{4mm}

 karl-search= End Diagram Valid paths vs. valid conclusions
\vspace{7mm}

 *************************************

\vspace{7mm}

\clearpage

$++++++++++++++++++++++++++++++++++++++++++++++++++++++$
\subsubsection{Diagram WeakTR}
\subsubsection{Diagram WeakTR}

 {\LARGE karl-search= Start Diagram WeakTR }

\vspace{10mm}

\begin{diagram}

{\xssc LABEL: {Diagram WeakTR}}
\label{Diagram WeakTR-b}
\index{Diagram WeakTR}

\centering
\setlength{\unitlength}{1mm}
{\renewcommand{\dashlinestretch}{30}
\begin{picture}(130,90)(0,0)

\put(5,50){\line(1,0){30}}
\put(35,50){\line(1,2){6}}
\put(35,50){\line(1,-2){6}}

\put(5,50){\circle*{1.5}}
\put(35,50){\circle*{1.5}}

\put(41,62){\circle*{1.5}}
\put(41,38){\circle*{1.5}}

\put(5,47){$x$}
\put(32,47){$y$}
\put(43,61){$a$}
\put(43,37){$b$}

\put(65,50){\line(1,0){35}}
\put(100,50){\line(1,2){12}}
\put(100,50){\line(1,-2){12}}

\put(65,50){\circle*{1.5}}
\put(100,50){\circle*{1.5}}

\put(112,74){\circle*{1.5}}
\put(112,26){\circle*{1.5}}

\put(65,47){$x$}
\put(97,47){$y$}
\put(114,73){$a'$}
\put(114,25){$b'$}

\put(29,10){Indiscernible by revision}

\end{picture}

}

\end{diagram}

\vspace{4mm}

 karl-search= End Diagram WeakTR
\vspace{7mm}

 *************************************

\vspace{7mm}

\clearpage

$++++++++++++++++++++++++++++++++++++++++++++++++++++++$
\subsubsection{Diagram CumSmall}
\subsubsection{Diagram CumSmall}

 {\LARGE karl-search= Start Diagram CumSmall }

\vspace{10mm}

\begin{diagram}

{\xssc LABEL: {Diagram CumSmall}}
\label{Diagram CumSmall}
\index{Diagram CumSmall}

\centering
\setlength{\unitlength}{1mm}
{\renewcommand{\dashlinestretch}{30}
\begin{picture}(100,100)(0,0)
\put(50,50){\circle{80}}
\path(24.3,80.6)(75.7,80.6)
\path(84.64,70)(84.64,30)

\put(50,50){{\xssc $B$}}
\put(50,85){{\xssc $A$}}
\put(85,50){{\xssc $A'$}}

\put(30,5) {{\rm\bf Cumulativity}}

\end{picture}
}
\end{diagram}

\vspace{4mm}

 karl-search= End Diagram CumSmall
\vspace{7mm}

 *************************************

\vspace{7mm}

\newpage

%
%
%

\end{document}